\documentstyle[11pt]{article}
\newtheorem{theorem}{Theorem}[section]
\newtheorem{lemma}[theorem]{Lemma}
\newtheorem{corollary}[theorem]{Corollary}

\newtheorem{proposition}[theorem]{Proposition}

\makeindex \makeglossary
\begin{document}                
%
%

\long\def\ig#1{\relax}
\ig{Thanks to Roberto Minio for this def'n.  Compare the def'n of
\comment in AMSTeX.}

\newcount \coefa
\newcount \coefb
\newcount \coefc
\newcount\tempcounta
\newcount\tempcountb
\newcount\tempcountc
\newcount\tempcountd
\newcount\xext
\newcount\yext
\newcount\xoff
\newcount\yoff
\newcount\gap%
\newcount\arrowtypea
\newcount\arrowtypeb
\newcount\arrowtypec
\newcount\arrowtyped
\newcount\arrowtypee
\newcount\height
\newcount\width
\newcount\xpos
\newcount\ypos
\newcount\run
\newcount\rise
\newcount\arrowlength
\newcount\halflength
\newcount\arrowtype
\newdimen\tempdimen
\newdimen\xlen
\newdimen\ylen
\newsavebox{\tempboxa}%
\newsavebox{\tempboxb}%
\newsavebox{\tempboxc}%

\makeatletter
\setlength{\unitlength}{.01em}%
\def\settypes(#1,#2,#3){\arrowtypea#1 \arrowtypeb#2 \arrowtypec#3}
\def\settoheight#1#2{\setbox\@tempboxa\hbox{#2}#1\ht\@tempboxa\relax}%
\def\settodepth#1#2{\setbox\@tempboxa\hbox{#2}#1\dp\@tempboxa\relax}%
\def\settokens[#1`#2`#3`#4]{%
     \def\tokena{#1}\def\tokenb{#2}\def\tokenc{#3}\def\tokend{#4}}
\def\setsqparms[#1`#2`#3`#4;#5`#6]{%
\arrowtypea #1
\arrowtypeb #2
\arrowtypec #3
\arrowtyped #4
\width #5
\height #6
}
\def\setpos(#1,#2){\xpos=#1 \ypos#2}

\def\bfig{\begin{picture}(\xext,\yext)(\xoff,\yoff)}
\def\efig{\end{picture}}

\def\putbox(#1,#2)#3{\put(#1,#2){\makebox(0,0){$#3$}}}

\def\settriparms[#1`#2`#3;#4]{\settripairparms[#1`#2`#3`1`1;#4]}%

\def\settripairparms[#1`#2`#3`#4`#5;#6]{%
\arrowtypea #1
\arrowtypeb #2
\arrowtypec #3
\arrowtyped #4
\arrowtypee #5
\width #6
\height #6
}

\def\resetparms{\settripairparms[1`1`1`1`1;500]\width 500}

\resetparms

\def\mvector(#1,#2)#3{
\put(0,0){\vector(#1,#2){#3}}%
\put(0,0){\vector(#1,#2){30}}%
}
\def\evector(#1,#2)#3{{
\arrowlength #3
\put(0,0){\vector(#1,#2){\arrowlength}}%
\advance \arrowlength by-30
\put(0,0){\vector(#1,#2){\arrowlength}}%
}}

\def\horsize#1#2{%
\settowidth{\tempdimen}{$#2$}%
#1=\tempdimen
\divide #1 by\unitlength
}

\def\vertsize#1#2{%
\settoheight{\tempdimen}{$#2$}%
#1=\tempdimen
\settodepth{\tempdimen}{$#2$}%
\advance #1 by\tempdimen
\divide #1 by\unitlength
}

\def\vertadjust[#1`#2`#3]{%
\vertsize{\tempcounta}{#1}%
\vertsize{\tempcountb}{#2}%
\ifnum \tempcounta<\tempcountb \tempcounta=\tempcountb \fi
\divide\tempcounta by2
\vertsize{\tempcountb}{#3}%
\ifnum \tempcountb>0 \advance \tempcountb by20 \fi
\ifnum \tempcounta<\tempcountb \tempcounta=\tempcountb \fi
}

\def\horadjust[#1`#2`#3]{%
\horsize{\tempcounta}{#1}%
\horsize{\tempcountb}{#2}%
\ifnum \tempcounta<\tempcountb \tempcounta=\tempcountb \fi
\divide\tempcounta by20
\horsize{\tempcountb}{#3}%
\ifnum \tempcountb>0 \advance \tempcountb by60 \fi
\ifnum \tempcounta<\tempcountb \tempcounta=\tempcountb \fi
}

\ig{ In this procedure, #1 is the paramater that sticks out all the way,
#2 sticks out the least and #3 is a label sticking out half way.  #4 is
the amount of the offset.}

\def\sladjust[#1`#2`#3]#4{%
\tempcountc=#4
\horsize{\tempcounta}{#1}%
\divide \tempcounta by2
\horsize{\tempcountb}{#2}%
\divide \tempcountb by2
\advance \tempcountb by-\tempcountc
\ifnum \tempcounta<\tempcountb \tempcounta=\tempcountb\fi
\divide \tempcountc by2
\horsize{\tempcountb}{#3}%
\advance \tempcountb by-\tempcountc
\ifnum \tempcountb>0 \advance \tempcountb by80\fi
\ifnum \tempcounta<\tempcountb \tempcounta=\tempcountb\fi
\advance\tempcounta by20
}

\def\putvector(#1,#2)(#3,#4)#5#6{{%
\xpos=#1
\ypos=#2
\run=#3
\rise=#4
\arrowlength=#5
\arrowtype=#6
\ifnum \arrowtype<0
    \ifnum \run=0
        \advance \ypos by-\arrowlength
    \else
        \tempcounta \arrowlength
        \multiply \tempcounta by\rise
        \divide \tempcounta by\run
        \ifnum\run>0
            \advance \xpos by\arrowlength
            \advance \ypos by\tempcounta
        \else
            \advance \xpos by-\arrowlength
            \advance \ypos by-\tempcounta
        \fi
    \fi
    \multiply \arrowtype by-1
    \multiply \rise by-1
    \multiply \run by-1
\fi
\ifnum \arrowtype=1
    \put(\xpos,\ypos){\vector(\run,\rise){\arrowlength}}%
\else\ifnum \arrowtype=2
    \put(\xpos,\ypos){\mvector(\run,\rise)\arrowlength}%
\else\ifnum\arrowtype=3
    \put(\xpos,\ypos){\evector(\run,\rise){\arrowlength}}%
\fi\fi\fi
}}

\def\putsplitvector(#1,#2)#3#4{
\xpos #1
\ypos #2
\arrowtype #4
\halflength #3
\arrowlength #3
\gap 140
\advance \halflength by-\gap
\divide \halflength by2
\ifnum \arrowtype=1
    \put(\xpos,\ypos){\line(0,-1){\halflength}}%
    \advance\ypos by-\halflength
    \advance\ypos by-\gap
    \put(\xpos,\ypos){\vector(0,-1){\halflength}}%
\else\ifnum \arrowtype=2
    \put(\xpos,\ypos){\line(0,-1)\halflength}%
    \put(\xpos,\ypos){\vector(0,-1)3}%
    \advance\ypos by-\halflength
    \advance\ypos by-\gap
    \put(\xpos,\ypos){\vector(0,-1){\halflength}}%
\else\ifnum\arrowtype=3
    \put(\xpos,\ypos){\line(0,-1)\halflength}%
    \advance\ypos by-\halflength
    \advance\ypos by-\gap
    \put(\xpos,\ypos){\evector(0,-1){\halflength}}%
\else\ifnum \arrowtype=-1
    \advance \ypos by-\arrowlength
    \put(\xpos,\ypos){\line(0,1){\halflength}}%
    \advance\ypos by\halflength
    \advance\ypos by\gap
    \put(\xpos,\ypos){\vector(0,1){\halflength}}%
\else\ifnum \arrowtype=-2
    \advance \ypos by-\arrowlength
    \put(\xpos,\ypos){\line(0,1)\halflength}%
    \put(\xpos,\ypos){\vector(0,1)3}%
    \advance\ypos by\halflength
    \advance\ypos by\gap
    \put(\xpos,\ypos){\vector(0,1){\halflength}}%
\else\ifnum\arrowtype=-3
    \advance \ypos by-\arrowlength
    \put(\xpos,\ypos){\line(0,1)\halflength}%
    \advance\ypos by\halflength
    \advance\ypos by\gap
    \put(\xpos,\ypos){\evector(0,1){\halflength}}%
\fi\fi\fi\fi\fi\fi
}

\def\putmorphism(#1)(#2,#3)[#4`#5`#6]#7#8#9{{%
\run #2
\rise #3
\ifnum\rise=0
  \puthmorphism(#1)[#4`#5`#6]{#7}{#8}{#9}%
\else\ifnum\run=0
  \putvmorphism(#1)[#4`#5`#6]{#7}{#8}{#9}%
\else
\setpos(#1)%
\arrowlength #7
\arrowtype #8
\ifnum\run=0
\else\ifnum\rise=0
\else
\ifnum\run>0
    \coefa=1
\else
   \coefa=-1
\fi
\ifnum\arrowtype>0
   \coefb=0
   \coefc=-1
\else
   \coefb=\coefa
   \coefc=1
   \arrowtype=-\arrowtype
\fi
\width=2
\multiply \width by\run
\divide \width by\rise
\ifnum \width<0  \width=-\width\fi
\advance\width by60
\if l#9 \width=-\width\fi
\putbox(\xpos,\ypos){#4}
{\multiply \coefa by\arrowlength
\advance\xpos by\coefa
\multiply \coefa by\rise
\divide \coefa by\run
\advance \ypos by\coefa
\putbox(\xpos,\ypos){#5} }%
{\multiply \coefa by\arrowlength
\divide \coefa by2
\advance \xpos by\coefa
\advance \xpos by\width
\multiply \coefa by\rise
\divide \coefa by\run
\advance \ypos by\coefa
\if l#9%
   \put(\xpos,\ypos){\makebox(0,0)[r]{$#6$}}%
\else\if r#9%
   \put(\xpos,\ypos){\makebox(0,0)[l]{$#6$}}%
\fi\fi }%
{\multiply \rise by-\coefc
\multiply \run by-\coefc
\multiply \coefb by\arrowlength
\advance \xpos by\coefb
\multiply \coefb by\rise
\divide \coefb by\run
\advance \ypos by\coefb
\multiply \coefc by70
\advance \ypos by\coefc
\multiply \coefc by\run
\divide \coefc by\rise
\advance \xpos by\coefc
\multiply \coefa by140
\multiply \coefa by\run
\divide \coefa by\rise
\advance \arrowlength by\coefa
\ifnum \arrowtype=1
   \put(\xpos,\ypos){\vector(\run,\rise){\arrowlength}}%
\else\ifnum\arrowtype=2
   \put(\xpos,\ypos){\mvector(\run,\rise){\arrowlength}}%
\else\ifnum\arrowtype=3
   \put(\xpos,\ypos){\evector(\run,\rise){\arrowlength}}%
\fi\fi\fi}\fi\fi\fi\fi}}

\def\puthmorphism(#1,#2)[#3`#4`#5]#6#7#8{{%
\xpos #1
\ypos #2
\width #6
\arrowlength #6
\putbox(\xpos,\ypos){#3\vphantom{#4}}%
{\advance \xpos by\arrowlength
\putbox(\xpos,\ypos){\vphantom{#3}#4}}%
\horsize{\tempcounta}{#3}%
\horsize{\tempcountb}{#4}%
\divide \tempcounta by2
\divide \tempcountb by2
\advance \tempcounta by30
\advance \tempcountb by30
\advance \xpos by\tempcounta
\advance \arrowlength by-\tempcounta
\advance \arrowlength by-\tempcountb
\putvector(\xpos,\ypos)(1,0){\arrowlength}{#7}%
\divide \arrowlength by2
\advance \xpos by\arrowlength
\vertsize{\tempcounta}{#5}%
\divide\tempcounta by2
\advance \tempcounta by20
\if a#8 %
   \advance \ypos by\tempcounta
   \putbox(\xpos,\ypos){#5}%
\else
   \advance \ypos by-\tempcounta
   \putbox(\xpos,\ypos){#5}%
\fi}}

\def\putvmorphism(#1,#2)[#3`#4`#5]#6#7#8{{%
\xpos #1
\ypos #2
\arrowlength #6
\arrowtype #7
\settowidth{\xlen}{$#5$}%
\putbox(\xpos,\ypos){#3}%
{\advance \ypos by-\arrowlength
\putbox(\xpos,\ypos){#4}}%
{\advance\arrowlength by-140
\advance \ypos by-70
\ifdim\xlen>0pt
   \if m#8%
      \putsplitvector(\xpos,\ypos){\arrowlength}{\arrowtype}%
   \else
      \putvector(\xpos,\ypos)(0,-1){\arrowlength}{\arrowtype}%
   \fi
\else
   \putvector(\xpos,\ypos)(0,-1){\arrowlength}{\arrowtype}%
\fi}%
\ifdim\xlen>0pt
   \divide \arrowlength by2
   \advance\ypos by-\arrowlength
   \if l#8%
      \advance \xpos by-40
      \put(\xpos,\ypos){\makebox(0,0)[r]{$#5$}}%
   \else\if r#8%
      \advance \xpos by40
      \put(\xpos,\ypos){\makebox(0,0)[l]{$#5$}}%
   \else
      \putbox(\xpos,\ypos){#5}%
   \fi\fi
\fi
}}

\def\topadjust[#1`#2`#3]{%
\yoff=10
\vertadjust[#1`#2`{#3}]%
\advance \yext by\tempcounta
\advance \yext by 10
}
\def\botadjust[#1`#2`#3]{%
\vertadjust[#1`#2`{#3}]%
\advance \yext by\tempcounta
\advance \yoff by-\tempcounta
}
\def\leftadjust[#1`#2`#3]{%
\xoff=0
\horadjust[#1`#2`{#3}]%
\advance \xext by\tempcounta
\advance \xoff by-\tempcounta
}
\def\rightadjust[#1`#2`#3]{%
\horadjust[#1`#2`{#3}]%
\advance \xext by\tempcounta
}
\def\rightsladjust[#1`#2`#3]{%
\sladjust[#1`#2`{#3}]{\width}%
\advance \xext by\tempcounta
}
\def\leftsladjust[#1`#2`#3]{%
\xoff=0
\sladjust[#1`#2`{#3}]{\width}%
\advance \xext by\tempcounta
\advance \xoff by-\tempcounta
}
\def\adjust[#1`#2;#3`#4;#5`#6;#7`#8]{%
\topadjust[#1``{#2}]
\leftadjust[#3``{#4}]
\rightadjust[#5``{#6}]
\botadjust[#7``{#8}]}

\def\putsquarep<#1>(#2)[#3;#4`#5`#6`#7]{{%
\setsqparms[#1]%
\setpos(#2)%
\settokens[#3]%
\puthmorphism(\xpos,\ypos)[\tokenc`\tokend`{#7}]{\width}{\arrowtyped}b%
\advance\ypos by \height
\puthmorphism(\xpos,\ypos)[\tokena`\tokenb`{#4}]{\width}{\arrowtypea}a%
\putvmorphism(\xpos,\ypos)[``{#5}]{\height}{\arrowtypeb}l%
\advance\xpos by \width
\putvmorphism(\xpos,\ypos)[``{#6}]{\height}{\arrowtypec}r%
}}

\def\putsquare{\@ifnextchar <{\putsquarep}{\putsquarep%
   <\arrowtypea`\arrowtypeb`\arrowtypec`\arrowtyped;\width`\height>}}
\def\square{\@ifnextchar< {\squarep}{\squarep
   <\arrowtypea`\arrowtypeb`\arrowtypec`\arrowtyped;\width`\height>}}
\def\squarep<#1>[#2`#3`#4`#5;#6`#7`#8`#9]{{
\setsqparms[#1]
\xext=\width                                          
\yext=\height                                         
\topadjust[#2`#3`{#6}]
\botadjust[#4`#5`{#9}]
\leftadjust[#2`#4`{#7}]
\rightadjust[#3`#5`{#8}]
\begin{picture}(\xext,\yext)(\xoff,\yoff)
\putsquarep<\arrowtypea`\arrowtypeb`\arrowtypec`\arrowtyped;\width`\height>%
(0,0)[#2`#3`#4`#5;#6`#7`#8`{#9}]%
\end{picture}%
}}

\def\putptrianglep<#1>(#2,#3)[#4`#5`#6;#7`#8`#9]{{%
\settriparms[#1]%
\xpos=#2 \ypos=#3
\advance\ypos by \height
\puthmorphism(\xpos,\ypos)[#4`#5`{#7}]{\height}{\arrowtypea}a%
\putvmorphism(\xpos,\ypos)[`#6`{#8}]{\height}{\arrowtypeb}l%
\advance\xpos by\height
\putmorphism(\xpos,\ypos)(-1,-1)[``{#9}]{\height}{\arrowtypec}r%
}}

\def\putptriangle{\@ifnextchar <{\putptrianglep}{\putptrianglep
   <\arrowtypea`\arrowtypeb`\arrowtypec;\height>}}
\def\ptriangle{\@ifnextchar <{\ptrianglep}{\ptrianglep
   <\arrowtypea`\arrowtypeb`\arrowtypec;\height>}}

\def\ptrianglep<#1>[#2`#3`#4;#5`#6`#7]{{
\settriparms[#1]%
\width=\height                         
\xext=\width                           
\yext=\width                           
\topadjust[#2`#3`{#5}]
\botadjust[#3``]
\leftadjust[#2`#4`{#6}]
\rightsladjust[#3`#4`{#7}]
\begin{picture}(\xext,\yext)(\xoff,\yoff)
\putptrianglep<\arrowtypea`\arrowtypeb`\arrowtypec;\height>%
(0,0)[#2`#3`#4;#5`#6`{#7}]%
\end{picture}%
}}

\def\putqtrianglep<#1>(#2,#3)[#4`#5`#6;#7`#8`#9]{{%
\settriparms[#1]%
\xpos=#2 \ypos=#3
\advance\ypos by\height
\puthmorphism(\xpos,\ypos)[#4`#5`{#7}]{\height}{\arrowtypea}a%
\putmorphism(\xpos,\ypos)(1,-1)[``{#8}]{\height}{\arrowtypeb}l%
\advance\xpos by\height
\putvmorphism(\xpos,\ypos)[`#6`{#9}]{\height}{\arrowtypec}r%
}}

\def\putqtriangle{\@ifnextchar <{\putqtrianglep}{\putqtrianglep
   <\arrowtypea`\arrowtypeb`\arrowtypec;\height>}}
\def\qtriangle{\@ifnextchar <{\qtrianglep}{\qtrianglep
   <\arrowtypea`\arrowtypeb`\arrowtypec;\height>}}

\def\qtrianglep<#1>[#2`#3`#4;#5`#6`#7]{{
\settriparms[#1]
\width=\height                         
\xext=\width                           
\yext=\height                          
\topadjust[#2`#3`{#5}]
\botadjust[#4``]
\leftsladjust[#2`#4`{#6}]
\rightadjust[#3`#4`{#7}]
\begin{picture}(\xext,\yext)(\xoff,\yoff)
\putqtrianglep<\arrowtypea`\arrowtypeb`\arrowtypec;\height>%
(0,0)[#2`#3`#4;#5`#6`{#7}]%
\end{picture}%
}}

\def\putdtrianglep<#1>(#2,#3)[#4`#5`#6;#7`#8`#9]{{%
\settriparms[#1]%
\xpos=#2 \ypos=#3
\puthmorphism(\xpos,\ypos)[#5`#6`{#9}]{\height}{\arrowtypec}b%
\advance\xpos by \height \advance\ypos by\height
\putmorphism(\xpos,\ypos)(-1,-1)[``{#7}]{\height}{\arrowtypea}l%
\putvmorphism(\xpos,\ypos)[#4``{#8}]{\height}{\arrowtypeb}r%
}}

\def\putdtriangle{\@ifnextchar <{\putdtrianglep}{\putdtrianglep
   <\arrowtypea`\arrowtypeb`\arrowtypec;\height>}}
\def\dtriangle{\@ifnextchar <{\dtrianglep}{\dtrianglep
   <\arrowtypea`\arrowtypeb`\arrowtypec;\height>}}

\def\dtrianglep<#1>[#2`#3`#4;#5`#6`#7]{{
\settriparms[#1]
\width=\height                         
\xext=\width                           
\yext=\height                          
\topadjust[#2``]
\botadjust[#3`#4`{#7}]
\leftsladjust[#3`#2`{#5}]
\rightadjust[#2`#4`{#6}]
\begin{picture}(\xext,\yext)(\xoff,\yoff)
\putdtrianglep<\arrowtypea`\arrowtypeb`\arrowtypec;\height>%
(0,0)[#2`#3`#4;#5`#6`{#7}]%
\end{picture}%
}}

\def\putbtrianglep<#1>(#2,#3)[#4`#5`#6;#7`#8`#9]{{%
\settriparms[#1]%
\xpos=#2 \ypos=#3
\puthmorphism(\xpos,\ypos)[#5`#6`{#9}]{\height}{\arrowtypec}b%
\advance\ypos by\height
\putmorphism(\xpos,\ypos)(1,-1)[``{#8}]{\height}{\arrowtypeb}r%
\putvmorphism(\xpos,\ypos)[#4``{#7}]{\height}{\arrowtypea}l%
}}

\def\putbtriangle{\@ifnextchar <{\putbtrianglep}{\putbtrianglep
   <\arrowtypea`\arrowtypeb`\arrowtypec;\height>}}
\def\btriangle{\@ifnextchar <{\btrianglep}{\btrianglep
   <\arrowtypea`\arrowtypeb`\arrowtypec;\height>}}

\def\btrianglep<#1>[#2`#3`#4;#5`#6`#7]{{
\settriparms[#1]
\width=\height                         
\xext=\width                           
\yext=\height                          
\topadjust[#2``]
\botadjust[#3`#4`{#7}]
\leftadjust[#2`#3`{#5}]
\rightsladjust[#4`#2`{#6}]
\begin{picture}(\xext,\yext)(\xoff,\yoff)
\putbtrianglep<\arrowtypea`\arrowtypeb`\arrowtypec;\height>%
(0,0)[#2`#3`#4;#5`#6`{#7}]%
\end{picture}%
}}

\def\putAtrianglep<#1>(#2,#3)[#4`#5`#6;#7`#8`#9]{{%
\settriparms[#1]%
\xpos=#2 \ypos=#3
{\multiply \height by2
\puthmorphism(\xpos,\ypos)[#5`#6`{#9}]{\height}{\arrowtypec}b}%
\advance\xpos by\height \advance\ypos by\height
\putmorphism(\xpos,\ypos)(-1,-1)[#4``{#7}]{\height}{\arrowtypea}l%
\putmorphism(\xpos,\ypos)(1,-1)[``{#8}]{\height}{\arrowtypeb}r%
}}

\def\putAtriangle{\@ifnextchar <{\putAtrianglep}{\putAtrianglep
   <\arrowtypea`\arrowtypeb`\arrowtypec;\height>}}
\def\Atriangle{\@ifnextchar <{\Atrianglep}{\Atrianglep
   <\arrowtypea`\arrowtypeb`\arrowtypec;\height>}}

\def\Atrianglep<#1>[#2`#3`#4;#5`#6`#7]{{
\settriparms[#1]
\width=\height                         
\xext=\width                           
\yext=\height                          
\topadjust[#2``]
\botadjust[#3`#4`{#7}]
\multiply \xext by2 
\leftsladjust[#3`#2`{#5}]
\rightsladjust[#4`#2`{#6}]
\begin{picture}(\xext,\yext)(\xoff,\yoff)%
\putAtrianglep<\arrowtypea`\arrowtypeb`\arrowtypec;\height>%
(0,0)[#2`#3`#4;#5`#6`{#7}]%
\end{picture}%
}}

\def\putAtrianglepairp<#1>(#2)[#3;#4`#5`#6`#7`#8]{{
\settripairparms[#1]%
\setpos(#2)%
\settokens[#3]%
\puthmorphism(\xpos,\ypos)[\tokenb`\tokenc`{#7}]{\height}{\arrowtyped}b%
\advance\xpos by\height
\advance\ypos by\height
\putmorphism(\xpos,\ypos)(-1,-1)[\tokena``{#4}]{\height}{\arrowtypea}l%
\putvmorphism(\xpos,\ypos)[``{#5}]{\height}{\arrowtypeb}m%
\putmorphism(\xpos,\ypos)(1,-1)[``{#6}]{\height}{\arrowtypec}r%
}}

\def\putAtrianglepair{\@ifnextchar <{\putAtrianglepairp}{\putAtrianglepairp%
   <\arrowtypea`\arrowtypeb`\arrowtypec`\arrowtyped`\arrowtypee;\height>}}
\def\Atrianglepair{\@ifnextchar <{\Atrianglepairp}{\Atrianglepairp%
   <\arrowtypea`\arrowtypeb`\arrowtypec`\arrowtyped`\arrowtypee;\height>}}

\def\Atrianglepairp<#1>[#2;#3`#4`#5`#6`#7]{{%
\settripairparms[#1]%
\settokens[#2]%
\width=\height
\xext=\width
\yext=\height
\topadjust[\tokena``]%
\vertadjust[\tokenb`\tokenc`{#6}]
\tempcountd=\tempcounta                       
\vertadjust[\tokenc`\tokend`{#7}]
\ifnum\tempcounta<\tempcountd                 
\tempcounta=\tempcountd\fi                    
\advance \yext by\tempcounta                  
\advance \yoff by-\tempcounta                 %
\multiply \xext by2 
\leftsladjust[\tokenb`\tokena`{#3}]
\rightsladjust[\tokend`\tokena`{#5}]%
\begin{picture}(\xext,\yext)(\xoff,\yoff)%
\putAtrianglepairp
<\arrowtypea`\arrowtypeb`\arrowtypec`\arrowtyped`\arrowtypee;\height>%
(0,0)[#2;#3`#4`#5`#6`{#7}]%
\end{picture}%
}}

\def\putVtrianglep<#1>(#2,#3)[#4`#5`#6;#7`#8`#9]{{%
\settriparms[#1]%
\xpos=#2 \ypos=#3
\advance\ypos by\height
{\multiply\height by2
\puthmorphism(\xpos,\ypos)[#4`#5`{#7}]{\height}{\arrowtypea}a}%
\putmorphism(\xpos,\ypos)(1,-1)[`#6`{#8}]{\height}{\arrowtypeb}l%
\advance\xpos by\height
\advance\xpos by\height
\putmorphism(\xpos,\ypos)(-1,-1)[``{#9}]{\height}{\arrowtypec}r%
}}

\def\putVtriangle{\@ifnextchar <{\putVtrianglep}{\putVtrianglep
   <\arrowtypea`\arrowtypeb`\arrowtypec;\height>}}
\def\Vtriangle{\@ifnextchar <{\Vtrianglep}{\Vtrianglep
   <\arrowtypea`\arrowtypeb`\arrowtypec;\height>}}

\def\Vtrianglep<#1>[#2`#3`#4;#5`#6`#7]{{
\settriparms[#1]
\width=\height                         
\xext=\width                           
\yext=\height                          
\topadjust[#2`#3`{#5}]
\botadjust[#4``]
\multiply \xext by2 
\leftsladjust[#2`#3`{#6}]
\rightsladjust[#3`#4`{#7}]
\begin{picture}(\xext,\yext)(\xoff,\yoff)%
\putVtrianglep<\arrowtypea`\arrowtypeb`\arrowtypec;\height>%
(0,0)[#2`#3`#4;#5`#6`{#7}]%
\end{picture}%
}}

\def\putVtrianglepairp<#1>(#2)[#3;#4`#5`#6`#7`#8]{{
\settripairparms[#1]%
\setpos(#2)%
\settokens[#3]%
\advance\ypos by\height
\putmorphism(\xpos,\ypos)(1,-1)[`\tokend`{#6}]{\height}{\arrowtypec}l%
\puthmorphism(\xpos,\ypos)[\tokena`\tokenb`{#4}]{\height}{\arrowtypea}a%
\advance\xpos by\height
\putvmorphism(\xpos,\ypos)[``{#7}]{\height}{\arrowtyped}m%
\advance\xpos by\height
\putmorphism(\xpos,\ypos)(-1,-1)[``{#8}]{\height}{\arrowtypee}r%
}}

\def\putVtrianglepair{\@ifnextchar <{\putVtrianglepairp}{\putVtrianglepairp%
    <\arrowtypea`\arrowtypeb`\arrowtypec`\arrowtyped`\arrowtypee;\height>}}
\def\Vtrianglepair{\@ifnextchar <{\Vtrianglepairp}{\Vtrianglepairp%
    <\arrowtypea`\arrowtypeb`\arrowtypec`\arrowtyped`\arrowtypee;\height>}}

\def\Vtrianglepairp<#1>[#2;#3`#4`#5`#6`#7]{{%
\settripairparms[#1]%
\settokens[#2]
\xext=\height                  
\width=\height                 
\yext=\height                  
\vertadjust[\tokena`\tokenb`{#4}]
\tempcountd=\tempcounta        
\vertadjust[\tokenb`\tokenc`{#5}]
\ifnum\tempcounta<\tempcountd%
\tempcounta=\tempcountd\fi
\advance \yext by\tempcounta
\botadjust[\tokend``]%
\multiply \xext by2
\leftsladjust[\tokena`\tokend`{#6}]%
\rightsladjust[\tokenc`\tokend`{#7}]%
\begin{picture}(\xext,\yext)(\xoff,\yoff)%
\putVtrianglepairp
<\arrowtypea`\arrowtypeb`\arrowtypec`\arrowtyped`\arrowtypee;\height>%
(0,0)[#2;#3`#4`#5`#6`{#7}]%
\end{picture}%
}}

\def\putCtrianglep<#1>(#2,#3)[#4`#5`#6;#7`#8`#9]{{%
\settriparms[#1]%
\xpos=#2 \ypos=#3
\advance\ypos by\height
\putmorphism(\xpos,\ypos)(1,-1)[``{#9}]{\height}{\arrowtypec}l%
\advance\xpos by\height
\advance\ypos by\height
\putmorphism(\xpos,\ypos)(-1,-1)[#4`#5`{#7}]{\height}{\arrowtypea}l%
{\multiply\height by 2
\putvmorphism(\xpos,\ypos)[`#6`{#8}]{\height}{\arrowtypeb}r}%
}}

\def\putCtriangle{\@ifnextchar <{\putCtrianglep}{\putCtrianglep
    <\arrowtypea`\arrowtypeb`\arrowtypec;\height>}}
\def\Ctriangle{\@ifnextchar <{\Ctrianglep}{\Ctrianglep
    <\arrowtypea`\arrowtypeb`\arrowtypec;\height>}}

\def\Ctrianglep<#1>[#2`#3`#4;#5`#6`#7]{{
\settriparms[#1]
\width=\height                          
\xext=\width                            
\yext=\height                           
\multiply \yext by2 
\topadjust[#2``]
\botadjust[#4``]
\sladjust[#3`#2`{#5}]{\width}
\tempcountd=\tempcounta                 
\sladjust[#3`#4`{#7}]{\width}
\ifnum \tempcounta<\tempcountd          
\tempcounta=\tempcountd\fi              
\advance \xext by\tempcounta            
\advance \xoff by-\tempcounta           %
\rightadjust[#2`#4`{#6}]
\begin{picture}(\xext,\yext)(\xoff,\yoff)%
\putCtrianglep<\arrowtypea`\arrowtypeb`\arrowtypec;\height>%
(0,0)[#2`#3`#4;#5`#6`{#7}]%
\end{picture}%
}}

\def\putDtrianglep<#1>(#2,#3)[#4`#5`#6;#7`#8`#9]{{%
\settriparms[#1]%
\xpos=#2 \ypos=#3
\advance\xpos by\height \advance\ypos by\height
\putmorphism(\xpos,\ypos)(-1,-1)[``{#9}]{\height}{\arrowtypec}r%
\advance\xpos by-\height \advance\ypos by\height
\putmorphism(\xpos,\ypos)(1,-1)[`#5`{#8}]{\height}{\arrowtypeb}r%
{\multiply\height by 2
\putvmorphism(\xpos,\ypos)[#4`#6`{#7}]{\height}{\arrowtypea}l}%
}}

\def\putDtriangle{\@ifnextchar <{\putDtrianglep}{\putDtrianglep
    <\arrowtypea`\arrowtypeb`\arrowtypec;\height>}}
\def\Dtriangle{\@ifnextchar <{\Dtrianglep}{\Dtrianglep
   <\arrowtypea`\arrowtypeb`\arrowtypec;\height>}}

\def\Dtrianglep<#1>[#2`#3`#4;#5`#6`#7]{{
\settriparms[#1]
\width=\height                         
\xext=\height                          
\yext=\height                          
\multiply \yext by2 
\topadjust[#2``]
\botadjust[#4``]
\leftadjust[#2`#4`{#5}]
\sladjust[#3`#2`{#5}]{\height}
\tempcountd=\tempcountd                
\sladjust[#3`#4`{#7}]{\height}
\ifnum \tempcounta<\tempcountd         
\tempcounta=\tempcountd\fi             
\advance \xext by\tempcounta           %
\begin{picture}(\xext,\yext)(\xoff,\yoff)
\putDtrianglep<\arrowtypea`\arrowtypeb`\arrowtypec;\height>%
(0,0)[#2`#3`#4;#5`#6`{#7}]%
\end{picture}%
}}

\def\setrecparms[#1`#2]{\width=#1 \height=#2}%
%

\def\recursep<#1`#2>[#3;#4`#5`#6`#7`#8]{{%
\width=#1 \height=#2
\settokens[#3]
\settowidth{\tempdimen}{$\tokena$}
\ifdim\tempdimen=0pt
  \savebox{\tempboxa}{\hbox{$\tokenb$}}%
  \savebox{\tempboxb}{\hbox{$\tokend$}}%
  \savebox{\tempboxc}{\hbox{$#6$}}%
\else
  \savebox{\tempboxa}{\hbox{$\hbox{$\tokena$}\times\hbox{$\tokenb$}$}}%
  \savebox{\tempboxb}{\hbox{$\hbox{$\tokena$}\times\hbox{$\tokend$}$}}%
  \savebox{\tempboxc}{\hbox{$\hbox{$\tokena$}\times\hbox{$#6$}$}}%
\fi
\ypos=\height
\divide\ypos by 2
\xpos=\ypos
\advance\xpos by \width
\xext=\xpos \yext=\height
\topadjust[#3`\usebox{\tempboxa}`{#4}]%
\botadjust[#5`\usebox{\tempboxb}`{#8}]%
\sladjust[\tokenc`\tokenb`{#5}]{\ypos}%
\tempcountd=\tempcounta
\sladjust[\tokenc`\tokend`{#5}]{\ypos}%
\ifnum \tempcounta<\tempcountd
\tempcounta=\tempcountd\fi
\advance \xext by\tempcounta
\advance \xoff by-\tempcounta
\rightadjust[\usebox{\tempboxa}`\usebox{\tempboxb}`\usebox{\tempboxc}]%
\bfig
\putCtrianglep<-1`1`1;\ypos>(0,0)[`\tokenc`;#5`#6`{#7}]%
\puthmorphism(\ypos,0)[\tokend`\usebox{\tempboxb}`{#8}]{\width}{-1}b%
\puthmorphism(\ypos,\height)[\tokenb`\usebox{\tempboxa}`{#4}]{\width}{-1}a%
\advance\ypos by \width
\putvmorphism(\ypos,\height)[``\usebox{\tempboxc}]{\height}1r%
\efig
}}

\def\recurse{\@ifnextchar <{\recursep}{\recursep<\width`\height>}}

\def\puttwohmorphisms(#1,#2)[#3`#4;#5`#6]#7#8#9{{%
%
\puthmorphism(#1,#2)[#3`#4`]{#7}0a
\ypos=#2
\advance\ypos by 20
\puthmorphism(#1,\ypos)[\phantom{#3}`\phantom{#4}`#5]{#7}{#8}a
\advance\ypos by -40
\puthmorphism(#1,\ypos)[\phantom{#3}`\phantom{#4}`#6]{#7}{#9}b
}}

\def\puttwovmorphisms(#1,#2)[#3`#4;#5`#6]#7#8#9{{%
%
%
%
\putvmorphism(#1,#2)[#3`#4`]{#7}0a
\xpos=#1
\advance\xpos by -20
\putvmorphism(\xpos,#2)[\phantom{#3}`\phantom{#4}`#5]{#7}{#8}l
\advance\xpos by 40
\putvmorphism(\xpos,#2)[\phantom{#3}`\phantom{#4}`#6]{#7}{#9}r
}}

\def\puthcoequalizer(#1)[#2`#3`#4;#5`#6`#7]#8#9{{%
%
\setpos(#1)%
\puttwohmorphisms(\xpos,\ypos)[#2`#3;#5`#6]{#8}11%
\advance\xpos by #8
\puthmorphism(\xpos,\ypos)[\phantom{#3}`#4`#7]{#8}1{#9}
}}

\def\putvcoequalizer(#1)[#2`#3`#4;#5`#6`#7]#8#9{{%
%
%
%
%
\setpos(#1)%
\puttwovmorphisms(\xpos,\ypos)[#2`#3;#5`#6]{#8}11%
\advance\ypos by -#8
\putvmorphism(\xpos,\ypos)[\phantom{#3}`#4`#7]{#8}1{#9}
}}

\def\putthreehmorphisms(#1)[#2`#3;#4`#5`#6]#7(#8)#9{{%
\setpos(#1) \settypes(#8)
\if a#9 %
     \vertsize{\tempcounta}{#5}%
     \vertsize{\tempcountb}{#6}%
     \ifnum \tempcounta<\tempcountb \tempcounta=\tempcountb \fi
\else
     \vertsize{\tempcounta}{#4}%
     \vertsize{\tempcountb}{#5}%
     \ifnum \tempcounta<\tempcountb \tempcounta=\tempcountb \fi
\fi
\advance \tempcounta by 60
\puthmorphism(\xpos,\ypos)[#2`#3`#5]{#7}{\arrowtypeb}{#9}
\advance\ypos by \tempcounta
\puthmorphism(\xpos,\ypos)[\phantom{#2}`\phantom{#3}`#4]{#7}{\arrowtypea}{#9}
\advance\ypos by -\tempcounta \advance\ypos by -\tempcounta
\puthmorphism(\xpos,\ypos)[\phantom{#2}`\phantom{#3}`#6]{#7}{\arrowtypec}{#9}
}}

\def\putarc(#1,#2)[#3`#4`#5]#6#7#8{{%
\xpos #1
\ypos #2
\width #6
\arrowlength #6
\putbox(\xpos,\ypos){#3\vphantom{#4}}%
{\advance \xpos by\arrowlength
\putbox(\xpos,\ypos){\vphantom{#3}#4}}%
\horsize{\tempcounta}{#3}%
\horsize{\tempcountb}{#4}%
\divide \tempcounta by2
\divide \tempcountb by2
\advance \tempcounta by30
\advance \tempcountb by30
\advance \xpos by\tempcounta
\advance \arrowlength by-\tempcounta
\advance \arrowlength by-\tempcountb
\halflength=\arrowlength \divide\halflength by 2
\divide\arrowlength by 5
\put(\xpos,\ypos){\bezier{\arrowlength}(0,0)(50,50)(\halflength,50)}
\ifnum #7=-1 \put(\xpos,\ypos){\vector(-3,-2)0} \fi
\advance\xpos by \halflength
\put(\xpos,\ypos){\xpos=\halflength \advance\xpos by -50
   \bezier{\arrowlength}(0,50)(\xpos,50)(\halflength,0)}
\ifnum #7=1 {\advance \xpos by
   \halflength \put(\xpos,\ypos){\vector(3,-2)0}} \fi
\advance\ypos by 50
\vertsize{\tempcounta}{#5}%
\divide\tempcounta by2
\advance \tempcounta by20
\if a#8 %
   \advance \ypos by\tempcounta
   \putbox(\xpos,\ypos){#5}%
\else
   \advance \ypos by-\tempcounta
   \putbox(\xpos,\ypos){#5}%
\fi
}}

\makeatother

\sloppy

\newcommand{\nl}{\hspace{2cm}\\ }

\def\nec{\Box}
\def\pos{\Diamond}
\def\diam{{\tiny\Diamond}}

\def\lc{\lceil}
\def\rc{\rceil}
\def\lf{\lfloor}
\def\rf{\rfloor}
\def\lk{\langle}
\def\rk{\rangle}
\def\lse{[\!|}
\def\rse{|\!]}
\def\le{(\!|}
\def\re{|\!)}

\def\homo{{\approx\!\! >}}
\def\inn{\in\!\!\!\!\ra}
\def\dsum{\stackrel{\cdot}{\sqcup}}
\def\dsum{\sqcup\!\!\!\!\cdot\;}

\def\tl{\triangleleft}
\def\tr{\triangleright}

\def\lhb{\lhd \hspace {-1mm}\bullet}

\newcommand{\pa}{\parallel}
\newcommand{\lra}{\longrightarrow}
\newcommand{\hra}{\hookrightarrow}
\newcommand{\hla}{\hookleftarrow}
\newcommand{\ra}{\rightarrow}
\newcommand{\lla}{\longleftarrow}
\newcommand{\da}{\downarrow}
\newcommand{\ua}{\uparrow}
\newcommand{\dA}{\downarrow\!\!\!^\bullet}
\newcommand{\uA}{\uparrow\!\!\!_\bullet}
\newcommand{\Da}{\Downarrow}
\newcommand{\DA}{\Downarrow\!\!\!^\bullet}
\newcommand{\UA}{\Uparrow\!\!\!_\bullet}
\newcommand{\Ua}{\Uparrow}
\newcommand{\Lra}{\Longrightarrow}
\newcommand{\Ra}{\Rightarrow}
\newcommand{\Lla}{\Longleftarrow}
\newcommand{\La}{\Leftarrow}
\newcommand{\nperp}{\perp\!\!\!\!\!\setminus\;\;}
\newcommand{\pq}{\preceq}

\newcommand{\lms}{\longmapsto}
\newcommand{\ms}{\mapsto}
\newcommand{\subseteqnot}{\subseteq\hskip-4 mm_\not\hskip3 mm}

\newcommand{\bth}{\begin{theorem}}
\newcommand{\eth}{\end{theorem}}

\def\phi{\varphi}
\def\ve{\varepsilon}
\def\o{{\omega}}

\def\bA{{\bf A}}
\def\bM{{\bf M}}
\def\bC{{\bf C}}
\def\bI{{\bf I}}
\def\bL{{\bf L}}
\def\bT{{\bf T}}
\def\bS{{\bf S}}
\def\bD{{\bf D}}
\def\bB{{\bf B}}
\def\bW{{\bf W}}
\def\bP{{\bf P}}
\def\bX{{\bf X}}
\def\bY{{\bf Y}}
\def\ba{{\bf a}}
\def\bb{{\bf b}}
\def\bc{{\bf c}}
\def\bd{{\bf d}}
\def\bh{{\bf h}}
\def\bi{{\bf i}}
\def\bj{{\bf j}}
\def\bk{{\bf k}}
\def\bm{{\bf m}}
\def\bn{{\bf n}}
\def\bp{{\bf p}}
\def\bq{{\bf q}}
\def\be{{\bf e}}
\def\br{{\bf r}}
\def\bi{{\bf i}}
\def\bs{{\bf s}}
\def\bt{{\bf t}}
\def\b1{{\bf 1}}

\def\bio{{\mbox{\boldmath $\iota$}}}
\def\bkappa{{\mbox{\boldmath $\kappa$}}}
\def\bpi{{\mbox{\boldmath $\pi$}}}
\def\brho{{\mbox{\boldmath $\rho$}}}
\def\bvrho{{\mbox{\boldmath $\varrho$}}}
\def\bbeta{{\mbox{\boldmath $\beta$}}}
\def\btheta{{\mbox{\boldmath $\theta$}}}
\def\bmu{{\mbox{\boldmath $\mu$}}}
\def\bcomp{{\mbox{\boldmath $;$}}}

\def\cBL{{\cal BL}}
\def\cB{{\cal B}}
\def\cA{{\cal A}}
\def\cC{{\cal C}}
\def\cD{{\cal D}}
\def\cE{{\cal E}}
\def\cF{{\cal F}}
\def\cG{{\cal G}}
\def\cI{{\cal I}}
\def\cJ{{\cal J}}
\def\cK{{\cal K}}
\def\cL{{\cal L}}
\def\cN{{\cal N}}
\def\cM{{\cal M}}
\def\cP{{\cal P}}
\def\cQ{{\cal Q}}
\def\cR{{\cal R}}
\def\cS{{\cal S}}
\def\cT{{\cal T}}
\def\cU{{\cal U}}
\def\cV{{\cal V}}
\def\cX{{\cal X}}
\def\cY{{\cal Y}}

\def\cat{{\bf Cat}}

\def\oC{{{\omega}Cat}}
\def\kC{{kCat}}
\def\nC{{n{\bf Cat}}}
\def\njC{{(n+1){\bf Cat}}}
\def\oG{{{\omega}Gr}}
\def\mts{{MltSet}}

\def\hg{{\bf Hg}}

\def\ofs{{\bf oFs}}
\def\lfs{{\bf lFs}}

\def\onfs{{\bf oFs}_n}
\def\okfs{{\bf oFs}_k}
\def\nofs{{\bf noFs}}
\def\onjfs{{\bf oFs}_{n+1}}
\def\onmjfs{{\bf oFs}_{n-1}}
\def\pfs{{\bf pFs}}
\def\nfs{{\bf nFs}}
\def\cfs{{\bf cFs}}
\def\posfs{{\bf Fs}^{+/1}}

\def\ctop{{{\cal C}topes}^{m/1}}
\def\ctyp{{{\cal C}types}^{m/1}}
\def\nctyp{{n{\cal C}types}^{m/1}}
\def\cctyp{{c{\cal C}types}^{m/1}}

\def\comp{{\bf Comp}} 


\def\mnComp{{{\bf Comp}^{m/1}_n}} 
\def\mnjComp{{{\bf Comp}^{m/1}_{n+1}}} 
\def\mnmjComp{{{\bf Comp}^{m/1}_{n-1}}} 
\def\mComp{{\bf Comp}^{m/1}} 

\def\nComma{{\bf Comma}_n}        
\def\njComma{{\bf Comma}_{n+1}}        

\def\mnComma{{\bf Comma}^{m/1}_n} 
\def\mkComma{{\bf Comma}^{m/1}_k} 
\def\mnjComma{{\bf Comma}^{m/1}_{n+1}} 
\def\mnmjComma{{\bf Comma}^{m/1}_{n-1}} 

\pagenumbering{arabic} \setcounter{page}{1}

\title{On ordered face structures and many-to-one computads}
\author{Marek Zawadowski
\\
Instytut Matematyki, Uniwersytet Warszawski\\
ul. S.Banacha 2, 00-913 Warszawa, Poland\\
zawado@mimuw.edu.pl\\}
\maketitle

\begin{abstract}  We introduce the notion of an ordered face
structure.  The ordered face structures to many-to-one computads
are like positive face structures, c.f. \cite{Z}, to
positive-to-one computads. This allow us to give an explicit combinatorial
description of many-to-one computads in terms of ordered face structures.
\end{abstract}

\tableofcontents

\section{Introduction}

The definition of multitopic categories the weak $\o$-categories
in the sense of Makkai contains two ingredients. The first
constitutes a description of shapes of cells that are considered
(this includes the relation between cells and their domains and
codomains), c.f. \cite{HMP} and the second constitutes a mechanism
of composition, c.f. \cite{M}. This paper is a contribution to a
better understanding of the first ingredient of the M.Makkai's
definition of multitopic categories, and we provide a relatively
simple combinatorial description of
the category many-to-one computads. The paper goes much along with
\cite{Z} except it deals with all many-to-one computads rather
than positive-to-one computads. This generates some substantial
complications and the structure of cells turns out to be much
richer.

\subsection*{Ordered face structures}

Our main combinatorial device introduced and studied in this paper
is the {\em ordered face structure}. The ordered face structures
correspond to all possible 'shapes' of cells (not only
indeterminates) in many-to-one computads\footnote{For the
definition of many-to-one computad see the appendix.}. In order to
relate them to our previous work \cite{MZ}, \cite{Z} we can draw
an analogy in the following table.
 \[ \begin{array}{|c|c|c|c|c|c|} \hline
              & \multicolumn{4}{c|}{shapes\; of}   \\ \cline{2-5}
     type\; of   &  \multicolumn{2}{c|}{indeterminates} &  \multicolumn{2}{c|}{arbitrary\; cells} \\ \cline{2-5}
     computads   & \multicolumn{4}{c|}{described\; in\; terms\; of }   \\ \cline{2-5}
                 &  graph\!-\! like  & computads &  graph\!-\!like  & computads \\
                 &   structures  &  &   structures  & \\ \hline

    one\!-\! to\!-\! one  & \alpha^n & (\alpha^n)^*& simple    & simple  \\
     {\rm \cite{MZ}}      &&  & \o\!-\!graphs  & categories  \\ \hline

    positive\!-\! to\!-\! one & principal & positive & positive\; face & positive  \\
     {\rm \cite{Z}}       & positive\; face & computopes & structures  & computypes  \\
               & structures       &   &  & \\ \hline

   many\!-\! to\!-\! one & principal & computopes & ordered\; face & pointed  \\
   {\rm [this\; paper]}& ordered\; face & & structures  & computypes \\
               & structures       &   &  & \\ \hline
  \end{array} \]
Now are going to explain it in an intuitive way. In the table we
describe cells in computads of three kinds. The later being
strictly more general than the former. The one-to-one computads
are the simplest.  They are free $\o$-categories over
$\o$-graphs\footnote{In the literature $\o$-graphs are sometimes
called globular sets.}. The positive-to-one computads are
computads in which the indeterminates (or indets) on the higher
dimension have as codomains indeterminates and as domains cells
that are not identities. Finally, the many-to-one computads are
computads in which the indets have as codomains indets again but
there is no specific restriction for the domains (other than that
they must be parallel to codomains).

Fix $n\in\o$. The $\o$-graph (also called globular set)
$\alpha^n$, has one $n$-face and exactly two faces of lower
dimensions than $n$, i.e.
\[ \alpha^n_l \;\;= \;\; \left\{ \begin{array}{ll}
                \emptyset    & \mbox{ if  $l> n$} \\
                \{ 2n \}     & \mbox{ if  $l= n$} \\
                \{ 2l+1,\, 2l  \}     & \mbox{ if  $0\leq l< n$}
                                    \end{array}
                            \right.\label{n.alpha0} \]
with domain and codomain given by $d,\, c: \alpha^n_l \lra
\alpha^n_{l-1}$, $d(x)=\{ 2l-1 \}$, $c(x)=2l-2$ for
$x\in\alpha^n_l$, and $1\leq l\leq n$. For example
$\alpha^4$\label{n.alpha1} can be pictured as follows:
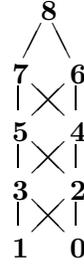
\begin{center}
\begin{picture}(1700,900)
\put(700, 0){${\bf 1}$} \put(900, 0){${\bf 0}$}
\put(720,110){\line(0,1){80}} \put(770,90){\line(1,1){110}}
\put(770,200){\line(1,-1){110}} \put(920,110){\line(0,1){80}}
\put(700,200){${\bf 3}$} \put(900, 200){${\bf 2}$}

\put(720,310){\line(0,1){80}} \put(770,290){\line(1,1){110}}
\put(770,400){\line(1,-1){110}} \put(920,310){\line(0,1){80}}

\put(700,400){${\bf 5}$} \put(900, 400){${\bf 4}$}

\put(720,510){\line(0,1){80}} \put(770,490){\line(1,1){110}}
\put(770,600){\line(1,-1){110}} \put(920,510){\line(0,1){80}}
\put(700,600){${\bf 7}$} \put(900,600){${\bf 6}$}

\put(740,690){\line(1,2){60}}
 \put(840,810){\line(1,-2){60}}
 \put(800,820){${\bf 8}$}
\end{picture}
\end{center}
i.e. $8$ is the unique face of dimension $4$ in $\alpha^4$ that
has $7$ as its domain and $6$ as its codomain, $7$ and $6$ have
$5$ as its domain and $4$ as its codomain, and so on. More
visually we can draw $\alpha^4$ as follows
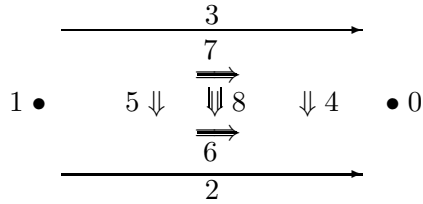
\begin{center} \xext=1300 \yext=680
\begin{picture}(\xext,\yext)(\xoff,\yoff)
 \put(-100,320){$1\;\bullet$}
  \put(1200,320){$\bullet\; 0$}

 \putmorphism(50,600)(1,0)[``3]{1100}{1}a
 \putmorphism(50,100)(1,0)[``2]{1100}{1}b

 \put(300,320){$5\Da$}
  \put(900,320){$\Da 4$}

 \put(540,420){$\Longrightarrow$}
  \put(540,220){$\Longrightarrow$}
  \put(550,445){\line(1,0){135}}
  \put(550,245){\line(1,0){135}}
  \put(570,145){6}
   \put(570,500){7}
 \put(580,320){$\Da 8$}
 \put(590,330){\line(0,1){70}}
 \put(630,330){\line(0,1){70}}
\end{picture}
\end{center}
The free category $(\alpha^n)^*$ generated by $\alpha^n$ has the
property that for any $\o$-category $C$, the set
$\oC((\alpha^n)^*,C)$ of $\o$-functors from $(\alpha^n)^*$ to $C$
correspond naturally to the set $C_n$ of $n$-cells of $C$. Thus in
one-to-one computads the shapes of indets are particularly simple
and this is why the $\o$-graphs describing them are called {\em
simple}. Simple $\o$-graphs are some 'special' pushouts of
$\alpha$'s. Instead of trying to repeat the definition from
\cite{MZ} we rather show an example:
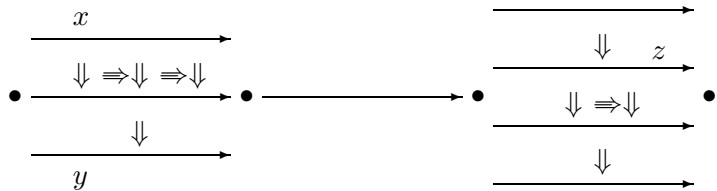
\begin{center}
\xext=2400 \yext=600
\begin{picture}(\xext,\yext)(\xoff,\yoff)
\putmorphism(0,300)(1,0)[\bullet`\bullet`]{800}{1}a
\putmorphism(0,500)(1,0)[\phantom{\bullet}`\phantom{\bullet}`]{800}{1}a
\putmorphism(0,100)(1,0)[\phantom{\bullet}`\phantom{\bullet}`]{800}{1}a
\put(200,350){\makebox{$\Downarrow$}}
\put(300,350){\makebox{$\Rightarrow\!\!\!\!\!\!\!\! -$ }}
\put(400,350){\makebox{$\Downarrow$}}
\put(500,350){\makebox{$\Rightarrow\!\!\!\!\!\!\!\! -$ }}
\put(600,350){\makebox{$\Downarrow$}}
\put(400,150){\makebox{$\Downarrow$}}
\putmorphism(800,300)(1,0)[\phantom{\bullet}`\bullet`]{800}{1}a
\putmorphism(1600,300)(1,0)[\phantom{\bullet}`\bullet`]{800}{0}a
\putmorphism(1600,600)(1,0)[\phantom{\bullet}`\phantom{\bullet}`]{800}{1}a
\putmorphism(1600,400)(1,0)[\phantom{\bullet}`\phantom{\bullet}`]{800}{1}a
\putmorphism(1600,200)(1,0)[\phantom{\bullet}`\phantom{\bullet}`]{800}{1}a
\putmorphism(1600,0)(1,0)[\phantom{\bullet}`\phantom{\bullet}`]{800}{1}a

\put(2000,450){\makebox{$\Downarrow$}}
\put(1900,250){\makebox{$\Downarrow$}}
\put(2000,250){\makebox{$\Rightarrow\!\!\!\!\!\!\!\! -$ }}
\put(2100,250){\makebox{$\Downarrow$}}
\put(2000,50){\makebox{$\Downarrow$}}

\put(200,550){\makebox{$x$}} \put(200,0){\makebox{$y$}}
\put(2200,430){\makebox{$z$}}
\end{picture}
\end{center}
 So indets have still indets as domains and codomains and
even if there is no one indet that generates all the $\o$-graph,
as in $\alpha^n$'s, the domains and codomains of indets so fit
together that they could be (uniquely) composed 'if they were
placed in an $\o$-category'.  Simple $\o$-categories, c.f.
\cite{MZ}, are $\o$-categories generated by such $\o$-graphs. The
category of simple $\o$-categories is dual to the category of
disks introduced in \cite{Joyal} . Note that there are two
definite ways the indets of the same dimension can be compared.
The face $x$ is smaller from $y$ in one way and from $z$ the other
way. We write $x<^+y$ and $x<^-z$.  The first order\footnote{Here
and later by order we mean {\em strict order} i.e. irreflexive and
transitive relation.} is called {\em upper} and the second is
called {\em lower}. More formally, the upper order on cells of
dimension $n$ is the least transitive relation such that
$d(a)<^+c(a)$ for any face $a$ of dimension $n+1$ ($d$ and $c$ are
operations of domain and codomain, respectively). Similarly, the
lower order on cells of dimension $n$ is the least transitive
relation such that if $d(x)=c(z)$ then $x<^-z$. In this case both
orders are definable using $d$ and $c$. For more on this see
\cite{MZ}.

The shapes of indeterminates in positive-to-one face structures
are more complicated.  We again use drawing to explain what
principal positive face structures are. The one below has
dimension $3$.
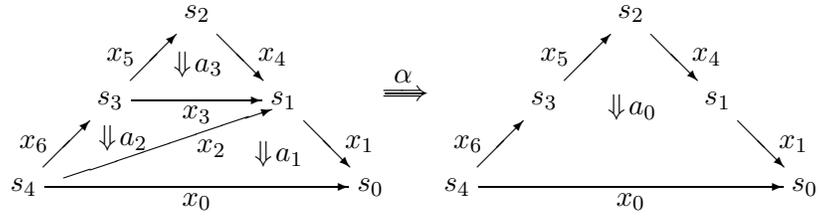
\begin{center} \xext=2400 \yext=680
\begin{picture}(\xext,\yext)(\xoff,\yoff)
 \label{figure ppof}
 \settriparms[-1`1`1;300]
 \putAtriangle(300,350)[s_2`s_3`s_1;x_5`x_4`x_3]
 \put(520,450){$\Da\! a_3$}
 \settriparms[-1`0`0;300]
 \putAtriangle(0,50)[\phantom{s_3}`s_4`;x_6``]
 \settriparms[0`1`0;300]
 \putAtriangle(600,50)[\phantom{s_1}``s_0;`x_1`]
 \putmorphism(0,50)(1,0)[\phantom{s_4}`\phantom{s_0}`x_0]{1200}{1}b
 \putmorphism(350,150)(3,1)[\phantom{s_4}`\phantom{s_0}`]{300}{1}a
 \put(600,170){$x_2$}
 \put(260,195){$\Da\! a_2$}
 \put(800,140){$\Da\! a_1$}

\put(1240,350){$\Longrightarrow$}
  \put(1250,375){\line(1,0){135}}
   \put(1280,410){$\alpha$}

\settriparms[-1`1`0;300]
 \putAtriangle(1800,350)[s_2`s_3`s_1;x_5`x_4`]
 \put(2020,300){$\Da\! a_0$}
 \settriparms[-1`0`0;300]
 \putAtriangle(1500,50)[\phantom{s_3}`s_4`;x_6``]
 \settriparms[0`1`0;300]
 \putAtriangle(2100,50)[\phantom{s_1}``s_0;`x_1`]
 \putmorphism(1500,50)(1,0)[\phantom{s_4}`\phantom{s_0}`x_0]{1200}{1}b
\end{picture}
\end{center}
Thus in positive face structures the codomains of indets are still
indets but the domains are so called {\em pasting diagrams} of
indets, i.e. domains contains indets that 'suitably fit together
so that we could compose them'. In these structures we have the
usual operation of taking codomain but the 'operation' of taking
domain of a face returns a non-empty set of faces rather than a
single face. To emphasize this change we use for these operations
the Greek letters $\gamma$ and $\delta$ instead of $c$ and $d$.
Thus $\gamma(\alpha)=a_0$, $\gamma(a_3)=x_3$, $\delta(\alpha)=\{
a_1,a_2,a_3\}$, $\delta(a_0)=\{ x_1,x_4,x_5,x_6\}$,
$\delta(a_2)=\{ x_3,x_6\}$. From the table we have that positive
face structures to principal positive face structures are like
simple $\o$-graphs to $\o$-graphs of form $\alpha^n$, for some
$n$. Thus it should be not surprising that positive face
structures looks like this:
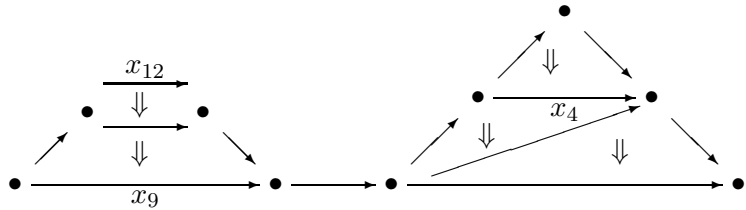
\begin{center} \xext=2400 \yext=640
\begin{picture}(\xext,\yext)(\xoff,\yoff)
\settriparms[-1`0`0;250]
 \putAtriangle(0,0)[\bullet`\bullet`;``]
 \settriparms[0`1`0;250]
 \putAtriangle(400,0)[\bullet``\bullet;``]
 \putmorphism(0,0)(1,0)[\phantom{\bullet}`\phantom{\bullet}`x_9]{900}{1}b
  \put(400,80){$\Da$}
 \putmorphism(250,200)(1,0)[\phantom{\bullet}`\phantom{\bullet}`]{400}{1}b
 \putmorphism(250,350)(1,0)[\phantom{\bullet}`\phantom{\bullet}`x_{12}]{400}{1}a
 \put(400,250){$\Da$}

\putmorphism(900,0)(1,0)[\phantom{\bullet}`\phantom{\bullet}`]{400}{1}b

 \settriparms[-1`1`1;300]
 \putAtriangle(1600,300)[\bullet`\bullet`\bullet;``x_4]
 \put(1820,400){$\Da$}
 \settriparms[-1`0`0;300]
 \putAtriangle(1300,0)[\phantom{\bullet}`\bullet`;``]
 \settriparms[0`1`0;300]
 \putAtriangle(1900,0)[\phantom{\bullet}``\bullet;``]
 \putmorphism(1300,0)(1,0)[\phantom{\bullet}`\phantom{\bullet}`]{1200}{1}b
 \putmorphism(1650,100)(3,1)[\phantom{\bullet}`\phantom{\bullet}`]{300}{1}a
 \put(1600,145){$\Da$}
 \put(2060,90){$\Da$}
\end{picture}
\end{center}
Different points, arrows etc. denote necessarily different cells,
 and if we omit their names in  figures it is for making it
less baroque. Note that in this case the indets of the same
dimension can be compared much the same way as indets in simple
$\o$-graphs in two definite ways. The face $x_{12}$ is smaller
than $x_9$ in one way and than $x_4$ the other way, and again we
write $x_{12}<^+x_9$ and $x_{12}<^-x_4$. Again the first order is
called {\em upper} and the second is called {\em lower}. More
formally, the upper order on faces of dimension $n$ is the least
transitive relation such that $x<^+y$ whenever there is a face $a$
of dimension $n+1$ such that $x\in\delta(a)$ and $\gamma(a)=y$.
Similarly, the lower order on faces of dimension $n$ is the least
transitive relation such that $x<^-y$ whenever
$\gamma(x)\in\delta(y)$. Any positive face structure $T$ generates
a computads $T^*$. The cells of dimension $n$ of such a computad
are positive face substructures of $T$ of dimension at most $n$.
These computads are called {\em positive computypes}. If $T$ is a
principal positive face structure then $T^*$ is a {\em positive
computope}\footnote{The word 'positive' is used here more like a
shorthand and in presence of 'other positive' notions this one
should be named properly as 'positive-to-one'.}. In this case
$T^*$ determines $T$ up to an isomorphism. For more on this see
\cite{Z}.

The shapes of indeterminates in many-to-one face structures are
even more complicated as this time the domains of indets might be
identities (='empty on something'). This generates a lot of
complications as we have three new kinds of faces. Apart from {\em
positive faces} like in previous case we have {\em empty-domain
faces} and then as a consequence we have {\em loops} (=faces with
domain equal codomain) and we also need to deal with {\em empty
faces}. The last kind of faces is not indicated in the pictures.
On each face $x$ of dimension $n$ there is an empty face $1_x$ of
dimension $n+1$.  They are much like with identities whose role
they play. We again use drawing to explain intuitively what
principal ordered face structures are:
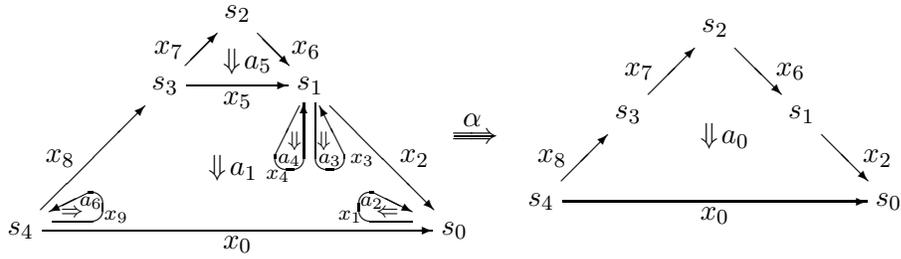
\begin{center} \xext=3000 \yext=800
\begin{picture}(\xext,\yext)(\xoff,\yoff)
 \label{figure pofs}
 \settriparms[-1`1`1;250]
 \putAtriangle(500,500)[s_2`s_3`s_1;x_7`x_6`x_5]
 \put(700,560){$\Da\! a_5$}
 \settriparms[-1`0`0;500]
 \putAtriangle(0,0)[\phantom{s_2}`s_4`;x_8``]
 \settriparms[0`1`0;500]
 \putAtriangle(500,0)[\phantom{s_1}``s_0;`x_2`]
 \putmorphism(0,0)(1,0)[\phantom{s_4}`\phantom{s_0}`x_0]{1500}{1}b
 \put(650,190){$\Da\! a_1$}

 \settriparms[-1`1`0;300]
 \putAtriangle(2100,400)[s_2`s_3`s_1;x_7`x_6`]
 \settriparms[-1`0`0;300]
 \putAtriangle(1800,100)[\phantom{s_3}`s_4`;x_8``]
 \settriparms[0`1`0;300]
 \putAtriangle(2400,100)[\phantom{s_4}``s_0;`x_2`]
 \putmorphism(1800,100)(1,0)[\phantom{s_4}`\phantom{s_0}`x_0]{1200}{1}b
 \put(2350,300){$\Da\! a_0$}

\put(1490,300){$\Longrightarrow$}
  \put(1500,325){\line(1,0){135}}
   \put(1530,360){$\alpha$}

  \put(930,260){\oval(100,100)[b]}
  \put(880,260){\line(1,2){85}}
  \put(980,260){\vector(0,1){180}}
  \put(925,255){$^\Da$}
    \put(890,195){$^{a_4}$}
     \put(850,140){$^{x_4}$}

  \put(1070,260){\oval(100,100)[b]}
   \put(1020,260){\line(0,1){180}}
  \put(1120,260){\vector(-1,2){85}}
   \put(1025,255){$^\Da$}
   \put(1030,195){$^{a_3}$}
   \put(1140,195){$^{x_3}$}

\put(140,10){$^\Rightarrow$}
 \put(235,80){\oval(100,100)[r]}
  \put(110,30){\line(1,0){135}}
  \put(235,130){\vector(-2,-1){135}}
    \put(205,50){$^{a_6}$}
    \put(290,0){$^{x_9}$}
\put(1230,10){$^\Leftarrow$}
 \put(1220,80){\oval(100,100)[l]}
  \put(1220,30){\line(1,0){135}}
  \put(1220,130){\vector(2,-1){135}}
    \put(1175,50){$^{a_2}$}
    \put(1100,0){$^{x_1}$}
\end{picture}
\end{center}
and a bit more fancy
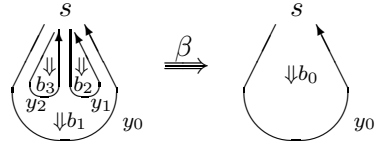
\begin{center} \xext=1200 \yext=480
\begin{picture}(\xext,\yext)(\xoff,\yoff)
\label{figure pofs2}
 \put(200,200){\oval(360,360)[b]}
 \put(160,30){$^{\Da b_1}$}
 \put(20,200){\line(1,2){110}}
 \put(380,200){\vector(-1,2){110}}
  \put(180,450){$s$}
   \put(400,30){$^{y_0}$}

 \put(1000,200){\oval(360,360)[b]}
 \put(960,190){$^{\Da b_0}$}
 \put(820,200){\line(1,2){110}}
 \put(1180,200){\vector(-1,2){110}}
  \put(980,450){$s$}
   \put(1200,30){$^{y_0}$}

\put(540,250){$\Longrightarrow$}
  \put(550,275){\line(1,0){135}}
   \put(580,310){$\beta$}

  \put(130,220){\oval(100,100)[b]}
  \put(80,220){\line(1,2){85}}
  \put(180,220){\vector(0,1){180}}
  \put(125,220){$^\Da$}
    \put(100,155){$^{b_3}$}
     \put(65,95){$^{y_2}$}

  \put(270,220){\oval(100,100)[b]}
   \put(220,220){\line(0,1){180}}
  \put(320,220){\vector(-1,2){85}}
   \put(235,220){$^\Da$}
   \put(230,155){$^{b_2}$}
   \put(290,95){$^{y_1}$}
\end{picture}
\end{center}
In these structures we also use the Greek letters $\gamma$ and
$\delta$ for domains and codomains, respectively.  Similarly as in
positive face structures the codomain is an operation associating
faces to faces. But the domain operation is still more involved as
it may associate to a face a non-empty set of faces or a single
empty face. Thus we have $\delta(\alpha)=\{ a_1,\ldots,a_6 \}$,
$\gamma(\alpha)=a_0$ but $\delta(a_2)=1_{s_0}$, $\delta(b_0)=1_s$.
Note that we should write $\delta(x_1)=\{ s_0 \}$ instead of
$\delta(x_0)=s_0$ but we will, as we did in \cite{Z}, mix
singletons with elements when dealing with faces or sets of faces
e.g. both conditions $\gamma(x_0)\in\delta(x_0)$ and
$\gamma(x_0)=\delta(x_0)$ are meaningful in this convention and in
fact, as we will see later, due to this 'double meaning' they are
equivalent in all ordered face structures saying that $x_0$ is a
loop. This time the relations between faces and their domains and
codomains does not encode all the needed data. The upper order
$<^+$ can be defined like in positive face structures from
$\gamma$ and $\delta$. However, due to existence of loops, the
relation $<^-$ defined as before is not a strict order in general.
In the above examples we have $x_3<^-x_4$, $x_4<^-x_3$ and
similarly $y_2<^-y_1$, $y_1<^-y_2$. But we definitely need to know
that $x_4$ comes before $x_3$ and that $y_2$ comes before $y_1$.
This is why we need as a separate additional data a strict order
$<^\sim$ that is contained in $<^-$ telling us that $x_4<^\sim
x_3$ and $y_2<^\sim y_1$ but not that $x_3<^\sim x_3$ and
$y_1<^\sim y_2$. As we need to have the strict order $<^\sim$ as
an additional piece of data we call those face structures {\em
ordered}. Note however that in the above cases we could solve our
problem of ordering the faces locally that is having just
restriction of the order $<^\sim$ to sets that are domains of
other faces. But to describe all the cells of many-to-one
computads we need more than just that. Below we have some examples
of ordered face structure
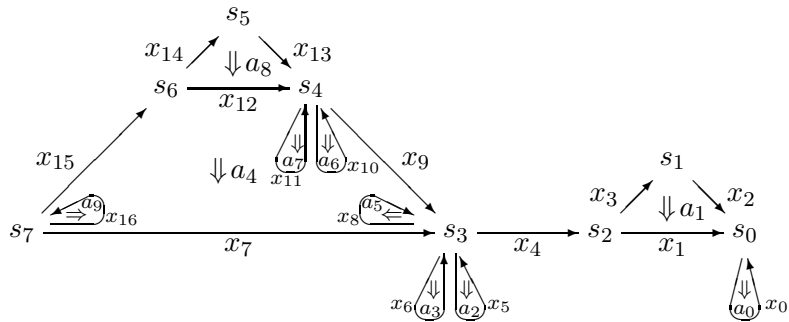
\begin{center} \xext=2320 \yext=1050
\begin{picture}(\xext,\yext)(\xoff,\yoff)
 \settriparms[-1`1`1;250]
 \putAtriangle(500,800)[s_5`s_6`s_4;x_{14}`x_{13}`x_{12}]
 \put(700,860){$\Da\! a_8$}
 \settriparms[-1`0`0;500]
 \putAtriangle(0,300)[\phantom{s_2}`s_7`;x_{15}``]
 \settriparms[0`1`0;500]
 \putAtriangle(500,300)[\phantom{s_1}``s_3;`x_9`]
 \putmorphism(0,300)(1,0)[\phantom{s_4}`\phantom{s_3}`x_7]{1500}{1}b
 \put(650,490){$\Da\! a_4$}

  \put(930,560){\oval(100,100)[b]}
  \put(880,560){\line(1,2){85}}
  \put(980,560){\vector(0,1){184}}
  \put(930,555){$^\Da$}
    \put(900,495){$^{a_7}$}
     \put(860,435){$^{x_{11}}$}

  \put(1070,560){\oval(100,100)[b]}
   \put(1020,560){\line(0,1){180}}
  \put(1120,560){\vector(-1,2){85}}
   \put(1035,555){$^\Da$}
   \put(1025,495){$^{a_6}$}
   \put(1125,480){$^{x_{10}}$}

\put(150,310){$^\Rightarrow$}
 \put(235,380){\oval(100,100)[r]}
  \put(100,330){\line(1,0){135}}
  \put(235,430){\vector(-2,-1){135}}
    \put(205,350){$^{a_9}$}
    \put(290,300){$^{x_{16}}$}
\put(1240,310){$^\Leftarrow$}
 \put(1220,380){\oval(100,100)[l]}
  \put(1220,330){\line(1,0){135}}
  \put(1220,430){\vector(2,-1){135}}
    \put(1175,350){$^{a_5}$}
    \put(1090,300){$^{x_8}$}

  \put(1410,50){\oval(100,100)[b]}
  \put(1360,50){\line(1,2){85}}
  \put(1460,50){\vector(0,1){180}}
  \put(1400,50){$^\Da$}
    \put(1375,-15){$^{a_3}$}
     \put(1275,0){$^{x_6}$}

  \put(1550,50){\oval(100,100)[b]}
   \put(1500,50){\line(0,1){180}}
  \put(1600,50){\vector(-1,2){85}}
   \put(1515,50){$^\Da$}
   \put(1510,-15){$^{a_2}$}
   \put(1610,0){$^{x_5}$}

     \put(2500,50){\oval(100,100)[b]}
   \put(2450,50){\line(1,4){40}}
  \put(2550,50){\vector(-1,4){40}}
   \put(2475,50){$^\Da$}
   \put(2460,-15){$^{a_0}$}
   \put(2570,0){$^{x_0}$}

\putmorphism(1500,300)(1,0)[\phantom{s_0}`\phantom{s_0}`x_4]{500}{1}b

 \settriparms[-1`1`1;250]
 \putAtriangle(2000,300)[s_1`s_2`s_0;x_3`x_2`x_1]
 \put(2200,360){$\Da\! a_1$}
\end{picture}
\end{center}
and
\begin{center} \xext=300 \yext=1050
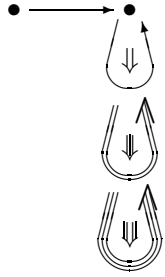
\begin{picture}(\xext,\yext)(\xoff,\yoff)
  \put(110,110){\oval(220,220)[b]}
  \put(110,110){\oval(190,190)[b]}
  \put(110,110){\oval(160,160)[b]}
  \put(0,110){\line(1,4){40}}
  \put(15,110){\line(1,4){40}}
  \put(30,110){\line(1,4){40}}
  \put(220,110){\line(-1,4){40}}
  \put(205,110){\line(-1,4){40}}
  \put(190,110){\line(-1,4){35}}
  \put(135,230){$\bigwedge$}
  \put(80,100){$\Da$}
  \put(90,110){\line(0,1){60}}
  \put(130,110){\line(0,1){60}}

  \put(110,410){\oval(190,190)[b]}
  \put(110,410){\oval(160,160)[b]}
  \put(15,410){\line(1,4){40}}
  \put(30,410){\line(1,4){40}}
  \put(205,410){\line(-1,4){44}}
  \put(190,410){\line(-1,4){38}}
  \put(125,530){$\bigwedge$}
  \put(80,400){$\Da$}
  \put(110,400){\line(0,1){70}}

  \put(110,710){\oval(160,160)[b]}
  \put(30,710){\line(1,4){40}}
  \put(190,710){\vector(-1,4){38}}
  \put(80,700){$\Da$}

\putmorphism(-290,900)(1,0)[\bullet`\bullet`]{400}{1}b
\end{picture}
\end{center}
and
\begin{center} \xext=320 \yext=390
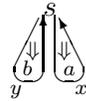
\begin{picture}(\xext,\yext)(\xoff,\yoff)
\put(180,310){$s$}
  \put(130,130){\oval(100,100)[b]}
  \put(80,130){\line(1,2){85}}
  \put(180,130){\vector(0,1){180}}
  \put(120,130){$^\Da$}
    \put(110,65){$^{b}$}
     \put(65,0){$^{y}$}

  \put(270,130){\oval(100,100)[b]}
   \put(220,130){\line(0,1){180}}
  \put(320,130){\vector(-1,2){85}}
   \put(230,130){$^\Da$}
   \put(250,65){$^{a}$}
   \put(290,0){$^{x}$}
\end{picture}
\end{center}
We see that faces $x_6$ and $x_5$ must be comparable via $<^\sim$
but they are not in domain of any other face. Thus a kind of
global order $<^\sim$ is needed. Note however that the fact that
$x_{11}$ comes before both $x_6$ and $x_4$ and that $x_0$ comes
after all of them could be deduced in a different way. The way the
ordered face structure $T$ generate a many-to-one computad $T^*$
is more involved then in case of positive face structures. An
$n$-cell in $T^*_n$ is a {\em local morphism} $\varphi:X\ra T$ where
$X$ is an ordered face structure of dimension at most $n$ and
$\varphi$ is a map that preserves $\gamma$, $\delta$ but the order
$<^\sim$ is preserved only locally i.e. for $a\in X$ $\varphi: (
\delta(a),<^\sim_a )\ra(\delta(\varphi(a)),<^\sim_{\varphi(a)})$
is an order isomorphism, where $<^\sim_a$,  $<^\sim_{\varphi(a)}$
are restrictions of orders $<^\sim$ to $\delta(a)$ and
$\delta(\varphi(a))$, respectively. Thus we have a cell $\varphi$:
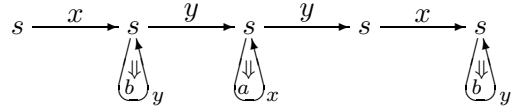
\begin{center} \xext=1650 \yext=320
\begin{picture}(\xext,\yext)(\xoff,\yoff)
\putmorphism(0,300)(1,0)[s`s`x]{400}{1}a
\putmorphism(400,300)(1,0)[\phantom{s}`s`y]{400}{1}a
\putmorphism(800,300)(1,0)[\phantom{s}`s`y]{400}{1}a
\putmorphism(1200,300)(1,0)[\phantom{s}`s`x]{400}{1}a

  \put(400,100){\oval(100,100)[b]}
  \put(350,100){\line(1,4){40}}
  \put(450,100){\vector(-1,4){40}}
  \put(380,100){$^\Da$}
    \put(370,35){$^{b}$}
     \put(465,10){$^{y}$}

  \put(800,100){\oval(100,100)[b]}
   \put(750,100){\line(1,4){40}}
  \put(850,100){\vector(-1,4){40}}
   \put(775,100){$^\Da$}
   \put(760,35){$^{a}$}
   \put(860,10){$^{x}$}

     \put(1600,100){\oval(100,100)[b]}
  \put(1550,100){\line(1,4){40}}
  \put(1650,100){\vector(-1,4){40}}
  \put(1580,100){$^\Da$}
    \put(1570,35){$^{b}$}
     \put(1665,10){$^{y}$}
\end{picture}
\end{center}
in the computad generated by the ordered face structure $T^*$
(where $T$ is the last example of an ordered face structure
above). The faces of the above ordered face structures are
labelled by the faces they are sent to by the local morphism $\varphi$.
Clearly, in this case the local preservation of the order $<^\sim$
does not impose any restriction on the map $\varphi:X\ra T$ other
than preservation of $\gamma$ and $\delta$. From this it should be
clear that we cannot in general determine $T$ having just $T^*$.
For example the ordered face structures
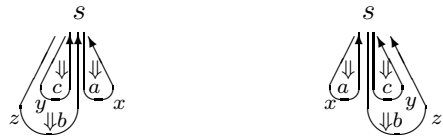
\begin{center} \xext=1320 \yext=520
\begin{picture}(\xext,\yext)(\xoff,\yoff)
 \put(100,150){\oval(200,200)[b]}
 \put(80,40){$^{\Da b}$}
 \put(0,150){\line(1,2){120}}
 \put(200,150){\vector(0,1){250}}
  \put(180,450){$s$}
   \put(-40,50){$^{z}$}

  \put(120,220){\oval(100,100)[b]}
  \put(70,220){\line(1,2){85}}
  \put(170,220){\vector(0,1){180}}
  \put(120,220){$^\Da$}
    \put(110,155){$^{c}$}
     \put(50,100){$^{y}$}
    \put(270,220){\oval(100,100)[b]}
  \put(220,220){\line(0,1){180}}
  \put(320,220){\vector(-1,2){85}}
  \put(235,220){$^\Da$}
    \put(235,155){$^{a}$}
     \put(320,100){$^{x}$}

 \put(1300,150){\oval(200,200)[b]}
 \put(1240,40){$^{\Da b}$}
 \put(1200,150){\line(0,1){250}}
 \put(1400,150){\vector(-1,2){120}}
  \put(1180,450){$s$}
  \put(1420,50){$^{z}$}

  \put(1270,220){\oval(100,100)[b]}
  \put(1220,220){\line(0,1){180}}
  \put(1320,220){\vector(-1,2){85}}
  \put(1240,220){$^\Da$}
    \put(1250,155){$^{c}$}
     \put(1330,110){$^{y}$}
    \put(1120,220){\oval(100,100)[b]}
  \put(1070,220){\line(1,2){85}}
  \put(1170,220){\vector(0,1){180}}
  \put(1110,220){$^\Da$}
    \put(1095,155){$^{a}$}
     \put(1050,100){$^{x}$}
\end{picture}
\end{center}
are not isomorphic, as $y<^\sim x$ and $z<^\sim x$ in the left one
and $x<^\sim y$ and  $x<^\sim z$ in the right one, but they
generate isomorphic computads. In other word passing from $T$ to
$T^*$ we are loosing part of data and this is why $T^*$ is not
sufficient, in general (unlike $T$), to determine the shape of a
cell in a many-to-one computad. To keep this information we need
to choose one cell in $T^*$ with the natural choice being the
identity on $T$, $id_T:T\ra T$. The $\o$-categories $T^*$ together
with a distinguished cells $id_T$ are pointed computypes which are
the computad-like descriptions of types of all cells in
many-to-one computads. The pointed computypes can be defined
abstractly but we are going to explain it elsewhere. If $T$ is a
principal ordered face structure then this distinguished cell can
be chosen in a unique way and hance it does not to be chosen at
all as we know anyway which one we were to choose. This is why the
computopes, the computad-like descriptions of types of indets in
many-to-one computads are the $\o$-categories generated by
principal ordered face structures (without an additional cell
chosen).

\subsection*{Primitive notions and axioms}

Thus we related ordered face structures to simple $\o$-graphs and
positive face structures and we have described the primitive
notions $\gamma$, $\delta$, $<^\sim$ that we had chosen to
axiomatize them. Now we shall describe some intuitions behind the
axioms of ordered faces structures. Even if they are more involved
they are quite close in the spirit to the axioms of positive face
structures.

As in case of positive face structures, the most important axiom
is the axiom of {\em globularity}. In case of $\o$-graphs it is
just $cc=cd$ and $dc=dd$ which, if we rebaptize $c$ as $\gamma$
and $d$ as $\delta$, take form
\begin{equation}\label{glob-0} \gamma\gamma(\alpha)=\gamma\delta(\alpha),\;\;\;\;
 \delta\gamma(\alpha)=\delta\delta(\alpha).\end{equation}
As it was pointed out in \cite{Z} this equations cannot hold even
for positive-to-one faces as the right hand sides might be much
bigger the left hand sides. In the example of a principal positive
computad from page \pageref{figure ppof}, we have
\[ \gamma\gamma(\alpha)=x_0 \neq \{ x_0,x_2,x_3\} =\gamma\delta(\alpha),\]
\[\delta\gamma(\alpha)=\{x_1,x_4,x_5,x_6 \}\subseteqnot
\{x_1,x_2,x_3,x_4,x_5,x_6 \}= \delta\delta(\alpha).\]
 Thus we corrected the formula (\ref{glob-0}) by subtracting some
 faces from the right side getting
\begin{equation}\label{glob-1}
\gamma\gamma(\alpha)=\gamma\delta(\alpha)-\delta\delta(\alpha),\;\;\;\;
 \delta\gamma(\alpha)=\delta\delta(\alpha)-\gamma\delta(\alpha).\end{equation}
Now it works for positive-to-one faces but if we allow loops in
the domains of faces, and we must if we allow empty-domain faces,
these formulas still doesn't work as we can see for the face $a_1$
in positive ordered face structure on page \pageref{figure pofs}.
We have
\[ \gamma\gamma(a_1)=s_0\neq\emptyset=
\gamma\delta(a_1)-\delta\delta(a_1) ,\;\;\;\;
\delta\gamma(a_1)=s_4\neq\emptyset=
\delta\delta(a_1)-\gamma\delta(a_1) \] Thus we see that we
subtracted too much. Correcting this we drop these loops and we
get
\begin{equation}\label{glob-2}
\gamma\gamma(\alpha)=\gamma\delta(\alpha)-\delta\delta^{-\lambda}(\alpha),\;\;\;\;
 \delta\gamma(\alpha)=\delta\delta(\alpha)-\gamma\delta^{-\lambda}(\alpha).\end{equation}
where $\delta^{-\lambda}(\alpha)$ means the set of those faces in
$\delta(\alpha)$ that are not loops. Now the formula
(\ref{glob-2}) works for the face $a_1$ and even for the face
$\alpha$ on page \pageref{figure pofs2}. But there is still a
problem with empty-domain faces, as we have for $b_0$ in the same
ordered face structure.
\[ \gamma\gamma(b_0)=s\neq\emptyset=\gamma\delta(b_0)-\delta\delta^{-\lambda}(b_0),\]
\[ \delta\gamma(b_0)=s\neq\emptyset=\delta\delta(b_0)-\gamma\delta^{-\lambda}(b_0).\]
As a remedy for this we shall still diminish the set that we
subtract by dropping empty faces which might be there. So we drop
these empty-faces and we get
\begin{equation}\label{glob-3}
\gamma\gamma(\alpha)=\gamma\delta(\alpha)-\delta\dot{\delta}^{-\lambda}(\alpha),\;\;\;\;
 \delta\gamma(\alpha)=\delta\delta(\alpha)-\gamma\dot{\delta}^{-\lambda}(\alpha).\end{equation}
where $\dot{\delta}^{-\lambda}(\alpha)$ means the set of those
faces in $\delta(\alpha)^{-\lambda}$ that are not empty
faces\footnote{This means that this set is either empty if
$\delta(\alpha)$ is an empty face or it is
$\delta(\alpha)^{-\lambda}$.}. We are almost there but if in the
domain $\delta(\alpha)$ of a face $\alpha$ we have both
empty-domain faces and faces with positive domains as we have in
$\delta(\alpha)$ in on page \pageref{figure pofs}, then the set
$\delta\delta(\alpha)$ may contain both empty and non-empty faces
whereas $\delta\gamma(\alpha)$ definitely contain just one kind of
faces either one single empty faces or a non-empty set of
non-empty faces.  However if we have faces of both kinds in
$\delta\delta(\alpha)$ the empty faces must be empty-faces on
domains or codomains of the non-empty faces in this set. And this
is the final modification that we do to our equation:
\begin{equation}\label{glob-4}
\gamma\gamma(\alpha)=\gamma\delta(\alpha)-\delta\dot{\delta}^{-\lambda}(\alpha),\;\;\;\;
 \delta\gamma(\alpha)\equiv_1\delta\delta(\alpha)-\gamma\dot{\delta}^{-\lambda}(\alpha).\end{equation}
where $A \equiv_1 B$ is equivalence of two set of faces of the
same dimension modulo empty faces which means that
\begin{enumerate}
  \item  if one set contains only empty faces then the other also
  contains only empty faces and these sets are equal,
  \item  or else both sets contain the same non-empty
  faces and any empty face in either
  set is an empty face on domain or codomain of a non-empty face is those
  sets.
\end{enumerate}
In other words if $\dot{A}$ denote non-empty faces in $A$,
$\ddot{A}$ denote empty faces in $A$ we have $A\equiv_1B$ iff
$\dot{A}=\dot{B}$ and $\ddot{A}\subseteq\ddot{B}\cup
1_{\gamma(\dot{B})\cup\delta(\dot{B})}$ and
 $\ddot{B}\subseteq\ddot{A}\cup
 1_{\gamma(\dot{A})\cup\delta(\dot{A})}$. Still in other words $A$
 and $B$ are sets of faces that generates, via $\gamma$ and $\delta$,
 the same substructures.

The last axiom, {\em loop filling} is the only other axiom that
does not mention order explicitly, it says that there are no empty
loops, i.e. if there is a loop it must be a codomain of at least
one face which is not  a loop.

The remaining four axioms talk about orders $<^+$ and $<^\sim$.
{\em Local discreteness} says that faces in a domain of any other
face cannot be comparable via the upper order $<^+$. The {\em
strictness}, {\em disjointness} together with {\em pencil
linearity} say in a sense that $<^\sim$ is the maximal strict
order order relation that is contained in the relation $<^-$ and
disjoint from $<^+$.

Note that as $<^+$, $<^-$ are the transitive closures of
elementary relations so these axioms are not first order axiom and
in fact they are expressed in the transitive closure logic.

\subsection*{Future work}

This paper covers only part of the program developed in \cite{Z}
for positive-to-one computads. We end this here as it is already
very long paper.  But the remaining parts of the program from
\cite{Z} for many-to-one computads and the application of this to
the multitopic categories, c.f. [HPM],[M], will be presented soon.

\subsection*{Content of the paper}

Section 2 contains the definition of a hypergraph and some
notation needed to introduce the notion of an ordered face
structure. In section 3 we introduce the main notion of this paper
the notion of an ordered face structure. In section 4 we develop
most of the needed elementary theory of ordered face structures.
This section should be more consulted when needed than read
through. The monotone morphism, is the stricter of two kind of morphisms
between ordered face structures, it preserves the order $<^\sim$
globally. The image of such a morphism is a convex set. In section
5 we show that from the convex set we can recover the whole
morphism up to an isomorphism. The domain of such a morphism is
recovered via cuts of empty loops in the convex subset. The next
two sections show the connection between positive and ordered face
structures. In section 6 we describe how we can divide a positive
face structure by an ideal to get an ordered face structure. In
section 7 we show that for any ordered face structure $S$ there is
a positive face structure $S^\dag$ and ideal so that $S^\dag$
divided by this ideal is isomorphic to $S$. The positive face
structure $S^\dag$ is defined with the help of cuts of so called
initial faces.
In sections 8 and 9 we describe some abstract structure of the category
$\ofs$ and show some of their properties.
This allow us to define in section 10 a free functor $(-)^*:\lfs\lra \oC$ from the
category of local face structures $\lfs$ to the category of
$\o$-categories. Local face structures are structures that have
operations $\gamma$ and $\delta$ as in ordered face structures but
with the order $<^\sim$ (in fact a binary relation) restricted to
domains of faces only. The section 11 discusses basic properties
of principal and normal ordered face structures the face
structures that are generated by a single face and such that can
be domains (in the sense of monoidal globular category $\ofs$) of
such structures. In section 12 we study decompositions of ordered
face structures. In \cite{Z} we have defined the decompositions of
positive face structures along some faces. Here we decompose
ordered face structure $S$ along a cut of initial faces
$\check{a}$ rather than a face $a$ as this decomposition in more
like a decomposition of the positive cover $S^\dag$ and then after
decomposition divided to get the decomposition of $S$. Doing it
this way we can deduce most of the properties of this
decomposition from the corresponding decomposition of positive
face structures. In section 13 we show that the $\o$-category
$T^*$, for $T$ being an ordered face structure, is in fact a
many-to-one computad. The next two sections 14 and 15 describe
with the help of ordered face structures the terminal many-to-one
computad and all the cells in an arbitrary many-to-one computad.

\subsection*{Notation and conventions}

As we already indicated we will intensionally confuse singletons
with elements when dealing with faces in ordered face structures.
In the paper we often will be using cells of different but
neighboring dimensions. As it is a bit confusing anyway we try to
make it a bit easier to follow by a careful use of the following
convention. $\alpha,\beta$ are faces of the same dimension, say
$n$, then $a,b$ are the faces of the same dimension $n-1$, $x,y,z$
are the faces of the same dimension $n-2$, $t,s$ are the faces of
the same dimension $n-3$, $u,v$ are the faces of the same
dimension $n-4$. We may use occasionally $A$, $B$ to denote faces
of dimension $n+1$. These faces may appear with indices but these
letter should be a direct hint which dimension we are working on.
The above examples were already using this convention. Last but
not least the composition of two morphisms
\begin{center}
\xext=1000 \yext=100
\begin{picture}(\xext,\yext)(\xoff,\yoff)
\putmorphism(0,0)(1,0)[X`Y`f]{500}{1}a
 \putmorphism(500,0)(1,0)[\phantom{Y}`Z`g]{500}{1}a
\end{picture}
\end{center}
may be denoted as either $g\circ f$ or more often $f;g$. In any
case we will write which way we mean the composition.


\subsection*{Acknowledgements} I like to thank Bill Boshuck, Victor
Harnik, and Mihaly Makkai for the conversation we had at the early
stage of this work.

\section{Hypergraphs}
A {\em hypergraph}\index{hypergraph} $S$ is
\begin{enumerate}
  \item a family $\{ S_k\}_{k\in \o}$ of finite sets of faces; only finitely many among these sets are
  non-empty;
  \item a family of
functions $\{ \gamma^S_k : S_{k+1}\dsum \b1_{S_k}\ra S_k \}_{k
\in\o }$; where $\b1_{S_k}=\{ 1_u : u\in S_k\}$ is the set of
empty faces of dimension $k$; the face $1_u$ is the empty
$k$-dimensional face on a non-empty face $u$ of dimension $k-1$.
  \item a family of total relations $\{ \delta^S_k :
S_{k+1}\dsum \b1_{S_k}\ra S_k\dsum \b1_{S_{k-1}}\}_{k\in\o}$; for
$a\in S_{k+1}$ we denote $\delta^S_{k}(a)=\{ x\in S_k\dsum
\b1_{S_{k-1}} \, :\, (a,x)\in\delta^S_{k} \}$; $\delta^S_{k}(a)$
is either singleton or it is non-empty subset of $S_k$.
Moreover $\delta^S_0 : S_{1}\dsum \b1_{S_0}\ra S_0\dsum
\b1_{S_{-1}}$ is a function (for this expression to make sense we
 assume that $S_{-1}=\emptyset$). We put
$\dot{\delta}(a)=\delta(a)\cap S$ and
$\ddot{\delta}(a)=\delta(a)\cap \b1_S$.
\end{enumerate}
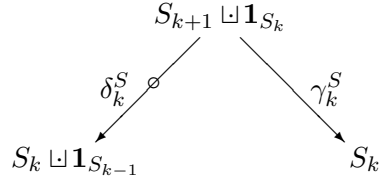
\begin{center} \xext=600 \yext=550
\begin{picture}(\xext,\yext)(\xoff,\yoff)
\settriparms[1`1`0;500]
 \putAtriangle(0,50)[S_{k+1}\dsum \b1_{S_k}`S_k\dsum \b1_{S_{k-1}}`S_k;\delta^S_k`\gamma^S_k`]
  \put(220,270){\makebox(100,100){$\circ$}}
\end{picture}
\end{center}
A {\em morphism of
hypergraphs}\index{morphism! of
hypergraphs}\index{hypergraph!morphism} $f:S\lra T$ is a family of
functions $f_k : S_k \lra T_k$, for $k \in\o$, such that the
diagrams
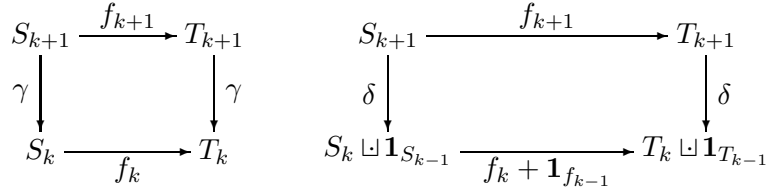
\begin{center}
\xext=2300 \yext=500 \adjust[`I;I`;I`;`I]
\begin{picture}(\xext,\yext)(\xoff,\yoff)
 \setsqparms[1`1`1`1;600`400]
 \putsquare(0,0)[S_{k+1}`T_{k+1}`S_k`T_k;
f_{k+1}`\gamma`\gamma`f_k]
 \setsqparms[1`1`1`1;1100`400]
 \putsquare(1200,0)[S_{k+1}`T_{k+1}`S_k\dsum \b1_{S_{k-1}}`T_k\dsum \b1_{T_{k-1}};
f_{k+1}`\delta`\delta`f_k+\b1_{f_{k-1}}]
\end{picture}
\end{center}
commute (where $\b1_{f_{k-1}}(1_x)=1_{f_{k-1}(x)}$, for $x\in
S_{k-2}$), for $k\in\o$.

The commutation of the left hand square is the commutation of the
diagram of sets an functions but in case of the right hand square
we mean more than commutation of a diagram of relations, i.e. we
demand that for any $a\in S_{\geq 1}$, $f_a:\delta(a)\lra
\delta(f(a))$ be a bijection, where $f_a$ is the restriction of
$f$ to $\dot{\delta}(a)$ (if $\delta(a)=1_u$ we mean by that
$\delta(f(a))=1_{f(u)}$).
The category of \index{category!$\hg$} hypergraphs is denoted by
$\hg$.

{\em Convention.} If $a\in S_k$ we treat $\gamma(a)$ sometimes as
an element of $S_{k-1}$ and sometimes as a subset $\{ \gamma(a)
\}$ of $S_{k-1}$. Similarly $\delta(a)$ is treated sometimes as a
set of faces or as a single face if this set of faces is a
singleton. In particular, when we write $\gamma(a)=\delta(b)$ we
mean rather $\{ \gamma(a) \}=\delta(b)$ or in other words that
$\delta(b)$ has one element this element is a face (not an empty
face) and that this face is $\gamma(a)$. We can also write
$\gamma(a)\in\delta(b)$ to mean that $\delta(b)\subset S$ and that
$\gamma(a)$ is one (of possibly many) elements of $\delta(b)$.
This convention simplifies the formulas considerably.

{\em Notation.} Before we go on, we need some notation.  Let $S$
be an ordered hypergraph.
\begin{enumerate}
   \item The dimension of $S$ is $\max\{ k\in\o : S_k\neq\emptyset
   \}$, and it is denoted by $dim(S)$.
  \item The sets of faces of different dimensions are assumed to be disjoint
  (i.e. $S_k\cap S_l=\emptyset$, for $k\neq l$);
   $S$ is also used to mean the set of all faces of $S$ i.e. $\bigcup_{k=0}^nS_k$;
   the notation $A\subseteq S$ mean
  that $A$ is a set of some faces of $S$; $A_k=A\cap S_k$, for
  $k\in\o$.
  \item If $a\in S_k$ then the face $a$ has dimension $k$ and we
   write $dim(a)=k$.
  \item $S_{\geq k} = \bigcup_{i\geq k} S_i $,
  $\;\;\;\;\; S_{\leq k} = \bigcup_{i\leq k} S_i $.
  The set $S_{\leq k} = \bigcup_{i\leq k} S_i $ is closed under
  $\delta$ and $\gamma$ so it is a sub-hypergraph of $S$, called
  $k$-{\em truncation}\index{truncation} of $S$.
  \item $\delta(A)=\bigcup_{a\in A}\delta(a)$  is the image of
  $A\subseteq S$ under $\delta$;
  \\ $\gamma(A)=\{ \gamma(a) \,:\, a\in A\}$ is the image of
  $A$ under $\gamma$. \\ Following the convention mentioned above if
  either $\gamma(A)$ or $\delta(A)$ is a singleton we may treat
  them as a (possibly empty) single face.
  \item For $a\in S_{\geq 1}$, the set $\theta(a)=\delta(a)\cup\gamma(a)$
  is {\em the set of codimension 1 faces} in $a$. We put $\dot{\theta}(a)=\theta(a)\cap S$.
\item Let $x,\alpha\in S$. We define the following subsets of faces of $S$:
\begin{enumerate}
  \item {\em empty domain faces}\index{face!empty domain -}: $S^\varepsilon = \{ a\in S :
  \delta(a)\in\b1_{S}\}$;
   \item {\em non-empty domain faces}\index{face!non-empty domain}: $S^{-\varepsilon} = S-S^\varepsilon$;
  we write $\delta^{-\varepsilon}(A)$ for $\delta(A)\cap S^{-\varepsilon}$;

  \item {\em loops}\index{face!loop}: $S^\lambda = \{ a\in S : \delta(a)=\gamma(a)\}$;
  \item {\em non-loops}\index{face!non-loop}: $S^{-\lambda} = S-S^\lambda$;
  we also write $\delta^{-\lambda}(A)$ for $\delta(A)\cap S^{-\lambda}$;
   \item {\em unary faces}\index{face!unary}: $S^u=\{ a\in S : |\dot{\delta}(a)|=1\}$;
   \item for $\alpha\in S_{\geq 2}$ we define the set of {\em internal
  faces} of $\alpha$;
 \[ \iota(\alpha)=\{ x\in S : \exists a,b\in \dot{\delta}^{-\lambda}(\alpha) :
 \gamma(a)=x\in\delta(b) \}=\gamma\dot{\delta}^{-\lambda}(\alpha)\cap\delta\dot{\delta}^{-\lambda}(\alpha) \]
   \item {\em internal faces}: $\iota(S)$;
  \item {\em initial faces}: $\cI =\cI^S =S^\varepsilon-\gamma(S^{-\lambda})$;
  \item $x-${\em cluster (of initial faces)}: $\cI_x =\cI^S_x = \{ \alpha\in\cI^S :\gamma\gamma(\alpha)=x\}$;
  \item {\em initial faces over} $\alpha$: $\cI^{\leq^+\alpha} =\{ \beta\in \cI : \beta\leq^+ \alpha \}$;
   \item $x-${\em cluster (of initial faces) over} $\alpha$: $\cI_x^{\leq^+\alpha} =\cI^{\leq^+\alpha}\cap\cI_x$.
\end{enumerate}

\item On each set $S_k$ we introduce two binary relations $<^{S_k,-}$ and
$<^{S_k,+}$. We usually omit $k$ in the superscript and sometimes
even $S$.
\begin{enumerate}
  \item $<^{S_0,-}$ is the empty relation. For $k>0$, the relation $<^{S_k,-}$ is the transitive closure of the relation $\lhd^{S_k,-}$ on $S_k$,
such that $a \lhd^{S_k,-} b$ iff $\gamma(a)\in\delta(b)$.  We
write $a\perp^{S_k,-}b$ if either $a <^{S_k,-} b$ or $b <^{S_k,-}
a$, and we write $a \leq^- b$ iff $a=b$ or $a <^- b$;

  \item $<^{S_k,+}$ is the transitive closure of the relation
  $\lhd^{S_k,+}$ on $S_k$,
such that $a \lhd^{S_k,+} b$ iff there is $\alpha\in
S_{k+1}^{-\lambda}$, such that $a\in \delta (\alpha)$ and
$\gamma(\alpha)=b$. We write $a\perp^{S_k,+}b$ if either $a
<^{S_k,+} b$ or $b <^{S_k,+} a$, and we write $a \leq^+ b$ if
either $a=b$ or $a <^+ b$.
  \item  $a\not\perp b$ if  both conditions $a\not\perp^+ b$
and $a\not\perp^- b$ hold.
\end{enumerate}


\item   Let $a,b\in S_k$. A {\em lower path from}\index{path!lower -}
$a$ to $b$ in $S$ is a sequence of faces $a_0, \ldots , a_m$ in
$S_k$ such that $a=a_0$, $b=a_m$ and,
$\gamma(a_{i-1})\in\delta(a_i)$,  for $i=1,\ldots ,m$.

A lower path is a {\em flat lower path}\index{path!flat lower -}
if it contains no loops other than $a$ or $b$.

\item   Let $x,y\in S_k$. An
{\em upper path from}\index{path!upper -} $x$ to $y$ in $S$ is a
sequence $a_0, \ldots , a_m$ in $S_{k+1}$ such that
$x\in\delta(a_0)$, $y=\gamma(a_m)$ and,
$\gamma(a_{i-1})\in\delta(a_i)$,  for $i=1,\ldots ,m$.

An upper path is a {\em flat upper path}\index{path!flat upper -}
if it contains no loops.

 \item The iterations of $\gamma$, $\delta$ and $\theta$ will be denoted in two different
 ways. By $\gamma^k$, $\delta^k$ and $\theta^k$ we mean $k$ applications of $\gamma$
 and $\delta$, respectively. By $\gamma^{(k)}$, $\delta^{(k)}$ and $\theta^{(k)}$ we mean the application as many times
 $\gamma$, $\delta$
 and $\theta$, respectively, to get faces of dimension $k$.  For
 example if $a\in S_5$ then $\delta^3(a)=\delta\delta\delta(a)\subseteq S_2$
 and $\delta^{(3)}(a)=\delta\delta(a)\subseteq S_3$.
\end{enumerate}

\section{Face structures}\label{face stuctures}

To simplify the notation, we treat both $\delta$ and $\gamma$ as
functions acting on faces as well as on sets of faces, which means
that sometimes we confuse elements with singletons. Clearly, both
$\delta$ and $\gamma$ when considered as functions on sets are
monotone.

We need the following relation. Let $S$ be a hypergraph. We
introduce an 'equality' relations  $\equiv_1$ on subsets of
$S_k\cup \b1_{S_{k-1}}$, for $k\in\o$, that may ignore the
$\b1_S$-part of the sets in presence of faces from $S$.  Let
$A,B\subseteq S_k\cup \b1_{S_{k-1}}$. We set that $A$ is 1-{\em
equal} $B$, notation $A\equiv_1B$, iff $A\cup 1_{\theta(A\cap S)}
= B\cup 1_{\theta(B\cap S)}$.

An {\em ordered face structure}\index{face structure!ordered}
$(S,<^{S_k,\sim})_{k\in\o}$ (also denoted $S$) is a hypergraph $S$
together with a family of $\{<^{S_k,\sim}\}_{k\in\o}$ of binary
relations ($<^{S_k,\sim}$ is a relation on $S_k$),   if it is
non-empty, i.e. $S_0\neq\emptyset$ and

\begin{enumerate}
 \item {\em Globularity:}\index{globularity}  for  $a\in S_{\geq 2}$:
\[ \gamma\gamma(a) =\gamma\delta(a)-\delta\dot{\delta}^{-\lambda}(a),\hskip 15mm
  \delta\gamma(a) \equiv_1\delta\delta(a)-\gamma\dot{\delta}^{-\lambda}(a).\]
  and for any $a\in S$:
  \[ \delta(1_a)=a=\gamma(1_a). \]
  \item {\em Local discreteness}:\index{local discreteness} if
  $x,y\in\delta(a)$ then $x\not\perp^+y$.
  \item {\em Strictness}:\index{strictness} for $k\in\o$, the
  relations $<^{S_k,+}$ and $<^{S_{k+1},\sim}$ are  strict orders\footnote{By {\em strict order}
  we mean an irreflexive and transitive relation.};
  $<^{S_0,+}$ is linear; (i.e. no flat path is a cycle).
  \item {\em Disjointness}:\index{disjointness} for $k\in\o$, the relation $<^{S_k,\sim}$ is
  a maximal strict order relation contained in $<^{S_k,-}$ that is disjoint from $<^{S_k,+}$, i.e. for $k>0$,
  \[\perp^{S_k,\sim}\cap \perp^{S_k,+}=\emptyset\]
  for any $a,b\in S_k$:
  \[ \mathrm{if}\;\; a<^\sim b\;\; \mathrm{then}\;\;  a<^- b \]
 \[ \mathrm{if}\;\; \theta(a)\cap\theta(b)= \emptyset\;\;\mathrm{then}\;\;
  a<^\sim b\;\; \mathrm{iff}\;\; a<^- b \]
  (i.e. if faces are not incident then $<^\sim$ is the same as $<^-$).
  \item {\em Pencil linearity}:
  \index{pencil linearity}\index{linearity!pencil -} for any $a,b\in S_{\geq 1}$, $a\neq b$,
  \[ \mathrm{if}\;\; \dot{\theta}(a)\cap\dot{\theta}(b)\neq \emptyset\;\;\mathrm{then}\;\;
  \mathrm{either}\;\; a\perp^\sim b\;\; \mathrm{or}\;\; a\perp^+ b \]
  for any $a\in S^{\varepsilon}_{\geq 2}$, $b\in S_{\geq 2}$,
  \[ \mathrm{if}\;\; \gamma\gamma(a)\in\iota(b)\;\;\mathrm{then}\;\;
  \mathrm{either}\;\; a<^\sim b\;\; \mathrm{or}\;\; a<^+ b \]
  (i.e. if faces are incident then they are comparable).
  \item {\em Loop-filling}:
  \index{loop-filling} $S^\lambda\subseteq \gamma(S^{-\lambda})$; (i.e. no
  empty loops).
\end{enumerate}
The relation $<^+$ is called the {\em upper order} and $<^\sim$ is
called {\em lower order}.

The morphism of ordered face structures, the {\em monotone
morphism}\index{morphism!monotone}, $f:S\lra T$ is a hypergraph
morphism that preserves the order $<^\sim$. The category of
ordered face structures, \index{category!$\ofs$} is denoted by
$\ofs$.

The relation $<^\sim$ in  an ordered face structure $S$ induces a
binary relation $(\dot{\delta}(a),<^{\sim}_a) $ for each $a\in
S_{>0}$ (where $<^{\sim}_a$ is the restriction of $<^\sim$ to the
set $\dot{\delta}(a)$). In the construction of the free
$\o$-categories over an ordered face structure we need to consider
hypergraph morphisms that preserves only this induced  structure
(not the whole relation $<^\sim$).  This is why we introduce the
category of local face structures.

A {\em local face structure}\index{face structure!local -}
$(S,<^{S_k,\sim}_a)_{a\in S}$ is a hypergraph $S$ together with a
family of $\{(\delta(a),<^{S_k,\sim}_a)\}_{a\in S}$ of binary
relations. The morphism of local face structures, the {\em local
morphism}\index{morphism!local}, $f:S\lra T$ is a hypergraph
morphism that is a local isomorphism  i.e. for $a\in S_{>1}$ the
restricted map
$f_a:(\dot{\delta}(a),<^\sim_a)\lra(\dot{\delta}(f(a)),<^\sim_{f(a)})$
is an order isomorphism, where $f_a$ is the restriction of $f$ to
$\dot{\delta}(a)$. The category of local face structures,
\index{category!$\lfs$} is denoted by $\lfs$.

Clearly we have a 'forgetful' functor:
\[ |-|: \ofs \lra \lfs \]
sending $(S,<^{S_k,\sim})_{k\in\o}$ to $(S,<^{\sim}_a)_{a\in
S_{>1}}$, where $<^{\sim}_a$ is the restriction of $<^\sim$ to
$\dot{\delta}(a)$, for $a\in S_{>1}$.

{\em Remarks.} Before we go on, we shall comment on the notions
introduced above.
\begin{enumerate}
\item The reason why we call the first condition 'globularity' is
that it will imply the usual globularity condition in the
$\o$-categories generated by ordered face structures. The word
'local' in 'Local discreetness' as anywhere else in the paper
refers to the fact that this property concerns sets of faces
constituting the domain of a face rather than the set of all
faces.

The property of 'pencil linearity' is strongly connected with the
property of positive face structures with the same name, c.f. [Z].
There it means that the set of faces with a fixed codomain $x$,
$\gamma$-{\em pencil}, (as well as the set of faces whose domains
contain a fixed face $x$, $\delta$-{\em pencil},) are linearly
ordered by $<^+$. For ordered face structures the same is true
about the faces that are not loops. The full condition also has
some implications for loops in pencils.

\item  The relation $\equiv_1$, needed to express the
$\delta$-globularity, is a way to say that two sets of faces, that
may contain empty faces, are essentially equal, even if they
differ by some empty faces. We identify via $\equiv_1$ two such
sets if those empty faces are morally there anyway.
$A,B\subseteq S_k\cup \b1_{S_{k-1}}$. Then the following
conditions are equivalent
\begin{enumerate}
  \item $A\equiv_1B$;
  \item the subhypergraphs of $S$ generated by $A$ and $B$ are equal;
  \item $\dot{A}=\dot{B}$, and $\ddot{A}\cup\b1_{\theta(\dot{A})}=\ddot{B}\cup\b1_{\theta(\dot{B})}$.
\end{enumerate}

\item We shall analyze in details $\gamma$-globularity and
$\delta$-globularity but some easier observations first:

\begin{enumerate}
  \item $\delta\delta^{-\varepsilon}(a)=\dot{\delta}\delta(a)$,
  $\;\delta\delta^{\varepsilon}(a)=\ddot{\delta}\delta(a)$.
  \item If $x\in T^\varepsilon$ then $\dot{\delta}^{-\lambda}(x)=\emptyset$ and
  $\gamma\gamma(x)=\gamma\delta(x)=\gamma(1_u)=\delta(1_u)=\delta\delta(x)=\gamma\delta(x)$.
  In particular, $\gamma(x)$ is a loop and
  $\delta(x)=1_{\gamma\gamma(x)}$, (i.e. $u=\gamma\gamma(x)$).
  \item If $\gamma(a)\in T^\varepsilon$ then
  $\delta\gamma(a)=1_{\gamma\gamma\gamma(a)}$.
\end{enumerate}

For $\delta$-globularity we distinguish two cases $\gamma(a)\in
T^\varepsilon$ and $\gamma(a)\in T^{-\varepsilon}$, and each has
two parts, for faces, and for empty-faces (the condition for empty
faces is translated to the condition about faces one dimension
lower).

\begin{description}
  \item[Case  $\gamma(a)\in T^{-\varepsilon}$]:
\begin{description}
  \item[faces:] $\delta\gamma(a)=\dot{\delta}\delta(a)-\gamma\dot{\delta}^{-\lambda}(a)$;
  \item[e-faces:]$\gamma\gamma\delta^\varepsilon(a)\subseteq\theta\delta\gamma(a)$;\\
  this is because we must have
  $\delta\delta^\varepsilon(a)\subseteq\b1_{\theta\delta\gamma(a)}$.

\end{description}
  \item[Case $\gamma(a)\in T^{\varepsilon}$]:
  \begin{description}
  \item[faces:]$\dot{\delta}\delta(a)\subseteq\gamma\dot{\delta}^{-\lambda}(a)$;
  \item[e-faces:]
  $\gamma\gamma\gamma(a)=\gamma\gamma\delta^\varepsilon(a)$;\\
   this is because we must have
   $1_{\gamma\gamma\gamma(a)}=\delta\gamma(a)=\delta\delta^{\varepsilon}(a)=
   1_{\gamma\gamma\delta^\varepsilon(a)}$.
\end{description}
\end{description}

The $\gamma$-globularity is much easier. We notice that if $a\in
T^\varepsilon$ then the condition is still slightly simpler, empty
faces play no role. We consider again two cases:

\begin{description}
  \item[Case  $a\in
  T^{\varepsilon}$:]$\gamma\gamma(a)=\gamma\delta(a)$.
  \item[Case $a\in T^{-\varepsilon}$:]  $\gamma\gamma(a) =\gamma\delta(a)-\delta\dot{\delta}^{-\lambda}(a)$; \\
  i.e. all elements of $\gamma\delta(a)$ but $\gamma\gamma(a)$ are
  in $\delta\dot{\delta}^{-\lambda}(a)$.  So we have
  $x_0\in\dot{\delta}(a)$ $\sim$-maximal in $\dot{\delta}(a)$ such
  that $\gamma\gamma(a)=\gamma(x_0)$.  $x_0$ might be a loop in
  which case, if $\gamma(a)$ is not a loop, there is another (unique) element
  $x_1\in\dot{\delta}^{-\lambda}(a)$, such that
  $\gamma(x_1)=\gamma\gamma(a)$.
\end{description}

\item If $S$ has dimension $n$, as a hypergraph, then we say that $S$
is {\em ordered $n$-face structure}\index{face
structure!ordered!$n$ -}.

\item A $k$-truncation of an ordered $n$-face structure $S$ is not
in general an an ordered $k$-face structure.  However
$k$-truncation of a local $n$-face structure is a local $k$-face
structure. This will be important later, in the description of the
many-to-one computads.

\item The {\em size of an ordered face structure}\index{size!of
ordered face structure} $S$ is the sequence natural numbers
$size(S)=\{ | S_n - \delta (S_{n+1}^{-\lambda})| \}_{n\in\o}$,
with almost all being equal $0$. We have an order $<$ on such
sequences, so that $\{ x_n \}_{n\in\o} < \{ y_n \}_{n\in\o}$ iff
there is $k\in\o$ such that $x_k< y_k$ and for all $l>k$, $x_l =
y_l$. This order is well founded and many facts about ordered face
structures will be proven by induction on the size.

\item Let $S$ be an ordered face structure. $S$ is $k$-normal\index{face structure!ordered!normal -} iff
$dim(S)\leq k$ and $size(S)_l=1$, for $l<k$. $S$ is
$k$-principal\index{face structure!ordered!principal -} iff
$size(S)_l=1$, for $l\leq k$. $S$ is principal iff
$size(S)_l\leq1$, for $l\in\o$. $S$ is principal of dimension $k$
iff $S$ is principal and $dim(S)=k$. By
$\pfs$\index{category!pfs@$\pfs$}
($\nfs$\index{category!nfs@$\nfs$}) we denote full subcategories
of $\ofs$ whose objects are principal (normal) ordered face
structures, respectively.
\end{enumerate}

\section{Combinatorial properties of ordered face structures}
\subsection*{Local properties}
\begin{lemma}
\label{fact1} Let $S$ be an ordered face structure, $x,a\in
S$. Then
\begin{enumerate}
   \item if $\delta(a)=1_x$ then  $x=\gamma\gamma(a)$;
  \item if $a\in S^\varepsilon$ then $\gamma(a)\in S^\lambda$;
  \item if $\gamma(a)\in S^\varepsilon$ then
  $\delta^\varepsilon(a)\neq\emptyset$;
  \item if $a\in S^{-\lambda}$ then $\gamma(a)\not\in\delta(a)$;
  \item $\ddot{\theta}\theta(a)=\delta\delta^\varepsilon(a)$;
  \item if $x<^+y$ then $y\not\in \cI$;
  \item if $x<^\sim y$ then $y\not S^{-\varepsilon}$.
\end{enumerate}
\end{lemma}
 {\em Proof.} Ad 1.  Assume that $\delta(a)=1_x$ for some $x\in S$. Then
$\dot{\delta}^{-\lambda}(a)=\emptyset$ and using
$\gamma$-globularity we get
\[  \gamma\gamma(a)=\gamma\delta(a)-\delta\dot{\delta}^{-\lambda}(a)
 =\gamma\delta(a)=\gamma(1_x)=x \]

 Ad 2. Suppose $a\in S^\varepsilon$. By 1. we have
 $\delta(a)=1_{\gamma\gamma(a)}$. Then
 $\dot{\delta}^{-\lambda}(a)=\emptyset$. So using $\delta$-globularity we
 have
 \[ \delta\gamma(a)=\delta\delta(a)-\gamma\dot{\delta}^{-\lambda}(a)
 = \delta(1_{\gamma\gamma(a)})=\gamma\gamma(a) \]
 i.e. $\gamma(a)$ is a loop.

Ad 3. Assume $\gamma(a)\in S^\varepsilon$. The using 2. and
globularity we obtain
\[1_{\gamma\gamma\gamma(a)}=\delta\gamma(a)=\delta\delta(a)-\gamma\dot{\delta}^{-\lambda}(a).\]
Thus there is $x\in\delta(a)$ such that
$\delta(x)=1_{\gamma\gamma\gamma(a)}$, i.e.
$x\in\delta^\varepsilon(a)$.

Ad 4. If we were to have $\gamma(a)\in\delta(a)$ then we would
have $\gamma(a)<^+\gamma(a)$ contradicting strictness of $<^+$.

Ad 5. As we have $\gamma\gamma(a),\gamma\delta(a)\subseteq S$ and
$\delta\gamma(a)\subseteq \delta\delta(a)$ we have
\[ \ddot{\theta}\theta(a)= (\gamma\gamma(a)\cap\delta\gamma(a)\cap\gamma\delta(a)\cap\delta\delta(a) )\cap
1_S = \ddot{\delta}\delta(a) = \delta\delta^\varepsilon(a).
\]

6. and 7. are obvious.
 $~\Box$

\begin{lemma}
\label{fact2} Let $S$ be an ordered face structure, $t,a,b,\alpha\in
S$.
\begin{enumerate}
  \item If $a\neq b$,  $a,b\in \dot{S}^{-\lambda}$, and either $\gamma(a)=\gamma(b)$ or
  $\dot{\delta}(a)\cap\dot{\delta}(b)\neq\emptyset$ then $a\perp^+b$.
  \item If $a,b\in \dot{\delta}^{-\lambda}(\alpha)$,  and either $\gamma(a)=\gamma(b)$
  or $\dot{\delta}(a)\cap\dot{\delta}(b)\neq\emptyset$ then $a=b$.
  \item Let $t\in \dot{\delta}\delta(a)$. Then there is a unique
  flat upper $\dot{\delta}^{-\lambda}(a)$-path from $t$ to
  $\gamma\gamma(\alpha)$.
  \item If $\alpha\in S^{-\varepsilon}$ then there is the
  $\sim$-largest element $a\in \dot{\delta}(\alpha)$. For this
  $a$ we have $\gamma(a)=\gamma\gamma(\alpha)$. All other elements
  of $\dot{\delta}(\alpha)$ have a well-defined $\sim$-successor.
  \item If $\gamma(\alpha)\in S^{-\lambda}$ then there is the
  $\sim$-largest element $a\in \dot{\delta}^{-\lambda}(\alpha)$. For this
  $a$ we have $\gamma(a)=\gamma\gamma(\alpha)$. All other elements
  of $\dot{\delta}^{-\lambda}(\alpha)$ have a well-defined $\sim$-successor in $\dot{\delta}^{-\lambda}(\alpha)$.
  \item If $\gamma(\alpha)\in S^{-\lambda}$ and $x\in
  \dot{\delta}\gamma(\alpha)$ then there is $a\in \dot{\delta}^{-\lambda}(\alpha)$ such that $x\in\delta(a)$.
  \item If $a<^+b$ then $\gamma(a)\leq^+\gamma(b)$.
\end{enumerate}
\end{lemma}
{\em Proof.} Ad 1. Let $a, b$ be as in the Lemma. By pencil
linearity, as $\dot{\theta}(a)\cap\dot{\theta}(b)\neq\emptyset$,
it is enough to show that $a\perp^\sim b$ does not hold. Suppose
contrary that $a<^\sim b$. Then, by disjointness, $a<^-b$. Thus we
have a flat lower path $a=a_0,\ldots a_k=b$, with $k>0$.  Now if
$\gamma(a)=\gamma(b)$ then we have a flat upper path
$\gamma(a),a_1,\ldots, a_k,\gamma(a)$ showing that
$\gamma(a)<^+\gamma(a)$. This contradicts strictness. on the other
hand, if some $x\in\delta(a)\cap\delta(b)$ then we have a flat
upper path $x,a_0,\ldots, a_{k-1},\gamma(a_{k-1})$. Hence
$x<^+\gamma(a_{k-1})$. As $x,\gamma(a_{k-1})\in\delta(b)$ we get a
contradiction with local discreetness.

Ad 2.  This is immediate consequence of 1. and local discreetness.

Ad 3. Fix $t\in \dot{\delta}\delta(a)$. Let
$t,x_1,\ldots,x_k,\gamma(x_k)$ be the longest flat upper
$\dot{\delta}\delta(a)$-path starting from $t$ (it might be
empty). Such a path exists by strictness. If
$\gamma(x_k)=\gamma\gamma(a)$ we are done. So assume that
$\gamma(x_k)\neq\gamma\gamma(a)$. We have
$\gamma(x_k)\in\gamma\delta(a)$.  So by globularity
$\gamma(x_k)\in\delta\dot{\delta}^{-\lambda}$ and hence there is
$x_{k+1}\in\dot{\delta}^{-\lambda}(a)$ such that
$\gamma(x_k)\in\delta(x_{k+1})$, i.e. $t,x_1,\ldots,
x_k,x_{k+1},\gamma(x_{k+1})$ is a longer flat path starting from
$t$ and we get a contradiction. Thus $\gamma(x_k)=\gamma\gamma(a)$
and we have $\dot{\delta}\delta(a)$-path from $t$ to
$\gamma\gamma(a)$.  The uniqueness of this path follows from 2.

Ad 4. As $\alpha\in S^{-\varepsilon}$, by globularity, there is
$a\in\delta(\alpha)$ such that $\gamma(a)=\gamma\gamma(\alpha)$.
If there is a loop $a\in\delta(a)$ such that
$\gamma(a)=\gamma\gamma(a)$ then by local discreetness there is
$\sim$-largest such loop.  Let $a_0$ be the $\sim$-largest loop in
$\delta(\alpha)$ such that $\gamma(a_0)=\gamma\gamma(\alpha)$ if
such a loop exists or else the unique
$a_0\in\delta^{-\lambda}(\alpha)$ such that
$\gamma(a_0)=\gamma\gamma(\alpha)$. We shall show that $a_0$ is
the $\sim$-largest element in $\delta(\alpha)$.

We consider two cases. If $\gamma(b)=\gamma\gamma(\alpha)$ by
pencil linearity, 2. and definition of $a_0$ we have $b\leq^\sim
a_0$. If $\gamma(b)\neq\gamma\gamma(\alpha)$ then by 1. there is a
flat upper $\delta(\alpha)$-path
$\gamma(b),b_1,\ldots,b_k,\gamma\gamma(\alpha)$ with $k>0$. By the
previous argument $b_k\leq^\sim a_0$ and hence $b<^\sim a_0$. Thus
in either case $a_0$ is the $\sim$-largest element in
$\delta(\alpha)$. For $a\in \delta(\alpha)-\{ a_0 \}$ we define
the successor in $\delta(\alpha)$ as follows:
\[ suc_{\alpha}(a)= \left\{ \begin{array}{ll}
           \inf_\sim(A)&
           \mbox{if  $A=\{ a'\in\delta^\lambda(\alpha):\gamma(a)=\gamma(a')\,{\rm and\;} a'<^\sim a\}\neq\emptyset$,}\\
           a''& \mbox{such that $a''\in\delta^{-\lambda}(\alpha)$,  $\gamma(a)\in\delta(a'')$, otherwise.}
                                    \end{array}
                \right. \]
The verification that it is a well defined successor is left for
the reader.

Ad 5. Assume that $\gamma(\alpha)\in S^{-\lambda}$. First assume
we have $a_0\in\dot{\delta}^{-\lambda}(\alpha)$ such that
$\gamma(a_0)=\gamma\gamma(\alpha)$. Then by 2. such $a_0$ is
unique and by an argument similar to the one above $a_0$ is
$\sim$-largest in $\dot{\delta}^{-\lambda}(\alpha)$ and all other
elements in $\dot{\delta}^{-\lambda}(\alpha)$ have a successor
there. Thus it remains to find $a_0$.

Note that to find $a_0$ it is enough to find
$x\in\dot{\delta}\delta(\alpha)$ such that
$x\neq\gamma\gamma(\alpha)$.  Having such $x$, by 3., we have a
flat upper $\dot{\delta}^{-\lambda}(\alpha)$-path
$x,b_1,\ldots,b_k,\gamma\gamma(\alpha)$ with $k>0$. We put
$a_0=b_k$.

To find $x$ we consider two cases. If $\gamma\in S^{-\varepsilon}$
then
$\gamma\gamma(\alpha)\not\in\delta\gamma(\alpha)\subseteq\dot{\delta}\delta(\alpha)$
and $\delta\gamma(\alpha)\neq\emptyset$. Then any element of
$\delta\gamma(\alpha)$ can be taken as $x$.

If $\gamma\in S^{\varepsilon}$ then
$1_{\gamma\gamma\gamma(\alpha)}=\delta\gamma(\alpha)\subseteq\delta\delta(\alpha)$.
So there is $a\in\delta(\alpha)$ such that $\delta(\alpha)=
1_{\gamma\gamma\gamma(\alpha)}$.  If
$\gamma(a)=\gamma\gamma(\alpha)$ then we found $a_0=a$ directly.
Otherwise
$\gamma\gamma(\alpha)\neq\gamma(a)\in\gamma\delta(\alpha)$. So
$\gamma(a)\in\dot{\delta}\delta(a)$ and we put $x=\gamma(a)$.

Ad 6. Suppose $\gamma(\alpha)\in S^{-\lambda}$ and
$x\in\dot{\delta}\gamma(\alpha)$. Then $\gamma(\alpha)\in
S^{-\varepsilon}$ and we have
$\gamma\gamma(\alpha)\not\in\delta\gamma(\alpha)\subseteq\dot{\delta}\delta(\alpha)$.
By 3. there is a flat upper $\dot{\delta}^{-\lambda}(\alpha)$-path
$x,a_1,\ldots,a_k,\gamma\gamma(\alpha)$. Then $x\in\delta(a_1)$
and $a_1\in\dot{\delta}^{-\lambda}(\alpha)$, as required.

Ad 7.  The essential case $a\lhd^+b$ follows from 3.  Then use
induction. $~\Box$

{\em Notation}. Having \ref{fact2}.4 we can introduce farther
notation. If $\alpha\in S^{-\varepsilon}$ then the $\sim$-largest
element in $\dot{\delta}(\alpha)$ will be denoted
$\varrho(\alpha)$.

If $\gamma(\alpha)\in S^{-\lambda}$ then the $\sim$-largest
element in $\dot{\delta}^{-\lambda}(\alpha)$ will be denoted
$\varrho^{-\lambda}(\alpha)$.

Clearly, whenever the formulas make sense, we have
\[ \gamma\gamma(\alpha)=\gamma(\varrho(\alpha))=\gamma(\varrho^{-\lambda}(\alpha)).\]

\begin{lemma}
\label{fact3} Let $S$ be an ordered face structure,
$a,b,\alpha\in S$.
\begin{enumerate}
  \item $\gamma(a)\in S^\lambda$ iff $\delta(a)\subseteq
  S^\lambda$ or $a\in S^\varepsilon$.
  \item If $b\in S^\lambda$ and $a<^+b$ then $a\in S^\lambda$.
  \item If $a\in S^\lambda$ then there is $\alpha\in
  \cI_{\gamma(a)}$ such that $\gamma(\alpha)\leq^+a$ and
  $\delta(\alpha)=1_{\gamma(a)}$.
  \item If $a\in S^\varepsilon$ then there is $b\in \cI$ such that $b\leq^+a$.
  In that case $\delta(b)=\delta(a)$.
  \item If $x\in S^\lambda$ then there is $b\in \cI$ such that $\gamma(b)\leq^+x$.
  In that case $\delta(b)=1_{\gamma(x)}$.
\end{enumerate}
\end{lemma}

{\em Proof.} Ad 1. $\La$. By Lemma \ref{fact1}.2 if $a\in
S^\varepsilon$ then $\gamma(a)$ is a loop. So assume that
$\delta(a)\subseteq S^\lambda$.  Then
$\dot{\delta}^{-\lambda}(a)=\emptyset(=\delta\dot{\delta}^{-\lambda}(a)=\gamma\dot{\delta}^{-\lambda}(a)$.
As for $x\in\delta(a)$ we have $\gamma(x)=\delta(x)$, we also have
$\gamma\delta(a)=\delta\delta(a)$. Thus by globularity, we have
\[ \gamma\gamma(a)=\gamma\delta(a)-\delta\dot{\delta}^{-\lambda}(a)=
\gamma\delta(a)=\delta\delta(a)=\delta\delta(a)-\gamma\dot{\delta}^{-\lambda}(a)=\delta\gamma(a)\]
i.e. $\gamma(a)$ is a loop, as required.

$\Ra$. Suppose now that $\gamma(a)\in S^\lambda$. If
$\delta(a)=1_{\gamma\gamma(a)}$ then $a\in S^\varepsilon$.  So
assume that $\delta(a)\subseteq S$.  By globularity, we have
\begin{equation}\label{f3}
\gamma\delta(a)-\delta\dot{\delta}^{-\lambda}(a)=\gamma\gamma(a)=\delta\gamma(a)
=\delta\delta(a)-\gamma\dot{\delta}^{-\lambda}(a).
\end{equation}
If $\dot{\delta}^{-\lambda}(a)=\emptyset$ then
$\delta\dot{\delta}^{-\lambda}(a)=\gamma\dot{\delta}^{-\lambda}(a)=\emptyset$
and $\gamma\delta(a)=\delta\delta(a)=\gamma\gamma(a)$ i.e.
$\delta(a)\subseteq S^\lambda$, as required.

So suppose now that $\dot{\delta}^{-\lambda}(a)\neq\emptyset$. Let
$x\in\dot{\delta}^{-\lambda}(a)$ be the $\sim$-largest element in
$\dot{\delta}^{-\lambda}(a)$.  Such an $x$ exists by Lemma
\ref{fact2}.  Then
$\gamma(x)\in\gamma\dot{\delta}^{-\lambda}(a)\subseteq
\gamma\delta(a)$. By (\ref{f3}) we have
$\gamma(x)\in\delta\dot{\delta}^{-\lambda}(a)$, and hence we have
$x'\in\dot{\delta}^{-\lambda}(a)$ such that
$\gamma(x)\in\delta(x')$. As $x,x'\in\delta(a)$ we have
$x\not\perp^+x'$. Moreover if we were to have
$\gamma(x')\in\delta(x)$ we would have $\gamma(x)<^+\gamma(x)$
contradicting strictness.  Thus, by pencil linearity, we have
$x<^\sim x'$, i.e. $x$ is not $\sim$-largest element in
$\dot{\delta}^{-\lambda}(a)$ contrary to the supposition. Hence
$\dot{\delta}^{-\lambda}(a)=\emptyset$ indeed, and
$\delta(a)\subseteq S^\lambda$, as required.

Ad 2. Use 1. and then induction.

Ad 3. Use loop-filling, pencil linearity, strictness, and 1. to
get a maximal upper $S-\gamma(S^{-\lambda})$-path
$\alpha_1,\ldots,\alpha_k,a$ ending at $a$ with $k>0$. Then
$\alpha_1\in\cI_{\gamma(a)}$ and $\gamma(\alpha_1)\leq^+a$.

Ad 4. Suppose $\gamma(\alpha)\in S^\varepsilon$. By globularity we
have
$1_{\gamma\gamma\gamma(\alpha)}=\delta\gamma(\alpha)\subseteq\delta\delta(\alpha)$.
Then there is $a\in\delta^\varepsilon(\alpha)\neq\emptyset$. The
thesis follows form the above observation= , strictness, and
inductions.

Ad 5. Fix $x\in S^\lambda$. By loop-filling and strictness there
is $a\in S^\varepsilon$ such that $\gamma(a)\leq^+x$. The rest
follows from 4. $~\Box$

\begin{lemma}
\label{fact4} Let $S$ be an ordered face structure $\alpha,a,b\in S$.
\begin{enumerate}
  \item $\iota\delta(\alpha)=\iota\gamma(\alpha)$.
  \item If $a<^+b$ then $\iota(a)\subseteq \iota(b)$.
\end{enumerate}
\end{lemma}

{\em Proof.} Ad 1. First we prove
$\iota\gamma(\alpha)\subseteq\iota\delta(\alpha)$. Fix
$u\in\iota\gamma(\alpha)$, i.e. we have
$x,y\in\dot{\delta}^{-\lambda}\gamma(\alpha)$ such that
$\gamma(x)=u\in\delta(y)$.  Let $x,a_1,\ldots, a_k,x'$ be a
maximal flat upper ${\delta}^{-\lambda}\gamma(\alpha)$-path
starting from $x$ such that $\gamma\gamma(a_i)=u$ for
$i=1,\ldots,k$. It might be empty in which case $x=x'$.  As $x\in
S^{-\lambda}$ by Lemma \ref{fact3}.2, $x'\in S^{-\lambda}$. Since
$\gamma(x')=u\in\iota\gamma(\alpha)$ it follows that
$\gamma(x')\neq\gamma\gamma\gamma(\alpha)$.  Thus
$x'\neq\gamma\gamma(\alpha)$. By Lemma \ref{fact2}.3 there is
$a\in\dot{\delta}^{-\lambda}(\alpha)$ such that $x'\in\delta(a)$.
By maximality of $x,a_1,\ldots, a_k,x'$ we have
$\gamma\gamma(a)\neq\gamma(x')$. Again by  Lemma \ref{fact2}.3
there is $y'\in\dot{\delta}^{-\lambda}(a)$ such that
$\gamma(x')=u\in\delta(y')$. But then
$u\in\iota(a)\subseteq\iota\delta(\alpha)$.

Now we shall show
$\iota\delta(\alpha)\subseteq\iota\gamma(\alpha)$. Let
$a\in\delta(\alpha)$ and $u\in\iota(a)$, i.e. there are
$x',y'\in\delta^{-\lambda}(a)$ such that
$\gamma(x')=u\in\delta(y')$.  We shall construct
$x,y\in{\delta}^{-\lambda}\gamma(\alpha)$ such that
$\gamma(x)=u\in\delta(y)$, i.e. $u\in\iota(\alpha)$.

Construction of $x$. Let $a_l,\ldots,a_1,x'$ be the maximal flat
$\delta^{-\lambda\varepsilon}(\alpha)$-path (possibly empty)
ending at $x'$ such that $\gamma\gamma(a_i)=\gamma(x')$ and
$\gamma(a_i)\in S^{-\lambda}$, for $i=1,\ldots,l$. Thus there is
$x\in\delta^{-\lambda}(a_l)$ such that
$\gamma(x)=\gamma\gamma(a_l)$. If $x\in\delta\gamma(\alpha)$ we
have $x$ with the required property (if the sequence is empty
$x=x'$). So suppose contrary that $x\not\in\delta\gamma(\alpha)$.
Since $x\in\delta\delta(\alpha)$, by globularity, it follows that
$x\in\gamma\dot{\delta}^{-\lambda}(\alpha)$. So there is
$a_{l+1}\in\dot{\delta}^{-\lambda}(\alpha)$ such that
$\gamma(a_{l+1})=x$. Since $x\in S^{-\lambda}$, we have
$a_{l+1}\in S^{-\varepsilon}$.  But then the path
$a_{l+1},a_l,\ldots,a_1,x'$ is longer then the maximal one and we
get a contradiction.

Construction of $y$. Let $b_k,\ldots,b_1,x'$ be the maximal flat
$\delta^{-\lambda\varepsilon}(\alpha)$-path (possibly empty)
ending at $y'$ such that
$u\in\delta\dot{\delta}^{-\lambda}(b_i)=\gamma(x')$ and
$\gamma(b_i)\in S^{-\lambda}$, for $i=1,\ldots,k$. By Lemma
\ref{fact2}.6 there is $y\in\delta^{-\lambda}(b_k)$ such that
$u\in\delta(y)$ (if the sequence is empty $y=y'$). If
$y\in\delta\gamma(\alpha)$ we have $y$ as required.  So suppose
that $y\not\in\delta\gamma(\alpha)$.  As
$y\in\delta\delta(\alpha)$, by globularity, we have that
$y\in\gamma\dot{\delta}^{-\lambda}(\alpha)$.  So there is
$b_k+1\in\dot{\delta}^{-\lambda}(\alpha)$ such that
$\gamma(b_{k+1})=y$. Since $y\in S^{-\lambda}$, we have
$b_{k+1}\in S^{-\varepsilon}$. But then the path
$b_{k+1},b_b,\ldots,b_1,y'$ is longer then the maximal one and we
get a contradiction again.

Ad 2. Use 1. and induction. $~\Box$

\begin{lemma}
\label{fact5} Let $S$ be an ordered face structure $\alpha,a,b\in S$.
\begin{enumerate}
  \item We have inclusions
  \begin{center} \xext=600 \yext=500
\begin{picture}(\xext,\yext)(\xoff,\yoff)
 \setsqparms[1`1`1`1;600`400]
 \putsquare(0,50)[\dot{\delta}\gamma\gamma(\alpha)`\dot{\delta}\gamma\delta(\alpha)`
 \dot{\delta}\delta\gamma(\alpha)`\dot{\delta}\delta\delta(\alpha);```]
 \end{picture}
\end{center}
  \item $\dot{\theta}\theta(a)=\gamma\gamma(a)\cup\dot{\delta}\dot{\delta}^{-\lambda}(a)$,
  $\gamma\gamma(a)\cap\dot{\delta}\dot{\delta}^{-\lambda}(a)=\emptyset$.
  \item $\dot{\theta}\theta(a)=\dot{\delta}\gamma(a)\cup\gamma\dot{\delta}^{-\lambda}(a)$,
  $\dot{\delta}\gamma(a)\cap\gamma\dot{\delta}^{-\lambda}(a)=\emptyset$.
   \item $\dot{\theta}\theta(a)=\gamma\gamma(a)\cup\iota(a)\cup\delta\gamma^{-\lambda\varepsilon}(a)$ (disjoint sum).
   \item $\dot{\theta}\theta(a)=\dot{\theta}\dot{\theta}(a)=\dot{\theta}\delta(a)$.
   \item $\dot{\theta}\theta\theta(\alpha)=\dot{\theta}\theta\gamma(\alpha)$.
   \item If $a<^+b$ then $\dot{\theta}\theta(a)\subseteq\dot{\theta}\theta(b)$.
   \item $\dot{\theta}\theta^{(k+1)}(a)=\dot{\delta}\dot{\delta}^{-\lambda}\gamma^{(k+2)}(a)\cup\gamma^{(k)}(a)$.
   \item $\gamma\ddot{\theta}\theta^{(k+2)}(a)\subseteq
   \dot{\theta}^{(k)}(a)$. (don't bother with $1_x$'s)
   \item $\gamma\dot{\delta}^{-\lambda}\delta(\alpha)=\gamma\dot{\delta}^{-\lambda}\gamma(\alpha)$.
   \item $\gamma\dot{\delta}^{-\lambda}\theta^{(k+2)}
   (\alpha)=\gamma\dot{\delta}^{-\lambda}\gamma^{(k+2)}(\alpha)$.
\end{enumerate}
\end{lemma}
{\em Proof.} Ad 1.  This is an easy consequence of
$\delta\gamma(a)\subseteq\delta\delta(a)$ and
$\gamma\gamma(a)\subseteq\gamma\delta(a)$.

Ad 2. Let
$A=\gamma\gamma(a)\cup\dot{\delta}\dot{\delta}^{-\lambda}(a)$.
Clearly $A\subseteq\dot{\theta}\theta(a)$.  We shall show the
converse inclusion. From globularity we have
$\gamma\gamma(a)\in\gamma\delta(a)\subseteq A$ and
$\dot{\delta}\gamma\subseteq\dot{\delta}\delta(a)=
\dot{\delta}\dot{\delta}^{-\lambda}(a)\cup\dot{\delta}\delta^\lambda(a)\cup\dot{\delta}\ddot{\delta}(a)$.
Moreover $\dot{\delta}\delta^\lambda(a)=
\gamma\delta^\lambda(a)\subseteq\gamma\delta(a)$. Finally, if
$\ddot{\delta}(a)\neq\emptyset$ then
$\ddot{\delta}(a)=1_{\gamma\gamma(a)}$. So
$\dot{\delta}\ddot{\delta}(a)=\delta(1_{\gamma\gamma(a)})=\gamma\gamma(a)\in
A$. Thus the other inclusion holds as well. The second part
follows directly from $\gamma$-globularity.

Ad 3. Let
$B=\dot{\delta}\gamma(a)\cup\gamma\dot{\delta}^{-\lambda}(a)$.
Clearly $B\subseteq\dot{\theta}\theta(a)$.  We shall show the
converse inclusion. From globularity we have
$\dot{\delta}\delta(a)\subseteq B$ and
$\gamma\gamma\in\gamma\delta(a)=\gamma\dot{\delta}^{-\lambda}(a)
\cup\gamma\delta^\lambda(a)\cup\gamma\ddot{\delta}(a)$ Moreover
$\gamma\delta^\lambda(a)=
\delta\delta^\lambda(a)\subseteq\delta\delta(a)$. Finally, if
$\ddot{\delta}(a)\neq\emptyset$ then
$\ddot{\delta}(a)=1_{\gamma\gamma(a)}$. So
$\gamma\ddot{\delta}(a)=\gamma(1_{\gamma\gamma(a)})=\gamma\gamma(a)\in
A$. Thus the other inclusion holds as well. The second part
follows directly from $\delta$-globularity.

Ad 4. Using 2. and 3. we have
\[ \iota(a)\cap(\gamma\gamma(a)\cup\delta\gamma(a))=
(\gamma\dot{\delta}^{-\lambda}(a)\cap\delta\dot{\delta}^{-\lambda}(a))\cap(\gamma\gamma(a)\cup\delta\gamma(a))
\subseteq \]
\[ \subseteq
(\gamma\dot{\delta}^{-\lambda}(a)\cap\delta\gamma(a))\cup(\gamma\gamma(a)\cap\delta\dot{\delta}^{-\lambda}(a))=
\emptyset\cup\emptyset=\emptyset
\]

\[ \iota(a)\cup(\gamma\gamma(a)\cup\delta\gamma(a))=
(\gamma\dot{\delta}^{-\lambda}(a)\cap\delta\dot{\delta}^{-\lambda}(a))\cup(\gamma\gamma(a)\cup\delta\gamma(a))
\supseteq \]
\[ \supseteq
(\gamma\dot{\delta}^{-\lambda}(a)\cup\delta\gamma(a))\cap(\gamma\gamma(a)\cup\delta\dot{\delta}^{-\lambda}(a))=
\dot{\theta}\theta(a)\cap\dot{\theta}\theta(a)=\dot{\theta}\theta(a)
\]
Note that $\dot{\delta}\gamma(a)=\delta\gamma^{-\varepsilon}(a)$
and if $\gamma(a)\in S^\lambda$ then
$\dot{\delta}\gamma(a)=\gamma\gamma(a)$.  From these observations
the rest follows.

5. Easy application of 2. and 3.

6. We need to show that
$\dot{\theta}\theta\delta(\alpha)\subseteq\dot{\theta}\theta\gamma(\alpha)$.
By 4. it is enough to show the following three inclusions:
$\iota\delta(\alpha)\subseteq\dot{\theta}\theta\gamma(\alpha)$,
$\gamma\gamma\delta(\alpha)\subseteq\dot{\theta}\theta\gamma(\alpha)$,
$\dot{\delta}\gamma\delta(\alpha)\subseteq\dot{\theta}\theta\gamma(\alpha)$.
The first follows immediately from 4. If $\gamma(\alpha)\in
S^\lambda$ then, as in this case
$\gamma\delta(\alpha)=\gamma\gamma(\alpha$  the second and third
inclusions hold as well.  Thus we shall show the second and third
inclusion in case $\gamma(\alpha)\in S^{-\lambda}$.

Assume $t\in\gamma\gamma\delta(\alpha)$.  Pick $\sim$-minimal
$a\in\delta(\alpha)$ such that $t=\gamma\gamma(a)$. Then either
$a\in S^\varepsilon$ or $a\in S^{-\varepsilon}$. In the former
case $\delta(a)=1_{\gamma\gamma(a)}\in\delta\delta(\alpha)$, and
by $\delta$-globularity (see definition of $\equiv_1$) we have
that $t=\gamma\gamma(a)\in\dot{\theta}\delta\gamma(\alpha)$. In
the later case by Lemma \ref{fact2}.4 (see also notation after the
proof) $x=\varrho^{-\lambda}(a)\in\delta^{-\lambda}(a)$ is well
defined and we have $\gamma(x)=\gamma\gamma(a)$. By
$\sim$-minimality of $a$, we get
$x\not\in\gamma\dot{\delta}^{-\lambda}(\alpha)$.  As
$x\in\delta\delta(\alpha)$ again by $\delta$-globularity we have
that $x\in\delta\gamma(\alpha)$. Thus
$t=\gamma(x)\in\gamma\delta\gamma(\alpha)$. Thus in either case
$t\in\dot{\theta}\theta\gamma(\alpha)$.  This end the proof of the
second inclusion.

Now assume $t\in\delta\gamma\delta(\alpha)$. Pick $\sim$-minimal
$a\in\delta(\alpha)$ such that $t\in\delta\gamma(a)$. Then either
$a\in S^\varepsilon$ or $a\in S^{-\varepsilon}$. In the former
case $\gamma(a)\in S^\lambda$. Then using the second inclusion we
get
\[ t\in\delta\gamma(a)=\gamma\gamma(a)\in\gamma\gamma\delta(\alpha)\subseteq
\dot{\theta}\theta\gamma(\alpha) \]

In the later case by Lemma \ref{fact2}.6 there is $x\in\delta(a)$
(i.e. $x\in\delta\delta(\alpha)$) such that $t\in\delta(x)$. By
$\sim$-minimality of $a$, $x\not\in\gamma\dot{\delta}^{-\lambda}$.
So by $\delta$-globularity we have $x\in\delta\gamma(\alpha)$.
Thus
$t\in\dot{\delta}(x)\subseteq\dot{\delta}\delta\gamma(\alpha)$.
Thus in either case $t\in\dot{\theta}\theta\gamma(\alpha)$. This
end the proof of the third inclusion and the whole statement 6.

For 7.  and 8. Use 1., 5.,  6., and induction.

9. Exercise.

Ad 10. $\supseteq$. As $\delta\gamma(\alpha)\subseteq
\delta\delta(\alpha)$ we have
$\dot{\delta}^{-\lambda}\gamma(\alpha)\subseteq
\dot{\delta}^{-\lambda}\delta(\alpha)$.  So
$\gamma\dot{\delta}^{-\lambda}\gamma(\alpha)\subseteq
\gamma\dot{\delta}^{-\lambda}\delta(\alpha)$.

$\subseteq$. Let
$t\in\gamma\dot{\delta}^{-\lambda}\delta(\alpha)$. Pick
$\sim$-minimal $a\in\delta(\alpha)$ such that there is
$x\in\dot{\delta}^{-\lambda}(a)$ so that  $t=\gamma(x)$. By
$\sim$-minimality of  $a$
$x\not\in\gamma\dot{\delta}^{-\lambda}(\alpha)$ and hence by
$\delta$-globularity $x\in\dot{\delta}^{-\lambda}\gamma(\alpha)$.
Thus $t=\gamma(x)\in\gamma\dot{\delta}^{-\lambda}\gamma(\alpha)$,
as required.

11. follows from 10. by induction. $~\Box$

\begin{lemma}[$\dot{\theta}\theta$ induction]\label{theta induction}
Let $S$ be an ordered face structure $\alpha,a,b\in S$.
\begin{enumerate}
  \item if $a\in\dot{\delta}^{-\lambda\varepsilon}(\alpha)$ then
  $\gamma\gamma(\alpha)\not\in\delta(a)$;
  \item if $a,b\in\dot{\delta}^{-\lambda\varepsilon}(\alpha)$ then
  $\delta(a)\cap\delta(b)=\emptyset$;
  \item $\dot{\theta}\theta(\alpha)=\gamma\gamma(\alpha)\cup\delta\dot{\delta}^{-\lambda\varepsilon}
  (\alpha)$;
  \item $\dot{\theta}\theta$-induction.  Whenever
\begin{enumerate}
  \item if $A\subset \dot{\theta}\theta(\alpha)$;
  \item $\gamma\gamma(\alpha)\in A$;
  \item for all $a\in \dot{\delta}^{-\lambda\varepsilon}(\alpha)$, if
  $\gamma(a)\in A$ then $\dot{\delta}(a)\subseteq A$
\end{enumerate}
 we have $A=\dot{\theta}\theta(\alpha)$.
\end{enumerate}
\end{lemma}

{\em Proof.} 1. and 2. follows from pencil linearity. 3. follows
from Lemma \ref{fact5}.2 and that
$\delta\dot{\delta^{-\lambda\varepsilon}}=\dot{\delta}\dot{\delta^{-\lambda}}$.
The $\dot{\theta}\theta$-induction follows from 3. $~\Box$

\begin{lemma}
\label{fact6} Let $S$ be an ordered face structure $a,a'\in
S$.
\begin{enumerate}
  \item If $a,a'\in S^{-\lambda}$ and $\gamma(a)\in\delta(a')$ then $a<^\sim a'$.
  \item If $a<^-a'$ and $\theta(a)\cap\theta(a')\neq\emptyset$ then
  $\gamma(a)\in\delta(a')$.
\end{enumerate}
\end{lemma}
{\em Proof.} Ad 1. Let $a,b\in S^{-\lambda}$ such that
$\gamma(a)\in\gamma(b)$. By strictness we cannot have $b<^-a$.
Thus by pencil linearity it is enough to to show that
$a\not\perp^+b$.

Suppose $a<^+b$, i.e. there is a flat upper path
$a,\alpha_1,\ldots,\alpha_r,b$. As $a\in S^{-\lambda}$, by Lemma
\ref{fact3}, $\gamma(\alpha_i)\in S^{-\lambda}$, for
$i=1,\ldots,r$. Now either
$\gamma(a)=\gamma\gamma(\alpha_r)=\gamma(b)$ or there is $1\leq
i\leq r$ such that $\gamma(a)\neq\gamma\gamma(\alpha_i)$. As $b\in
S^{-\lambda}$ and $\gamma(a)\in\delta(b)$ the former is
impossible. Fix minimal $i_0$ such that
$\gamma(a)\neq\gamma\gamma(\alpha_{i_0})$. Then, by Lemma
\ref{fact5}.2, $\gamma(a)\in\delta\dot{\delta}(\alpha_{i_0})$. As
$\gamma(\alpha_{i_0-1})\in S^{-\lambda}$ (or $a\in S^{-\lambda}$
if $i_0=1$) using Lemma \ref{fact4} we get
$\gamma(a)\in\iota(\alpha_{i_0})\subseteq\iota(\alpha_r)$. On the
other hand $\gamma(a)\in\delta(a)\subseteq\delta\gamma(\alpha_r)$
and $\delta\gamma(\alpha_r)\cap\iota(\alpha_r)=\emptyset$. But
this contradicts Lemma \ref{fact5}.4. Thus $a<^+b$ cannot hold.

Now suppose that $b<^+a$, i.e. there is a flat upper path
$b,\beta_1,\ldots,\beta_r,b$. As $a\in S^{-\lambda}$, by Lemma
\ref{fact3}, $\gamma(\beta_i)\in S^{-\lambda}$, for
$i=1,\ldots,r$. Now either $\gamma(a)\in\delta\gamma(\beta_r)$ or
there is $1\leq i\leq r$ such that
$\gamma(a)\in\delta\gamma(\beta_i)$. As $a\in S^{-\lambda}$ and
$\gamma(\beta_r)=a$ the former is impossible. Fix minimal $i_1$
such that $\gamma(a)\not\in\delta\gamma(\beta_{i_1})$. By Lemma
\ref{fact5}.3 we have
$\gamma(a)\in\gamma\dot{\delta}(\beta_{i_1})$. As $b\in
S^{-\lambda}$, if $i_1>1$ then $\gamma(\beta_{i_1-1})\in
S^{-\lambda}$, as well. Thus
$\gamma(a)\in\iota(\beta_{i_1})\subseteq\iota(\beta_r)$. But
$\gamma(a)=\gamma\gamma(\beta_r)$.  But this contradicts Lemma
\ref{fact5}.4. again and hence $b<^+a$ cannot hold either.

Therefore $a\not\perp^+b$ and then $a<^\sim b$.

Ad 2. If $a<^-b$ then we have a lower path $a=a_0,a_1,\ldots,
a_k,a_{k+1}=a'$, with $k\geq 0$, such that $a_1,\ldots, a_k$ is
flat. If $k=0$ then $\gamma(a)\in\delta(a')$ and we are done. We
shall show, using $\theta(a)\cap\theta(a')\neq\emptyset$, that
$k>0$ is impossible.

If $\gamma(a)=\gamma(a')$ then
$\gamma(a),a_1,\ldots,a_k,(a_{k+1}),\gamma(a')$ is a flat upper
path, where the face $(a_{k+1})$ in parenthesis $()$  is optional
i.e. it is in the path iff it is not a loop. If $k>0$ then $<^+$
is not strict.

If $x\in\delta(a)\cap\delta(a')$ then we have a flat upper
$x,(a_0),a_1,\ldots,a_k,\gamma(a_k)$, with $a_0$ optional. If
$k>0$ then $x<^+\gamma(a_k)\in\delta(a')$. If $x=\gamma(a_k)$ then
$<^+$ is not strict and if $x\neq\gamma(a_k)$ then, as
$x,\gamma(a_k)\in\delta(a')$, we get contradiction with local
discreetness.

If $\gamma(a')=\delta(a)$ then
$\gamma(a'),(a_0),a_1,\ldots,a_k,(a_{k+1}),\gamma(a')$ is a flat
upper path, and again if $k>0$, we get contradiction with
strictness of $<^+$.
 $~\Box$

 \begin{lemma}
 \label{fact9} Let $S$ be an ordered face structure $\alpha,a,b,x\in S$.
\begin{enumerate}
  \item If $\theta(a)\cap\iota(\alpha)\neq\emptyset$ then
  $a<^+\gamma(\alpha)$.

\item If $\alpha\in S-\gamma(S^{-\lambda})$ and $x\in\theta(a)\cap\iota(\alpha)$
  then there is $b\in \dot{\delta}(\alpha)$ such that $a\leq^+b$ and $x\in\theta(b)$.

\end{enumerate}
\end{lemma}
{\em Proof.} 2. can be easily deduced from the proof of 1.

 Ad 1.  First we show that if $x\in\iota(\alpha)$ then there is
$\alpha'\leq^+\alpha$ such that $x\in\iota(\alpha')$ and
$\alpha'\in S-\gamma(S^{-\lambda})$. Take as $\alpha'$ a
$+$-minimal face such that $x\in\iota(\alpha')$ and
$\alpha'\leq^+\alpha$. If $\alpha'\in\gamma(S^{-\lambda})$ then
there is $\xi\in S^{-\lambda}$ such that $\gamma(\xi)=\alpha'$. As
$x\in\iota(\alpha')=\iota\gamma(\xi)=\iota\delta(\xi)$ there is
$\alpha''\in\delta(\xi)$ such that $x\in\iota(\alpha'')$. As
$\alpha''<^+\alpha'$, $\alpha'$ was not minimal contrary to the
supposition.  Thus $\alpha'\in S-\gamma(S^{-\lambda})$.

Next we show that it is enough to show 1. in case $\alpha\in
S-\gamma(S^{-\lambda})$. Suppose that
$x\in\theta(a)\cap\iota(\alpha)$ for some $x\in S$. By the above
there is $\alpha'\leq^+\alpha$ such that $x\in\iota(\alpha')$ and
$\alpha'\in S-\gamma(S^{-\lambda})$. So by Lemma \ref{fact2}.7 and
the above we have $a<^+\gamma(\alpha')\leq^+\gamma(\alpha)$, as
required.

So assume that $\alpha\in S-\gamma(S^{-\lambda})$ and
$x\in\theta(a)\cap\iota(\alpha)$ for some $x\in S$. We consider
three cases:
\begin{enumerate}
  \item $a\in S^{-\lambda}$ and $\gamma(a)=x\in\iota(\alpha)$;
  \item $a\in S^{-\lambda}$ and
  $x\in\delta(a)\cap\iota(\alpha)$;
  \item $a\in S^\lambda$.
\end{enumerate}
We fix $b,c\in\delta^{-\lambda}(\alpha)$ such that
$\gamma(b)=x\in\delta(c)$ for the rest of the argument.

Case 1.  By Lemma \ref{fact2}.1 we have $a\perp^+b$ or $a=b$. If
$a\leq^+b$ then we are done. So we need to show that $b\not<^+a$.
Suppose contrary that $b<^+a$ and $b,\beta_1,\ldots,\beta_r,a$ is
a flat upper $(S-\gamma(S^{-\lambda}))$-path from $b$ to $a$.  As
$b\in\delta(\beta_1)\cap\delta(\alpha)$ and $\beta_1,\alpha\in
S-\gamma(S^{-\lambda})$ we have $\beta_1=\alpha$. Thus
$x=\gamma(a)\in\iota(\alpha)=\iota(\beta_1)\subseteq\iota(\beta_r)$
and $x=\gamma(a)=\gamma\gamma(\beta_r)\not\in\iota(\beta_r)$ and
we get a contradiction.

Case 2. Again by Lemma \ref{fact2}.1 we get that $a\perp^+c$ or
$a=c$. If $a\leq^+c$ we are done. We shall show that $c\not<^+a$.
Suppose contrary that $c<^+a$ and $c,\beta_1,\ldots,\beta_r,a$ is
a flat upper $(S-\gamma(S^{-\lambda}))$-path from $c$ to $a$.  As
$c\in\delta(\beta_1)\cap\delta(\alpha)$ and $\beta_1,\alpha\in
S-\gamma(S^{-\lambda})$ we have $\beta_1=\alpha$. Hence
$x=\gamma(a)\in\iota(\alpha)=\iota(\beta_1)\subseteq\iota(\beta_r)$
and $x\in\delta(a)=\delta\gamma(\beta_r)\not\in\iota(\beta_r)$ and
we get a contradiction again.

Case 3. By loop-filling, pencil linearity, and strictness we have
a flat $(S-\gamma(S^{-\lambda}))$-path $\alpha_0,\ldots,\alpha_k$
ending at $a$ such that $\alpha_0\in S^\varepsilon$. As $a\in
S^\lambda$ we have $\gamma(\alpha_i)\in S^\lambda$ ,for
$i=0,\ldots,k$, and $\gamma\gamma(\alpha_0)=\gamma(a)$. Thus by
pencil linearity we have either $\alpha_0<^\sim\alpha$ or
$\alpha_0<^+\alpha$. As $\alpha_0,\alpha\in
S-\gamma(S^{-\lambda})$ the later is impossible. It remains to
show that if $\alpha_0<^\sim\alpha$ then $a<^+\gamma(\alpha)$. Let
$\alpha_0,\beta_1,\ldots,\beta_r=\alpha$ be a flat lower
$S-\gamma(S^{-\lambda})$-path.  Since $\alpha_i,\beta_j\in
S-\gamma(S^{-\lambda})$ we have $\alpha_i=\beta_i$ for
$i=1,\ldots,\min(k,r)$. If $r\leq k$ then $\gamma(\alpha)\leq^+a$.
But $\gamma(\alpha)\not\in S^\lambda\ni a$. This is a
contradiction with Lemma \ref{fact3}.1. So $k<r$ and hence
$a=\gamma(\alpha_k)\leq^+\gamma(\beta_r)=\gamma(\alpha)$. Since
$a\in S^\lambda\ni\gamma(\alpha)$, we have in fact that
$a<^+\gamma(a)$, as required. $~\Box$

 \begin{lemma}
 \label{fact9.5} Let $S$ be an ordered face structure $a,b,c,d\in S$.
\begin{enumerate}
  \item If $a<^\sim b <^\sim c$ and $a,c<^+d$ then $b<^+d$.
\end{enumerate}
\end{lemma}
{\em Proof.} 1. is easy. $~\Box$

\subsection*{Global properties}

 Let $S$ be an ordered face structure, $n,i \in\o$, $a,a_i\in S_n$, for $i=1,\ldots,k$. The
 {\em weight}\index{weight}\index{face!weight of -} of a face $a$ is the number
 \[ wt(a)=|\{b\in S^{-\lambda}\, : \, b<^+a \}| \]
 The {\em weight of a flat path}  $\vec{a}= a_1,\ldots a_k$ is is
 the sum of weights of its faces
 $wt(\vec{a})=\sum_{i=1}^k wt(a_i)$.

\begin{lemma}
\label{fact7} Let $S$ be an ordered face structure $\alpha,a\in S$, and
$a_1,\ldots,a_k$ flat lower $\dot{\delta}^{-\lambda}(\alpha)$-path
with $k>0$. Then
\[ wt(\gamma(\alpha))>\sum_{i=1}^k wt(a_i). \]
Moreover, $wt(a)=0$ iff $a\in S^{\lambda}$ or
$a\not\in\gamma(S^{-\lambda})$.
\end{lemma}
{\em Proof.} Let $\alpha,a_i$ be as in assumptions of Lemma. Then
for $b\in S$, if $b<^+a_i$ then $b<^+\gamma(\alpha)$. As
$a_i\not<^+a_j$ for any $1\leq i,j\leq k$,  to prove the
inequality we need to show that the faces on the right hand side
are calculated at most once, i.e. for $b\in S^{-\lambda}$, if
$b<^+a_i$ then $b\not<^+a_j$ for $j\neq i$. Suppose contrary, that
$b<^+a_i$ and $b<^+a_j$, with $i<j$. Then, by Lemma \ref{fact2}.7,
$\gamma(b)\leq^+\gamma(a_i)$ and hence $b<^-a_j$. So by Lemma
\ref{fact6} $b<^\sim a_j$. But $b<^+a_j$ and we get a
contradiction with disjointness.

The last statement of Lemma is left for the reader. $~\Box$

\begin{lemma}
\label{fact8} Let $S$ be an ordered face structure,
$X$ convex subset of $S$, $x,y\in X$. If $x<^+y$ then there is a
unique flat upper $(X-\gamma(X^{-\lambda}))$-path from $x$ to $y$.
\end{lemma}

{\em Proof.} Let $X$, $x$, and $y$ be as in assumptions of Lemma.
As $X$ is convex there is a flat upper $X$-path $x,a_1,\ldots ,
a_k,y$. Assume that it is a path with the smallest weight. Suppose
that there is $i\leq k$ such that $a_i\in\gamma(X^{-\lambda})$.
Hence there is $\alpha\in X^{-\lambda}$ such that
$\gamma(\alpha)=a_i$. Let
\[ x'= \left\{ \begin{array}{ll}
           x & \mbox{if  $i=1$,}  \\
           \gamma(a_{i-1}) & \mbox{otherwise.}
                                    \end{array}
                \right. \]
Then
$x'\in\delta(a_i)\cap\delta\gamma(\alpha)\subseteq\delta\delta(a_i)$.
By Lemma \ref{fact2}.3 there is a flat upper
$\delta^{-\lambda}(\alpha)$-path $x',b_1,\ldots,b_r,\gamma(a_i)$.
By Lemma \ref{fact7}, $wt(b_1,\ldots,b_r)<wt(a_i)$, and hence
$wt(a_1,\ldots ,
a_{i-1},b_1,\ldots,b_r,a_{i+1},\ldots,a_k)<wt(a_1,\ldots,a_k)$,
contrary to the supposition that the weight of the path
$x,a_1,\ldots , a_k,y$ is minimal. Thus $x,a_1,\ldots , a_k,y$ is
a $(X-\gamma(X^{-\lambda}))$-path, as required.  The uniqueness of
the path follows from Lemma \ref{fact6} and pencil linearity.
 $~\Box$

A lower flat path $a_0,\ldots, a_k$ is a {\em maximal
path}\index{path!maximal} if
$\delta(a_1)\subseteq\delta(S)-\gamma(S^{-\lambda})$ and
$\gamma(a_k)\in\gamma(S)-\delta(S^{-\lambda})$, i.e. if it can't
be extended either way.

\begin{lemma}[Path Lemma]
\label{fact10} Let $k\geq 0$, $a_0,\ldots, a_k $ be a
maximal lower flat path in an ordered face structure $S$,  $b\in
S$, $0\leq s\leq k$, $a_s<^+b$. Then there are  $0\leq l\leq s\leq
p \leq k$ such that
\begin{enumerate}
   \item $a_i<^+b$ for $i= l, \ldots , p$;
   \item $\gamma(a_p)=\gamma(b)$;
   \item either $l>0$ and $\gamma(a_{l-1})\in\delta(b)$ \\
   or $l=0$ and either $a_0\in S^\varepsilon$ and
   $\gamma\gamma(a_0)\in\theta\delta(b)$
    or $a_0\in S^{-\varepsilon}$ and
   $\delta(a_0)\subseteq\delta(b)$;
    \item $a_i<^\sim b<^\sim a_j$, for $i= 1, \ldots , l-1$ and $j= p+1, \ldots , k$;
  \item $\gamma(a_i)\in\iota(S)$, for $l\leq i <p$.
\end{enumerate}
\end{lemma}
 {\em Proof.} We put
\[ l=\min\{l'\leq s: \forall_{l'\leq i\leq s}\; a_i<^+b \} \hskip 1cm
p=\max\{p'\geq s: \forall_{s\leq i\leq p'}\; a_i<^+b \}.  \] Then
1. holds by definition.

Ad 2. Suppose contrary that $\gamma(a_p)\neq\gamma(b)$. Let
$a_p,\beta_0,\ldots\beta_r,b$ be a flat upper path from $a_p$ to
$b$. As $a\in S^{-\lambda}$ we have $\gamma(\beta_r)\in
S^{-\lambda}$ for $i=1,\ldots,r$. Let
$i_0=\min\{i:\gamma(a_p)\neq\gamma(\beta_i) \}$. Then
$\gamma(a_p)\in\iota(\beta_{i_0})$ and hence, by maximality of the
path, $p< k$ and
$\delta(a_{p+1})\cap\iota(\beta_{i_0})\neq\emptyset$. Thus, by
Lemma \ref{fact9},
$a_{p+1}<^+\gamma(\beta_{i_0})\leq\gamma(\beta_r)=b$ contrary to
the definition of $p$.

Ad 3.  Let $a_l,\beta_1,\ldots,\beta_r,b$ be a flat upper path. We
consider two cases: $l>0$ and $l=0$.

Case $l>0$. Suppose contrary that
$\gamma(a_{l-1})\not\in\delta\gamma(\beta_r)$.  Let
$i_1=\min\{i:\gamma(a_{l-1})\not\in\delta\gamma(\beta_i) \}$.
Then $\gamma(a_{l-1})\in\iota(\beta_{i_1})$ and hence, by Lemma
\ref{fact9},
$a_{l-1}<^+\gamma(\beta_{i_1})\leq^+\gamma(\beta_r)=b$ contrary to
the definition of $l$.

Case $l=0$. If $a_0\in S^\varepsilon$ then, using Lemma
\ref{fact5}, we have
\[ \gamma\gamma(a_0)\in
\theta\theta\delta(\beta_0)\subseteq\theta\theta\gamma(\beta_r)=\theta\theta(\alpha)=\theta\delta(\alpha),\]
as required in this case.

So now assume that $a_0\in S^{-\varepsilon}$.  As, by maximality
of the path, there is no face $a\in S^{-\lambda}$ such that
$\gamma(a)\in\delta(a_0)$, we have
$\delta(a_0)\cap\gamma\dot{\delta}^{-\lambda}(\beta_i)=\emptyset$
for $i=1,\ldots,r$. Clearly
$\delta(a_0)\subseteq\delta\delta(\beta_0)$. Suppose that
$\delta(a_0)\subseteq\delta\delta(\beta_i)$ with $i\leq r$.  Then

\[ \delta(a_0)\subseteq\delta\delta(\beta_i)-\gamma\dot{\delta}^{-\lambda}(\beta_i)=
\delta\gamma(\beta_i)\subseteq \delta\delta(\beta_{i+1}) \]
 (last $\subseteq$ make sense only for $i<r$). Thus
 $\delta(a_0)\subseteq\delta\gamma(\beta_r)=\delta(b)$, as
 required.

4. follows easily from Lemma \ref{fact3}.7.

Ad 5. Fix $l\leq i\leq p$. Let $a_i,\beta_1,\ldots,\beta_r,b$ be a
flat upper path.  If we were to have
$\gamma(a_i)\not\in\bigcup_{i=0}^r\iota(\beta_r)$, by an argument
similar as in 2., we would have $\gamma(a_i)=\gamma(b)$
contradicting strictness.
 $~\Box$

\begin{lemma}[Second Path Lemma]
\label{fact11} Let $k\geq 0$, $a_0,\ldots, a_k $ be
 a flat lower path in an ordered face structure $S$, $x\in\delta(a_0)-\gamma(S^{-\lambda})$,   $b\in S$,
 $a_k<^+b$. Then either $x\in \delta(b)$ or there is $0\leq i <k$, such that $\gamma(a_i)\in\delta(b)$,
 and hence $x \leq^+ y$ for some $y\in\delta(b)$,  ($y=\gamma(a_i$)).
\end{lemma}
{\em Proof.} This is an easy consequence of Path Lemma. $~\Box$

\subsection*{Convex sets}

 Let $S$ be an ordered face structure, $n,i \in\o$, $a,a_i\in S_n$, for $i=1,\ldots,k$. The
 {\em height}\index{height}\index{face!height of -} of a face $a$ in $S$ is the length of the longest flat upper
 $(S-\gamma(S^{-\lambda}))$-path ending at $a$. The height of $a$ is
 denoted by $ht_S(a)$ or if it does not lead to confusions by $ht(a)$.
 The {\em height of a flat path}  $\vec{a}= a_1,\ldots a_k$ is is
 the sum of heights of its faces
 $ht(\vec{a})=\sum_{i=1}^k ht(a_i)$.

 The {\em depth}\index{depth}\index{face!depth of -} of a face $a$ in $S$ is the length of the longest flat upper
 $(S-\gamma(S^{-\lambda}))$-path starting from $a$. The depth of $a$ is
 denoted by $dh_S(a)$ or if it does not lead to confusions by $dh(a)$.
 The {\em depth of a flat path}  $\vec{a}= a_1,\ldots a_k$ is is
 the sum of depths of its faces
 $dh(\vec{a})=\sum_{i=1}^k dh(a_i)$.

 Let $X$ be a subhypergraph of $S$. We say that $X$ is a {\em
convex subset}\index{hypergraph!convex sub-}\index{convex subset}
in $S$ if it is non-empty and the relation $<^{X,+}$ is the
restriction of $<^{S,+}$ to $X$.

Let $X$ be a convex subset of $S$, $a\in X$. The $X$-{\em
depth}\index{depth}\index{face!depth of -} of a face $a$ is the
length of the longest flat upper  $(X-\gamma(X^{-\lambda}))$-path
starting from $a$. The $X$-depth of $a$ is  denoted by $dh_X(a)$.
If $X=S$ and it does not lead to confusions we write $dh(a)$. The
$X$-{\em depth of a flat path}  $\vec{a}= a_1,\ldots a_k$ is is
the sum of depths of its faces $dh_X(\vec{a})=\sum_{i=1}^k
dh_X(a_i)$.

\begin{lemma}\label{convex} Let $X$ be a convex subset of an ordered face structure
$S$.  Then $X$ satisfy  all the axioms of ordered face structures
but loop-filling, where as $<^{X,\sim}$ we take $<^{S,\sim}$
restricted to $X$.
\end{lemma}
{\em Proof.} The only fact that needs a comment is that if
$a<^{X,\sim}b$ then $a<^{X,-}b$.  But this follows from the
observation that $a=a_0,a_1,\ldots,a_{k-1},a_k=b$, ($k>0$) is a
lower path iff $\gamma(a_0),a_1,\ldots,a_{k-1},\gamma(a_{k-1})$ is
a (possibly empty) upper path. $~\Box$

\begin{lemma}
\label{fact12} Let $X$ be a convex subset of an ordered face structure $S$,
and  $\alpha,a,b\in X$.
\begin{enumerate}
  \item $dh_X(a)=0$ iff $a\not\in\dot{\delta}^{-\lambda}(X)$.
  \item If $\alpha\in X^{-\lambda}-\gamma(X^{-\lambda})$ and
  $b\in\dot{\delta}(\alpha)$ then $dh_X(b)=dh_X(\gamma(\alpha))+1$.
  \item $a<^+b$ then $dh_X(a)>dh_X(b))$.
\end{enumerate}
\end{lemma}
{\em Proof.} Easy. $~\Box$

\begin{lemma}
\label{fact13}
Let $S$ be an ordered face structure, $X$ convex subset of $S$,
$x,y\in X-\iota(X)$ and $x<^+y$. Then there is a flat upper
$(X-\delta(X^{-\lambda}))$-path from $x$ to $y$.
\end{lemma}
{\em Proof.} Let $x,y\in X-\iota(X)$ and $x<^+y$.  Let
$x,a_1,\ldots,a_m,y$ be a flat upper $X$-path of least $X$-depth.
Suppose that it is not $(X-\delta(X^{-\lambda}))$-path. Thus by
Lemma \ref{fact12} $dh_X(a_1,\ldots,a_m)>0$.  Pick $a_s$ of
maximal $X$-depth in $a_1,\ldots,a_m$. Let $\alpha\in
X-\gamma(X^{-\lambda})$ such that $a_s\in\delta(\alpha)$. Let
\[ l=\min\{l'\leq s: \forall_{l'\leq i\leq s}\; a_i\in\delta(\alpha) \} \hskip 1cm
p=\max\{p'\geq s: \forall_{s\leq i\leq p'}\; a_i\in\delta(\alpha)
\}.  \]  Since $x,y\in X-\iota(X)$, by an argument similar to the
one given in Path Lemma, we get that
$\gamma(a_p)=\gamma\gamma(\alpha)$ and with
\[ x'= \left\{ \begin{array}{ll}
           x & \mbox{if  $l=1$,}  \\
           \gamma(a_{l-1}) & \mbox{otherwise.}
                                    \end{array}
                \right. \]
$x'\in \delta\gamma(\alpha)$.  Thus
$x,a_1,\ldots,a_{l-1},\gamma(\alpha),a_{p+1}\ldots,a_m,y$ is a
path of a smaller $X$-depth then $x,a_1,\ldots,a_m,y$, contrary to
the assumption. Therefore $x,a_1,\ldots,a_m,y$ is in fact a
$(X-\delta(X^{-\lambda}))$-path. $~\Box$

\subsection*{Order}

\begin{lemma}\label{le_linearity}
Let $S$ be an ordered face structure, $a\in S$. Then the set
\[ \{ b\in S\, :\, a\leq^+ b   \} \]
is linearly ordered by $\leq^+ $.
\end{lemma}
{\em Proof.} Let $a,\alpha_1,\ldots,\alpha_k$ be a maximal flat
upper $S-\gamma(S^{-\lambda})$-path starting from $a$. Then the
set
\[ \{ a \}\cup\{ \gamma(\alpha_i):i=1,\ldots,k\} = \{ b\in S\, :\, a\leq^+ b
\}\] is obviously linearly ordered. $~\Box$

\begin{lemma}
\label{fact15.5} Let $S$ be an ordered face structure $a,b,c\in S$.
\begin{enumerate}
  \item If $a<^+b$ and $b<^\sim c$ then $a<^\sim c$. 
  \item If $a<^\sim b$ and $b<^+ c$ then either $a<^\sim c$ or $a<^+ c$.
\end{enumerate}
\end{lemma}
{\em Proof.} Ad 1. Assume $a<^+b$ and $b<^\sim c$.  Let
$a,\alpha_1,\ldots,\alpha_k,b$ be a flat upper path from $a$ to
$b$ and  $b=b_0,b_1,\ldots, b_l=c$ a lower path from from $b$ to
$c$. Using Lemma \ref{fact2} we get a flat upper
$\bigcup_i\delta(\alpha_i)$-path
$\gamma(a),a_1,\ldots,a_r\gamma(b)$. Thus we have a lower path
$a,a_1,\ldots,a_r,b_1,\ldots,b_l=c$, i.e. $a<^-c$. If
$\theta(a)\cap\theta(c)=\emptyset$ then clearly $a<^\sim c$.

If $\theta(a)\cap\theta(c)\neq\emptyset$ then by pencil linearity
we have $a\perp^+c$ or $a\perp^\sim c$. We shall show that the
only condition that does not lead to a contradiction is $a<^\sim
c$.

If $a<^+c$, then, as $a<^+b$, by Lemma \ref{le_linearity} we have
$b\perp^+c$. Contradiction.

If $c<^+a$, then, as $a<^+b$, we have $c\perp^+b$. Contradiction.

If $c<^\sim a$, then, as $b<^\sim c$, we have $b\perp^\sim a$.
Contradiction.

Ad 2. Assume  $a<^\sim b$ and $b<^+ c$. First note that  by an
argument as above we can show that if $a$ and $c$ are comparable
at all then either $a<^\sim c$ or $a<^+ c$. Thus it is enough to
show that $a$ and $c$ are comparable.  Let $a,a_1,\ldots,a_k,b$ be
a lower path with $k\geq 0$.  We can assume that $a_1,\ldots,a_k$
is a flat path. As, $b<^+c$ by Path Lemma either $a_0<^\sim b$ or
$a_0<^+b$.  In the former case we have
$\theta(a)\cap\theta(c)=\emptyset$ and that $a<^-c$.  Thus
$a<^\sim c$.  In the later either $\gamma(a)\in\delta(c)$ and we
are done or there is a flat upper path
$a,\alpha_1,\ldots,\alpha_r,c$ and $i\leq r$ such that
$\gamma(a)\in\iota(\alpha_i)$. Then by Lemma \ref{fact9} we have
$a<^+\gamma(\alpha_i)\leq^+\gamma(\alpha_r)=c$, as required.
$~\Box$

\begin{lemma}
\label{tech_lemma3} Let $S$ be an ordered face structure, $a,b\in
S$. Then we have
\begin{enumerate}
 \item\label{tl4} If $a<^+b$ then $\gamma(a)\leq^+\gamma(b)$;
  \item\label{tl5} If $a<^\sim b$ then $\gamma(a)\leq ^+\gamma(b)$;
  \item\label{tl1} If $\gamma(a)=\gamma(b)$ then either $a=b$ or
  $a\perp^+b$ or $a\perp^\sim b$;
  \item\label{tl2} If $\gamma(a)<^+\gamma(b)$ then either $a<^+b$ or
  $a<^\sim b$;
  \item\label{tl3} If $\gamma(a)\perp^\sim \gamma(b)$ then  $a\not\perp^\sim b$ and
  $a\not\perp^+b$.
\end{enumerate}
\end{lemma}
{\em Proof.}~ 1. is repeated from Lemma \ref{fact2}.7.

Ad 2. If $a<^\sim b$ then there is a lower path
$a=a_0,a_1,\ldots,a_m=b$. Hence
$\gamma(a),a_1,\ldots,a_m,\gamma(b)$ is an upper path.  So either
$\gamma(a)=\gamma(b)$ or after dropping loops from the sequence
$a_1,\ldots,a_m$ we get a flat upper path from $\gamma(a)$ to
$\gamma(b)$, as required.

Ad 3. This is an immediate consequence of pencil linearity.

Ad 4.  Suppose that $\gamma(a)<^+\gamma(b)$. If
$\theta(a)\cap\theta(b)\neq\emptyset$ then the thesis is obvious.
So assume that $\theta(a)\cap\theta(b)=\emptyset$. Thus, by
disjointness, it is enough to show that either $a<^+b$ or $a<^-b$.
There is a flat upper path $\gamma(a),a_1,\ldots,a_m,\gamma(b)$,
with $m\geq 1$.

Now we argue by cases. If $b$ is a loop the clearly $a<^-b$.
Similarly, if $b=a_m$ then $m>1$ and hence $a<^-b$. Finally,
assume that $a_m<^+b$. If $a_1<^\sim b$ then $a<^b$. If $a_1<^+ b$
then we have a flat path $a_1,\alpha_1,\ldots,\alpha_r,b$. Using
our assumptions  we find $i\leq r$ such that
$\gamma(a)\in\iota(\alpha_i)$. Then by Lemma \ref{fact9}.1 we get
that $a<^+\gamma(\alpha_i)\leq^+b$, as required.

Ad 5. It is an immediate consequence of \ref{tl4}., \ref{tl5}. and
disjointness. $\Box$

\begin{proposition}\label{comparison} Let $S$ be an ordered face
structure, $a,b\in S$. Let $\{ a_i\}_{0\leq i\leq n}$, $\{
b_i\}_{0\leq i\leq n}$ be two sequences of codomains of $a$ and
$b$, respectively, so that
\[ a_i=\gamma^{(i)}(a),  \hskip 20mm b_i=\gamma^{(i)}(b)\]
(i.e. $dim(a_i)=i$), for $i=0,\ldots ,n$. Then, there are two
numbers $l$ and $k$ such that $0\leq l\leq k \leq n$, $1\leq k$
and either
\begin{enumerate}
  \item $a_i=b_i$ for $i<l$,
  \item $a_i<^+b_i$ for $l\leq i < k$,
  \item $a_i<^\sim b_i$ for $k = i \leq n$,
  \item $a_i \not\perp b_i$ for $k < i \leq n$,
\end{enumerate}
or 1.-4. holds with the roles of $a$ and $b$ interchanged.
\end{proposition}
{\em Proof.}~ The above conditions  we can present more visually
as:
\[ a_0=b_0,\ldots , a_{l-1}=b_{l-1},a_{l}<^+b_{l}, \ldots  a_{k}<^+b_{k}, \]
\[  a_{k+1}<^-b_{k+1}, a_{k+2}\not\perp b_{k+2}, \ldots , a_{n}\not\perp b_{n}. \]
These conditions we will verify from the bottom up. Note that by
strictness $<^{S_0,+}$ is a linear order. So either $a_0=b_0$ or
$a_0\perp^+b_0$. In the later case $l=0$. As $a\neq b$ then there
is $i\leq n$ such that $a_i\neq b_i$. Let $l$ be minimal such,
i.e.  $l=\min\{ i : a_i\neq b_i \}$. By Lemma
\ref{tech_lemma3}.\ref{tl1}, $a_l\perp^+b_l$ or $a_l\perp^\sim
b_l$. We put $k = \max \{ i\leq n : a_i\perp ^+b_i\; {\rm or }\;
i=l\}$. If $k=n$ we are done. If $k<n$ then by Lemma
\ref{tech_lemma3}.\ref{tl2}, we have $a_{k+1}\perp^\sim b_{k+1}$.
Then if $k+1<n$, by Lemma \ref{tech_lemma3}.\ref{tl3}, $a_i
\not\perp b_i$ for $k+2\leq i \leq n$. Finally, by Lemma
\ref{tech_lemma3}.\ref{tl4} and .\ref{tl5} all the inequalities
head in the same direction. This ends the proof. $\Box$

For $a,b\in S$ we define
 $a<^\sim_lb$  iff   $\gamma^{(l)}(a)<^\sim \gamma^{(l)}(b)$.

\begin{corollary}\label{order}
 Let $S$ be an ordered face structure, $a,b\in S_n$,
$a\neq b$. Then either $a\perp^+b$ or there is a unique $0\leq
l\leq n$ such that $a\perp^\sim_lb$, but not both.
\end{corollary}

The above Corollary allows us to define an  order $<^S$ (also
denoted $<$) on all cells of $S$ as follows. For $a,b\in S_n$,

\[ a<^S b \;\;\; {\rm iff}\;\;\; a<^+b \;\; {\rm or}\;\; \exists_l \;\; a<^\sim_l b. \]

\begin{corollary}
\label{linord} For any ordered face structure $S$, and $k\in\o$,
the relation $<^S$ restricted to $S_k$ is a linear order.
\end{corollary}
{\em Proof.}~ In the proof we use Lemma \ref{fact15.5} without
mention. We need to verify that $<^S$ is transitive. Let $a,b,c\in
S_n$, $l,k\leq n\in\o$.  We argue by cases.

Case $a<^+b<^+c$. Then by transitivity of $<^+$ we have $a<^+c$.

Case $a<^+b<^\sim_lc$. Then
$\gamma^{(l)}(a)<^+\gamma^{(l)}(b)<^\sim\gamma^{(l)}(c)$.
Therefore $\gamma^{(l)}(a)<^\sim\gamma^{(l)}(c)$ and hence
$a<^\sim_lc$.

Case $a<^\sim_lb<^+c$. Then
$\gamma^{(l)}(a)<^\sim\gamma^{(l)}(b)<^+\gamma^{(l)}(c)$. Thus
either $\gamma^{(l)}(a)<^\sim\gamma^{(l)}(c)$ and hence
$a<^\sim_lc$ or $\gamma^{(l)}(a)<^+\gamma^{(l)}(c)$. In the later
case by Proposition \ref{comparison} we have either $a<^+c$ or
there is $l'$ such that $l\leq l'\leq n$ and $a<^\sim_{l'}c$.

Case $a<^\sim_kb<^\sim_lc$. If $k=l$ then by transitivity of
$<^\sim$ we have $a<^\sim_kc$.

If $k>l$ then
$\gamma^{(l)}(a)\leq^+\gamma^{(l)}(b)<^\sim\gamma^{(l)}(c)$.
Therefore $\gamma^{(l)}(a)<^\sim\gamma^{(l)}(c)$ and hence
$a<^\sim_lc$.

If $k<l$ then
$\gamma^{(k)}(a)\leq^\sim\gamma^{(k)}(b)<^+\gamma^{(k)}(c)$.
Therefore either $\gamma^{(k)}(a)\leq^\sim\gamma^{(k)}(c)$ and
hence $a<^\sim_kb$ or $\gamma^{(k)}(a)\leq^+\gamma^{(k)}(c)$. In
the later case again by Proposition \ref{comparison} we have
either $a<^+c$ or there is $k'$ such that $k\leq k'\leq n$ and
$a<^\sim_{l'}c$. $\Box$

\subsection*{Monotone morphisms}

From Corollary \ref{linord} we also get

\begin{corollary}
\label{pres and refl of order} Let $f:S\ra T$ be a monotone morphism of
ordered face structures, and $l,k\in\o$, $l\leq k$, $x,y\in S_k$.
Then
\begin{enumerate}
  \item $x<^\sim_l y$ iff $f(x)<^\sim_l f(y)$;
  \item $x\leq^+ y$ iff $f(x)\leq^+ f(y)$.
\end{enumerate}
\end{corollary}
{\em Proof.}~ Obvious. $\Box$

{\em Remak.} Note however that monotone morphisms do not preserve the
relation $<^+$ in general.

\begin{corollary}
\label{monotone bi are iso} Any monotone morphism of ordered face
structures which is a bijection is an isomorphism.
\end{corollary}
{\em Proof.}~ If $f$ is a monotone bijection of ordered face
structures it is clearly a local isomorphism. But by Lemma
\ref{pres and refl of order} it reflects $<^\sim$ as well, i.e. it
is a monotone isomorphism. $\Box$

\begin{lemma}
\label{fiber} Let $f:S\ra T$ be a  monotone morphism of ordered face
structures. If $a\in T$ then if $f^{-1}(a)\neq\emptyset$ there are
$b,c\in f^{-1}(a)$ and a flat upper
$S^u-\gamma(S^{-\lambda})$-path $b,\alpha_1,\ldots,\alpha_r,c$,
with $r\geq 0$, such that $f^{-1}(a)=\{ b \}\cup\{\gamma(\alpha_i)
\}_{1\leq i\leq r}$. In particular all faces in $f^{-1}(a)$ are
parallel to each other and the whole set is linearly ordered by
$<^+$.
\end{lemma}
{\em Proof.} Suppose $b,c\in S$ such that $f(b)=f(c)=a\in T$. Then
for any $l$, $f(b)\not<^\sim_l f(c)$. Therefore for any $l$,
$b\not<^\sim_l c$. Thus, by Lemma \ref{comparison}, $b<^+c$. Hence
There is a flat upper $S-\gamma(S^{-\lambda})$-path
$b,\alpha_1,\ldots,\alpha_r,c$. As $f(b)=f(c)$, we have
$f(\alpha_i)\in T^\lambda$, for $i=1,\ldots,r$. In particular,
$\alpha_i\in S^u$, for $i=1,\ldots,r$. So we have shown that
between any two different elements of $f^{-1}(a)$ there is a flat
upper $S^u-\gamma(S^{-\lambda})$-path. This clearly imply all the
remaining parts of Lemma. $~\Box$

\begin{corollary}\label{mono}
 Any endomorphism of an ordered
face structures is an identity.
\end{corollary}
{\em Proof.}  Let $f:S\ra S$ be a monotone morphism. First note that as
$S_0$ is linearly ordered by $<^+$, if $x,y\in S_0$ then
$ht(x)+dh(x)=ht(y)+dh(y)$, and $x=y$ iff $ht(x)=ht(y)$. As $f$
preserves $<^\sim$, using Lemma \ref{fact6}, we get that
$f_0:S_0\ra S_0$ is an identity.

In order to get a contradiction we suppose that $f$ is not
identity. Let $k$ be the minimal such that $f_k\neq 1_{S_k}$ and
let $a\in S_k$ be $<^S$-minimal such that $f(a)\neq a$. By
minimality of $k$ we have $f(a)\| a$. We shall show that
$f(a)\perp^+ a$. By previous observation and pencil linearity we
have that either $f(a)\perp^+a$ or $f(a)\perp^\sim a$. If $a<^\sim
f(a)$ then we get an infinite sequence $a<^\sim f(a)<^\sim
ff(a)<^\sim\ldots$ contradicting strictness of $<^\sim$. The
condition $f(a)<^\sim a$ cannot hold for the similar reasons.
Suppose that $a<^+f(a)$ (the argument for the case $f(a)<^+a$ is
similar). Let $a,\alpha_1,\ldots,\alpha_r,f(a)$ be a flat upper
path with $r>0$. Then, by Lemma \ref{fact6}, we have
$\alpha_i<^\sim \alpha_{i+1}$ for $i=1,\ldots,r-1$. Hence we get
an infinite sequence
\[ \alpha_1<^\sim f(\alpha_1)<^\sim ff(\alpha_1)<^\sim\ldots \]
again contradicting strictness of $<^\sim$.  Thus $f$ must be an
identity indeed. $~\Box$

\begin{proposition}\label{loc map + are mono} Let $f:S\ra T$
be a local morphism of ordered face structures that preserves $<^\sim$
on sets $S_k-\delta(S_{k+1}^{-\lambda})$, for $k\in\o$. Then $f$
is a monotone morphism. In particular if $S$ is $n$-normal then $f$ is
a monotone morphism iff $f_n:S_n\lra T_n$ preserves $<^\sim$.
\end{proposition}
{\em Proof.} Let $\prec$ denote the relation $<^{S,\sim}$
restricted to such pairs of elements $a,b$ that either
$a,b\in\delta(\alpha)$ for some $\alpha\in S$ or $a,b\in
S-\delta(S^{-\lambda})$. Thus we must show that if a hypergraph
morphism $f:S\ra T$ between ordered face structures preserves
$\prec$, i.e. is such that for any $a,b\in S$ if $a\prec b$ then
$f(a)<^{T,\sim}b$ then it preserves $<^\sim$, i.e. for any $a,b\in
S$ if $a<^{S,\sim}b$ then $f(a)<^{T,\sim}b$.

So fix $f:S\ra T$ preserving $\prec$ and $a,b\in S$ such that
$a<^\sim b$. Then by disjointness and few other facts, there is a
lower path $a=a_0,\ldots, a_k=b$ such that $a_i<^\sim a_{i+1}$,
for $i=0,\ldots,k-1$. By transitivity of $<^\sim$ it is enough to
show that $f(a)<^\sim f(b)$ only in case $\gamma(a)\in\delta(b)$.
We shall prove by induction of the sum of depth of $a$ and $b$,
$s=dh(a)+dh(b)$ that if $\gamma(a)\in\delta(b)$ and $a<^\sim b$
then $f(a)<^\sim f(b)$.

If $s=0$ then $a,b\in S-\delta(S^{-\lambda})$ and hence, by
assumption on $f$, $f(a)<^\sim f(b)$.

So assume that $s>0$ and that for $s'$ smaller than $s$ the
inductive hypothesis holds. We consider two cases.

Case $dh(a)\geq dh(b)$. So we have $\alpha\in
S^{-\lambda}-\gamma(S^{-\lambda})$ such that $a\in\delta(\alpha)$.
Hence $\delta(b)\cup\theta\delta(\alpha)\neq\emptyset$.  If
$\delta(b)\cap\iota(\alpha)\neq\emptyset$ then, by Lemma
\ref{fact9}, $b<^+\gamma(\alpha)$. If
$\delta(b)\cap\theta\gamma(\alpha)\neq\emptyset$ then either
$b<^+\gamma(\alpha)$ or $\gamma(\alpha)<^\sim b$, as two other
cases easily lead to a contradiction. If $b<^\sim\gamma(\alpha)$
then $a<^\sim\gamma(\alpha)$, and if $\gamma(\alpha)\leq^+ b$ then
$a<^+b$. So in both cases we get a contradiction. Now if
$b<^+\gamma(\alpha)$, as $\alpha\in
S^{-\lambda}-\gamma(S^{-\lambda})$ and $dh(a)\geq dh(b)$ we have
$b\in\delta(\alpha)$. Hence $a\prec b$ and then by assumption on
$f$ we get $f(a)<^\sim f(b)$, as required. If
$\gamma(\alpha)<^\sim b$ then, as
$\gamma\gamma(\alpha)\in\delta(b)$ by strictness and induction
hypothesis, we have $f(a)\leq^+f(\gamma(\alpha))<^\sim f(b)$. So
by Lemma \ref{fact15.5} $f(a)<^\sim f(b)$, as well.

Case $dh(a)< dh(b)$. So we have $\beta\in
S^{-\lambda}-\gamma(S^{-\lambda})$ such that $b\in\delta(\beta)$.
Then $a\not\in\delta(\alpha)$ and $a\not\perp^+\gamma(\beta)$. We
shall show that $\gamma(a)\in\delta\gamma(\beta)$. Clearly we have
$\gamma(a)\in\delta\delta(\beta)$. If we were to have
$\gamma(a)\in\iota(\beta)$ then $a<^+\gamma(\beta)$ and hence
$dh(a)\geq dh(\gamma(\beta))+1\geq dh(b)$ contradicting our
assumption. Therefore $\gamma(a)\in\delta\gamma(\beta)$.  Now, to
get a contradiction, we assume that
$\gamma(a)=\gamma\gamma(\beta)\not\in\delta\gamma(\beta)$. By
$\delta$-globularity we have $c\in\dot{\delta}^{-\lambda}(\beta)$
such that $\gamma(c)=\gamma(a)$. Thus, by pencil linearity, either
$c\perp^\sim a$ or $c\perp^+ a$. As if $a<^+c$ then
$a<^+\gamma(\beta)$, if $c<^+a$ then $b<^+a$, and if $a<^\sim c$
then $\gamma(a)\neq\gamma(c)$ the only non-trivial case, we have
to consider, is $c<^\sim a$. Thus by Lemma \ref{fact9.5} we have
$a<^+\gamma(\beta)$, and we get a contradiction again. Therefore
$\gamma(a)\in\delta\gamma(\beta)$ as claimed.

As $b<^+\gamma(\beta)<^\sim a$ would lead to $b<^\sim a$ and
contradiction, we must have $a<^\sim\gamma(\beta)$. Thus by
induction hypothesis we have $f(a)<^\sim f(\gamma(\beta))$.
Clearly, we also have $f(b)\leq^+f(\gamma(\beta))$. As
$\gamma(a)\in\delta(b)$ we have $\gamma(f(a))\in\delta(f(b))$ and
hence, by pencil linearity, either $f(a)\perp^\sim f(b)$ or
$f(a)\perp^+ f(b)$. We shall show that the only condition that
does not lead to a contradiction is $f(a)<^\sim f(b)$. If $f(a)<^+
f(b)$ then  $f(a)<^+ f(\gamma(\beta))$ and contradiction. If
$f(b)<^\sim f(a)$ then,  by Lemma \ref{le_linearity},
$f(a)\perp^+f(\gamma(\beta))$ and contradiction. If $f(b)<^\sim
f(a)$ then  $f(b)<^\sim f(\gamma(\beta))$ and again we get a
contradiction. Thus $f(a)<^\sim f(b)$, as required. $~\Box$

\begin{corollary}
\label{local to global} The ordered face structure
$S$ is uniquely determined by its local face structure part $|S|$
and the order $<^{\sim}$ restricted to the sets
$S_k-\delta(S^{-\lambda}_{k+1})$, for $k\in \o$.
\end{corollary}

{\em Proof.} Let $S$ and $S'$ be ordered face structures such that
their local parts are equal, i.e. $|S|=|S'|$ and that the relation
$<^{\sim,S}$ agree with $<^{\sim,S'}$ on the set
$S-\delta(S^{-\lambda})$.  Then the identity morphism is preserving
$<^\sim$ on the set $S-\delta(S^{-\lambda})$. Thus by Proposition
\ref{loc map + are mono} it is a monotone morphism, considered as a map
either way, i.e. $S=S'$.
 $~\Box$

In general, there are more local than monotone morphisms between
ordered face structures.  However if we restrict our attention to
the principal ordered face structures those two notions agree. We
have

\begin{corollary}\label{embb} The embedding $|-|:\pfs\lra\lfs$ is full and faithful.
\end{corollary}
{\em Proof.} Fix a local morphism $f:S\lra T$ between ordered face
structures, with $S$ being principal.   Then for $k\in\o$ the sets
$S_k-\delta(S_{k+1}^{-\lambda})$ has at most one element.  So
$<^\sim$ is obviously preserved on these sets, i.e. by Lemma
\ref{loc map + are mono} $f$ is monotone. $~\Box$

The limits and colimits in $\ofs$ are rather rare and in $\lfs$
also do not always exist. For example if we take a local face
structure $S$ such that $S_0=\{ u\}$, $S_1=\{ x,y,z\}$, $S_2=\{
a\}$, and with
$\gamma(x)=\gamma(y)=\gamma(z)=\delta(x)=\delta(y)=\delta(z)=u$,
$\delta(a)=\{ x,y\}$, $\gamma(a)=z$, and $<_a$ being empty
relation, then we have two local isomorphisms from $S$ to itself
identity $1_S$ and a map $\sigma$ switching $x$ and $y$. Clearly
in the coequalizer $Q$ of $1_S$ and $\sigma$ the faces $x$ and $y$
should be identified but then the map $q:S\ra Q$ in to it cannot
be local as there cannot be a  bijection from $\delta(a)$ to
$\delta(q(a))$.  Thus the coequalizer $1_S$ and $\sigma$ does not
exists in $\lfs$. However we have

\begin{proposition}\label{lims ofs} The colimits and connected limits of
diagrams from $\ofs$ exists in $\lfs$ and are calculated in $Set$
in each dimension separately.
\end{proposition}
{\em Proof.} The main property of the monotone morphisms of order
face structures that allow calculations of the above limits and
colimits is the following. If $f,g:S\ra T$ are monotone morphisms and
$a\in S_{>0}$ such that $f(a)=g(a)=b$ then the functions
$f_a,g_a:\delta(a)\ra \delta(b)$ are equal. This is an immediate
consequence of Corollary \ref{pres and refl of order}. $~\Box$

\section{Stretching the convex subhypergraphs}\label{convex subhyp}

From Corollary \ref{pres and refl of order} follows immediately
that the (hypergraph) image of a monotone morphism is a convex subset
of the codomain.  In this section we shall show that the converse
is also true and it is true in an essentially unique way, i.e. if
$X$ is a convex subset of $T$ then there is a monotone morphism
$\nu_X:[X]\ra T$ such that image of $\nu_X$ is $X$, i.e. we can
cover a convex set by an ordered face structure. Moreover, if
$f_i:S_i\ra T$, $i=0,1$, are monotone morphisms such that their images
are equal, $im(f_0)=im(f_1)$, then there is a monotone isomorphism
$g:S_0\ra S_1$ making the triangle
\begin{center}
\xext=600 \yext=350
\begin{picture}(\xext,\yext)(\xoff,\yoff)
\settriparms[1`1`1;300] \putVtriangle(0,0)[S_0`S_1`T;g`f_0`f_1]
\end{picture}
\end{center}
commutes, i.e. the covering is essentially unique. As the title of
the section suggests, the construction of $[X]$ is done by
stretching $X$. The stretching means in this case the splitting
all the empty loops in the convex set $X$.

Let $T$ be an ordered face structure, and $X\subseteq T$ a
subhypergraph. Recall that $X$ is {\em
convex}\index{hypergraph!convex sub-} in $T$ if it is non-empty
and the relation $<^{X,+}$ is the restriction of $<^{T,+}$ to $X$.
For the rest of the section  assume that $X$ is a convex
subhypergraph of $T$. We shall define an ordered face structure
$[X]$ and a monotone morphism $\nu_X:[X]\lra T$. But first we need to
explain what are cuts of empty loops.

We define the set of empty loops in $X$ as
\[ \cE^X=X^\lambda -\gamma(X^{-\lambda}) \]
 and the set of empty loops in $X$ over $x\in X$ as
\[ \cE^X_x=\{ a\in \cE : \gamma(a)=x \}. \]
As $X$ is convex, the relation $<^{T,\sim}$ restricted to
$\cE^X_x$ is a linear order. We say that a triple $(x,L,U)$ is an
{\em $x$-cut}\index{cut} i.e. a cut of $\cE^X_x$ iff  $L\cup
U=\cE^X_x$ and for $l\in L$ and $l'\in U$, $l<^\sim l'$. Note that
both $L$ and $U$ might be empty and hence, there is an $x$-cut for
any $x\in X$ (e.g. $(x,\emptyset, \cE^X_x)$). Let $\cC(\cE^X_x)$
be the set of $x$-cuts.

We need some notation for cuts in $X$. If $(x,L,U)$ is a $x$-cut
then $L$ determines $U$ and vice versa ($L=\cE^X_x-U$ and
$U=\cE^X_x-L$).

 Therefore we sometimes denote this cut by
describing only the lower cut $(x,L,-)$ or only the upper cut
$(x,-,U)$, whichever is easier to define. Let $a$ be an arbitrary
face in $X$,  $y\in\dot{\delta}(a)$. We define two sets
\[ \ua a =  \{ b\in \cE^X_{\gamma(a)} : a<^\sim b\}, \hskip 10mm \da_y a =\{
b\in \cE^X_y : b<^\sim a\}\] that determine the cuts $(x, -, \ua
a)$ and $(y,\da_y a,-)$. We often omit subscript $y$ inside the
cuts, i.e. we usually write $(y,\da a,-)$ when we mean $(y,\da_y
a,-)$. If $\gamma(a)=$ then we sometimes write $\ua_xa$ instead of
$\ua a$ to emphasis that the cut is over$x$.

Now we are ready to define $[X]$. We put for $k\in\o$
\[ [X]_k= \bigcup_{x\in X_k} \cC(\cE^X_x) \]
We put, for $(x,L,U)\in [X]_l$
\[ \gamma(x,L,U)= (\gamma(x),-,\ua x)\]
and
\[ \delta(x,L,U)= \left\{ \begin{array}{ll}
          1_{(\gamma\gamma(x),-,\ua \gamma(x))} & \mbox{if  $\delta(x)=1_{\gamma\gamma(x)}$,}  \\
          \{ (t,\da x,-) : t\in\dot{\delta}(x) \} & \mbox{otherwise.}
                                    \end{array}
                \right. \]

We have a hypergraph map
\[ \nu_X : [X] \lra T \]
such that $\nu_X(x,L,U)=x$, for $(x,L,U)\in [X]$. The order
$<^{[X],\sim}$ is inherited from $T$ via $\nu_X$, i.e.
$(x,L,U)<^{[X],\sim}(x',L',U')$ iff $x<^{T,\sim}x'$, for
$(x,L,U),(x',L',U')\in [X]$.

Let $Z\subseteq T$. By $<Z>$ we mean the least subhypergraph of
$T$ containing $Z$. We call $Z$ {\em convex set} if $<Z>$ is a
convex subhypergraph.  If $Z$ is a convex set we write $\cE^Z$,
$[Z]$ and $\nu_Z$ instead of $\cE^{<Z>}$, $[<Z>]$ and $\nu_{<Z>}$.
Moreover, if $Z=\{\alpha\}$ we write $\cE^{\alpha }$, $[\alpha ]$
and $\nu_\alpha$ instead of $\cE^{\{\alpha\} }$, $[\{\alpha\} ]$
and $\nu_{\{\alpha\} }$.

{\em Example.} An ordered face structure $T$
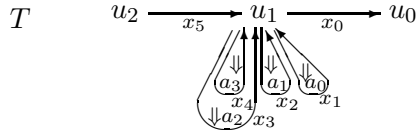
\begin{center} \xext=1320 \yext=520
\begin{picture}(\xext,\yext)(\xoff,\yoff)
\putmorphism(0,450)(1,0)[u_2`u_1`_{x_5}]{480}{1}b
\putmorphism(480,450)(1,0)[\phantom{u_1}`u_0`_{x_0}]{480}{1}b
 \put(-400,400){$T$}

 \put(350,150){\oval(200,200)[b]}
 \put(285,40){$^{\Da a_2}$}
 \put(250,150){\line(1,2){120}}
 \put(450,150){\vector(0,1){250}}
  \put(440,30){$^{x_3}$}

  \put(360,220){\oval(100,100)[b]}
  \put(310,220){\line(1,2){85}}
  \put(410,220){\vector(0,1){180}}
  \put(360,220){$^\Da$}
    \put(325,155){$^{a_3}$}
     \put(365,90){$^{x_4}$}
    \put(520,220){\oval(100,100)[b]}
  \put(470,220){\line(0,1){180}}
  \put(570,220){\vector(-1,2){85}}
  \put(480,220){$^\Da$}
    \put(490,155){$^{a_1}$}
     \put(520,90){$^{x_2}$}
     \put(650,220){\oval(100,100)[b]}
  \put(600,220){\line(-1,2){85}}
  \put(700,220){\vector(-1,1){180}}
  \put(600,190){$^\Da$}
    \put(620,155){$^{a_0}$}
     \put(670,100){$^{x_1}$}
\end{picture}
\end{center}
has a convex subset $X$
\begin{center} \xext=1320 \yext=520
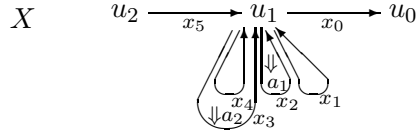
\begin{picture}(\xext,\yext)(\xoff,\yoff)
\putmorphism(0,450)(1,0)[u_2`u_1`_{x_5}]{480}{1}b
\putmorphism(480,450)(1,0)[\phantom{u_1}`u_0`_{x_0}]{480}{1}b
 \put(-400,400){$X$}
 \put(350,150){\oval(200,200)[b]}
 \put(285,40){$^{\Da a_2}$}
 \put(250,150){\line(1,2){120}}
 \put(450,150){\vector(0,1){250}}
  \put(440,30){$^{x_3}$}

  \put(360,220){\oval(100,100)[b]}
  \put(310,220){\line(1,2){85}}
  \put(410,220){\vector(0,1){180}}
  \put(365,90){$^{x_4}$}
    \put(520,220){\oval(100,100)[b]}
  \put(470,220){\line(0,1){180}}
  \put(570,220){\vector(-1,2){85}}
  \put(480,220){$^\Da$}
    \put(490,155){$^{a_1}$}
     \put(520,90){$^{x_2}$}
     \put(650,220){\oval(100,100)[b]}
  \put(600,220){\line(-1,2){85}}
  \put(700,220){\vector(-1,1){180}}
     \put(670,100){$^{x_1}$}

\end{picture}
\end{center}
with $\cE^X=\{ x_4,x_1 \}$ whose stretching is the following
ordered face structure $[X]$:
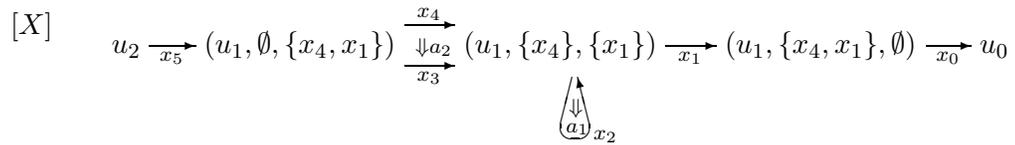
\begin{center} \xext=2820 \yext=520
\begin{picture}(\xext,\yext)(\xoff,\yoff)
 \put(-400,400){$[X]$}
 \putmorphism(0,360)(1,0)[u_2`(u_1,\emptyset,\{x_4,x_1\})`_{x_5}]{600}{1}b
 \putmorphism(600,430)(1,0)[\phantom{(u_1,\emptyset,\{x_4,x_1\})}`\phantom{(u_1,\{x_4\},\{x_1\})}`_{x_4}]{900}{1}a
 \putmorphism(600,360)(1,0)[\phantom{(u_1,\emptyset,\{x_4,x_1\})}`(u_1,\{x_4\},\{ x_1\})`]{900}{0}b
 \putmorphism(600,290)(1,0)[\phantom{(u_1,\emptyset,\{x_4,x_1\})}`\phantom{(u_1,\{x_4\},\{ x_1\})}`^{x_3}]{900}{1}b
  \put(1000,300){$^{\Da a_2}$}
  \putmorphism(1500,360)(1,0)[\phantom{(u_1,\{x_4\},\{ x_1\})}`(u_1,\{x_4,x_1\},\emptyset)`_{x_1}]{900}{1}b
  \putmorphism(2400,360)(1,0)[\phantom{(u_1,\{x_4\},\{ x_1\})}`u_0`_{x_0}]{600}{1}b
  \put(1550,90){\oval(100,100)[b]}
  \put(1500,90){\line(1,4){40}}
  \put(1600,90){\vector(-1,4){40}}
  \put(1530,90){$^\Da$}
    \put(1520,25){$^{a_1}$}
     \put(1615,0){$^{x_2}$}
\end{picture}
\end{center}
We adopt the convention that the empty cut in $[X]$, say
$(x_5,\emptyset,\emptyset)$, is identified with the corresponding
face in $X$, $x_5$ in this case.

\begin{lemma}\label{subfs0} Let $T$ be an ordered face
  structure, $X$ be a convex subset of $T$, $(a,L,U),(a',L',U')\in [X]$. Then
\begin{enumerate}
  \item $(a,L,U)$ is an empty domain face in $[X]$ iff $a$ is an empty domain face in $T$;
  \item $(a,L,U)$ is a loop in $[X]$ iff $a$ is a loop in $T$ and
  there is no empty loop $l\in \cE^X$ such that $l\leq^{T,+}a$;
  \item $(a,L,U)<^+(a',L',U')$ iff $a<^+a'$ or ($a=a'$ and $L\subseteqnot
  L'$);
  \item $(a,L,U)<^-(a',L',U')$ iff $a<^-a'$ and  (if $\gamma(a)\in\delta(a')$ then
  $\cE_{\gamma(a)}=\ua_{\gamma(a)}a\cup\da_{\gamma(a)}a' $);
   \item  If $(a,L,U)\in [X]^{-\varepsilon}$ then
 '$\delta(a,L,U)$ is $\delta(a)$ with colors': if
  $x,y\in\delta(a)$ such that $y$ is $<^\sim$-successor of $x$ in
  $\delta(a)$ then $(y,\da a,-)$ is $<^\sim$-successor of $(x,\da a,-)$ in
  $\delta(a,L,U)$; in particular $\gamma(x,\da a,-)\in\delta(y,\da
  a,-)$.
  \item If $(a,L,U)\in [X]^{\varepsilon}$ then $\delta(a,L,U)=1_{(\gamma\gamma(a),\da
  \gamma(a),-)}$.
\end{enumerate}
\end{lemma}{\it Proof.}~ 1., 5. and 6. are obvious.

Ad 2. If $(a,L,U)$ is a loop in $[X]$ so must be $a$ in $T$. So
fix a cut $(a,L,U)$ in $[X]$ such that $a$ is a loop. Let us
denote $\gamma(a,L,U)=(\gamma(a),-,\ua a)$ and
$\delta(a,L,U)=(\gamma(a),\da a,-)$.

If there is $l\in \cE^X_{\gamma(a)}$ such that $l\leq^+a$, then
$a\not\perp^\sim l$. Hence $l\in L'$ and $l\in U'$. So
$\gamma(a,L,U)\neq\delta(a,L,U)$ and $(a,L,U)$ is not a loop.

If there is no $l\in \cE^X_{\gamma(a)}$ such that $l\leq^+a$ then
any empty loop $l\in \cE^X_{\gamma(a)}$ is $<^\sim$-comparable
with $a$. Thus $\gamma(a,L,U)=(\gamma(a),-,\ua a)=(\gamma(a),\da
a,-)=\delta(a,L,U)$.

Ad 3. Fix $(a,L,U)$, $(a',L',U')$ in $[X]$.

First we shall show that the condition is necessary. Suppose that
$(a,L,U),(\alpha_1,L_1,U_1),\ldots, (\alpha_k,L_k,U_k),
(a',L',U')$ is a flat upper path in $[X]$.

Suppose that $a\neq a'$. Clearly $a,\alpha_1,\ldots,\alpha_k,a'$
is an upper path in $X$. So after deleting loops we get a flat
upper $X$-path from $a$ to $a'$, i.e. $a<^+a'$.

Suppose now that $a=a'$ and $L\neq L'$. As
$(a,L,U)\in\delta(\alpha_1,L_1,U_1)\neq\emptyset$ we have $L=\da_a
\alpha_1$. Moreover, $\gamma(\alpha_k,L_k,U_k)=(a',L',U')$ implies
that $\ua_a\alpha_k=U'$. Since
\[(a,-,\ua_a\alpha_i) =\gamma(\alpha_i,L_i,U_i)\in\delta(\alpha_{i+1},L_{i+1},U_{i+1})=(a,\da_a\alpha_{i+1},-) \]
we have that $(a,\da_a\alpha_{i+1},\ua_a\alpha_i)$ is a cut, for
$i=1,\ldots,k-1$. As $\ua_a\alpha_i\cap\da_a\alpha_i=\emptyset$,
we have that $\da_a \alpha_i\subseteq \da_a \alpha_{i+1}$, for
$i=1,\ldots,k-1$. Thus
\[ L=\da_a\alpha_1\subseteq \da_a\alpha_k\subseteq \cE^X_a-\ua_a\alpha_k=\cE^X_a-U'=L'\]
i.e. the condition is necessary.

Now we shall show that the condition is sufficient. First note
that if $l\in \cE^X_a$ then
\[ \delta(l,\emptyset,-)=(a,\da l,-)=(a,\da l,\{ l \}\cup\ua l),\]
\[ \gamma(l,\emptyset,-)=(a,-,\ua l)=(a,\da l\cup\{ l \},\ua l).\]
Thus if $a=a'$ we have that  $(a,L,U)<^+(a',L',U')$ iff
$L\subseteqnot L'$.

Assume now that $a<^+a'$. Let $a,\alpha_1,\ldots,\alpha_k,a'$ be a
flat upper $X$-path of minimal weight. We claim that it is
$X-\gamma(X^{-\lambda})$-path. Suppose contrary that there is
$A\in X^{-\lambda}$ such that $\gamma(A)=\alpha_i$, for some $i$.
Then we will have a flat upper $\delta(A)$-path
$\beta_1,\ldots,\beta_r$ from $\gamma(\alpha_{i-1})$ (or $a$ if
$i=0$) to $\gamma(\alpha_i)$. Replacing $\alpha_i$ by
$\beta_1,\ldots,\beta_r$ we get a flat upper $X$-path of smaller
weight than $\alpha_1,\ldots,\alpha_k$ contrary to the choice of
this path.  This $\alpha_1,\ldots,\alpha_k$ is a flat upper
$X-\gamma(X^{-\lambda})$-path indeed.

As $\alpha_i\in X-\gamma(X^{-\lambda})$ we have
\[\delta(\alpha_i,\emptyset,-)=\{ (b,\da
\alpha_i,\emptyset) : b\in\delta(\alpha_i) \},\;\;\;\;\;\;
\gamma(\alpha_i,\emptyset,-)=(\gamma(\alpha_i),\emptyset,\ua
\alpha_i). \]
 From this and previous we get that
\[ (a,L,U)\leq^+(a,\da \alpha_1,\emptyset)\in\delta(\alpha_{i+1},\emptyset,-),\;\;\; \]
\[ \gamma(\alpha_i,\emptyset,-)\leq^+(\gamma(\alpha_i),\da
\alpha_{i+1},\emptyset)\in \delta(\alpha_{1},\emptyset,-) \]
\[ \gamma(\alpha_k,\emptyset,-)\leq^+(a',L',U') \]
and this shows that $(a,L,U)<^+(a',L',U')$, as required.

 Ad 4. Using 3. we have the following equivalent statement

\[ \begin{array}[t]{c}
 (a,L,U)<^-(a',L',U') \\ \hline
 \exists_{x\in\delta(a')}\; (\gamma(a),-,\ua a)\leq^+(x,\da a',-) \\ \hline
 \exists_{x\in\delta(a')}\; \gamma(a)<^+x \;\mbox{or }\; (\gamma(a)=x \;\mbox{and }\;
  \cE_{\gamma(a)}-\ua_{\gamma(a)} a\subseteq\da_{\gamma(a)}a') \\ \hline
 a<^-a' \;\mbox{and}\;( \mbox{if }\; \gamma(a)\in\delta(a') \;\mbox{then}\;
 \cE_{\gamma(a)}=\ua_{\gamma(a)} a\cup\da_{\gamma(a)} a')
\end{array} \]
$~\Box$

\begin{lemma}
\label{fact16.5} Let $S$  be an ordered face structure,
$X$ a convex subset of $S$, $a\in X$, $u\in\theta\delta(a)$, $l\in
\cE^X_{u}$. Then
\begin{enumerate}
\item If $a\in S^{\varepsilon}$ then $\gamma(a)\perp^\sim l$.
\item If $a\in S^{-\varepsilon}$ then
\begin{enumerate}
  \item if $u=\gamma\gamma(a)$ then $\gamma(a)<^\sim l$ iff
  $\varrho(a)<^\sim l$;
  \item if $u\in\iota(a)$ then if $x,y\in\delta(a)$ and $x<^\sim l <^\sim
  y$ there is $z\in\delta(a)$ such that $l\leq^+z$ and $x<^\sim z <^\sim y$;
  \item if $u\in\delta\gamma(a)$ then $l<^\sim \gamma(a)$ iff
  $l<^\sim x_0$, where $x_0$ is the $\sim$-minimal element in
  $\delta(a)$ such that $u\in\delta(x)$.
\end{enumerate}
\end{enumerate}
\end{lemma}

{\em Proof.} Ad 1. By pencil linearity we have that either
$\gamma(a)\perp^\sim l$ or $\gamma(a)\perp^+ l$. As $l$ is an
empty loop we cannot have $\gamma(a)\leq^+ l$. We shall show that
$l<^+\gamma(a)$ is impossible, as well.

Suppose not, and that we have a flat upper
$X-\gamma(X^{-\lambda})$-path $l,a_1,\ldots,a_k,\gamma(a)$. As
$a_k\in X^{-\varepsilon}$ we have that $a_k<^+a$. But $a\in
X^\varepsilon$, so by Path Lemma, we have $a_i<^+a$, for
$i=1,\ldots,n$. Let $a_1,\alpha_1,\ldots,\alpha_r,a$ be a flat
upper $X$-path. As $l\in\delta(a)$ and
$\delta\gamma(\alpha_i)\cup\gamma\dot{\delta}^{-\lambda}(\alpha_i)=\theta\delta(\alpha_i)$,
either for all $i$, $l\in\delta\gamma(\alpha_i)$, or there is
$i_0$ such that $l\in\gamma\dot{\delta}^{-\lambda}(\alpha_{i_0})$.
In the former case $l\in\delta\gamma(\alpha_r)=\delta(a)$ and we
get a contradiction, as $a\in T^\varepsilon$. In the later case
$l$ is not an empty loop in $X$ contrary to the assumption.

Ad 2(a).  First, assume $\gamma(a)<^\sim l$. As
$\gamma\gamma(a)=\gamma(l)$, by pencil linearity we have either
$\varrho(a)\perp^\sim l$ or $\varrho(a)\perp^+ l$.  We shall show
that the other cases then $\varrho(a)<^\sim l$ lead quickly to a
contradiction. If $l<^\sim \varrho(a)$ then
$\gamma(a)<^\sim\varrho(a)$ and this is a contradiction. If $l<^+
\varrho(a)$ then $l<^+\gamma(a)$ and this is a contradiction. If
$\varrho(a)<^+ l$ then $l\perp^+\gamma(a)$ and this is again a
contradiction. Thus we must have $\varrho(a)<^\sim l$\footnote{In
the following, the similar simple arguments we will describe in a
shorter form as follows: as other cases are easily excluded, we
have $\varrho(a)<^\sim l$.}.

Next we assume $\varrho(a)<^\sim l$. By pencil linearity we have
either $\gamma(a)\perp^\sim l$ or $\gamma(a)\perp^+ l$. We need to
show that $\gamma(a)<^\sim l$.  We shall show that the condition
$l<^+\gamma(a)$ leads to a contradiction.  The other two are
easily excluded. Clearly $l\not\in\delta(a)$.

So suppose that $l<^+\gamma(a)$. Let $l,a_1,\ldots,a_k,\gamma(a)$
be a flat upper $X-\gamma(X^{-\lambda})$-path. By Path Lemma,
either there is $i<k$ such that $\gamma(a_i)\in\delta(a)$ or
$a_i<^+a$ for $i=1,\ldots,k$. In the former case we have
$l<^+\gamma(a_i)<^\sim\varrho(a)$. Thus, by Lemma \ref{fact15.5},
$l<^\sim \varrho(a)$ and this is a contradiction. In the later
case, there is a flat upper $X$-path
$a_1,\alpha_1,\ldots,\alpha_r,a$.  As $l\in\delta(a_1)$ and
$l\not\in\delta(a)$ there is $i$ such that $l\in\iota(\alpha_i)$.
In particular $l$ is not an empty loop in $X$ and we get a
contradiction again.

Ad 2(b). Fix $x,y\in\delta(a)$ such that
$\gamma(x)=u\in\delta(a)$, $l\in\cE^X_u$ such that $x<^\sim
l<^\sim y$. If $l\in\delta(a)$ 2(ii) obviously holds, so assume
that $l\not\in\delta(a)$. We have
$\gamma(l)\leq^+\gamma(y)\leq^+\gamma\gamma(a)$. Thus by
Proposition \ref{tech_lemma3}, either $l<^\sim\gamma(a)$ or
$l<^+\gamma(a)$. The former case gives immediately
$x<^\sim\gamma(a)$ and a contradiction.

Thus we have $l<^+\gamma(a)$. Fix a flat $X$-path
$l,a_1,\ldots,a_k,\gamma(a)$. By Path Lemma, either there is $i<k$
such that $\gamma(a_i)\in\delta(a)$ or $a_i<^+a$ for all $1\leq
i\leq k$. In the former case $\gamma(a_i)$ is the $z$ we are
looking for.  We shall show that the later case leads to a
contradiction. Take a flat upper $X$-path
$a_1,\alpha_1,\ldots,\alpha_r,a$. As $l\in\delta(a_1)$ and
$l\not\in\delta(a)$,  there is $1\leq i\leq r$ such that
$l\in\gamma\dot{\delta}^{-\lambda}(\alpha_i)$. In particular,
$l\not\in\cE_u$, contrary to the assumption.

Ad 2(c). Suppose $l<^\sim\gamma(a)$. Then, as other cases are
easily excluded, we have $l<x_0$ indeed.

On the other hand, if $l<^\sim x_0$ the only case not easily
excluded is $l<^+\gamma(a)$. Let $l,a_1,\ldots,a_k,\gamma(a)$ be a
flat upper $X-\gamma(X^{-\lambda})$-path.  By Path Lemma either
there is $i_0<k$ such that $\gamma(a_{i_0})\in\delta(a)$ or
$a_i<^+a$ for $i=1,\ldots,k$. In the later case, there is a flat
upper $X$-path $a_1,\alpha_1,\ldots, \alpha_k,a$. As
$l\not\in\delta(a)$ and $l\in\delta(a_1)$, there is $1\leq i\leq
k$ such that $l\in \iota(\alpha_i)$.  In particular
$l\not\in\cE^X_u$, contrary to the assumption.

In the former case we shall show that $u\in\delta\gamma(a_i)$, for
$i=1,\ldots,i_0$.  We have $u=\gamma(l)\in\delta\delta(a_1)$. Note
that if for some $i\leq i_0$, we would have that $u\in\iota(a_i)$,
then, as $u\in\delta(x_0)$, by Lemma \ref{fact9}, we would have
$x_0<^+\gamma(a_{i_0})$, contradicting local discreteness. Now
suppose contrary, that for some $i_1\leq i_0$, we have
$u\not\in\delta\gamma(a_{i_1})$. Then, by the previous
observation, we have
$u=\gamma\gamma(a_{i_1})\in\gamma\dot{\delta}^{-\lambda}(a_{i_1})$.
In particular, $\dot{\delta}^{-\lambda}(a_{i_1})\neq\emptyset$ and
$\gamma(a_i)$ is not a loop, for $i_1\leq i \leq i_0$. As
$u\not\in\iota(a_i)$ for $i\leq i_0$, we have
$u=\gamma\gamma(a_i)$ for $i\leq i_0$.  In particular
$u=\gamma\gamma(a_{i_0})\in\gamma\dot{\delta}^{-\lambda}(a)$. But
$u\in\delta\gamma(a)$ and
$\delta\gamma(a)\cap\gamma\dot{\delta}^{-\lambda}(a)=\emptyset$ so
we get a contradiction. $~\Box$

\begin{proposition}\label{subfs} Let $T$ be an ordered face
  structure, $X$ be a convex subset of $T$. Then
\begin{enumerate}
  \item $[X]$ is an ordered face structure, and  $\nu_X: [X]\ra T$ is a monotone morphism;
  \item if
$f_i:S_i\ra T$, $i=0,1$, are monotone morphisms such
$im(f_0)=im(f_1)$, then there is a monotone isomorphism $g:S_0\ra
S_1$ making the triangle
\begin{center}
\xext=600 \yext=350
\begin{picture}(\xext,\yext)(\xoff,\yoff)
\settriparms[1`1`1;300] \putVtriangle(0,0)[S_0`S_1`T;g`f_0`f_1]
\end{picture}
\end{center}
commutes;
  \item $\cE^X$ is empty iff $\nu_X$ is an
 embedding;
  \item if $X$ is a proper subset of $T$ then $size(T)>size([X])$.
\end{enumerate}
\end{proposition}

{\it Proof.}~ Ad 1. The fact that $\nu_X$ is a monotone morphism is
immediate from the definition of $[X]$. We need to check that
$[X]$ satisfies the axioms of ordered face structures. {\em Local
discreteness} and {\em Strictness}, are easy using Lemma
\ref{subfs0}.

{\em Disjointness}. We shall check that if
$\theta(a,L,U)\cap\theta(a',L',U')=\emptyset$ and
$(a,L,U)<^-(a',L',U')$ then $(a,L,U)<^\sim(a',L',U')$. The
remaining parts of the condition are easy.

Suppose that $\theta(a,L,U)\cap\theta(a',L',U')=\emptyset$ and
$(a,L,U)<^-(a',L',U')$. By Lemma \ref{subfs0} we have that
$a<^-a'$. If $\theta(a)\cap\theta(a')=\emptyset$, we get by
disjointness in $T$ that $a<^\sim a'$, and we are done. Assume
that $\theta(a)\cap\theta(a')\neq\emptyset$.  Thus by Lemma
\ref{fact6}.2 $\gamma(a)\in\delta(a')$. By characterization of
$<^-$ in $[X]$ we have
$\ua_{\gamma(a)}a\cup\da_{\gamma(a)}a'=\cE_{\gamma(a)}$. As
$(\gamma(a),-,\ua a)\in\theta(a,L,U)$ and $(\gamma(a),\da
a',-)\in\theta(a',L',U')$ we get that $(\gamma(a),-,\ua
a)\neq(\gamma(a),\da a',-)$. So
$\cE_{\gamma(a)}\,-\,\ua_{\gamma(a)}a\neq\da_{\gamma(a)}a'$. But
then there is $l\in\ua_{\gamma(a)}a\cap\da_{\gamma(a)}a'$. Hence
$a<^\sim l<^\sim a'$, i.e. $a<^\sim a'$ as required.

{\em Loop filling}. Suppose $(a,L,U)$ is a loop in $[X]$. If
$L\neq\emptyset$ then let $l=\max_\sim(L)$. By Lemma \ref{subfs0}
$(l,\emptyset,-)$ is not a loop in $[X]$. We have
$\gamma(l,\emptyset,-)=(a,-,\ua l)=(a,L,U)$.

Now consider the case $L=\emptyset$. $a$ is not an empty loop
since otherwise $(a,L,U)$ wouldn't be a loop in $[X]$. Thus there
is $\alpha\in X^{-\lambda}$ such that $\gamma(\alpha)=a$. Clearly,
we can choose such $\alpha$ in
$X^{-\lambda}-\gamma(X^{-\lambda})$. Then $(\alpha,\emptyset,-)$
is not a loop and
\[ \gamma(\alpha,\emptyset,-)=(\gamma(\alpha),-,\ua\alpha)=
(a,\emptyset,\ua\alpha)=(a,L,U).\]

{\em Pencil linearity}. Let $(a,L,U)\neq(a',L',U')$ be some faces
in $[X]$ such that
$\dot{\theta}(a,L,U)\cap\dot{\theta}(a',L',U')\neq\emptyset$. Then
either $a\perp^\sim b$ or $a\perp^+ b$ or ($a=b$ and $L\neq L'$).
In the first case we have $(a,L,U)\perp^\sim(a',L',U')$ and in the
remaining cases we have $(a,L,U)\perp^+(a',L',U')$.

To see the second part of the pencil linearity assume that
$\check{a}=(a,L,U)\in [X]^\varepsilon$, $\check{b}=(b,L',U')\in
[X]$, $x\in\delta(b)$, such that $(x,\da b,-),(y,\da b,-)\in
\delta^{-\lambda}(b,L',U')$, and
\[ \gamma\gamma(\check{a})=\gamma(x,\da b,-)\in\delta(y,\da b,-)\]
i.e. for some $t\in\delta(y)$
\begin{equation}\label{lleq}
 (\gamma\gamma(a),-,\ua \gamma(a))=(\gamma(x),-,\ua x)=(t,\da y,-)
\end{equation}
We need to show that either $\check{a}<^\sim \check{b}$ or
$\check{a}<^+ \check{b}$. From the characterization of $<^\sim$
and $<^+$ in $[X]$ it is enough to show that either $a<^\sim b$ or
$a<^+b$. And for that, by Lemma \ref{tech_lemma3}.4, it is enough
to show that $\gamma(a)<^+\gamma(b)$. We shall consider four cases
separately:
\begin{enumerate}
  \item $x,y\in T^{-\lambda}$;
  \item $x\in T^{-\lambda}$ and $y\in T^\lambda$ and there is
  $l_y\in\cE_{\gamma\gamma(a)}$ such that $l_y\leq^+ y$;
  \item $y\in T^{-\lambda}$ and $x\in T^\lambda$ and there is
  $l_x\in\cE_{\gamma\gamma(a)}$ such that $l_x\leq^+ x$;
  \item $x,y\in T^\lambda$ and there are
  $l_x,l_y\in\cE_{\gamma\gamma(a)}$ such that $l_x\leq^+ x$ and $l_y\leq^+ y$.
\end{enumerate}

Case 1. If $x,y\in T^{-\lambda}$ then
$\gamma\gamma(a)\in\iota(b)$. So, by pencil linearity in $T$, we
get that either $a<^\sim b$ or $a<^+b$.

Case 2. In this case we have $\gamma(a)<^\sim l_y$. As
$\gamma(x)=\gamma\gamma(a)$ we have either $x\perp^\sim\gamma(a)$
or $x\perp^+\gamma(a)$. We have $\gamma(a)\not<^-x$. Moreover, by
Lemma \ref{fact3}, $x\not<^+\gamma(a)$, as $x\in T^{-\lambda}$ and
$\gamma(a)\in T^\lambda$. So we have either either $\gamma(a)<^+x$
or  $x<^\sim\gamma(a)$. In the former case we get immediately
$\gamma(a)<^+\gamma(b)$. In the later case, we have
\begin{equation}\label{four_cells}
x<^\sim \gamma(a)<^\sim l_y, \;\;\;\;\; x<^+ \gamma(b)>^+ l_y.
\end{equation}
So we have
$\gamma\gamma(a)\leq^+\gamma(l_y)\leq^+\gamma\gamma(b)$. Thus, by
Lemma \ref{tech_lemma3}, we have $\gamma(a)\perp^+\gamma(b)$ or
$\gamma(a)\perp^\sim\gamma(b)$. Using (\ref{four_cells}) we see
that of four conditions only the $\gamma(a)<^+\gamma(b)$ does not
lead to a contradiction.

Case 3. First note that $\gamma(l_x)=\gamma(x)=\gamma\gamma(a)$
and hence $\gamma(a)\perp^+l_x$ or $\gamma(a)\perp^\sim l_x$. The
inequality $\gamma(a)<^\sim l_x$ is impossible as $l_x\not\in\ua
x=\ua\gamma(a)$.  $\gamma(a)\leq^+ l_x$ is impossible since $l_x$
is an empty loop and $\gamma(a)$ is not. Finally,
$l_x<^+\gamma(a)$ is impossible as $(\gamma(a),-,-)$ is a loop in
$[X]$, and cannot contain any empty loops. So we have shown that
$l_x<^\sim\gamma(a)$. As $\gamma\gamma(a)\in\delta(y)$ and $y\in
T^{-\lambda}$, we have $\gamma(a)<^-y$ and $y\not<^-\gamma(a)$. As
$y\in T^{-\lambda}$ and $\gamma(a)\in T^\lambda$, by Lemma
\ref{fact3}, we cannot have $y<^\sim\gamma(a)$. So we must have
either $\gamma(a)<^+ y$ or $\gamma(a)<^\sim y$. If $\gamma(a)<^+
y$ then clearly $\gamma(a)<^+\gamma(b)$. If
$l_x<^\sim\gamma(a)<^\sim y$ then having $l_x<^+\gamma(b)>^+y$ we
can easily verify, as before in (\ref{four_cells}), that we must
have $\gamma(a)<^+\gamma(b)$.

Case 4. As $l_y<^+y$ and $\da y=\cE_{\gamma\gamma(a)}-\ua
\gamma\gamma(a)$, we have $\gamma(a)<^\sim l_y$. As
$\gamma\gamma(a)=\gamma(l_x)$, we also have $\gamma(a)\perp^+l_x$
or $\gamma(a)\perp^\sim l_x$. It is easy to see that the only
inequality that does not lead to a contradiction in
$l_x<^\sim\gamma(a)$. So we have
\begin{equation}\label{four_cells1}
l_x<^\sim \gamma(a)<^\sim l_y, \;\;\;\;\; l_x<^+ \gamma(b)>^+ l_y.
\end{equation}
From (\ref{four_cells1}) we get, as before from
(\ref{four_cells}), that $\gamma(a)<^+\gamma(b)$.

\vskip 1mm

{\em Globularity}. Let us fix a face $\check{a}=(a,L,U)\in
[X]_{\geq 2}$. As different $a$-cuts are parallel, i.e. they have
the same domains and codomains, to verify Globularity condition in
$[X]$ we don't need know in fact the very cut over $a$ for which
we check the condition.  It is enough to know that it is a cut
over $a$.  So in the following $\check{a}$ will be treated as a
cut over $a$, for which we don't bother to specify exactly which
one it is. In the following $\gamma$ and $\delta$ when applied to
cuts are meant in $[X]$ and when applied to faces are meant in
$T$.

{\em $\gamma$-globularity}. If $\check{a}\in [X]^\varepsilon$ then
we have
\[ \gamma\gamma(\check{a})=(\gamma\gamma(a),-,\ua \gamma(a))=\gamma(1_{(\gamma\gamma(a),-,\ua
\gamma(a))})=\gamma\delta(\check{a}) \]

Now assume that $\check{a}\in [X]^{-\varepsilon}$. We need to
verify the following three conditions:

\begin{description}
  \item[{\em (i)}] $\gamma\gamma(\check{a})\in\gamma\delta(\check{a})$;
  \item[{\em (ii)}] $\gamma\gamma(\check{a})\not\in\delta\dot{\delta}^{-\lambda}(\check{a})$;
  \item[{\em (iii)}] $\gamma\delta(\check{a})\subseteq\gamma\gamma(\check{a})\cup
  \delta\dot{\delta}^{-\lambda}(\check{a})$.
\end{description}

Ad {\em (i)}.  We have $(\varrho(a),\da a,-)\in\delta(\check{a})$
and then using Lemma \ref{fact16.5}.2.(a) we have
\[ \gamma\gamma(\check{a}) = (\gamma\gamma(a),-,\ua\gamma(a))=
(\gamma\varrho(a),-,\ua\gamma(a))=\gamma(\varrho(a),\da a,-) \]

Ad {\em (ii)}.  Suppose that $x\in\delta(a)$ and $u\in\delta(x)$
so that $(x,\da a,-)\in\dot{\delta}^{-\lambda}(\check{a})$ and
$(u,\da x,-)\in\delta\dot{\delta}^{-\lambda}(\check{a})$. If
$u\neq\gamma\gamma(a)$ then clearly $(u,\da
x,-)\neq(\gamma\gamma(a),-,\ua \gamma(a))$. If $u=\gamma\gamma(a)$
then $x\in\delta^{\lambda}(a)$ and by characterization of loops in
$[X]$, Lemma \ref{subfs0}, there is $l\in\cE_{\gamma\gamma(a)}^X$
such that $l\leq^+ x$.  Then $l\leq^+\gamma(a)$. So
$\da_{\gamma\gamma(a)}x\not\ni l\not\in\ua \gamma(a)$ and hence
\[ \gamma\gamma(\check(a))=(\gamma\gamma(a),-,\ua \gamma(a))\neq(u,\da x,-)\in\delta(x,\da a,-). \]

Ad {\em (iii)}. Let $x\in\delta(a)$ so that
\[\gamma(x,\da a,-)=(\gamma(x),-,\ua x)\neq
(\gamma\gamma(a),-,\ua\gamma(a))=\gamma\gamma(\check{a}).\] Then
either $\gamma(x)\neq\gamma\gamma(a)$ or
$\gamma(x)=\gamma\gamma(a)$ and there is
$l\in\cE^X_{\gamma\gamma(a)}$ so that $x<^\sim l$ and
$\gamma(a)\not<^\sim l$. This implies that the face
\[ y_0 = \min_\sim\{ y\in\delta(a) : \gamma(x)\in\delta(y)\;
{\rm and\; either }\; y\in T^{-\lambda} \;{\rm or}\; \exists
_{l\in\cE_{\gamma(x)}^X}\, l\leq^+y \} \] is well defined, i.e.
the set over which the minimum is taken is not empty. Then we have
$(y_0,\da a, -)\in\dot{\delta}^{-\lambda}(\check{a})$ and
$(\gamma(x),-,\ua x)=(\gamma(x),\da y_0,-)\in
\delta\dot{\delta}^{-\lambda}(\check{a})$.

{\em $\delta$-globularity}. We consider separately two cases
$\gamma(\check{a})\in [X]^{\varepsilon}$ and $\gamma(\check{a})\in
[X]^{-\varepsilon}$. In the former case we need to verify two
conditions
\begin{description}
  \item[{\em (i)}] $\dot{\delta}\delta(\check{a})\subseteq\gamma\dot{\delta}^{-\lambda}(\check{a})$;
  \item[{\em (ii)}] $\gamma\gamma\gamma(\check{a})=\gamma\gamma\delta^{\varepsilon}(\check{a})$.
\end{description}
Ad {\em (i)}. Let $x\in\dot{\delta}(a)$ and $u\in\dot{\delta}(x)$
so that $(u,\da x,-)\in \dot{\delta}\delta(\check{a})$. As $T$ is
an ordered face structure,
$\dot{\delta}\delta(a)\subseteq\gamma\dot{\delta}^{-\lambda}(a)$
and hence $u\in\gamma\dot{\delta}^{-\lambda}(a)$. Thus there is a
$y\in\dot{\delta}^{-\lambda}(a)$ such that  $\gamma(y)=u$. From
this follows that the face
\[ y_1 = \max_\sim\{ y\in\delta(a) : \gamma(y)=u,\;  y<^\sim x\;
{\rm and\; either }\; y\in T^{-\lambda} \;{\rm or}\; \exists
_{l\in\cE_{\gamma(x)}^X}\, l\leq^+y \} \] is well defined. Then we
have $(y_1,\da a,-)\in\dot{\delta}^{-\lambda}(\check{a})$ and,
using Lemma \ref{fact16.5}.2.(b), we get
\[ \gamma(y_1,\da a,-)=(\gamma(y_1),-,\ua y_1)=(u,\da x,-). \]
This shows {\em (i)}.

Ad {\em (ii)}. First note that if
$\gamma(\check{a})\in[X]^\varepsilon$ then $\gamma(a)\in
X^\varepsilon$ and hence $\delta^\varepsilon(a)\neq\emptyset$.  So
$\delta^\varepsilon(\check{a})\neq\emptyset$, as well. Thus we
need to show that
$\gamma\gamma\delta^\varepsilon(\check{a})\subseteq\gamma\gamma\gamma(\check{a})$.

Fix $(x,\da a,-)\in\delta^{\varepsilon}(\check{a})$ and
$l\in\cE^X_{\gamma\gamma\gamma(a)}$. It is enough to show that
\begin{equation}\label{delta-convex}
\gamma(x)<^\sim l\;\;\;\;{\rm iff}\;\;\;\;\gamma\gamma(a)<^\sim l.
\end{equation}
Clearly $x\in\delta^\varepsilon(a)$. By Lemma \ref{fact16.5}.1
$\gamma\gamma(a)\perp^\sim l$. Since $l\in\cE^X$ and
$\gamma(x)\leq\gamma\gamma(a)$ we have $l\not\perp^\sim\gamma(x)$.

Thus ($\gamma\gamma(a)<^\sim l$ and $l<^\sim \gamma(x)$) or
($\gamma(x)<^\sim l$ and $l<^\sim \gamma\gamma(a)$) then
$\gamma(x)\perp^\sim \gamma\gamma(a)$ and this is a contradiction.
Therefore (\ref{delta-convex}) holds. This shows {\em (ii)} and
end up the case $\gamma(\check{a})\in [X]^{\varepsilon}$.

In case $\gamma(\check{a})\in [X]^{-\varepsilon}$ we need to
verify the following four conditions
\begin{description}
    \item[{\em (i)}] $\dot{\delta}\delta(\check{a})\subseteq\delta\gamma(\check{a})
    \cup\gamma\dot{\delta}^{-\lambda}(\check{a})$;
    \item[{\em (ii)}] $\delta\gamma(\check{a})\subseteq\dot{\delta}\delta(\check{a})$;
    \item[{\em (iii)}] $\delta\gamma(\check{a})\cap \gamma\dot{\delta}^{-\lambda}(\check{a})=\emptyset$;
    \item[{\em (iv)}] $\gamma\gamma\delta^\varepsilon(\check{a})\subseteq\theta\delta\gamma(\check{a})$.
\end{description}

Ad {\em (i)}. Fix $x\in\delta(a)$ and $u\in\delta(x)$ so that
$(u,\da x,-)\in\dot{\delta}\delta(\check{a})$.  Assume that
$(u,\da x,-)\not\in\delta\gamma(\check{a})$.  Then either
$u\not\in\delta\gamma(a)$ or $u\in\delta\gamma(a)$ and there is
$l\in\cE^X_{u}$ such that  $l<^\sim x$ and $l\not<^\sim\gamma(a)$.
Then, by Lemma \ref{fact16.5}.2.(c), the face
\[ y_2 = \max_\sim\{ y\in\delta(a) : \gamma(y)=u,\;
{\rm and\; either }\; y\in T^{-\lambda} \;{\rm or}\; \exists
_{l\in\cE_{u}^X}\, l\leq^+y \} \] is well defined. Clearly,
$(y_2,\da a,-)\in\dot{\delta}^{-\lambda}(\check{a})$.  By Lemma
\ref{fact16.5}.2.(b), we have
\[ (u,\da x,-)=(u,-,\ua y_2)=\gamma(y_2,\da a,-)\in\gamma\dot{\delta}(\check{a})\]

Ad {\em (ii)}. Fix $u\in\delta\gamma(a)$ so that $(u,\da
\gamma(x),-)\in\delta\gamma(\check{a})$. If $a\in X^\varepsilon$
then $\delta(\check{a})=1_{(u,\da \gamma(x),-)}$ and hence
$(u,\da\gamma(x),-)\in\delta\delta(\check{a})$. So suppose that
$a\in X^{-\varepsilon}$.  Then the face
\[ y_3=\min_\sim\{ y\in\dot{\delta}(a) :u\in\delta(y) \} \]
is well defined. By Lemma \ref{fact16.5}.2.(c), we have
\[ (u,\da\gamma(a),-)=(u,\da y_3,-)\in\delta\delta(\check{a})\]

Ad {\em (iii)}. Let $u\in\delta\gamma(a)$ so that $(u,\da
\gamma(a),-)\in\delta\gamma(\check{a})$. We shall show that if
$(u,\da \gamma(a),-)\in\gamma\delta(\check{a})$ then $(u,\da
\gamma(a),-)\in[X]^\lambda$.  So fix $z\in\delta(a)$ such that
$(u,\da a,-)=\gamma(z,\da a,-)=(\gamma(z),-,\ua z)$. As
$\gamma(z)=u\in\delta\gamma(a)$, $z$  is a loop. If we were to
have $l\in\cE_u$ such that $l\leq^+z$ then $l\leq^+\gamma(a)$ and
hence $l\not\in\da_u\gamma(a)$ and $l\not\in\ua z$.  Thus
\[ (u,\da \gamma(a),-)\neq(\gamma(z),-,\ua z) \]
contrary to the assumption.

Ad {\em (iv)}.  Let $x\in\delta^\varepsilon(a)$ so that $(x,\da
a,-)\in\delta^\varepsilon(\check{a})$. We shall show that
\begin{equation}\label{glob-empty}
\gamma\gamma(x,\ua a,-) = (\gamma\gamma(x),-,\ua
\gamma(x))\in\theta\delta\gamma(\check{a}).
\end{equation}
Note that as $T$ is an ordered face structure, we have
$\gamma\gamma(x)\in\gamma\gamma\delta^\varepsilon(a)\subseteq\theta\delta\gamma(a)$.

First we claim that for any $t\in\delta\gamma(a)$ we have
$t\not\perp^+\gamma(x)$. Fix $t\in\delta\gamma(a)$. As
$\gamma(x)\in\gamma\dot{\delta}^{-\lambda}(a)$, using
$\gamma\dot{\delta}^{-\lambda}\cap\delta\gamma(a)=\emptyset$ we
get that $\gamma(x)\not\leq^+t$. Now suppose contrary, that
$t<^+\gamma(x)$.  Thus there is a flat upper path
$t,x_1,\ldots,x_n,\gamma(x)$, with $r\geq 1$. As $\gamma(x)$ is a
loop and $x$ is an empty domain face, by Path Lemma, we have
$x_i<^+x$ for $i=1,\ldots,n$. In particular, there is a flat upper
path in $T$, $x_1,a_1,\ldots,a_m,x$.  As $x\in T^\varepsilon$ and
$x_1\in T^{-\varepsilon}$, for some $1\leq j_0\leq m$, we have
$t\in\iota(a_{j_0})$. Thus
$t\in\delta\gamma(a)\cap\iota(a_{j_0})$.  Hence by Lemma
\ref{fact9}, we have
$\gamma(a)<^+\gamma(a_{j_0})\leq^+\gamma(a_s)=x$. But
$x\in\delta(a)$ and we get a contradiction with strictness. This
ends the proof of the claim.

Now let $u=\gamma\gamma(x)$. Using the claim it is easy to see
that one of the following conditions holds:
\begin{description}
  \item[{\em(a)}] $u=\gamma\gamma\gamma(a)$ and
  $\varrho\gamma(a)<^\sim\gamma(x)$;
  \item[{\em(b)}] $u\in\delta\gamma\gamma(a)$ and with
  $s_1=\min_\sim\{s\in\delta\gamma(a):u\in\delta(s) \}$ we have
  $\gamma(x)<^\sim s_1$;
  \item[{\em(c)}] there are $s_0,s_1\in\delta\gamma(a)$ such that
  $\gamma(s_0)=u\in\delta(s_1)$ and $s_0<^\sim\gamma(x)<^\sim
  s_1$.
\end{description}
In each case we shall show (\ref{glob-empty}).

Ad {\em(a)}. Using Lemma \ref{fact16.5}.2.(a), we have
\[ (u,-,\ua \varrho\gamma(a))\leq^+(u,-,\ua
\gamma(x))\leq^+ (u,-,\ua \gamma\gamma(a))\leq^+(u,-,\ua
\varrho\gamma(a))\] Thus $\gamma\gamma(x,\da a,-)=(u,-,\ua
\varrho\gamma(a))\in \gamma\delta\gamma(\check{a})$.

Ad {\em(b)}. Using Lemma \ref{fact16.5}.1 and
\ref{fact16.5}.2.(c), we have
\[ (u,\da \gamma\gamma(a),-)\leq^+(u,\da
s_1,-)\leq^+(u,\da \gamma(x),-)\leq^+(u,\da \gamma\gamma(a),-)\]
Thus $\gamma\gamma(x,\da a,-)=(u,\da s_1,-)\in
\delta\delta\gamma(\check{a})$.

Ad {\em(c)}. Let $s_1$ be as above in $(c)$ and $s_0$ maximal such
as in $(c)$, i.e.
\[ s_0=\max_\sim\{ s\in\delta\gamma(a): \gamma(s)=u,\, s<^\sim\gamma(x) \}. \]
 Suppose $\ua s_0\neq \ua\gamma(x)$. Then there is $l\in\cE_u$
 such that $s_0<^\sim l <^\sim\gamma(x)$. Then by Lemma
 \ref{fact16.5}.2.(b) there is $t\in\delta\gamma(a)$ such that
 $l<^+t$ and $s_0<^\sim t<^\sim s_1$. We shall show that the
 existence of such a $t$ leads to a contradiction. Note that
 $t$  is a loop and that $\gamma(t)=u$. By the claim proven above
 it follows that $t\not\perp^+\gamma(x)$. So, by pencil linearity,
 we should have $t\perp^\sim\gamma(x)$. But if $\gamma(x)<^\sim t$
 then by transitivity we get $l<^\sim t$ and we get a
 contradiction with disjointness. On the other hand, it $t<^\sim
 \gamma(x)$ then, as other cases are easily excluded, we have that
 $s_0<^\sim t$. But this contradicts the choice of $s_0$.
Thus $\ua s_0= \ua\gamma(x)$ holds and we have $\gamma\gamma(x,\da
a,-)=(u,-,\ua s_0)\in \gamma\delta\gamma(\check{a})$. This ends
{\em (iv)} and 1.

Ad 2. By 1. it is enough to show that if $f:S\ra T$ is a monotone
morphism such that $f(S)=X$ then there is a monotone isomorphism
$g:S\ra [X]$ making the triangle
\begin{center}
\xext=600 \yext=350
\begin{picture}(\xext,\yext)(\xoff,\yoff)
\settriparms[1`1`1;300] \putVtriangle(0,0)[S`{[X]}`T;g`f`\nu_X]
\end{picture}
\end{center}
 commutes. We put, for $a\in S$
\[ g(a)= \left\{ \begin{array}{ll}
           (f(a),-,\ua f(a))  &  \mbox{if  $\alpha\in S^{-\lambda}-\gamma(S^{-\lambda})$ such that $a=\gamma(\alpha)$,}\\
           (f(a),\emptyset,-) & \mbox{such that $a\not\in\gamma(S^{-\lambda})$.}
                                    \end{array}
                \right. \]
Note that if $a\in\gamma(S^{-\lambda})$ then there is a unique
$\alpha\in S^{-\lambda}-\gamma(S^{-\lambda})$ such that
$a=\gamma(\alpha)$.  This shows that $g$ is a well defined
function. As monotone morphisms preserves and reflects $<^\sim$ (in
particular $f$ and $\nu_X$ does) it follows that $g$ preserves
$<^\sim$.

Before we verify the other properties of $g$ let us make one
observation.  Fix $x\in X$ and let
\[ y_{\min} =\min_{<^+}\{y'\in S : f(y)=x \}, \hskip 5mm y_{\max} =\max_{<^+}\{y'\in S : f(y)=x \}  \]
and $y_{\min}, l_1,\ldots,l_k,y_{\max}$ be a flat upper
$S-\gamma(S^{-\lambda})$-path from $y_{\min}$ to $y_{\max}$.  Then
$\cE^X_x=\{ f(l_i)\}_{1\leq i\leq k}$.

With this description it is easy to see that $g$ preserves both
$\gamma$ and $\delta$. Fix $a\in S_{\geq 1}$. Then if $b\in
S^{-\lambda}-\gamma(S^{-\lambda})$ and $\gamma(b)=\gamma(a)$ we
have
\[ g(\gamma(a)=(f(\gamma(a)),-,\ua f(b))=(\gamma(f(a)),-,\ua f(a))=\gamma (f(a)). \]
If $\gamma(a)\not\in\gamma(S^{-\lambda})$ we can show that
$g(\gamma(a)=(\gamma (f(a))$ in a similar way.

If $a\in S^{-\varepsilon}$ then we have
\[ g(\delta(a)=\{ (f(x),-,\ua f(b)) : \gamma(b)=x\in\delta(a),\;
b\in S^{-\lambda}-\gamma(S^{-\lambda}) \} \cup  \]
\[ \cup \{ (f(x),\emptyset,-) : x\in\delta(a),\; x\not\in\gamma(S^{-\lambda}) \}  = \]
\[ = \{ (f(x),\da f(a),-) : x\in\delta(a) \} =\delta(g(a)). \]
If $a\in S^{\varepsilon}$ we clearly have
$g(\delta(a)=\delta(g(a))$, as well.

It remains to show that $g$ is a bijection.  Suppose $g$ is not
one-to-one. Let $a,b\in S$ such that $g(a)=g(b)$ and $a\neq b$. In
particular  $f(a)=f(b)$.  By Lemma \ref{fiber} we can assume that
there is $\alpha\in S^u-\gamma(S^{-\lambda})$ such that
$a=\delta(\alpha)$ and $\gamma(\alpha)=b$. Then
 \[ g(a)=(f(a),L_a,U_a)=(f(b),L_b,U_b)=g(b).\]
But $U_a\ni f(\alpha)\not\in U_b$, and we get a contradiction,
i.e. $g$ is one-to-one.

To see that $g$ is onto fix $(a,L,U)\in [X]$.  First assume that
$L\neq\emptyset$.  Then let $\alpha=\max_{<^\sim}(L)\in X$ and let
$\alpha'=\min_{<^+}\{ \alpha''\in S : f(\alpha'')=\alpha \}$.
Clearly $\alpha'\in S^{-\lambda}-\gamma(S^{-\lambda})$.  Then
\[ g(\gamma(\alpha'))=(f(\gamma(\alpha')),-,\ua
\alpha')=(a,L,U).\] If $L=\emptyset$ then with
$b=\min_{<^+}\{b':f(b')=a \}$ we have
\[ g(b) = (f(b),\emptyset,-)=(a,L,U)\]
in this case $(a,L,U)$ is in the image of $g$, as well. Thus $g$
is onto and hence a bijection.

Ad 3. If $\cE^X$ is empty there is exactly one cut
$(x,\emptyset,\emptyset)$ over any face $x\in X$. So $\nu_X$ is an
embedding in that case.  If there is $l\in \cE^X$ then
$(\gamma(l),\da l,-)\neq (\gamma(l),-,\ua l)$ and
$\nu_X(\gamma(l),\da l,-)=\nu_X(\gamma(l),-,\ua l)$, i.e.  $\nu_X$
is not an embedding.

Ad 4. First note that for any $k\in\o$,
$[X]_k-\delta([X]^{-\lambda}_{k+1})=\{ (a,-,\emptyset) : a\in
X-\delta(X^{-\lambda}_{k+1}) \}$. In particular, we have
$size([X])_k=|X-\delta(X^{-\lambda}_{k+1}|$. Now fix $a$ and $k$
so that $a\in T_k-X_k$, $a$ is a face of the maximal dimension not
in $X$. Then $T_{k+1}=X_{k+1}$ and hence
$a\not\in\delta(T^{-\lambda}_{k+1})$ so we have
\[ a\in T_k-\delta(T^{-\lambda}_{k+1}) =  T_k-\delta(X^{-\lambda}_{k+1} )
\supset  X_k-\delta(X^{-\lambda}_{k+1})\not\ni a \] Thus
\[ size(T)_k=|T_k-\delta(T^{-\lambda}_{k+1})
|>|X_k-\delta(X^{-\lambda}_{k+1})|=size([X]_k. \]
 As $size(T)_l=size([X]_l$, for $l>k$, we have $size(T)>size([X])$.
$~\Box$

Even if the equivalence classes of objects of the comma category
$\ofs\da T$ corresponds to the elements of the poset $Convex((T)$
of convex subsets of $T$ we are not saying that $\ofs\da T$ and
$Convex((T)$ are equivalent as categories. In fact, if   $X\subset
Y$ are convex subsets of $T$, it does not mean that there is a
morphism from $[X]$ to $[Y]$ over $T$ as following example shows.

{\em Example 1.} Let $X\subset Y$ be convex subsets of an ordered
face structure $T$ as shown on the diagram below.

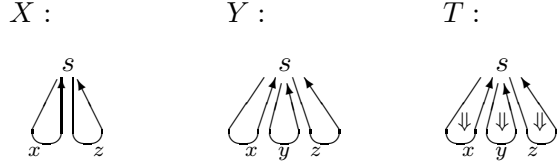
\begin{center} \xext=1500 \yext=660
\begin{picture}(\xext,\yext)(\xoff,\yoff)
\put(0,600){$X:$}
  \put(130,230){\oval(100,100)[b]}
  \put(80,230){\line(1,2){85}}
  \put(180,230){\vector(0,1){180}}
     \put(65,100){$^{x}$}
    \put(180,430){$s$}

  \put(270,230){\oval(100,100)[b]}
   \put(220,230){\line(0,1){180}}
  \put(320,230){\vector(-1,2){85}}
   \put(290,100){$^{z}$}

\put(750,600){$Y:$}
  \put(950,230){\oval(100,100)[b]}
  \put(900,230){\line(1,4){40}}
  \put(1000,230){\vector(-1,4){40}}
   \put(930,100){$^{y}$}

  \put(810,230){\oval(100,100)[b]}
  \put(760,230){\line(2,3){120}}
  \put(860,230){\vector(1,3){60}}
     \put(815,100){$^{x}$}
    \put(930,430){$s$}

  \put(1090,230){\oval(100,100)[b]}
   \put(1040,230){\line(-1,3){60}}
  \put(1140,230){\vector(-2,3){120}}
   \put(1040,100){$^{z}$}

\put(1500,600){$T:$}
  \put(1700,230){\oval(100,100)[b]}
  \put(1650,230){\line(1,4){40}}
  \put(1750,230){\vector(-1,4){40}}
  \put(1680,200){$^\Da$}
   \put(1680,100){$^{y}$}

  \put(1560,230){\oval(100,100)[b]}
  \put(1510,230){\line(2,3){120}}
  \put(1610,230){\vector(1,3){60}}
  \put(1550,200){$^\Da$}
     \put(1565,100){$^{x}$}
    \put(1680,430){$s$}

  \put(1840,230){\oval(100,100)[b]}
   \put(1790,230){\line(-1,3){60}}
  \put(1890,230){\vector(-2,3){120}}
   \put(1810,200){$^\Da$}
   \put(1790,100){$^{z}$}
\end{picture}
\end{center}
Clearly $X\subseteq Y$. And the stretching of $X$ and $Y$ gives
\begin{center}
\xext=1600 \yext=220
\begin{picture}(\xext,\yext)(\xoff,\yoff)
\put(0,200){\makebox(50,50){$[X]:$}}
\putmorphism(0,0)(1,0)[(s,\emptyset,-)`(s,\{ x \},-)`x]{800}{1}a
 \putmorphism(800,0)(1,0)[\phantom{(s,\{ x \},-)}`(s,\{ x,z\},-)`z]{800}{1}a
\end{picture}
\end{center}
and
\begin{center}
\xext=2800 \yext=220
\begin{picture}(\xext,\yext)(\xoff,\yoff)
\put(0,200){\makebox(50,50){$[Y]:$}}
\putmorphism(0,0)(1,0)[(s,\emptyset,-)`(s,\{ x \},-)`x]{800}{1}a
 \putmorphism(800,0)(1,0)[\phantom{(s,\{ x \},-)}`(s,\{ x,y\},-)`y]{900}{1}a
 \putmorphism(1700,0)(1,0)[\phantom{(s,\{ x,y \},-)}`(s,\{ x,y,z\},-)`z]{1000}{1}a
\end{picture}
\end{center}
respectively. Clearly there is no map from $[X]$ to $[Y]$ over
$T$.

Moreover is such a comparison map exists it does not need to be
unique, as the following example shows.

{\em Example 2.}
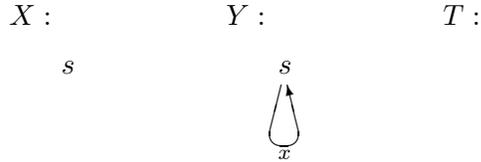
\begin{center} \xext=1500 \yext=660
\begin{picture}(\xext,\yext)(\xoff,\yoff)
\put(0,600){$X:$} \put(180,430){$s$}

\put(750,600){$Y:$}
  \put(950,230){\oval(100,100)[b]}
  \put(900,230){\line(1,4){40}}
  \put(1000,230){\vector(-1,4){40}}
   \put(930,100){$^{x}$}
    \put(930,430){$s$}

\put(1500,600){$T:$}
  \put(1700,230){\oval(100,100)[b]}
  \put(1650,230){\line(1,4){40}}
  \put(1750,230){\vector(-1,4){40}}
  \put(1680,200){$^\Da$}
   \put(1680,100){$^{x}$}
    \put(1680,430){$s$}
\end{picture}
\end{center}
Clearly $X\subseteq Y$. The stretching of $Y$ gives
\begin{center}
\xext=1300 \yext=120
\begin{picture}(\xext,\yext)(\xoff,\yoff)
\put(0,0){\makebox(50,50){$[Y]:$}}
\putmorphism(500,0)(1,0)[(s,\emptyset,-)`(s,\{ x \},-)`x]{800}{1}a
\end{picture}
\end{center}
Thus from $[X]=X$ to $[Y]$ there are two monotone morphisms and both of
them commutes over $T$.

\section{Quotients of positive face structures}\label{quotient of pfs}

Positive face structures can be thought of as ordered face
structures without empty-domain faces. If we collapse to
'identity' some domains of some unary faces which are not
codomains of any other face in a positive face structure we obtain
an ordered face structures which is not necessarily positive. In
this section we shall describe this construction of a quotient of
a positive face structure and prove its properties. In the next
section we shall show that, we can obtain any ordered face
structure in this way.


Let $T$, $S$ be ordered (or positive) face structures. We say that
$f : T \ra S$ is a {\em collapsing morphism}\index{morphism!collapsing -} if
$f=\{ f_k:T_k \ra S_k\cup 1_{S_{k-1}}\}_k\in\o$, $f$ preserves
$<^\sim$, $\gamma$ and $\delta$ i.e. for $k>0$, $a,b\in S_k$, we
have
\begin{enumerate}
  \item if $a<^\sim b$ and $f(a),f(b)\in T_k$ then
$f(a)<^\sim f(b)$
  \item $f(\gamma(a))=\gamma(f(a))$ and $f(\delta(a))\equiv_1\delta(f(a))$.
\end{enumerate}
The {\em kernel}\index{kernel} of the morphism $f$ is the set of faces
$ker(f)=f^{-1}(1_S)$.

{\em Example.}
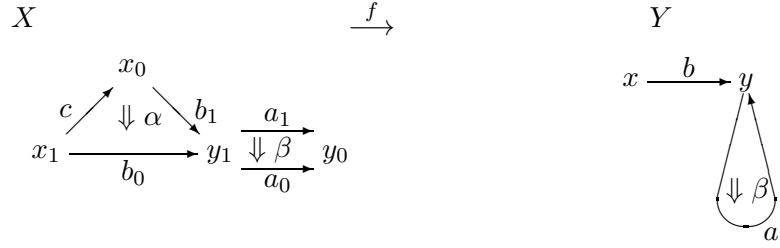
\begin{center} \xext=2600 \yext=650 \adjust[`I;I`;I`;`I]
\begin{picture}(\xext,\yext)(\xoff,\yoff)

\put(0,750){\makebox(50,50){$X$}}
\put(1200,750){\makebox(50,50){$\stackrel{f}{\lra}$}}

\settriparms[-1`1`1;300]
\putAtriangle(100,300)[x_0`x_1`\phantom{y_1};c`b_1`b_0]

\putmorphism(700,380)(1,0)[\phantom{y_1}`\phantom{y_0}`a_1]{400}{1}a
\putmorphism(700,300)(1,0)[y_1`y_0`]{400}{0}b
\putmorphism(700,250)(1,0)[\phantom{y_1}`\phantom{y_0}`a_0]{400}{1}b
 \put(850,290){\makebox(50,50){$\Downarrow \beta$}}
 \put(400,400){\makebox(50,50){$\Downarrow \alpha$}}

\put(2200,750){\makebox(50,50){$Y$}}
\putmorphism(2120,550)(1,0)[x`y`b]{400}{1}a

 \put(2520,150){\oval(200,200)[b]}
 \put(2500,150){\makebox(50,50){$\Downarrow \beta$}}
 \put(2620,150){\vector(-1,4){90}}
 \put(2420,150){\line(1,4){90}}
 \put(2580,0){\makebox(50,50){$a$}}
\end{picture}
\end{center}
\vskip -5mm where $f$ is given by: $x_i\mapsto x$, $y_i\mapsto y$,
$a_0\mapsto a$, $a_1\mapsto 1_y$, $b_i\mapsto b$, $c\mapsto 1_x$,
$\alpha\mapsto 1_b$, $\beta\mapsto \beta$. We have $ker(f)\{
c,\alpha,a_1 \}$.

{\em Remark.} The collapsing morphisms do not compose. If a map $f:X\ra
Y$ sends $\alpha$ to $1_a$ and  a map $g:Y\ra Z$ sends $a$ to
$1_x$ then $g\circ f$ should send $\alpha$ to $1_{1_x}$ but we
don't consider such faces in ordered face structures.

Let $T$ be a positive face structure, $T^u$ unary faces in $T$. A
set $\cJ\subseteq T^u$ is an {\em ideal}\index{ideal} iff
\begin{enumerate}
  \item $\cJ\cap\gamma(T)=\emptyset$;
  \item $\cJ\cap\delta(\cJ)=\emptyset$.
\end{enumerate}
$\sim_{\cJ_{k+1}}$ is the least equivalence relation on $T_k$
containing $\sim'_{\cJ_{k+1}}$; for $x,x'\in T_k$ we have
$x\sim'_{\cJ_{k+1}}x'$ iff there is $a\in \cJ_{k+1}$ such that
$x=\delta(a)$ and $\gamma(a)=x'$. The kernel of any collapsing morphism
is an ideal.

We define an ordered hypergraph $T_{/\cJ}$ the quotient of $T$ by
the ideal $\cJ$:
\begin{enumerate}
  \item $T_{/\cJ,k}=(T_k-\cJ_k)_{/\sim_{\cJ_{k+1}}}$,
  \item $\gamma_{/\cJ}: T_{/\cJ,k+1}\lra T_{/\cJ,k}$, $\;\;\;\delta_{/\cJ}: T_{/\cJ,k+1}\lra T_
  {/\cJ,k}\dsum \b1_{T_{/\cJ,k-1}}$,
   \[ \gamma_{/\cJ}([a])=[\gamma^T(a)],\hskip 10mm  \delta_{/\cJ}([a])= \left\{ \begin{array}{ll}
          1_{[\gamma^T\gamma^T(a)]} & \mbox{if  $\delta(a)\subseteq \cJ$,}  \\
          \{[x] : x\in\delta^T(a)-\cJ\} & \mbox{otherwise.}
                                    \end{array}
                \right. \]
   for $[a]\in T_{/\cJ,k+1}$,
  \item $[x]<^{T_{/\cJ},k,\sim}[x']$ iff $x<^{T_k,-}x'$, for
  $[x],[x']\in T_{/\cJ,k}$.
\end{enumerate}

We define $q_\cJ : T \lra T_{/ \cJ}$ by
\[ q_{\cJ}(a)= \left\{ \begin{array}{ll}
          1_{[\gamma(a)]} & \mbox{if  $a\in \cJ$,}  \\
          {[a]} & \mbox{otherwise.}
                                    \end{array}
                \right. \]
In the remaining part of the section we are going to prove

\begin{theorem}\label{quotient ofs} Let $T$ be a positive face structure, $\cJ\subseteq T^u$ is an ideal.
Then $T_{/\cJ}$ is an ordered face structure and $q_\cJ$ is a
collapsing morphism with kernel $\cJ$.
\end{theorem}

Before we prove this theorem we need some Lemmas.

The class $\cL_\cJ$ of $\cJ${\em -loops}\index{loop@$\cJ$-loop} in
$T$, is defined as the least set $X\subseteq T$ such that
\begin{center}
if $\alpha\in T$ and $\delta(\alpha)\subseteq \cJ\cup X$ then
$\gamma(\alpha)\in X$.
\end{center}
Note that $\cL_\cJ\cap \cJ=\emptyset$.

The following three lemmas concerns positive face structures and
their quotients.

\begin{lemma}\label{unary faces} Let $T$ be a positive face structure, $\cJ\subseteq T^u$ is an ideal.
Then
\begin{enumerate}
  \item If $a\in T$ then the following are equivalent:
\begin{enumerate}
  \item $\gamma(a)\in T^u$;
  \item $\delta(a)\subseteq T^u$;
  \item there is $v\in \delta\gamma(a)$ and a $\delta^u(a)$-path
  from $v$ to $\gamma\gamma(a)$.
\end{enumerate}
  \item $\cL_\cJ\subseteq T^u$.
  \item If $a<^+b$ and $b\in T^u$ then $a\in T^u$.
\end{enumerate}
\end{lemma}

{\em Proof.} 2. and 3. follows from 1. We shall prove 1.

Fix $a\in T$. Let $x_0=\max_{<^\sim}(\delta(a))$.  From
globularity we have
\[\delta\gamma(a)=\delta\delta(a)-(\gamma\delta(a)-\gamma\gamma(a))
\hskip 3mm {\rm and} \hskip 3mm
(\gamma\delta(a)-\gamma\gamma(a))\subseteq \delta\delta(a)
\]
Recall that if $x,y\in\delta(a)$ then
$\delta(x)\cap\delta(y)=\emptyset$. Using these observations, we
get
\[ |\delta\gamma(a)|=|\delta\delta(a)|-|\gamma\delta(a)-\gamma\gamma(a)|
=\bigcup_{x\in \delta(a)} |\delta(x)| -|\gamma(\delta(a)-x_0)|= \]
\[ = \bigcup_{x\in \delta(a)} |\delta(x)| -(|\delta(a)|-1) =
1+ \bigcup_{x\in \delta(a)} (|\delta(x)|-1) \]
 This shows that the set $\delta\gamma(a)$ is a singleton if and only if for
$x\in\delta(a)$ the sets $\delta(x)$ are singletons. This shows
that {\em (a)} is equivalent to {\em (b)}.

Clearly, {\em (b)} implies {\em (c)}.  We shall show the converse.
Let $v\in\delta\gamma(a)$ and $v,x_1,\ldots,x_k,\gamma\gamma(a)$
be an upper $\delta^u(a)$-path from $v$ to $\gamma\gamma(a)$. We
shall show that $\delta(a)=\{ x_1,\ldots,x_k \}$.   Suppose
contrary, that there $y\in \delta(a)-\{ x_1,\ldots,x_k \}$.  Let
$y=y_0,y_1,\ldots,y_r,\gamma\gamma(a)$ be an upper
$\delta(a)$-path to $\gamma\gamma(a)$. Hence
$\gamma(y_r)=\gamma(x_k)$ and then $y_r=x_k$. Let $r'=\min \{
i:y_i\not\in \{ x_1,\ldots,x_k \} \}$. Then $r'<r$ and
$y_{r'+1}=x_j$ for some $j$. If $j=1$ then
$\gamma(y_{r'})=v\in\delta(x_1)\subseteq\delta\gamma(a)$. But then
$v\in \delta\gamma(a)\cap\gamma\delta(a)$, which contradicts
globularity. If $j>1$ then, as $x_j\in T^u$, we have
$\gamma(y_{r'})=\gamma(x_{j-1})$. But, as $y_{r'}, x_{j-1}\in
\delta(a)$ we must have $y_{r'}=x_{j-1}$ contrary to the choice of
$r'$. Thus $\delta(a)=\{ x_1,\ldots,x_k \}$ and {\em (c)} implies
{\em (b)} as well.
 $~\Box$

The following lemma describes some basic properties $T_{/\cJ}$.

\begin{lemma}\label{quotient1} Let $T$ be a positive face structure,
$\cJ\subseteq T^u$ is an ideal, $a,x,y\in T-\cJ$. Then
\begin{enumerate}
   \item $[x]_{\sim_\cJ}=[y]_{\sim_\cJ}$ iff $x=y$ or there is an
  upper $\cJ$-path from $x$ to $y$ or from $y$ to $x$.
  \item The functions $\gamma_{/\cJ}$ and $\delta_{/\cJ}$ are well defined.
  \item $[a]\in T_{/\cJ}^\varepsilon$ if and only if $\delta(a)\subseteq\cJ$.
  \item $[a]\in T_{/\cJ}^\lambda$ if and only if $a\in\cL_\cJ$.
\end{enumerate}
\end{lemma}

{\em Proof.} Ad 1. It is enough to note that it $a,b\in \cJ$ then
$a,b\in T^u-\gamma(T)$ and therefore if $\gamma(a)=\gamma(b)$ or
$\delta(a)\cap\delta(b)\neq\emptyset$ then $a=b$.

Ad 2. Since for $a\in T$, the value $\gamma^{T_{/\cJ}}[a]$ and
$\delta_{/\cJ}[a]$ depend only on $\gamma(a)$ and $\delta(a)$ (and
not on $a$ itself) it is enough to show that if $a\sim'_\cJ b$
then $\gamma(a)=\gamma(b)$ and $\delta(a)=\delta(b)$.

So assume that there is $\alpha\in T$ such that $a=\delta(\alpha)$
and $\gamma(\alpha)=b$. Since $T$ is a positive face structures
$\gamma(a)\cap\delta(a)=\emptyset$. Thus using globularity (in
positive face structures), we have
\[ \gamma(b)=\gamma\gamma(\alpha)=\gamma\delta(\alpha)-\delta\delta(\alpha)=\gamma(a)-\delta(a)=\gamma(a)\]
and
 \[ \delta(b)=\delta\gamma(\alpha)=\delta\delta(\alpha)-\gamma\delta(\alpha)=\delta(a)-\gamma(a)=\delta(a) \]
as required.

Ad 3. This follows immediately from the definition of
$\delta_{/\cJ}([a])$.

Ad 4. We argue by induction on the height $ht(a)$. The inductive
assumption is:
\begin{center} $Ind_n$: for $a\in T-\cJ$ such that $ht(a)=n$ we have:
$a\in\cL_\cJ$ iff $[a]\in T^\lambda_{/\cJ}$. \end{center}

We can assume that $a\in T^u$, as each of the conditions
$a\in\cL_\cJ$ and $[a]\in T^\lambda_{/\cJ}$ implies it.

If $ht(a)=0$ then neither $a\in\cL_\cJ$ nor $[a]\in
T^\lambda_{/\cJ}$. Hence $Ind_0$ holds.

Assume that $ht(a)=1$. Let $\alpha\in T-\gamma(T)$ such that
$\gamma(\alpha)=a$.

Suppose that $a\in\cL_\cJ$. Then $\delta(\alpha)\subseteq \cJ$.
Hence $\delta(a)\sim_\cJ\gamma(a)$ and $[a]\in T^\lambda_{/\cJ}$.

On the other hand, if $[a]\in T^\lambda_{/\cJ}$, then
$[\gamma(a)]=\gamma_{/\cJ}([a])=\delta_{/\cJ}([a])=[\delta(a)]$.
So there is a $\cJ$-path from $\delta(a)$ to $\gamma(a)$. As
$ht(a)=1$ this must be a $\delta(\alpha)$-path. Since it is a
$T^u$-path, we have $\delta(\alpha)\subseteq \cJ$.  Thus
$a=\gamma(\alpha)\in\cL_\cJ$.

Finally, assume that $ht(a)=n>1$. Let $\alpha\in T-\gamma(T)$ such
that $\gamma(\alpha)=a$. As $a\in T^u$, so
$\delta(\alpha)\subseteq T^u$.  Let $a_1,\ldots,a_k$ be the lower
path containing all elements of $\delta(\alpha)$.

First suppose that $a\in\cL_\cJ$. Then $\delta(\alpha)\subseteq
\cL_\cJ\cup \cJ$. If $a_i\in \cJ$ then, by def
$\delta(a_i)\sim_\cJ\gamma(a_i)$. If $a_i\in \cL_\cJ$ then, as
$ht(a_i)<n$, by induction hypothesis $[a_i]$ is a loop. But this
means that $\delta(a_i)=[\gamma(a_i)]$ and in this case again we
have that $\delta(a_i)\sim_\cJ\gamma(a_i)$. By transitivity of
$\sim_\cJ$ we have
$\delta(a)=\delta(a_1)\sim_\cJ\gamma(a_k)=\gamma(a)$, i.e. $[a]$
is a loop in $T_{/\cJ}$, as required.

Now suppose that $[a]\in T^\lambda_{/\cJ}$. Thus there is an upper
$\cJ$-path $\delta(a), b_1,\ldots, b_m,\gamma(a)$.  We claim that
there are numbers $0=m_0<m_1<m_2<\ldots <m_k=m$ such that
\begin{description}
  \item[(i)] either $m_i=m_{i-1}+1$ and $a_i=b_{m_i}\in\cJ$
  \item[(ii)] or $b_{m_{i-1}+1},\ldots,b_{m_i}$ is a path from
  $\delta(a_i)$ to $\gamma(a_i)$ (i.e. $\delta(a_i)\in\delta(b_{m_{i-1}+1})$ and
   $\gamma(a_i)=\gamma(b_{m_i})$).
\end{description}

Having the above claim it follows that either $a_i\in\cJ$ or
$[a_i]\in T^\lambda_{\cJ}$, for $i=1,\ldots ,k$. As $ht(a_i)<n$,
by inductive hypothesis this means that $a_i\in \cL_\cJ\cup \cJ$,
for $i=1,\ldots, k$, i.e. $\delta(\alpha)\subseteq\cL_\cJ\cup
\cJ$. So by definition of $\cL_\cJ$, $a=\gamma(\alpha)\in
\cL_\cJ$, as required.

It remains to prove the claim. Suppose contrary, that the claim is
not true. Let $1\leq i_0\leq k$ be the least number for which it
does not hold, i.e. we have $m_0,m_1,\ldots, m_{i_0-1}$ satisfying
$\bf (i)$ or $\bf (ii)$. In particular, either $i_0=1$ or
$\gamma(b_{m_{i_0-1}})=\gamma(a_{{i_0-1}})=\delta(a_{{i_0}})$. So
$\delta(b_{m_{i_0}})=\delta(a_{{i_0}})$. As $\bf (i)$ does not
hold $a_{i_0}\neq b_{m_{i_0}}$. Since $b_{m_{i_0}}\in
T-\gamma(T)$, $b_{m_{i_0}}<^+a_{i_0}$.  As $\bf (ii)$ does not
hold, by Path lemma (for positive face structures),
$b_i<^+a_{i_0}$, for $i=m_{i_0-1}+1, \ldots ,m$ and
$\gamma(a_k)=\gamma(b_m)\neq\gamma(a_{i_0})$. Again by Path lemma,
the upper path $\delta(a_{i_0}),b_{m_{i_0-1}+1}, \ldots
,b_m,\gamma(a_k)$ can be extended to an upper path reaching
$\gamma(a_{i_0})$:
\[ \delta(a_{i_0}),b_{m_{i_0-1}+1}, \ldots
,b_m,\gamma(a_k),c_1,\ldots , c_r, \gamma(a_{i_0}) \]
 But this means that we have both $a_{i_0}\leq^-a_k$
 and $a_{i_0}<^+a_k$. This contradicts the disjointness and ends the
 proof of the claim and the Lemma.
 $~\Box$

The following lemma describes some relations between pathes in $T$
and $T_{/\cJ}$.

\begin{lemma} Let $T$ be a positive face structure, $\cJ\subseteq T^u$
is an ideal. $x,y,a\in T-\cJ$. Then
\begin{enumerate}
  \item If $a\not\in\gamma(T)$ then $[x]_{\sim_\cJ}\in\delta([a]_{\sim_\cJ})$ iff
  there are $y\in \delta(a)-\cJ$ and a $\cJ$-path (possibly empty) from $x$ to $y$.
  \item  If $a_1,\ldots,a_k$ is a flat path in $T-\gamma(T)$ then
  $<[a_i]_{\sim_\cJ} : 1\leq i \leq k,\: a_i\not\in\cJ\cup\cL_\cJ>$ is a flat path
  in $T_{/\cJ}-\gamma_{/\cJ}(T_{/\cJ}^{-\lambda})$.
  \item Assume $a_1,\ldots, a_k\in T-(\cJ\cup \cL_\cJ)$. Then
  $[a_1]_{\sim_\cJ}, \ldots, [a_k]_{\sim_\cJ}$ is a flat path
  in $T_{/\cJ}-\gamma_{/\cJ}(T_{/\cJ}^{-\lambda})$ iff for $1\leq j<k$
  there there is a $\cJ$-path $b_{j,1},\ldots , b_{j,l_j}$ so that
  \[ a_1,b_{1,1},\ldots , b_{1,l_1},a_2,b_{2,1},\ldots , b_{2,l_2},a_3,
  \ldots ,a_{k-1},b_{k-1,1},\ldots , b_{k-1,l_{k-1}}a_k \]
  is a path in $T$.
  \item $[a]_{\sim_\cJ}\in T_{/\cJ}-\gamma(T_{/\cJ}^{-\lambda})$ iff
  there is  $a'\in T-\gamma(T)$ such that $a'\sim_\cJ a$.
  \item $[x]_{\sim_\cJ}<^{T_{/ \cJ}, +}[y]_{\sim_\cJ}$ iff $x<^{T,+}y$ and
  the upper $(T-\gamma(T))$-path from $x$ to $y$ is not a $\cJ$-path.
\end{enumerate}
\end{lemma}

{\em Proof.} 1. follows easily from pencil linearity of positive
face structures,  2. is easy and 3. is a consequence of 1.

Ad 4. $\Ra$: Let $a' = \min_{<^{T,+}}([a])$, i.e. $a'$ is the
least element of $[a]$.  Suppose that $a'\in\gamma(T)$, i.e. there
is $\alpha\in T$ such that $\gamma(\alpha)=a'$. Clearly, we can
assume that $\alpha\not\in\gamma(T)$. If $\alpha\in\cJ\subseteq
T^u$ then $\delta(\alpha)\sim_{\cJ}a'$ and $\delta(\alpha)<^+a'$,
contrary to the choice of $a'$. If $\alpha\not\in\cJ$ then
$\gamma_{/\cJ}([\alpha])=[a]$ and by assumption $[\alpha]\in
T_{/\cJ}^{\lambda}$. But then by description of the loops is
$T_{/\cJ}$, $\alpha\in\gamma(T)$, contrary to the choice of
$\alpha$. The contradiction shows that $a'\not\in\gamma(T)$, as
required.

$\La$: Suppose $a\in T-\gamma(T)$. We need to show that $[a]\in
T_{/\cJ}-\gamma(T_{/\cJ}^{-\lambda})$.

Suppose contrary, that there is $\alpha\in T-(\cJ\cup\cL_\cJ)$
such that $\gamma_{/\cJ}([\alpha])=a$. Since $a\not\in\gamma(T)$
there is an upper $\cJ$-path $a,\alpha_1,\ldots ,
\alpha_k,\gamma(\alpha)$. Since $\alpha_k\in \cJ\subseteq
T-\gamma(T)$ and $\alpha\not\in\cJ$, by pencil linearity in
positive face structures, we have that $\alpha_k<^k\alpha$.  By
Path Lemma, (since $a\not\in\gamma(T)$) either
$a\in\delta(\alpha)$ or there is $i<k$ such that
$\gamma(\alpha_i)\in\delta(\alpha)$. Using the characterization of
loops in $T_{/\cJ}$ we get in either case that $[\alpha]$ is a
loop in $T_{/\cJ}$, i.e. $\alpha\in\cL_\cJ$ contrary to the
assumption.

Ad 5. The 'if' part is obvious.  We shall show the 'only if' part.
The essential argument consists of showing that in case the path
from $[x]$ to $[y]$ has length $1$ the conclusion hold. Then use
induction.

So assume that $[x],[a],[y]$ is an upper path in $T_{/\cJ}$ with
$[a]\in T_{/\cJ}^{-\lambda}$. Thus we have $x'\in\delta(a)-\cJ$ so
that one of the following four cases holds. There are $\cJ$-pathes
\begin{enumerate}
  \item from $x$ to $x'$ and from $\gamma(a)$ to $y$;
  \item from $x$ to $x'$ and from $y$ to $\gamma(a)$;
  \item from $x'$ to $x$ and from $\gamma(a)$ to $y$;
  \item from $x'$ to $x$ and from $y$ to $\gamma(a)$.
\end{enumerate}
In case 1. the conclusion follows immediately.  The case 4. is
most involved of 2., 3., and 4. and we will deal with this case
only.

Let $x',b_1,\ldots,b_k,x$ and $y,c_1,\ldots,c_r,\gamma(a)$ be
(non-empty) $\cJ$-pathes.  We have
$x'\in\delta(b_1)\cap\delta(a)$. As $b_1\in\cJ$ and $a\not\in\cJ$
we have $b_1<^+a$.  As $[a]\in T^{-\lambda}_{/\cJ}$, and
$b_1,\ldots, b_k$ is a $\cJ$-path we have
$\gamma(b_i)\neq\gamma(a)$, for $i=1,\ldots, k$. By Path Lemma, we
have a (non-empty) upper $T-\gamma(T)$-path $x,b_1',\ldots
,b_l',\gamma(a)$. As $[a]$ is not a loop $b_1',\ldots ,b_l'$ is
not a $\cJ$-path. Since $\gamma(c_r)=\gamma(b_l')$ and both
$c_1,\ldots,c_r$ and $b_1',\ldots ,b_l'$ are $T-\gamma(T)$-pathes,
it follows that one is the end-part of the other.  As the former
is a $\cJ$-path and the latter is not, $c_1,\ldots,c_r$ is the end
$b_1',\ldots ,b_l'$. Thus we have an upper path
$x,b_1',\ldots,b_{l-r}',y$ which is not a $\cJ$-path, as required.
$~\Box$

\vskip 1mm

{\em Proof of Theorem \ref{quotient ofs}.} We shall check that
$T_{/\cJ}$ satisfies all the conditions of the definition of an
ordered face structure.

{\em Local discreteness}, {\em Strictness}, and {\em Loop-filling}
are obvious from the Lemmas above.

{\em Globularity.} First we shall spell the definitions of the
sets involved. For $a\in T_{\geq 2}-\cJ$, we have

\[ \gamma_{/\cJ}\gamma_{/\cJ}([a])= [\gamma\gamma(a)] \]

\[ \delta_{/\cJ}\gamma_{/\cJ}([a])= \left\{ \begin{array}{ll}
          [\gamma\gamma(a)] & \mbox{if  $\delta(a)\subseteq \cJ$,}  \\
          1_{[\gamma\gamma\gamma(a)]} & \mbox{if  $\delta\gamma(a)\subseteq \cJ$,}  \\
          \{ [x] : x\in \delta\gamma(a)-\cJ \} & \mbox{otherwise.}
                                    \end{array}
                \right. \]

 \[ \gamma_{/\cJ}\delta_{/\cJ}([a])= \left\{ \begin{array}{ll}
          [\gamma\gamma(a)] & \mbox{if  $\delta(a)\subseteq \cJ$,}  \\
          \{ [\gamma(x)] : x\in \delta(a)-\cJ \} & \mbox{otherwise.}
                                    \end{array}
                \right. \]

\[ \delta_{/\cJ}\delta_{/\cJ}([a])= \left\{ \begin{array}{ll}
          [\gamma\gamma(a)] & \mbox{if  $\delta(a)\subseteq \cJ$,}  \\
          \{ [u] : \exists_x u\in \delta(x)-\cJ,\; x\in \delta(a)-\cJ \}\cup  & \mbox{}  \\
           \{ 1_{[\gamma\gamma(x)]} : \delta(x)\subseteq \cJ,\; x\in  \delta(a)-\cJ \} & \mbox{otherwise.}
                                    \end{array}
                \right. \]

 \[ \dot{\delta}^{-\lambda}_{/\cJ}([a])= \left\{ \begin{array}{ll}
          \emptyset & \mbox{if  $\delta(a)\subseteq \cJ$,}  \\
          \{ [x] : x\in \delta(a)-(\cL_\cJ\cup\cJ) \} & \mbox{otherwise.}
                                    \end{array}
                \right. \]

 \[ \gamma_{/\cJ}\dot{\delta}^{-\lambda}_{/\cJ}([a])= \left\{ \begin{array}{ll}
          \emptyset & \mbox{if  $\delta(a)\subseteq \cJ$,}  \\
          \{ [\gamma(x)] : x\in \delta(a)-(\cL_\cJ\cup\cJ) \} & \mbox{otherwise.}
                                    \end{array}
                \right. \]

 \[ \delta_{/\cJ}\dot{\delta}^{-\lambda}_{/\cJ}([a])= \left\{ \begin{array}{ll}
          \emptyset & \mbox{if  $\delta(a)\subseteq \cJ$,}  \\
          \{ [u] : \exists_x u\in \delta(x)-\cJ,\; x\in \delta(a)-(\cL_\cJ\cup\cJ) \}\cup  & \mbox{}  \\
           \{ 1_{[\gamma\gamma(x)]} : \delta(x)\subseteq \cJ,\; x\in  \delta(a)-\cJ \} & \mbox{otherwise.}
                                    \end{array}
                \right. \]

Thus if $\delta(a)\subseteq \cJ$ it is easy to see that
globularity holds. So we assume that $\delta(a)\not\subseteq \cJ$.

{\em $\gamma$-globularity}. We shall show:
\begin{enumerate}
  \item $[\gamma\gamma(a)]\in \gamma_{/\cJ}\delta_{/\cJ}([a])$
  \item $[\gamma\gamma(a)]\not\in \delta_{/\cJ}\dot{\delta}^{-\lambda}_{/\cJ}([a])$
  \item
  $\gamma_{/\cJ}\delta_{/\cJ}([a])-[\gamma\gamma(a)]\subseteq \delta_{/\cJ}\dot{\delta}^{-\lambda}_{/\cJ}([a])$
\end{enumerate}

Ad 1.  Let $x_0,x_1,\ldots,x_k$ be a lower $\delta(a)$-path  such
that $\gamma(x_k)=\gamma\gamma(a)$, $x_0\not\in\cJ$ and
$x_i\in\cJ$, for $i>0$ ($k$ is possibly equal $0$).  Such a
sequence exists since $\delta(a)\not\subseteq \cJ$ and is unique
since $x_i\in T^u$, for $i>0$. Then
\[ [\gamma\gamma(a)] = [\gamma(x_0)]\in
\gamma_{/\cJ}\delta_{/\cJ}([a]),\] as required.

Ad 2. Suppose contrary, that there is $x\in
\delta(a)-(\cL_\cJ\cup\cJ)$ (i.e. $[x]$ is not a loop) and
$u\in\delta(x)-\cJ$ such that $[u]=[\gamma\gamma(a)]$. Let
$u,x=x_0,x_1,\ldots, x_k,\gamma\gamma(a)$ be the upper
$\delta(a)$-path from $u$ to $\gamma\gamma(a)$. Since
$u\sim_\cJ\gamma\gamma(a)$ and $u<^+\gamma\gamma(a)$ there is an
upper $\cJ$-path $u,y_0,\ldots,y_l,\gamma\gamma(a)$. As
$x_0\not\in\cJ$ and $u\in\delta(x_0)\cap\delta(y_0)$, by pencil
linearity, we have $y_0<^+x_0$.  By Path Lemma, there is $0\leq
i\leq l$ such, that $\gamma(y_i)=\gamma(x_0)$.  Hence using the
characterization of the loops in $T_{/\cJ}$, $[x]$ is a loop
contrary to our assumption. From this contradiction if follows
that indeed $[\gamma\gamma(a)]\not\in
\delta_{/\cJ}\dot{\delta}^{-\lambda}_{/\cJ}([a])$.

Ad 3. Fix $x\in \delta(a)-\cJ$. Then
$[\gamma(x)]\in\gamma_{/\cJ}\delta_{/\cJ}([a])$. Let
$x_0,\ldots,x_k$ be the $\delta(a)$-path such that
$x_0\not\in\cL_\cJ\cup\cJ$ and $x_i\in\cL_\cJ\cup\cJ$ for $i>0$.
Clearly $x_i\in T^u$ if $i>0$ and possibly $k=0$. Then
$\gamma_{/\cJ}([x_0])=\gamma_{/\cJ}\gamma_{/\cJ}([a])$. So an
arbitrary element of
$\gamma_{/\cJ}\delta_{/\cJ}([a])-[\gamma\gamma(a)]$ is of form
$[\gamma(x)]$ for $x\in\delta(a)-(\cJ\cup\{x_0 \})$. Then we have
a lower $\delta(a)$-path $x=y_1,\ldots,y_l=x_0$ with $l>1$. Put
\[ l'=\max(\{ l'':\{ y_2,\ldots,y_{l''}\} \subseteq (\cL_\cJ\cup \cJ) \}\cup\{ 1\}) \]
As $x_0\not\in\cL_\cJ\cup\cJ$, we have $1\leq l'<l$, and
$y_{l'+1}\not\in\cL_\cJ\cup\cJ$, i.e.
$[y_{l'+1}]\in\dot{\delta}^{-\lambda}_{/\cJ}([a])$. Clearly,
$\gamma(y_{l'})\in\delta(y_{l'+1})-\cJ$ and hence
\[ [\gamma(x)]=[\gamma(y_{l'})]\in\delta_{/\cJ}\dot{\delta}^{-\lambda}_{/\cJ}([a]),\]
which ends the proof of $\gamma$-globularity.

{\em $\delta$-globularity.} We have there different cases:
\begin{description}
  \item[I] $\delta(a)\subseteq \cJ$,
  \item[II] $\delta(a)\not\subseteq \cJ$ and $\delta\gamma(a)\subseteq
  \cJ$,
  \item[III] $\delta(a)\not\subseteq \cJ$ and $\delta\gamma(a)\not\subseteq
  \cJ$.
\end{description}
Case {\bf I}, as we already mentioned, is obvious.

Case {\bf II}: $\delta(a)\not\subseteq \cJ$ and
$\delta\gamma(a)\subseteq \cJ$.

In this case we have:
\[ \delta_{/\cJ} \gamma_{/\cJ}([a])=1_{[\gamma\gamma\gamma(a)]}. \]

\[ \delta_{/\cJ}\delta_{/\cJ}([a])=  \{ [u] : \exists_{x\in\delta(a)-\cJ}\; u\in \delta(x)-\cJ  \}\cup
           \{ 1_{[\gamma\gamma(x)]} : \delta(x)\subseteq \cJ,\; x\in  \delta(a)-\cJ \} \]
and
\[ \gamma_{/\cJ}\dot{\delta}^{-\lambda}_{/\cJ}([a]) =\{ [\gamma(x)] : x\in \delta(a)-(\cL_\cJ\cup\cJ) \}.  \]

Let $u\in \delta(x)-\cJ$ and $x\in\delta(a)-\cJ$. As
$\delta\gamma(a)\subseteq \cJ$ and $u\not\in\cJ$, by globularity
(of positive face structures), $u\in\gamma\delta(a)$. Thus there
is $y_0\in\delta(a)$ such, that $\gamma(y_0)=u$. Since
$\delta(\cJ)\cap\cJ=\emptyset$ and $\delta\gamma(a)\subseteq \cJ$
there is a $\delta(a)$-path $y_k,\ldots,y_0$ such that
$y_k\not\in(\cJ\cup\cL_\cJ)$ and
$y_{k-1},\ldots,y_0\subseteq(\cJ\cup\cL_\cJ)$, $k\geq 0$.  Then
$[y_k]\in\dot{\delta}^{-\lambda}([a])$ and
$\gamma_{/\cJ}([y_k])=[u]$. Thus $\{ [u] :
\exists_{x\in\delta(a)-\cJ} u\in\delta(x)-\cJ \}\subseteq
\gamma_{/\cJ}\dot{\delta}_{/\cJ}^{-\lambda}([a])$.

It remains to show that
\begin{enumerate}
  \item there is $x\in\delta(a)-\cJ$, $\delta(x)\subseteq \cJ$,
  \item for any such $x$, we have
  $1_{[\gamma\gamma(x)]}=1_{[\gamma\gamma\gamma(a)]}$.
\end{enumerate}

Ad 1. The existence of such $x$ follows easily from Path Lemma.

Ad 2. Suppose $x\in\delta(a)-\cJ$, $\delta(x)\subseteq \cJ$.  As
$\cJ\cap\gamma(T)=\emptyset$, by globularity, we have that
$\delta(x)\subseteq\delta\gamma(a)$. Hence
$\gamma\gamma(x)\in\gamma\delta(x)\subseteq\gamma\delta\gamma(a)$,
and then there is an upper $\delta\gamma(a)$-path (possibly empty)
from $\gamma\gamma(x)$ to $\gamma\gamma\gamma(a)$. But
$\delta\gamma(a)\subseteq\cJ$, so this is a $\cJ$-path and this
means that $\gamma\gamma(x)\sim_\cJ\gamma\gamma\gamma(a)$, i.e.
$[\gamma\gamma(x)]=[\gamma\gamma\gamma(a)]$, as required.

Case {\bf III}: $\delta(a)\not\subseteq \cJ$ and
$\delta\gamma(a)\not\subseteq \cJ$. This is the only case, where
we do not have equality but $\equiv_1$ only.  We need to verify:
\begin{enumerate}
  \item $\delta_{/\cJ}\gamma_{/\cJ}([a])\subseteq \delta_{/\cJ}\delta_{/\cJ}([a])$;
   \item $\delta_{/\cJ}\gamma_{/\cJ}([a])\cap\gamma_{/\cJ}\dot{\delta}^{-\lambda}_{/\cJ}([a])=\emptyset$;
  \item $\delta_{/\cJ}\delta_{/\cJ}([a])-
  \delta_{/\cJ}\gamma_{/\cJ}([a])\subseteq\gamma_{/\cJ}\dot{\delta}^{-\lambda}_{/\cJ}([a])$;
  \item $\gamma_{/\cJ}\delta_{/\cJ}\delta^\varepsilon_{/\cJ}([a])\subseteq {\theta_{/\cJ}\delta_{/\cJ}\gamma_{/\cJ}([a])}$.
\end{enumerate}

Ad 1. We have $\delta_{/\cJ}\gamma_{/\cJ}([a])=\{ [u] :
u\in\delta\gamma(a)-\cJ \}$.

So, let $u\in\delta\gamma(a)-\cJ$ and
$u,x_1,\ldots,x_k,\gamma\gamma(a)$ be a $\delta(a)$-path, $k\geq
1$.  There is $1\leq l\leq k$ such that $x_i\in\cJ$, for $i<l$,
and $x_l\not\in\cJ$. Such $l$ exists, since
$\delta(a)\not\subseteq\cJ$. Let
\[ v= \left\{ \begin{array}{ll}
          u & \mbox{if  $l=1$,}  \\
          \gamma(x_{l-1}) & \mbox{otherwise.}
                                    \end{array}
                \right. \]
Then $[u]=[v]$. Moreover, $v\in\delta(x_{l+1})-\cJ$,
$x_{l+1}\in\delta(a)-\cJ$, i.e.
$[v]\in\delta_{/\cJ}\delta_{/\cJ}([a])$, as required.

Ad 2. Suppose contrary that there is $u\in\delta\gamma(a)-\cJ$ and
$x\in\delta(a)-(\cL_\cJ\cup\cJ)$ so that
$[u]=[\gamma(x)]\in\gamma_{/\cJ}\delta_{/\cJ}([a])$.

Thus we have a $\cJ$-path $u,x_1,\ldots,x_k,\gamma(x)$. As
$\delta\gamma(a)\cap\gamma\delta(a)=\emptyset$, $u\neq\gamma(x)$
and $k\geq 1$. Since $\gamma(x)=\gamma(x_k)$, $x_k\in\cJ$ and
$x\not\in\cJ$, by pencil linearity $x_k<^+x$. We have that
$\delta(x_l)\cap\delta(x)=\emptyset$ for $1\leq l\leq k$, since
otherwise $[x]$ would be a loop. Let
$y_1,\ldots,y_r,x_1,\ldots,x_k$ be a continuation of the path
$y_1,\ldots,y_r,x_1,\ldots,x_k$ through $u$ (i.e. $\gamma(y_r)=u$)
such that there is $v\in\delta(y_1)\cap\delta(x)$. Since
$x\in\delta(a)$, $v\in\delta\delta(a)$. So there is
$v'\in\delta\gamma(a)$ such that $v'\leq^+ v$.  But then $v'<^+u$
and $v',u\in\delta\gamma(a)$, which is impossible.  Thus 2. holds,
as well.

Ad 3. Let $u\in\delta(x)-\cJ$ and $x\in\delta(a)-\cJ$, i.e.
$[u]\in\delta_{/\cJ}\delta_{/\cJ}([a])$ and suppose that
$[u]\not\in \delta_{/\cJ}\gamma_{/\cJ}([a])$. Let
$x_0,\ldots,x_l,u$ be a $\delta(a)$-path, $l\geq 0$, such that
$x_1,\ldots,x_l\subseteq (\cL_\cJ\cup\cJ)$, and $x_0\not\in
(\cL_\cJ\cup\cJ)$. Such a path exists since $[u]\not\in
\delta_{/\cJ}\gamma_{/\cJ}([a])$.  Then
$[x_0]\in\dot{\delta}_{/\cJ}^{-\lambda}([a])$ and hence
$[u]=[\gamma(x_0)]\in
\gamma_{/\cJ}\dot{\delta}^{-\lambda}_{/\cJ}([a])$, as required.

Ad 4. We have
\[ \delta_{/\cJ}\delta^{\varepsilon}_{/\cJ}([a])
= \{ 1_{[\gamma\gamma(x)]} :
x\in\delta(a)-\cJ,\,\delta(x)\subseteq\cJ \}\]

Fix $x\in\delta(a)-\cJ$ such that $\delta(x)\subseteq\cJ$. We need
to show that $\gamma_{/\cJ}\gamma_{/\cJ}([x])=[\gamma\gamma(x)]\in
\theta_{/\cJ}\delta_{/\cJ}\gamma_{/\cJ}([a])$.

We have
$\gamma\gamma(x)\in\gamma\gamma\delta(a)\subseteq\gamma\gamma\gamma(a)\cup\delta\delta\gamma(a)$.
If $\gamma\gamma(x)=\gamma\gamma\gamma(a)$ then, using
$\gamma$-globularity of positive face structures, we have
\[
[\gamma\gamma(x)]=[\gamma\gamma\gamma(a)]=\gamma_{/\cJ}\gamma_{/\cJ}\gamma_{/\cJ}([a])\subseteq
\gamma_{/\cJ}\delta_{/\cJ}\gamma_{/\cJ}([a])\subseteq\theta_{/\cJ}\delta_{/\cJ}\gamma_{/\cJ}([a])\]
and 4. holds.

So now assume that $\gamma\gamma(x)\in\delta\delta\gamma(a)$. Thus
there is an upper $\delta\gamma(a)$-path
$\gamma\gamma(x),u_1,\ldots,u_k,\gamma\gamma\gamma(a)$. If it is a
$\cJ$-path then
$[\gamma\gamma(x)]=[\gamma\gamma\gamma(a)]\in\theta_{/\cJ}\delta_{/\cJ}\gamma_{/\cJ}([a])$.
If it is not a $\cJ$-path, then let
$i_0=\min\{i':u_{i'}\not\in\cJ\}$ and
\[ t= \left\{ \begin{array}{ll}
          \gamma\gamma(x) & \mbox{if  $i_0=1$,}  \\
          \gamma(u_{i_0-1}) & \mbox{otherwise.}
                                    \end{array}
                \right. \]
Then $u_1,\ldots,u_{i_0-1}$ is a $\cJ$-path and
$\gamma\gamma(x)\sim_\cJ t\in\delta(u_i)-\cJ$. Thus
$[\gamma\gamma(x)]\in\delta_{/\cJ}([u_{i_0}])$.  But
$u_{i_0}\in\delta\gamma(a)-\cJ$, so
$[u_{i_0}]\in\delta_{/\cJ}\gamma_{/\cJ}([a])$ and then
\[\gamma\gamma(x)]\in\delta_{/\cJ}\delta_{/\cJ}\gamma_{/\cJ}([a])\subseteq\theta_{/\cJ}\delta_{/\cJ}\gamma_{/\cJ}([a]),\]
as required. This ends the proof of globularity of $T_{/\cJ}$.

{\em Disjointness.} From the description of the orders we get
immediately
$\perp^{T_{/\cJ},+}\cap\perp^{T_{/\cJ},\sim}=\emptyset$. If
$a,b\in T-\cJ$ and $[a]<^\sim[b]$ the by definition $a<^{T,-}b$.
So we have a lower $T-\gamma(T)$-path $a=a_0,\ldots,a_k=b$. Let
$b_1,\ldots, b_l$ be the path (possibly empty) obtained from
$a_1,\ldots,a_{k-1}$ by dropping elements that belong to
$\cL_\cJ\cup\cJ$. Then $[a],[b_1],\ldots,[b_l],[b]$ is a lower
flat path in $T_{/\cJ}$. Hence $[a]<^\sim[b]$ implies $[a]<^-[b]$.

Now assume that
$\theta_{/\cJ}([a])\cap\theta_{/\cJ}([b])=\emptyset$ and
$[a]<^-[b]$. We need to show that $[a]<^\sim[b]$.

So we have a lower flat path $[a]=[a_0],\ldots,[a_k]=[b]$ in
$T_{/\cJ}$, with $k>1$. There are some cases to be considered. We
will deal with the one which is most involved: $k=2$, there are
$x\in \delta(a_1)-\cJ$, $y\in \delta(a_2)-\cJ$ there are upper
$\cJ$-pathes $x,b_1,\ldots,b_l,\gamma(a_0)$, and
$y,c_1,\ldots,c_r,\gamma(a_1)$.

As $[a_1]$ is not a loop $b_i<^+a_1$, for $i\leq l$ and
$c_i<^+a_1$, for $i\leq r$. Moreover they are
$T-\gamma(T)$-pathes.  It is easy to see that if we continue the
path $b_1,\ldots,b_l$ as $T-\gamma(T)$-path we shall get to
$c_1,\ldots,c_r$ (note that $\cJ$-faces are unary). Thus there is
a path $b_1,\ldots,b_l,d_1,\ldots,d_s ,c_1,\ldots,c_r$, with
$s\geq 0$. Then $a_0,d_1,\ldots,d_s ,a_2$ is a lower path showing
that $a=a_0<^-a_2=$, i.e. $[a]<^\sim[b]$, as required. This ends
the proof of disjointness.

{\em Pencil linearity.} Let $[a], [b]\in T_{/\cJ}$, and $[a]\neq
[b]$. First assume that
$\theta_{/\cJ}([a])\cap\theta_{/\cJ}([b])\neq\emptyset$. Thus we
have three cases to consider:
\begin{description}
  \item[I] $\gamma_{/\cJ}([a])=\gamma_{/\cJ}([b])$,
  \item[II] $\gamma_{/\cJ}([a])\in\delta_{/\cJ}([b])$,
  \item[III] $\delta_{/\cJ}([a])\cap\delta_{/\cJ}([b])\neq\emptyset$.
\end{description}

Case {\bf I}. Possibly changing the roles of $a$ and $b$ there is
a $\cJ$-path $\gamma(a),a_1,\ldots ,a_k,\gamma(b)$.  If
$\gamma(a)\in\delta(b)$ or there is $i<k$ such that
$\gamma(a_i)\in\delta(b)$ then $a<^{T,-}b$ and hence
$[a]<^{T_{/\cJ},-}[b]$. If it is not the case, that is
$\gamma(a)\not\in\delta(b)$ and for all $i<k$
$\gamma(a_i)\not\in\delta(b)$ then by Path Lemma there is
$y\in\delta(b)$ and an upper path
$y,b_1,\ldots,b_l,a,a_1,\ldots,a_k,\gamma(b)$ and $a<^+b$.
Therefore $[a]<^+[b]$.

Case {\bf II}. In this case we have $x\in\delta(b)-\cJ$ and a
$\cJ$-path $\gamma(a),a_1,\ldots,a_k,x$ or
$x,a_1,\ldots,a_k,\gamma(a)$.  In the former case we have
$a<^{T,-}b$ and hence $[a]<^{T_{/\cJ},-}[b]$. In the later case
either there is $1\leq i\leq k$ such that $\gamma(b)=\gamma(a_i)$
and then $\gamma_{/\cJ}([b])=\gamma_{/\cJ}([a])$, i.e. this case
is reduced to {\bf I} or, by Path Lemma, we have that $a_i<^+b$,
for $1\leq i\leq k$.  Let $a_k,\alpha_1,\ldots,\alpha_l,b$ be an
upper path in $T$. As $\gamma(a_k)\neq\gamma(b)$ there is $1\leq j
\leq l$ such that $\gamma(a_k)\in\iota(\alpha_j)$.  Since
$\gamma(a_k)=\gamma(a)$ we have that $\gamma(a)\in\iota(\alpha_j)$
and then $a<^+\gamma(\alpha_j)\leq^+\gamma(\alpha_l)=b$, i.e.
$[a]<^+[b]$.

Case {\bf III}. Possibly changing the roles of $a$ and $b$ there
are $x\in\delta(a)-\cJ$ and $y\in\delta(b)-\cJ$ and  $\cJ$-path
$x,a_1,\ldots ,a_k,y$. If there is $1\leq i\leq k$ such that
$\gamma(a_i)=\gamma(a)$ then $a<^-b$ and hence $[a]<^\sim [b]$. If
for all $i\leq n$ we have  $\gamma(a_i)\neq\gamma(a)$ then, by
Path Lemma, $b<^+a$ and hence $[b]<^+[a]$.

Next let assume that $[a]\in T_{/\cJ}^\varepsilon$, $[b]\in
T_{/\cJ}$ and
$\gamma_{/\cJ}\gamma_{/\cJ}([a])\in\iota_{/\cJ}([b])$. Thus there
are $x,y\in \delta(b)-(\cJ\cup\cL_\cJ)$ such that
$\gamma_{/\cJ}\gamma_{/\cJ}([a])=\gamma_{/\cJ}([x])
\in\delta_{/\cJ}([y])$. Hence there is $u\in\delta(y)-\cJ$ such
that $\gamma(x)\sim_\cJ u$. If we were to have a $\cJ$-path
$u,x_1\ldots,x_k,\gamma(x)$ then, as $u\in\delta(b)$ by Path
Lemma, either there is $1\leq i\leq k$ such that
$\gamma(x_i)=\gamma(y)$ or $\gamma(x_k)\neq\gamma(y)$ and
$x_k<^+y$. In the former case $[y]$ would be a loop in the later,
we would have $x<^+y$ and $x,y\in\delta(b)$. As none of the above
is possible, it follows that we cannot have a $\cJ$-path from $u$
to $\gamma(x)$. Hence we have a $\cJ$-path
$\gamma(x),x_1\ldots,x_k,u$.

 {\em Claim.} Exactly one of the following conditions holds:
\begin{description}
  \item[(i)] there is a $\cJ$-path (possibly empty) from
  $\gamma\gamma(a)$ to $\gamma(x)$;
  \item[(ii)] there is a $\cJ$-path (possibly empty) from
  u to $\delta\gamma(a)$;
  \item[(iii)] $\delta(a)\subseteq\{ x_1,\ldots,x_k\}$.
\end{description}

Clearly no two of the above three conditions can hold
simultaneously. We assume that {\bf (i)} and {\bf (ii)} does not
hold and we shall prove {\bf (iii)}. We can assume that $k\geq 1$.
As $\delta(a)\subseteq \cJ$ and $\{ x_1,\ldots,x_k\}\subseteq
\cJ$, it is enough to show:
\begin{description}
  \item[(a)] there is $1\leq i \leq k$ such that
  $\gamma(x_i)=\gamma\gamma(a)$;
  \item[(b)] either $\gamma(x)\in\delta\gamma(a)$ or there is $1\leq j < i$ such that
  $\gamma(x_j)=\delta\gamma(a)$;
\end{description}

Ad {\bf (a)}. Suppose that {\bf (a)} does not hold. Then, as {\bf
(i)} does not hold, we have an upper $\cJ$-path
$u,x_1,\ldots,x_l,\gamma\gamma(a)$, with $l>k$. As $x_l\in\cJ$ and
$\gamma(a)\not\in\cJ$, we have $x_l<^+\gamma(a)$. So by Path
Lemma, either $\gamma(x_{i_0})\in\delta\gamma(a)$, for some $k\leq
i_0<l$ or $u\in\iota(a)$. In the former case we get {\bf (ii)}
contrary to the supposition.  In the later case, on one hand, as
$u\in\delta(y)\cap\iota(a)$, we have that $y<^+\gamma(a)$. Thus
$\gamma(y)\leq^+\gamma\gamma(a)$. On the other hand, if we were to
have $k<i_1\leq l$ such that $\gamma(x_{i_1})=\gamma(y)$ then, as
$x_i\in\cJ$, we would have that $[y]$ is a loop. Hence, by Path
Lemma $x_i<^+y$, for $i=k+1,\ldots,l$, and
$\gamma\gamma(a)=\gamma(x_l)\neq\gamma(y)$.  But then, again by
Path Lemma, $\gamma\gamma(a)<^+\gamma(y)$.  Therefore we get a
contradiction once more.  This ends the proof of {\bf (a)}.

Ad {\bf (b)}. As we have established {\bf (a)}, let us fix $1\leq
i_1\leq k$ such that $\gamma(x_{i_1})=\gamma\gamma(a)$.  Suppose
that {\bf (b)} does not hold.  Then $\{
x_1,\ldots,x_{i_1}\}\subseteq \delta(a)$ and
$\gamma(x)\in\iota(a)$. Thus $x<^+\gamma(a)$. Let
$y_1,\ldots,y_r,x_1,\ldots,x_{i_1}$ be the lower $\cJ$-path
consisting of all the faces in $\delta(a)$.  Clearly,
$\delta(y_1)=\delta\gamma(a)$, $\gamma(y_r)=\gamma(x)$ and
$y_r<^+x$. If we were to have $1\leq i\leq r$ such that
$\delta(y_i)\cap\delta(x)\neq\emptyset$ then the face $[x]$ would
be a loop contrary to the supposition. Thus, by Path Lemma,
$y_i<^+x$, for $i=1,\ldots,r$ and there is $v\in\delta(x)$ such
that $v<^+\delta(y_1)=\delta\gamma(a)$ (both $\delta(y_1)$ and
$\delta\gamma(a)$ are singletons). On the other hand, as
$x<^+\gamma(a)$, we have, by Path Lemma, a lower path
$z_1,\ldots,z_s=x$, with $s\geq 1$, and
$w\in\delta(z_1)\cap\delta\gamma(a)$. Then, for
$w'=\gamma(z_{s-1})\in\delta(x)$ (or $w'=w$ if $s=1$) we have
$\delta\gamma(a)\leq^+w'$. Thus $v,w'\in\delta(x)$ and
$v<^+\delta\gamma(a)\leq^+w'$.  But this is impossible by
Proposition 5.1 of \cite{Z}. This ends the proof of {\bf (b)} and
of the Claim.

Having the Claim, it is easy to see, that in case {\bf (i)}
$\gamma(a)\leq x$ and in case {\bf (ii)} $\gamma(a)\leq y$. Thus
in both cases we have $[a]<^+[b]$.

Finally assume that {\bf (iii)} holds. Thus $x<^-\gamma(a)$ and
$\gamma\gamma(a)<^+\gamma\gamma(b)$. From the later and \cite{Z},
we have that either $\gamma(a)<^-\gamma(b)$ or
$\gamma(a)<^+\gamma(b)$. If $\gamma(a)<^-\gamma(b)$ then using the
former and transitivity of $<^-$ we would have $x<^-\gamma(b)$.
But $x<^+\gamma(b)$ and this contradicts disjointness.  Thus
$\gamma(a)<^+\gamma(b)$. Then, again by \cite{Z}, we have that
either $a<^-b$ or $a<^+b$, as required.

The fact that $q_\cJ$ is a collapsing morphism with the kernel $\cJ$ is
left for the reader. This ends the proof of the theorem. $~\Box$

\begin{proposition}\label{cover prop}
Let $T$ be a positive face structure, and $\cI$ be an ideal in
$T$. Then $size(T_{/\cI})=size(T)$.
\end{proposition}
{\it Proof.}~We have a quotient morphism $q:T\ra T_{/\cJ}$ such that,
for $a\in T-\cJ$, we have $q(a)=[a]_\cJ$. We shall show that, for
any $k\in\o$,  the restriction of this function
\[ \tilde{q}_k : T_k-\delta(T_{k+1})\lra
(T_{/\cJ})_k-\delta_{/\cJ}((T_{/\cJ}^{-\lambda})_{k+1}) \] is a
bijection. This is clearly sufficient to establish 2. To see that
$\tilde{q}_k$ is one-to-one, note that for $a,a'\in
T_k-\delta(T_{k+1})$, by Lemma \ref{quotient1}.1, we have
$a\not\sim_\cJ a'$.

We shall verify that $\tilde{q}_k$ is onto. Fix
$[a]_\cJ\in(T_{/\cJ})_k-\delta_{/\cJ}((T_{/\cJ}^{-\lambda})_{k+1})$
such that $a\in T-\cJ$ is $<^+$-maximal in its class $[a]_\cJ$.
Suppose that $a\in\delta(T_{k+1})$, and fix $\alpha\in
T_{k+1}-\gamma(T_{k+2})$ such that $a\in\delta(\alpha)$. Then if
$\alpha\in \cJ$, $a$ is not $<^+$-maximal in its class. If
$\alpha\not\in \cJ$, then $[\alpha]_\cJ\in T_{/\cJ}^{-\lambda}$
and $[a]_\cJ\in\delta_{/\cJ}([\alpha]_\cJ)$. In either case we get
a contradiction.  Thus there is no $\alpha\in T_{k+1}$ such that
$a\in\delta(\alpha)$, and $\tilde{q}_k$ is onto indeed. $~\Box$

\subsection*{Positive covers}
Recall that the kernel of a collapsing morphism $q:Y\ra X$ is the set
$ker(q)=q^{-1}(1_X)\subseteq Y$. In more concrete terms, as $q$
preserves codomains, we have $ker(q)=\{ a\in Y: q(a)=1_{q(a)}\}$.

We say that a collapsing morphism $q:Y\ra X$ is a {\em (positive)
cover}\index{positive cover} iff there is an ideal $\cJ$ in $Y$,
and a monotone isomorphism $h :Y/_\cJ\lra X$ such that the
triangle
\begin{center}
\xext=600 \yext=350
\begin{picture}(\xext,\yext)(\xoff,\yoff)
\settriparms[1`1`1;300]
 \putAtriangle(150,0)[Y`Y/_\cJ`X;p_\cJ`q`h]
\end{picture}
\end{center}
commutes.

\begin{proposition}\label{cover}
Let $q:Y\ra X$ be a collapsing morphism  and $\cJ$ an ideal in $Y$, and
$p_\cJ :Y\ra Y/_\cJ$ a positive cover.
\begin{enumerate}
   \item $ker(q)$ is an ideal iff $ker(q)\subseteq Y^u$.
  \item $q:Y\ra X$ is a positive cover iff $q$ is onto and $ker(q)$ is an ideal.
  \item If $ker(p_\cJ)\subseteq ker(q)$ then there is a unique collapsing morphism
$r : Y/_\cJ\ra X$ making the triangle
 \begin{center} \xext=800 \yext=350
 \begin{picture}(\xext,\yext)(\xoff,\yoff)
  \settriparms[1`1`1;300]
  \putAtriangle(0,50)[Y`Y/_\cJ`X; p_\cJ`q`r]
  \end{picture}
 \end{center}
 commutes.
\end{enumerate}
 \end{proposition}
{\it Proof.}~ Ad 1. Any ideal in $Y$ is contained in $Y^u$. Thus
we need to show that if $ker(q)\subseteq Y^u$ then
$ker(q)\cap\gamma(Y)=\emptyset = ker(q)\cap\delta(ker(q))$.

Suppose there is $a\in ker(q)\cap\gamma(Y)$. Let $\alpha\in Y$
such that $\gamma(\alpha)=a$. Then
\[ 1_{\gamma(q(a))}=q(a)=q(\gamma(\alpha))=\gamma(q(\alpha))\in X\]
and we get a contradiction. Thus $ker(q)\cap\gamma(Y)=\emptyset$.

Now suppose that $a\in ker(q)\cap\delta(ker(q))$. Fix $\alpha\in
ker(q)$ such that $a\in\delta(a)$.  As $ker(q)\subseteq Y^u$ we
have $a=\delta(\alpha)$. So we have
\[ 1_{\gamma(q(a))}=q(a)=q(\delta(\alpha))=\delta(q(a))=
\delta(1_{\gamma(q(\alpha))})=\gamma(q(\alpha))\in X \] and we get
a contradiction again.

Ad 2. Clearly the conditions are necessary.  To see that they are
sufficient it is enough to note that they  imply that the map
$h:Y_{/ker(q)}\ra X$ such that $h([a])=q(a)$, for $a\in Y-\cJ$ is
an isomorphism in $\ofs$.

Ad 3. We put $r([y])=q(y)$, for $y\in Y-\cJ$. Since
$ker(p_\cJ)\subseteq ker(q)$,  $r$ is well defined. As
\[  p_\cJ(y) = \left\{ \begin{array}{ll}
            [y]             & \mbox{if $y\in Y-\cJ$,}  \\
            1_{[\gamma(y)]} & \mbox{if $y\in\cJ$.}
                                    \end{array}
                \right. \]
we have $q=r\circ p_\cJ$.  It remains to verify that $r$ preserves
$\gamma$, $\delta$, and $<^\sim$.

Fix $y,y'\in Y-\cJ$. We have
\[ [y]<^{Y_{/\cJ},\sim}[y']\;\; \mbox{ iff }\;\; y<^{Y,\sim}y'
\;\;\mbox{ iff }\;\; q(y)<^{X,\sim}q(y') \;\;\mbox{ iff }\;\;
r([y])<^{X,\sim}r([y']) \] i.e. $r$ preserves $<^\sim$.

Now fix $y\in Y_{\geq 1}-\cJ$. We have
\[ r(\gamma([y])) = r([\gamma(y)]) = q(\gamma(y)) =
\gamma(q(y)) = \gamma(r([y])) \] i.e. $r$ preserves $\gamma$.

To see that $r$ preserves $\delta$ we consider two cases:
$\delta(y)\subseteq \cJ$ and $\delta(y)\not\subseteq \cJ$. If
$\delta(y)\subseteq \cJ$ then we have
\[ r(\delta([y])) = r(1_{[\gamma\gamma(y)]}) = 1_{r([\gamma\gamma(y)])} =
1_{q(\gamma\gamma(y))} = q(\delta(y)) \equiv_1 \delta(q(y)) =
\delta(r([y])) \] and if $\delta(y)\not\subseteq \cJ$ we have
\[ r(\delta([y])) = r(\{[u] : u\in \delta(y)-\cJ \}) =\{ q(u) : u\in \delta(y)-\cJ \}  \equiv_1 \]
\[ \equiv_1  \{ q(u) : u\in \delta(y) \} = q(\delta(y)) \equiv_1 \delta(q(y)) = \delta(r([y]). \]
Thus in both cases $\delta$ is preserved.
 $~\Box$

\section{Positive covers of ordered face structures}\label{positive covers}

In this section we describe a kind of inverse construction to the
quotient construction from previous section.  We shall show that
any ordered face structure $S$ can be covered by a positive one
$S^\dag$. We begin with some notation and the construction. Then
we shall prove few technical lemmas. Using these lemmas we shall
describe the properties of the construction, in particular that
$S^\dag$ is a positive face structure and that $q_S:S^\dag\ra S$
is a quotient morphism.  Finally, we will make some farther comments
about this construction.

\subsection*{The construction of $S^\dag$}

$S$ an ordered face structure fixed for the whole section. The
construction of $S^\dag$ is based on cuts, but this time we
consider the cuts of initial faces in $S$ not, as in section
\ref{convex subhyp}, of empty loops. We use essentially the same
notation for both cuts of empty loops and cuts of initial faces.
But, as we never use these different cuts in the same context so
there is no risk to mix them.

Recall from section \ref{face stuctures},  that
$\cI=\cI^S=S^\varepsilon-\gamma(S^{-\lambda})$ is the set of {\em
initial faces} in $S$, and $\cI_x=\cI^S_x=\{ a\in \cI:
\delta(a)=1_x\}$ is the set of {\em initial faces
based\index{face!-s based on $x$} on $x$}. $\cI_x$ is a linearly
ordered by $<$ (we have, for $a,b\in\cI_x$, that $a<b$ iff
$\gamma(a)<^\sim\gamma(b)$). An {\em $x$-cut}\index{cut} is a
triple $(x,L,U)$ such that $L\cup U=\cI_x$ and for $\alpha\in L$
and $\beta\in U$, $\alpha<\beta$. $\cC(\cI_x)$ is the set of all
$x$-cuts. $\cC(\cI_X)=\cC(\cI^T_X) = \bigcup \{ \cC(\cI_x) : x\in
X \}$, where $X\subseteq S$, is the set of all $X$-cuts, i.e. all
$x$-cuts with $x\in X$.

If $(x,L,U)$ is a $x$-cut then $L$ determines $U$ and vice versa
($L=\cI^X_x-U$ and $U=\cI^X_x-L$). Therefore we sometimes denote
this cut by describing only the lower cut $(x,L,-)$ {\em lower
description of the cut}\index{cut! lower description of -} or only
the upper cut $(x,-,U)$ {\em upper description of the
cut}\index{cut! upper description of -}, whichever is easier to
define.

For $a\in S$ and $x\in\dot{\delta}(a)$. We define the following
sets:
\[ \ua a\;=\{ \alpha\in \cI_{\gamma(a)}: a<^\sim \gamma(\alpha) \},\hskip 5mm
\da_x a\;=\{ \alpha\in \cI_x: \gamma(\alpha)<^\sim a \} \] and
cuts $(\gamma(a),-,\ua a)$ and $(x,\da_x a,-)$.  In order to save
the space we drop subscript $x$ in the notation $\da_x a$ inside
$x$-cuts, i.e. we often write $(x,\da a,-)$ instead of $(x,\da_x
a,-)$. Clearly, $(x,-,\ua b)=(x,\da_x a,-)$ iff $\da_x a\cup\ua
b=\cI_x$ and $\da_x a\cap\ua b=\emptyset$.

We describe below the positive hypergraph $S^\dag$. The set of
faces of dimension $k$ is
\[ S^\dag_k=\cC(\cI_{S_k})\cup \overline{\cI_{k+1}},\]
where $\overline{\cI_{k+1}}$ is another copy of the set
$\cI_{k+1}$ whose elements have bars on it, i.e.
$\overline{\cI_{k+1}}=\{ \overline{\alpha}:\alpha\in\cI_{k+1}\}$.
Thus the faces of each dimension are of two disjoint kinds: cuts
and bars. The  domains and codomains in $S^\dag$ we define
separately for bars and cuts. Fix $k>0$.  For
$\overline{\alpha}\in\overline{\cI_{k+1}}$ we have
 \[ \gamma^{\dag}(\overline{\alpha})=(\gamma\gamma(\alpha), -,\ua\gamma(\alpha)),\hskip 5mm
 \delta^{\dag}(\overline{\alpha})=(\gamma\gamma(\alpha),\da\gamma(\alpha), -),\]
  for $(a,L,U)\in \cC(\cI_a)$, with $a\in S_k$ we have
  \[ \gamma^{\dag}(a,L,U)= (\gamma(a), -, \ua a), \hskip 5mm
   \delta^{\dag}(a,L,U)= \overline{\cI^{\leq^+a}} \cup
  \{(x, \da a ,-) : x\in\dot{\delta}(a)  \}.  \]

We have a map $q_S:S^\dag\lra S$ such that
\[ q_S(z)= \left\{ \begin{array}{ll}
           a & \mbox{if  $z=(a,L,U)\in \cC(\cI_a)$,}\\
           1_{\gamma\gamma(\alpha)}& \mbox{if $z=\overline{\alpha}\in\cI$.}
                                    \end{array}
                \right. \]
i.e. it sends $a$-cuts to $a$, and any bar $\overline{\alpha}$ to
an empty-face  $1_{\gamma\gamma(\alpha)}$.

{\em Example.} The positive cover of the ordered face structure
$T$ as below
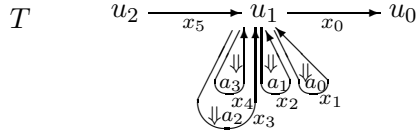
\begin{center} \xext=1320 \yext=520
\begin{picture}(\xext,\yext)(\xoff,\yoff)
\putmorphism(0,450)(1,0)[u_2`u_1`_{x_5}]{480}{1}b
\putmorphism(480,450)(1,0)[\phantom{u_1}`u_0`_{x_0}]{480}{1}b
 \put(-400,400){$T$}

 \put(350,150){\oval(200,200)[b]}
 \put(285,40){$^{\Da a_2}$}
 \put(250,150){\line(1,2){120}}
 \put(450,150){\vector(0,1){250}}
  \put(440,30){$^{x_3}$}

  \put(360,220){\oval(100,100)[b]}
  \put(310,220){\line(1,2){85}}
  \put(410,220){\vector(0,1){180}}
  \put(360,220){$^\Da$}
    \put(325,155){$^{a_3}$}
     \put(365,90){$^{x_4}$}
    \put(520,220){\oval(100,100)[b]}
  \put(470,220){\line(0,1){180}}
  \put(570,220){\vector(-1,2){85}}
  \put(480,220){$^\Da$}
    \put(490,155){$^{a_1}$}
     \put(520,90){$^{x_2}$}
     \put(650,220){\oval(100,100)[b]}
  \put(600,220){\line(-1,2){85}}
  \put(700,220){\vector(-1,1){180}}
  \put(600,190){$^\Da$}
    \put(620,155){$^{a_0}$}
     \put(670,100){$^{x_1}$}

\end{picture}
\end{center}
is the following positive face structure $T^\dag$
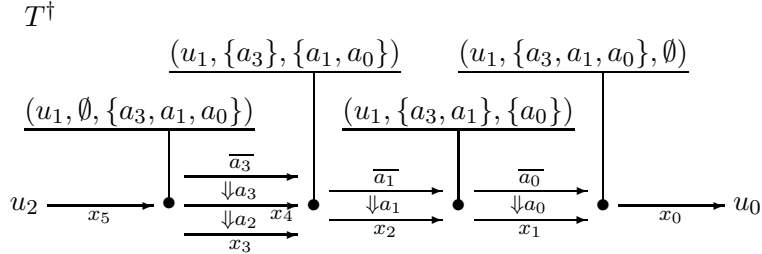
\begin{center} \xext=2820 \yext=820
\begin{picture}(\xext,\yext)(\xoff,\yoff)
  \put(0,770){$T^\dag$}
\putmorphism(0,150)(1,0)[u_2`\bullet`_{x_5}]{500}{1}b
 \putmorphism(500,250)(1,0)[\phantom{\bullet}`\phantom{\bullet}`_{\overline{a_3}}]{500}{1}a
 \putmorphism(500,150)(1,0)[\phantom{\bullet}`\bullet`]{500}{1}b
  \putmorphism(500,50)(1,0)[\phantom{\bullet}`\phantom{\bullet}`_{x_3}]{500}{1}b
   \put(680,200){$_{\Da a_3}$}
   \put(680,100){$_{\Da a_2}$}
   \put(850,120){$_{x_4}$}
   \putmorphism(1000,150)(1,0)[\phantom{\bullet}`\bullet`]{500}{0}b
   \putmorphism(1000,200)(1,0)[\phantom{\bullet}`\phantom{\bullet}`_{\overline{a_1}}]{500}{1}a
   \putmorphism(1000,100)(1,0)[\phantom{\bullet}`\phantom{\bullet}`^{x_2}]{500}{1}b
   \put(1180,150){$_{\Da a_1}$}
   \putmorphism(1500,150)(1,0)[\phantom{\bullet}`\bullet`]{500}{0}b
   \putmorphism(1500,200)(1,0)[\phantom{\bullet}`\phantom{\bullet}`_{\overline{a_0}}]{500}{1}a
   \putmorphism(1500,100)(1,0)[\phantom{\bullet}`\phantom{\bullet}`^{x_1}]{500}{1}b
   \put(1680,150){$_{\Da a_0}$}
    \putmorphism(2000,150)(1,0)[\phantom{\bullet}`u_0`_{x_0}]{500}{1}b

 \put(0,450){$\underline{(u_1,\emptyset,\{a_3,a_1,a_0\})}$}
   \put(500,410){\line(0,-1){240}}
 \put(500,650){$\underline{(u_1,\{a_3\},\{a_1,a_0\})}$}
  \put(1000,610){\line(0,-1){440}}
 \put(1100,450){$\underline{(u_1,\{a_3,a_1\},\{a_0\})}$}
 \put(1500,410){\line(0,-1){240}}
  \put(1500,650){$\underline{(u_1,\{a_3,a_1,a_0\},\emptyset)}$}
   \put(2000,610){\line(0,-1){440}}
\end{picture}
\end{center}
As before we use the convention that the empty cut in $T^\dag$,
say $(x_5,\emptyset,\emptyset)$, is identified with the
corresponding face in $T$, $x_5$ in this case. All bullets
$\bullet$ denote cuts and they are linked by a line to the
descriptions of the cuts they denote.

An ideal $\cI$ in an ordered face structure $S$ is an {\em unary
ideal}\index{ideal!unary -} iff $\cI\subseteq\delta(S^u)$. The
following is a kind of inverse of the Theorem \ref{quotient ofs}.

\begin{theorem}\label{posivive cover} Let $S$ be an ordered face structure. Then
$S^\dag$ is a positive face structure, $\overline{\cI}$ is an
unary ideal in $S^\dag$,  $q_S:S^\dag\lra S$ is a positive cover
with the kernel $\overline{\cI}$.
\end{theorem}

\subsection*{Some technical lemmas}

Since the Lemmas stated below are very technical we will comments
on them. Lemmas \ref{fact16new0}, \ref{fact16new1},
\ref{fact16new2} are there to be used in the proofs of Lemmas
\ref{fact17}, \ref{fact18}, \ref{fact19}. Lemma \ref{fact16new2}
is a suplement to the pencil linearity axiom and it says
intuitively that if some faces are incident then some (other)
faces are comparable.
Lemmas \ref{fact17}, \ref{fact18}, \ref{fact19} concern
$\theta\delta(x)$-cuts. They express cuts determined by some faces
in terms of cuts determined by some other faces.   Lemma
\ref{fact17}, is about the cuts determined by $\gamma(x)$, Lemma
\ref{fact18}, is about the cuts determined by a faces
$t\in\delta(x)$, and Lemma \ref{fact19}, is about the cuts
determined by a faces in $\gamma(\cI^{\leq^+x})$.

\begin{lemma}
\label{fact16new0}
Let $S$ be an ordered face structure $x,a\in S$.
 If $a\in\cI$ and $x\leq^+\gamma(a)$ then $x=\gamma(a)$
\end{lemma}
{\em Proof.}~Suppose $x<^+\gamma(a)$. Let
$x,a_1,\ldots,a_k,\gamma(a)$ be an upper
$(S-\gamma(S^{-\lambda}))$-path, $k\geq 1$. As
$\gamma(a)=\gamma(a_k)$ and $a,a_k\in S-S^{-\lambda})$, we have
$a=a_k$. But this is a contradiction, as $a\in S^\varepsilon$ and
$a_k\in S^{-\varepsilon}$. $~\Box$

\begin{lemma}
\label{fact16new1}
Let $S$ be an ordered face structure $t,t'\in S$,
$z\in\cI_{\gamma(t)}$, $t<^\sim t'$, $\gamma(t)\in\delta(t')$.
Then either $t<^\sim\gamma(z)$ or $\gamma(z)<^\sim t'$.
\end{lemma}
{\em Proof.}~Suppose contrary that $t\not<^\sim\gamma(z)\not<^\sim
t'$. Then, as $t<^\sim t'$, we also have
$t'\not<^\sim\gamma(z)\not<^\sim t$. So
$t\perp^+\gamma(z)\perp^+t'$. Thus, by Lemma \ref{fact16new0}, we
have $\gamma(z)\leq^+t,t'$.  Hence, by Lemma \ref{le_linearity},
$t\perp^+t'$ and we get a contradiction. $~\Box$

\begin{lemma}
\label{fact16new2}
Let $S$ be an ordered face structure $u,x,y,z\in S$, $z\in \cI_u$,
$y\in\cI_u^{\leq^+x}$.
\begin{enumerate}
  \item Let $x\in S^{-\varepsilon}$, $\gamma\gamma(x)=u$. If
  $\varrho(x)<^\sim\gamma(z)$ then either $z<^+x$ or $\gamma(x)<^\sim\gamma(z)$.
  \item Let $x\in S^{\varepsilon}$, $\gamma\gamma(x)=u$. If $\gamma(y)<^\sim\gamma(z)$ then either $z<^+x$ or
  $\gamma(x)<^\sim\gamma(z)$.
   \item  Let $x\in S^{-\varepsilon}$, $u\in\dot{\delta}\gamma(x)$,
   $t=\inf_\sim\{t'\in\dot{\delta}(x): u\in\delta(t') \}$. If
  $\gamma(z)<^\sim t$ then either $z<^+x$ or $\gamma(z)<^\sim\gamma(x)$.
  \item Let $x\in S^{\varepsilon}$, $u\in\dot{\delta}\gamma(x)$. If $\gamma(z)<^\sim\gamma(y)$
  then either $z<^+x$ or $\gamma(z)<^\sim\gamma(x)$.
\end{enumerate}
\end{lemma}

{\em Proof.} We use notation as above in the statement of Lemma.
Recall that if $z\in\cI$ then for no face $z'$ we have either
$z'<^\sim z$ or $z'<^+ z$.

Ad 1. As $\gamma\gamma(x)=\gamma\gamma(z)$ we have either
$\gamma(x)\perp^+\gamma(z)$ or $\gamma(x)\perp^\sim\gamma(z)$. In
the later case, as assumption $\gamma(z)<^\sim\gamma(x)$
immediately leads to contradiction, we get
$\gamma(x)<^\sim\gamma(z)$. In the former case we have either
$z<^+x$ or $z<^\sim x$. The later of these to is impossible, as we
would have $\gamma(z)\leq^+t\leq^\sim\varrho(x)$ and hence
$\gamma(z)\leq\varrho(x)$ contrary to the supposition. Thus we get
either $z<^+x$ or $\gamma(x)<^\sim\gamma(z)$.

Ad 2. In this case again we have $\gamma\gamma(x)=\gamma\gamma(z)$
and hence  either $\gamma(x)\perp^+\gamma(z)$ or
$\gamma(x)\perp^\sim\gamma(z)$. In the later case we again easily
get that $\gamma(x)<^\sim\gamma(z)$. In the former case, as
$z<^\sim x\in\cE$ is impossible, we get $z<^+x$. Thus again, we
get that either $z<^+x$ or $\gamma(x)<^\sim\gamma(z)$.

Ad 3. As $\gamma\gamma(z)=\dot{\delta}\gamma(x)$ we have either
$\gamma(x)\perp^+\gamma(z)$ or $\gamma(x)\perp^\sim\gamma(z)$. In
the later case we easily get (otherwise $\gamma(x)<^\sim t$) that
$\gamma(z)<^\sim\gamma(x)$. In the former case, as $z\in\cI$, we
get that either $z<^+x$ or $z<^\sim x$.  We shall show that
$z<^\sim x$ is impossible. Suppose contrary, then there is
$t'\in\delta(x)$ such that $\gamma(z)\leq^+t'$. If we were to have
$\gamma(z)=t'$ then, by definition of $t$, we would have
$t\leq^\sim t'=\gamma(z)$.  Thus $\gamma(z)<^+t'$ and there is a
flat upper path $\gamma(z),z_1,\ldots,z_k,t'$, with $k\geq 1$. If
$u\not\in\theta(t')$ then, as $u=\gamma\gamma(z)$, there is $1\leq
i\leq k$ such that $u\in\iota(z_i)$.  Hence
$t<^+\gamma(z_i)\leq^+\gamma(z_k)=t'$ and we get a contradiction
with local discreteness. If $u\in\theta(t')$ then, using the
definition of $t$, we easily get that $t<^\sim t'$. As
$\gamma(z)<^\sim t$ we get $\gamma(z)<^\sim t'$ contrary to the
definition of $t'$.  Thus the assumption $z<^\sim x$ leads to a
contradiction.

Ad 4. As $\gamma\gamma(z)=\dot{\delta}\gamma(x)$ we have either
$\gamma(x)\perp^+\gamma(z)$ or $\gamma(x)\perp^\sim\gamma(z)$. In
the later case we easily get that $\gamma(z)<^\sim\gamma(x)$. In
the former case, as $z,x\in\cE$ we cannot have $z\perp^\sim x$.
As, $x<^+z\in\cI$ is also impossible, we have $z<^+x$ in that
case. Thus we get either $z<^+x$ or $\gamma(z)<^\sim\gamma(x)$.
$~\Box$

The above Lemma had four parts with first two and second two
having the same conclusions.  The following Lemma contains in fact
four statement with essentially the same conclusion. This is why
we state it in a bit unusual form to emphasize it.

\begin{lemma}
\label{fact16new3}
Let $S$ be an ordered face structure $u,x,z\in S$, $z\in \cI_u$,
$u\in\theta\delta(x)$. Moreover, assume that one the following
four conditions
\begin{enumerate}
  \item $t,t''\in\dot{\delta}(x)$, $\gamma(t)=u\in\delta(t'')$,
  \item $y\in\cI_u^{\leq^+x}$,  $t''\in\dot{\delta}(x)$, $\gamma(y)=t$, $u\in\delta(t'')$,
  \item  $t\in\dot{\delta}(x)$, $y''\in\cI_u^{\leq^+x}$, $\gamma(t)=u$, $\gamma(y'')=t''$,
  \item $y,y''\in\cI_u^{\leq^+x}$, $\gamma(y)=t$,
  $\gamma(y'')=t''$,
\end{enumerate}
 holds. If $t<^\sim \gamma(z)<^\sim t''$ then either $z<^+x$ or there is
$t'\in\delta^\lambda(x)$ such that $t<^\sim t'<^\sim t''$,
$\gamma(z)\leq^+t'$ and $\gamma(t')=u$.
\end{lemma}

{\em Proof.} We use notation as above in the statement of Lemma.
Note that if $\gamma(z)\leq^+t'$ and $t'\in S^\lambda$ then
$\gamma(t')=u$.

First we shall show that any of the above four assumptions imply
the claim: either $z<^+x$ or $z<^\sim x$. Note that
$\gamma\gamma(z)\in\theta\delta(x)$. If
$\gamma\gamma(z)\in\iota(x)$ then the claim follows immediately
from pencil linearity. If $\gamma\gamma(z)\in\theta\gamma(x)$ then
by pencil linearity we get that either $\gamma(z)\perp^+\gamma(x)$
or $\gamma(z)\perp^\sim\gamma(x)$. In the former case we get,
again by pencil linearity, the claim.  In the later case, as
$t<^\sim \gamma(z)<^\sim t''$, we get either $t\perp^\sim
\gamma(x)$ or $t''\perp^\sim \gamma(x)$, i.e. a contradiction, as
$t,t''\leq^+\gamma(x)$ under each of the four assumptions above.
Thus we have the claim.

Now it remains to show that each of the following four assumptions
imply that if $z<^\sim x$ then there is $t'\in\delta^\lambda(x)$
such that $t<^\sim t'<^\sim t''$, $\gamma(z)\leq^+t'$. As all the
arguments are very similar we shall show this for the assumption
1.

 Assume $z<^\sim x$. Then there is $t'\in\delta(x)$ such that
$\gamma(z)\leq^+t'$. We need to show that $t<^\sim t'<^\sim t''$,
and $t'\in S^\lambda$.

If $\gamma(z)=t'$ we are done.  So assume that $\gamma(z)<^+t'$
and then we have a flat upper path $\gamma(z),z_1,\ldots,z_k,t'$,
with $k\geq 1$. If $\gamma\gamma(z)=u\not\in\theta(t')$ then there
is $1\leq i\leq k$ that $u\in\iota(z_i)$. So
$t,t''<^+\gamma(z_i)\leq^+\gamma(z_k)=t'$ and we get a
contradiction with local discreetness. Thus $u\in\theta(t')$ and
we have $t\perp^\sim t'\perp^\sim t''$. If we were to have
$t'\leq^\sim t$ then we would have $t'\leq^\sim \gamma(z)$ and if
we were to have $t''\leq^\sim t'$ then we would have
$\gamma(z)<^\sim t'$.  Thus we must have $t<^\sim t'<^\sim t''$.
Therefore there are $u'\in\delta(t')$ and $u''\in\delta(t'')$ such
that $u=\gamma(t)\leq^+u'\leq^+\gamma(t')\leq^+u''$. As
$u,u''\in\delta(t'')$ and $u\leq^+u''$ we have $u=u''$. Hence
$\gamma(t')\in\delta(t')$, i.e. $t'\in S^\lambda$. $~\Box$

\begin{lemma}
\label{fact17}
Let $S$ be an ordered face structure $u,x\in S$,
$u\in\dot{\delta}\gamma(x)$. We put
 \[ t_{\sup} =\sup_\sim(\{ \varrho(x) \} \cup\gamma(\cI_{\gamma\gamma(x)}^{\leq^+x})), \hskip 8mm
 t_{\inf} =\inf_\sim(\{ t\in\dot{\delta}(x):u\in\delta(t) \}\cup\gamma(\cI_u^{\leq^+x})). \]
 The elements $t_{\sup}$, $t_{\inf}$ are well defined and
\[1.\hskip 15mm (\gamma\gamma(x),-,\ua \gamma(x))=(\gamma\gamma(x),-,\ua t_{\sup}),\hskip 50mm \]
\[2.\hskip 15mm (u,\da \gamma(x),-)=(u,\da t_{\inf},-).\hskip 65mm\]
\end{lemma}
{\em Proof.}~Ad 1. We consider two cases depending on whether
$t_{\sup} =\varrho(x)$ or $ t_{\sup}
=\sup_\sim(\gamma(\cI_{\gamma\gamma(x)}^{\leq^+x}))$. Fix
$z\in\cI_{\gamma\gamma(x)}$.

Case $t_{\sup} =\varrho(x)$. Assume $\gamma(x)<^\sim\gamma(z)$. As
$\varrho(x)<^+\gamma(x)$, we have $\varrho(x)<\gamma(z)$. But if
we were to have $\varrho(x)<^+\gamma(z)$ we would have
$\gamma(x)\perp^+\gamma(z)$. Thus, as
$\gamma\varrho(x)=\gamma\gamma(z)$, we have $\varrho(x)<^\sim
\gamma(z)$. To see the converse, assume that
$\varrho(x)<^\sim\gamma(z)$. So by Lemma \ref{fact16new2}.1 we
have that either $\gamma(x)<^\sim\gamma(z)$ or $z<^+x$. But
$z<^+x$ would contradict the choice of $t_{\sup}$. Thus
$\gamma(x)<^\sim\gamma(z)$, and hence the equation 1. holds in
this case.

Case $t_{\sup}
=\sup_\sim(\gamma(\cI_{\gamma\gamma(x)}^{\leq^+x}))$. Fix
$y_{\sup}\in\cI_{\gamma\gamma(x)}^{\leq^+x}$ such that $t_{\sup}
=\gamma(y_{\sup})$. Assume that $\gamma(x)<^\sim\gamma(z)$.  As
$\gamma(y_{\sup})\leq^+\gamma(x)$ we have
$\gamma(y_{\sup})<\gamma(z)$. But
$\gamma(y_{\sup})\not\perp^+\gamma(z)$, so
$\gamma(y_{\sup})<^\sim\gamma(z)$. For converse, assume that
$\gamma(y_{\sup})<^\sim\gamma(z)$. If $x\in S^{-\varepsilon}$ then
$\varrho(x)<^\sim\gamma(y_{\sup})$ and again by Lemma
\ref{fact16new2}.1 we have that either $\gamma(x)<^\sim\gamma(z)$
or $z<^+x$. If $x\in S^{\varepsilon}$ then by Lemma
\ref{fact16new2}.2 we get once more that  either
$\gamma(x)<^\sim\gamma(z)$ or $z<^+x$. But $z<^+x$ would
contradict the choice of $t_{\sup}$. Thus
$\gamma(x)<^\sim\gamma(z)$, and hence the equation 1. holds in
this case as well.

Ad 2. We consider again two cases depending on the set $ t_{\inf}$
is in.  Fix $z\in\cI_u$.

Case $t_{\inf} =\inf_\sim(\{ t\in\dot{\delta}(x):u\in\delta(t)
\})$. Assume $\gamma(z)<^\sim t_{\inf}$. Then by  Lemma
\ref{fact16new2}.3 we have that either $\gamma(z)<^\sim\gamma(x)$
or $z<^+x$. But $z<^+x$ would contradict the choice of $t_{\inf}$.
Thus $\gamma(z)<^\sim\gamma(x)$. To see the converse, assume that
$\gamma(z)<^\sim\gamma(x)$. But then as other cases are easily
excluded we must have $\gamma(z)<^\sim t_{\inf}$ , and hence the
equation 2. holds in this case.

Case $t_{\inf} =\inf_\sim(\gamma(\cI_u^{\leq^+x}))$. Fix
$y_{\inf}\in\cI_u^{\leq^+x}$ such that $t_{\inf}
=\gamma(y_{\inf})$. Assume that $\gamma(z)<^\sim\gamma(x)$.  As
$\gamma(y_{\inf}),\gamma(z)\in\cI_u^{\leq^+x}$ we have
$\gamma(y_{\inf})\perp^\sim\gamma(z)$. As
$\gamma(y_{\inf})<^\sim\gamma(z)$ leads immediately to a
contradiction we have $\gamma(z)<^\sim\gamma(y_{\inf})$. To see
the converse assume that $\gamma(z)<^\sim\gamma(y_{\inf})$. If
$x\in S^{-\varepsilon}$ then by definition of $y_{\sup}$ we have
$\gamma(y_{\sup})<^\sim
\inf_\sim\{t'\in\dot{\delta}(x):u\in\delta(t') \}$ and again by
Lemma \ref{fact16new2}.3 we have that either
$\gamma(z)<^\sim\gamma(x)$ or $z<^+x$. If $x\in S^{\varepsilon}$
then by Lemma \ref{fact16new2}.4 we get once more that  either
$\gamma(z)<^\sim\gamma(x)$ or $z<^+x$. But $z<^+x$ would
contradict the choice of $t_{\inf}$. Thus
$\gamma(z)<^\sim\gamma(x)$, and hence the equation 2. holds in
this case as well. $~\Box$

\begin{lemma}
\label{fact18}
Let $S$ be an ordered face structure $u,t,x\in S$,
$u\in\dot{\delta}(t)$ and $t\in\dot{\delta}(x)$. We put
\[t_{\sup} =\sup_\sim(\{ t'\in\dot{\delta}(x):t'<^\sim t,\,\gamma(t')=u \} \cup
\gamma(\{y\in \cI_u^{\leq^+x}: \gamma(y)<^\sim t \})), \]
\[ t_{\inf} =\inf_\sim(\{ t'\in\dot{\delta}(x):t<^\sim t' \} \cup
\gamma(\{y\in \cI_{\gamma(t)}^{\leq^+x}: t<^\sim \gamma(y)\})).\]
The elements $t_{\sup}$, $t_{\inf}$ are not necessarily well
defined, due to the fact that these set might be empty, but we
have
\[3.\hskip 15mm (u,\da t,-)= \left\{ \begin{array}{ll}
           (u,\da \gamma(x),-) & \mbox{if  $t_{\sup}$ is undefined,}  \\
          (u,-,\ua t_{\sup}) & \mbox{otherwise.}
                                    \end{array}
                \right. \hskip 3cm \]
\[4.\hskip 15mm (\gamma(t),-,\ua t)= \left\{ \begin{array}{ll}
          (\gamma\gamma(x),-,\ua \gamma(x)) & \mbox{if  $t_{\inf}$ is undefined,}  \\
          (\gamma(t),\da t_{\inf},-) & \mbox{otherwise.}
                                    \end{array}
                             \right. \hskip 18mm \]
\end{lemma}
{\em Proof.}~Ad 3. First note that if $t_{\sup}$ is undefined then
$t$ is $t_{\inf}$ from Lemma \ref{fact17}. Thus the equation 3.
holds in this case by Lemma \ref{fact17}.2.  If $t_{\sup}$ is
defined then we consider two cases depending on the set $t_{\sup}$
is in. However in either case, by Lemma \ref{fact16new1}, we have
that $\da_ut\cup\ua t_{\sup}=\cI_u$.

Case $t_{\sup} =\sup_\sim(\{ t'\in\dot{\delta}(x):t'<^\sim
t,\,\gamma(t')=u \})$. Suppose there is $z\in\cI_u$ such that
$t_{\sup}<^\sim\gamma(z)<^\sim t$. Then, as the assumption 1. of
Lemma \ref{fact16new3} holds, we have that  either $z<^+x$ or
there is $t'\in\delta^\lambda(x)$ such that $t_{\sup}<^\sim
t'<^\sim t$, $\gamma(z)\leq^+t'$ and $\gamma(t')=u$. Both cases
contradict the choice $t_{\sup}$. Thus $\da_ut\cap\ua
t_{\sup}=\emptyset$, and hence $(u,\da t,-)=(u,-,\ua t_{\sup})$
i.e. the equation 3. holds in this case.

Case $t_{\sup} =\sup_\sim(\gamma(\{y\in \cI_u^{\leq^+x}:
\gamma(y)<^\sim t \}))$. Fix $y_{\sup}\in\cI_u^{\leq^+x}$ such
that $\gamma(y_{\sup})=t_{\sup}$. Suppose there is $z\in\cI_u$
such that $\gamma(y_{\sup})<^\sim\gamma(z)<^\sim t$. Then, as the
assumption 2. of Lemma \ref{fact16new3} holds, we have that either
$z<^+x$ or there is $t'\in\delta^\lambda(x)$ such that
$\gamma(y_{\sup})<^\sim t'<^\sim t$, $\gamma(z)\leq^+t'$ and
$\gamma(t')=u$. Both cases contradict the choice $y_{\sup}$. Thus
$\da_ut\cap\ua \gamma(y_{\sup})=\emptyset$, i.e. the equation 3.
holds in this case, as well.

Ad 4. First note that if $t_{\inf}$ is undefined then $t$ is
$t_{\sup}$ from Lemma \ref{fact17}. Thus the equation 4. holds in
this case by Lemma \ref{fact17}.1.  If $t_{\inf}$ is defined then
we consider two cases depending on the set $t_{\inf}$ is in.
However in either case, by Lemma \ref{fact16new1}, we have that
$\da_{\gamma(t)}t_{\inf}\cup\ua t=\cI_{\gamma(t)}$.

Case $t_{\inf} =\inf_\sim(\{ t'\in\dot{\delta}(x):t<^\sim t' \})$.
Suppose there is $z\in\cI_{\gamma(t)}$ such that
$t<^\sim\gamma(z)<^\sim t_{\inf}$. Then, as the assumption 1. of
Lemma \ref{fact16new3} holds, we have that  either $z<^+x$ or
there is $t'\in\delta^\lambda(x)$ such that $t<^\sim t'<^\sim
t_{\inf}$, $\gamma(z)\leq^+t'$ and $\gamma(t')=\gamma(t)$. Both
cases contradict the choice $t_{\inf}$. Thus
$\da_{\gamma(t)}t_{\inf}\cap\ua t=\emptyset$, and hence
$(\gamma(t),-,\ua t)=(\gamma(t),\da t_{\inf},-)$, i.e. the
equation 4. holds in this case.

Case $t_{\inf} =\inf_\sim(\gamma(\{y\in \cI_{\gamma(t)}^{\leq^+x}:
t<^\sim \gamma(y)\}))$. Fix $y_{\inf}\in\cI_{\gamma(t)}^{\leq^+x}$
such that $\gamma(y_{\inf})=t_{\inf}$. Suppose there is
$z\in\cI_u$ such that $t<^\sim\gamma(z)<^\sim\gamma(y_{\inf})$.
Then, as the assumption 3. of Lemma \ref{fact16new3} holds, we
have that either $z<^+x$ or there is $t'\in\delta^\lambda(x)$ such
that $\gamma(y_{\sup})<^\sim t'<^\sim t$, $\gamma(z)\leq^+t'$ and
$\gamma(t')=\gamma(t)$. Both cases contradict the choice
$y_{\inf}$. Thus $\da_{\gamma(t)}t\cap\da
\gamma(y_{\inf})=\emptyset$, i.e. the equation 4. holds in this
case, as well. $~\Box$

\begin{lemma}
\label{fact19}
Let $S$ be an ordered face structure $y,x\in S$,
$y\in\cI_{\gamma\gamma(y)}^{\leq^+x}$. Then
$\gamma\gamma(y)\in\theta\delta(x)$. We put
\[  t_{\sup} =\sup_\sim(\{ t\in\dot{\delta}(x):t<^\sim \gamma(y)\} \cup
\gamma(\{y'\in \cI_{\gamma\gamma(y)}^{\leq^+x}:
\gamma(y')<^\sim\gamma(y) \})),\]
\[ t_{\inf} =\inf_\sim(\{ t\in\dot{\delta}(x):\gamma(y)<^\sim t \} \cup
\gamma(\{y'\in \cI_{\gamma\gamma(y)}^{\leq^+x}:
\gamma(y)<^\sim\gamma(y') \})).\]
 The elements $t_{\sup}$,$t_{\inf}$ are not necessarily well defined,
 due to the fact that these set might be empty, but we have
\[{ 5.}\hskip 15mm (\gamma\gamma(y),\da \gamma(y),-)= \left\{ \begin{array}{ll}
           (\gamma\gamma(y),\da \gamma(x),-) & \mbox{if  $t_{\sup}$ is undefined,}  \\
          (\gamma\gamma(y),-,\ua t_{\sup}) & \mbox{otherwise.}
                                    \end{array}
                \right. \hskip 10mm\]
\[{ 6.}\hskip 15mm (\gamma\gamma(y),-,\ua \gamma(y))= \left\{ \begin{array}{ll}
          (\gamma\gamma(x),-,\ua \gamma(x)) & \mbox{if  $t_{\inf}$ is undefined,}  \\
          (\gamma\gamma(y),\da t_{\inf},-) & \mbox{otherwise.}
                                    \end{array}
                             \right. \hskip 10mm\]
\end{lemma}
{\em Proof.}~Ad 5. First note that if $t_{\sup}$ is undefined
then, with $u=\gamma\gamma(y)$, $y$ is $y_{\inf}$ from (the proof
of) Lemma \ref{fact17}. Thus the equation 5. holds in this case by
Lemma \ref{fact17}.2. If $t_{\sup}$ is defined then we consider
two cases depending on the set $t_{\sup}$ is in. However in either
case, by Lemma \ref{fact16new1}, we have that
$\da_{\gamma\gamma(y)}\gamma(y)\cup\ua
t_{\sup}=\cI_{\gamma\gamma(y)}$.

Case $t_{\sup} =\sup_\sim(\{ t\in\dot{\delta}(x):t<^\sim \gamma(y)
\})$. Suppose there is $z\in\cI_{\gamma\gamma(y)}$ such that
$t_{\sup}<^\sim\gamma(z)<^\sim\gamma(y)$. Then, as the assumption
3. of Lemma \ref{fact16new3} holds, we have that  either $z<^+x$
or there is $t'\in\delta^\lambda(x)$ such that $t_{\sup}<^\sim
t'<^\sim t$, $\gamma(z)\leq^+t'$ and $\gamma(t')=\gamma\gamma(y)$.
Both cases contradict the choice $t_{\sup}$. Thus
$\da_{\gamma\gamma(y)}\gamma(y)\cap\ua t_{\sup}=\emptyset$, i.e.
the equation 5. holds in this case.

Case $t_{\sup} =\sup_\sim(\gamma(\{y'\in
\cI_{\gamma\gamma(y)}^{\leq^+x}: \gamma(y')<^\sim \gamma(y) \}))$.
Fix $y_{\sup}\in\cI_{\gamma\gamma(y)}^{\leq^+x}$ such that
$\gamma(y_{\sup})=t_{\sup}$. Suppose there is
$z\in\cI_{\gamma\gamma(y)}$ such that
$\gamma(y_{\sup})<^\sim\gamma(z)<^\sim \gamma(y)$. Then, as the
assumption 4. of Lemma \ref{fact16new3} holds, we have that either
$z<^+x$ or there is $t'\in\delta^\lambda(x)$ such that
$\gamma(y_{\sup})<^\sim t'<^\sim t$, $\gamma(z)\leq^+t'$ and
$\gamma(t')={\gamma\gamma(y)}$. Both cases contradict the choice
$y_{\sup}$. Thus $\da_{\gamma\gamma(y)}\gamma(y)\cap\ua
\gamma(y_{\sup})=\emptyset$, i.e. the equation 5. holds in this
case, as well.

Ad 6. First note that if $t_{\inf}$ is undefined then,
$\gamma\gamma(y)=\gamma\gamma(x)$ and $y$ is $y_{\sup}$ from (the
proof of) Lemma \ref{fact17}. Thus the equation 6. holds in this
case by Lemma \ref{fact17}.1.  If $t_{\inf}$ is defined then
consider two cases depending on the set $t_{\inf}$ is in. However
in either case, by Lemma \ref{fact16new1}, we have that
$\da_{\gamma\gamma(y)}t_{\inf}\cup\ua
\gamma(y)=\cI_{\gamma\gamma(y)}$.

Case $t_{\inf} =\inf_\sim(\{ t\in\dot{\delta}(x):\gamma(y)<^\sim t
\})$. Suppose there is $z\in\cI_{\gamma\gamma(y)}$ such that
$\gamma(y)<^\sim\gamma(z)<^\sim t_{\inf}$. Then, as the assumption
2. of Lemma \ref{fact16new3} holds, we have that  either $z<^+x$
or there is $t'\in\delta^\lambda(x)$ such that $t<^\sim t'<^\sim
t_{\inf}$, $\gamma(z)\leq^+t'$ and $\gamma(t')=\gamma\gamma(y)$.
Both cases contradict the choice $t_{\inf}$. Thus
$\da_{\gamma\gamma(y)}t_{\inf}\cap\ua \gamma(y)=\emptyset$, hence
$(\gamma\gamma(y),-,\ua \gamma(y))=(\gamma\gamma(y),\da
t_{\inf},-)$, i.e. the equation 6. holds in this case.

Case $t_{\inf} =\inf_\sim(\gamma(\{y'\in
\cI_{\gamma\gamma(y)}^{\leq^+x}: \gamma(y)<^\sim \gamma(y')\}))$.
Fix $y_{\inf}\in\cI_{\gamma\gamma(y)}^{\leq^+x}$ such that
$\gamma(y_{\inf})=t_{\inf}$. Suppose there is
$z\in\cI_{\gamma\gamma(y)}$ such that
$\gamma(y)<^\sim\gamma(z)<^\sim\gamma(y_{\inf})$. Then, as the
assumption 4. of Lemma \ref{fact16new3} holds, we have that either
$z<^+x$ or there is $t'\in\delta^\lambda(x)$ such that
$\gamma(y)<^\sim t'<^\sim \gamma(y_{\inf})$, $\gamma(z)\leq^+t'$
and $\gamma(t')=\gamma\gamma(y)$. Both cases contradict the choice
$y_{\inf}$. Thus $\ua_{\gamma\gamma(y)}t\cap\da
\gamma(y_{\inf})=\emptyset$, i.e. the equation 6. holds in this
case, as well. $~\Box$

\subsection*{The Proof}

 {\it Proof of Theorem \ref{posivive cover}.}~ Fix an ordered face structure $S$. Clearly for
$a\in\cI$, $\delta^\dag(\overline{a})\neq\emptyset$. Suppose that
$(x,L,U)$ is a cut in $S^\dag_k$, with $k>0$. Then either
$\dot{\delta}(x)\neq\emptyset$ or $\dot{\delta}(x)=\emptyset$ and
then by Lemma \ref{fact3} we have that there is $y\in
\cI^{\leq^+x}$. In either case $\delta^\dag(x,L,U)\neq\emptyset$.
Thus $S^\dag$ is a positive hypergraph. We shall check that
$S^\dag$ satisfies all four positive face structure axioms.

{\em Globularity.} We need to verify globularity for both kinds of
faces in $S^\dag$: bars and cuts. First we shall check globularity
for bars. Fix $\alpha\in \cI$. We have
\[ \gamma^{\dag}\gamma^{\dag}(\overline{\alpha})=(\gamma\gamma\gamma(\alpha), -,
\ua\gamma\gamma(\alpha))=\gamma^{\dag}\delta^{\dag}(\overline{\alpha}),\]
\[ \delta^{\dag}\gamma^{\dag}(\overline{\alpha})= \overline{\cI^{\leq^+\gamma(\alpha)}} \cup
  \{(t, \da \gamma\gamma(\alpha) ,-) : t\in\dot{\delta}\gamma\gamma(\alpha) \}=
  \delta^{\dag}\delta^{\dag}(\overline{\alpha}).  \]
We need to show that
\[ \gamma^{\dag}\gamma^{\dag}(\overline{\alpha})\not\in
\delta^{\dag}\gamma^{\dag}(\overline{\alpha}).\]

Suppose contrary that
$\gamma^{\dag}\gamma^{\dag}(\overline{\alpha})\in
\delta^{\dag}\gamma^{\dag}(\overline{\alpha})$. Then, as
$\gamma^{\dag}\gamma^{\dag}(\overline{\alpha})$ is a cut, we would
have $\gamma\gamma\gamma(\alpha)\in\delta\gamma\gamma(\alpha)$,
i.e. $\gamma\gamma(\alpha)$ is a loop. Thus by Lemma \ref{fact3}
there is $a\in\cI_{\gamma\gamma\gamma(\alpha)}$ such that
$\gamma(a)\leq \gamma\gamma(\alpha)$. But then
$\ua\gamma\gamma(\alpha)\not\ni a \not\in \da\gamma\gamma(\alpha)$
and hence
\[ (\gamma\gamma\gamma(\alpha),-,\ua\gamma\gamma(\alpha))\neq
(\gamma\gamma\gamma(\alpha),\da\gamma\gamma(\alpha),-), \] which
means that $\gamma^{\dag}\gamma^{\dag}(\overline{\alpha})\not\in
\delta^{\dag}\gamma^{\dag}(\overline{\alpha})$ after all. From
this the globularity for bars follow easily.

Now we shall check globularity for cuts. Fix a cut $(x,L,U)$ in
$S^\dag_k$, with $k>1$. The sets involved in the globularity
equations are sums of some sets. We shall spell these sets below
giving names to their summands. We have
\[\gamma^{\dag}\gamma^{\dag}(x,L,U)=(\gamma\gamma(x),-,\ua\gamma(x))=\psi,\]
\[\delta^{\dag}\gamma^{\dag}(x,L,U)=\overline{\cI^{\leq^+\gamma(x)}}\cup
 \{ (u,\da\gamma(x),-):u\in\dot{\delta}\gamma(x)\}=Z_1\cup Z_2,\]
\[\gamma^{\dag}\delta^{\dag}(x,L,U)=\{
(\gamma\gamma(y),-,\ua\gamma(y)):y\in \cI^{\leq^+x} \} \cup
\{(\gamma(t),-,\ua t): t\in \dot{\delta}(x) \} =Z_3\cup Z_4 \]
 \[\delta^{\dag}\delta^{\dag}(x,L,U)=\{ (\gamma\gamma(y),\da\gamma(y),-):y\in \cI^{\leq^+x} \} \cup \{
 \overline{s} : s\in\cI^{\leq^+t }, \; t\in\dot{\delta}(x)\} \cup \]
 \[ \cup \{(u,\da t,-):  t\in \dot{\delta}(x),\; u\in \dot{\delta}(t) \} =Z_5\cup Z_6\cup Z_7\]

In order to verify $\gamma$-globularity, i.e.
\[ \gamma^{\dag}\gamma^{\dag}(x,L,U) = \gamma^{\dag}\delta^{\dag}(x,L,U)-\delta^{\dag}\delta^{\dag}(x,L,U),  \]
we shall show:
\begin{description}
  \item[{\rm (A)}] $\psi\in\gamma^{\dag}\delta^{\dag}(x,L,U)$,
  \item[{\rm (B)}] $\psi\not\in\delta^{\dag}\delta^{\dag}(x,L,U)$,
  \item[{\rm (C)}] $\gamma^{\dag}\delta^{\dag}(x,L,U)-\psi\subseteq
  \delta^{\dag}\delta^{\dag}(x,L,U)$.
\end{description}

Ad A. By Lemma \ref{fact17} ($\gamma(x)$-cuts) either $x\in
S^{-\varepsilon}$ and $\psi=(\gamma\varrho(x),-,\ua \varrho(x))\in
Z_4$ or there is $y\in\in\cI^{\leq^+x}_{\gamma\gamma(x)}$ such
that $\psi=(\gamma\gamma(y),-,\ua \gamma(y))\in Z_3$. In either
case $\psi\in\gamma^{\dag}\delta^{\dag}(x,L,U)$.

Ad B. As $\psi$ is not a bar, we have $\psi\not\in Z_6$.

Suppose $\psi\in Z_5$. Then there is
$y\in\cI^{\leq^+x}_{\gamma\gamma(x)}$ such that
$(\gamma\gamma(y),\da \gamma(y),-)=\psi$. So $y\in\ua \gamma(x)$,
i.e. $\gamma(x)<^\sim\gamma(y)$.  But $y\leq^+x$, so
$\gamma(y)\leq^+\gamma(x)$ and we have a contradiction with the
disjointness. Thus $\psi\not\in Z_5$.

Suppose now that $\psi\in Z_7$. So there is $t\in\dot{\delta}(x)$
such that $\gamma\gamma(x)\in\dot{\delta}(t)$ and
$(\gamma\gamma(x),\da t,-)=\psi$. As $t\in\dot{\delta}(x)$ we have
$\gamma(t)\leq^+\gamma\gamma(x)$. So $t$ is a loop. Then, by Lemma
\ref{fact3}, there is $y\in\cI_{\gamma\gamma(x)}$ such that
$\gamma(y)\leq^+t$. As $y\in S^{-\lambda}$, we have $y<^\sim x$,
and hence $\gamma(y)\leq^+\gamma(x)$. Thus
$y\not\in\da_{\gamma\gamma(x)}t$ and $y\not\in\ua\gamma(x)$, i.e.
$(\gamma\gamma(x),\da t,-)\neq \psi$ after all. Thus $\psi\not\in
Z_5$ and hence $\psi\not\in \delta^{\dag}\delta^{\dag}(x,L,U)$.

Ad C. Fix $\xi\in\gamma^{\dag}\delta^{\dag}(x,L,U)$, such that
$\xi\neq\psi$. If $\xi\in Z_4$ then there is $t\in\dot{\delta}(x)$
such that $\xi=(\gamma(t),-,\ua t)$. We shall use Lemma
\ref{fact18} (t-cuts). As $\xi\neq\psi$
\[ t_{\inf} =\inf_\sim(\{ t'\in\delta(x):t<^\sim t' \} \cup
 \gamma(\{y\in \cI_{\gamma(t)}^{\leq^+x}: t<^\sim \gamma(y)\})).\]
is well defined and then $\xi=(\gamma(t),\da t_{\inf},-)$. Now, if
$t_{\inf} =\inf_\sim(\{ t'\in\delta(x):t<^\sim t' \})$ then
$\xi\in Z_7$ and if $t_{\inf} =\inf_\sim( \gamma(\{y\in
\cI_{\gamma(t)}^{\leq^+x}: t<^\sim \gamma(y)\}))$ then $\xi\in
Z_5$.

If $\xi\in Z_3$ then there is $y\in\cI^{\leq^+x}$, so that $\xi
=(\gamma\gamma(y),-,\ua \gamma(y))$. We shall use Lemma
\ref{fact19} ($\gamma(y)$-cuts). As $\xi\neq\psi$ then
\[ t_{\inf} =\inf_\sim(\{ t\in\delta(x):\gamma(y)<^\sim t \} \cup
\gamma(\{y'\in \cI_{\gamma\gamma(y)}^{\leq^+x}:
\gamma(y)<^\sim\gamma(y') \}))\] is well defined and
$\xi=(\gamma\gamma(y),\da t_{\inf},-)$. Again, if $t_{\inf}
=\inf_\sim(\{ t\in\delta(x):\gamma(y)<^\sim t \})$ then $\xi\in
Z_7$ and if $t_{\inf} =\inf_\sim(\gamma(\{y'\in
\cI_{\gamma\gamma(y)}^{\leq^+x}: \gamma(y)<^\sim\gamma(y') \}))$
then $\xi\in Z_5$. Thus C. holds. This ends verification of
$\gamma$-globularity for $S^\dag$.

Now we shall check $\delta$-globularity for $S^\dag$, i.e.
\[ \delta^{\dag}\gamma^{\dag}(x,L,U) =
\delta^{\dag}\delta^{\dag}(x,L,U)-\gamma^{\dag}\delta^{\dag}(x,L,U).
\] Both sides of this equation contains both bars and cuts. We
show equalities for them separately.

First we shall show the equality for bars. We need to show that
$Z_6=Z_1$. Clearly, $Z_6\subseteq Z_1$. We shall verify that
$Z_1\subseteq Z_6$. Let $\overline{t} \in Z_1$, i.e. $t\in \cI$
and $t\leq^+\gamma(x)$. As $t\in\cI$ and
$\cI\cap\gamma(S^{-\lambda})=\emptyset$ we have that either
$t=\gamma(x)=\delta(x)$ or $t\neq\gamma(x)$. In the former case
clearly $\overline{t} \in Z_6$. In the later case there is an
upper $(S-\gamma(S^{-\lambda}))$-path $t,x_1,\ldots,x_k,\gamma(x)$
with $k\geq 1$. By pencil linearity $x_k\leq^+ x$. As
$t\in\cI\subseteq S-\gamma(S^{-\lambda})$ by Second Path Lemma,
either $t\in\delta(x)$ or there is $1\leq i<k$ such that
$\gamma(x_i)\in\delta(x)$. In either case there is $s\in\delta(x)$
($s=t$ or $s=\gamma(x_i$)) such that $t\leq^+s$, i.e.
$\overline{t} \in Z_6$. Thus the $\delta$-globularity for bars
holds.

Now we will show the $\delta$-globularity for cuts. Clearly, it is
enough to restrict ourself to cuts over $\dot{\theta}\delta(x)$ as
other cuts cannot appear in the equation. Moreover, by Lemma
\ref{fact6} (atlas), we have
$\dot{\theta}\delta(x)=\dot{\delta}\gamma(x)\cup\gamma\dot{\delta}^{-\lambda}(x)$.
As in $\delta^{\dag}\gamma^{\dag}(x,L,U)$ can appear only cuts
over $\dot{\delta}\gamma(x)$ we will split our proof farther by
considering these two case separately. Let
$u\in\dot{\theta}\delta(x)$. By $X_{u-cuts}$ we mean $u$-cuts in
the set $X$. To end the proof of $\delta$-globularity we need to
show:
\begin{description}
 \item[{\rm (I)}] if $u\in\gamma\dot{\delta}^{-\lambda}(x)$ then
 $\delta^{\dag}\delta^{\dag}(x,L,U)_{u-cuts}\subseteq \gamma^{\dag}\delta^{\dag}(x,L,U)$,
\item[{\rm (II)}] if $u\in\dot{\delta}\gamma(x)$ then the cut $\psi_u=(u,\da\gamma(x),-)$
is the only $u$-cut in $\delta^{\dag}\gamma^{\dag}(x,L,U)$.
Moreover we have:
\begin{description}
  \item[{\rm (A$_u$)}] $\psi_u\in\delta^{\dag}\delta^{\dag}(x,L,U)$,
  \item[{\rm (B$_u$)}] $\psi_u\not\in\gamma^{\dag}\delta^{\dag}(x,L,U)$,
  \item[{\rm (C$_u$)}] $\delta^{\dag}\delta^{\dag}(x,L,U)_{u-cuts}-\psi_u\subseteq
  \gamma^{\dag}\delta^{\dag}(x,L,U)$.
  \end{description}
\end{description}
Note the similarity of the conditions (A), (B), (C) with (A$_u$),
(B$_u$), (C$_u$).

Ad I. Fix $u\in \gamma\dot{\delta}^{-\lambda}(x)$ and
$t_u\in\dot{\delta}^{-\lambda}(x)$ such that $\gamma(t_u)=u$. Let
$\phi=(u,L',U')\in\delta^{\dag}\delta^{\dag}(x,L,U)$. We put
\[ t_\phi= \left\{ \begin{array}{ll}
          t & \mbox{if  $\phi\in Z_7$ and $\,t\in\dot{\delta}(x)$ such that $L'=\da_ut$,}  \\
          \gamma(y) & \mbox{if $\phi\in Z_5$, and $\,y\in\cI^{\leq^+x}_u$ such that $L'=\da_u\gamma(y)$.}
                                    \end{array}
                \right. \]
Thus $\phi=(u,\da t_\phi, -)$. Put
\[t_{\sup} =\sup_\sim(\{ t'\in\dot{\delta}(x):t'<^\sim t_\phi,\,\gamma(t')=u \} \cup
\gamma(\{y'\in \cI_u^{\leq^+x}: \gamma(y')<^\sim t_\phi \})). \]

As $t_u\in \{ t'\in\delta(x):t'<^\sim t_\phi,\,\gamma(t')=u
\}\neq\emptyset$ the face $t_{\sup}$ is well defined. Then, by
Lemmas \ref{fact18} and \ref{fact19} we have that $\phi=(u,-,\ua
t_{\sup})\in\gamma^{\dag}\gamma^{\dag}(x,L,U)$.

Ad II. The fact that $\psi_u$ is the only $u$-cut in
$\delta^{\dag}\gamma^{\dag}(x,L,U)$ is obvious from our
description of this set as sum $Z_1\cup Z_2$.

Ad A$_u$. Let $t_{\inf}=\inf_\sim(\{
t\in\dot{\delta}(x):u\in\delta(t) \} \cup
\gamma(\cI^{\leq^+}_u))$. If $\{ t\in\dot{\delta}(x):u\in\delta(t)
\}=\emptyset$ then $x\in S^\varepsilon$ and hence, by Lemma
\ref{fact3}, $\cI^{\leq^+}_u\neq\emptyset$.  Thus $t_{\inf}$ is
well defined. By Lemma \ref{fact17}, we have $\psi_u=(u,\da
t_{\inf},-)\in\delta^{\dag}\delta^{\dag}(x,L,U)$, as required.

Ad B$_u$. Suppose $\psi_u\in Z_3$. Then there is
$y\in\cI^{\leq^+x}_u$ such that $\psi_u=(u,-,\ua \gamma(y))$. As
$y\leq^+x$ we have $\gamma(y)\leq^+\gamma(x)$. Thus
$\gamma(y)\not\leq^\sim\gamma(x)$. This means that
$y\not\in\da_u\gamma(x)$. Clearly $y\not\in\ua \gamma(y)$. Thus
$\psi_u=(u,\da \gamma(x),-)\neq(u,-,\ua \gamma(y))$, after all.
This shows that $\psi_u\not\in Z_3$.

Suppose now that $\psi_u\in Z_4$. So there is
$t\in\dot{\delta}(x)$ such that $\psi_u=(u,-,\ua t)$. As
$\gamma(t)=u\in\delta\gamma(x)$, $t$ is a loop. Then, by Lemma
\ref{fact3}, there is $y\in\cI_{u}$ such that $\gamma(y)\leq^+t$
and, by transitivity of $<^+$, $\gamma(y)\leq^+\gamma(x)$. Thus
$y\not\in\ua t$ and $y\not\in\da_u\gamma(x)$. Then $\psi_u=(u,\da
\gamma(x),-)\neq(u,-,\ua t)$. So $\psi\not\in Z_4$, and hence
$\psi\not\in \gamma^{\dag}\delta^{\dag}(x,L,U)$.

Ad C$_u$.  Fix
$\xi=(u,L',U')\in\delta^{\dag}\delta^{\dag}(x,L,U)$, such that
$\xi\neq\psi_u$. If $\xi\in Z_7$ then there is
$t\in\dot{\delta}(x)$ such that $u\in\dot{\delta}(t)$ and
$\xi=(u,\da t,-)$.  We shall use Lemma \ref{fact18} (t-cuts). As
$\xi\neq\psi_u$  the face
\[t_{\sup} =\sup_\sim(\{ t'\in\dot{\delta}(x):t'<^\sim t,\,\gamma(t')=u \} \cup
\gamma(\{y\in \cI_u^{\leq^+x}: \gamma(y)<^\sim t \})), \]
 is well defined and then $\xi=(u,-,\ua t_{\sup})$. Now, if $t_{\sup}
=\sup_\sim(\{ t'\in\delta(x):t'<^\sim t,\,\gamma(t')=u \})$ then
$\xi\in Z_4$ and if $t_{\sup} =\sup_\sim( \gamma(\{y\in
\cI_u^{\leq^+x}: \gamma(y)<^\sim t \}))$ then $\xi\in Z_3$.

If $\xi\in Z_5$ then there is $y\in\cI^{\leq^+x}$, so that $\xi
=(u,\da \gamma(y),-)$. We shall use Lemma \ref{fact19}
($\gamma(y)$-cuts).

As $\xi\neq\psi_u$ the face
\[  t_{\sup} =\sup_\sim(\{ t\in\dot{\delta}(x):t<^\sim \gamma(y)\} \cup
\gamma(\{y'\in \cI_{\gamma\gamma(y)}^{\leq^+x}:
\gamma(y')<^\sim\gamma(y) \})),\] is well defined and
$\xi=(u,-,\ua t_{\sup})$. Again, if $t_{\sup} =\sup_\sim(\{
t\in\dot{\delta}(x):t<^\sim \gamma(y) \})$ then $\xi\in Z_4$ and
if $t_{\sup} =\sup_\sim(\gamma(\{y'\in
\cI_{\gamma\gamma(y)}^{\leq^+x}: \gamma(y')<^\sim\gamma(y) \}))$
then $\xi\in Z_3$. Thus C$_u$. holds. This ends verification of
$\delta$-globularity for $S^\dag$.

 {\em Orders in $S^\dag$.} Before we verify the remaining axioms
 of positive face structures we shall describe the order in
 $S^\dag$.  Let $(x,L,U), (y,L',U')\in \cC(\cI_{S_k})$ be two cuts in
 $S^\dag_k$ and $a,b\in\cI_{k+1}$ so that $\overline{a}, \overline{b}$ are
 two bars in  $S^\dag_k$.

 For $k\geq 0$, the upper order $<^{\dag,+}$ in  $S^\dag_k$ can be
characterized as follows ($<^+$ is the upper order in $S$):
\begin{enumerate}
  \item (cut,cut): $(x,L,U)<^{\dag,+}(y,L',U')$ iff either $x<^+y$ or $x=y$ and $L\subseteqnot L'$;
  \item (bar,cut): $\overline{a}<^{\dag,+}(y,L',U')$ iff $\gamma(a)\leq^+ y$;
  \item (cut,bar): $(x,L,U)<^{\dag,+}\overline{b}$ never holds true;
  \item (bar,bar): $\overline{a}<^{\dag,+}\overline{b}$ never holds true.
\end{enumerate}

For $k\geq 1$, the lower order $<^{\dag,-}$ in  $S^\dag_k$ can be
characterized as follows ($<^\sim$ is the lower order in $S$):
\begin{enumerate}
  \item (cut,cut): $(x,L,U)<^{\dag,-}(y,L',U')$ iff $x<^\sim y$;
  \item (bar,cut): $\overline{a}<^{\dag,-}(y,L',U')$ iff $\gamma(a)<^\sim y$;
  \item (cut,bar): $(x,L,U)<^{\dag,-}\overline{b}$ iff $x<^\sim\gamma(b)$;
  \item (bar,bar): $\overline{a}<^{\dag,-}\overline{b}$ iff $\gamma(a)<^\sim\gamma(b)$.
\end{enumerate}

 {\em Strictness.} The strictness is obvious from the above description of
 $<^{\dag,+}$.  Note that all faces in $S^\dag_0$ are cuts. So
 $<^{\dag,+}$ on $S^\dag_0$ is a linear order since $<^+$ is.

{\em Disjointness.} With the description of $<^{\dag,+}$ and
$<^{\dag,-}$ above the disjointness is a matter of a simple check
using disjointness of $<^{+}$ and $<^{\sim}$.

{\em Pencil linearity.} Let $\overline{a}$, $\overline{b}$ be two
different bars in $S^\dag$ and $(x,L,U)$, $(y,L,U)$ be two
different cuts in $S^\dag$. To show $\gamma$-linearity we need to
consider three cases:
\begin{enumerate}
  \item $\gamma(x,L,U)=\gamma(y,L',U')$,
  \item $\gamma(\overline{a})=\gamma(x,L,U)$,
  \item $\gamma(\overline{a})=\gamma(\overline{b})$.
\end{enumerate}

Ad 1. We have $(\gamma(x),-,\ua x)=(\gamma(y),-,\ua y)$.  If $x=y$
then either $L\subseteqnot L'$ or $L'\subseteqnot L$. Thus
$x\perp^+y$. If $x\neq y$ and $\gamma(x)=\gamma(y)$ then either
$x\perp^+y$ or $x\perp^\sim y$.  In case $x\perp^+y$ we have
$(x,L,U)\perp^+(y,L,U)$. We shall show that $x\perp^\sim y$ is
impossible. Suppose $x<^\sim y$. As $\gamma(x)=\gamma(y)$, it
follows that $y$ is a loop.  Let $c\in\cI_{\gamma(y)}$ be an
initial face such that $\gamma(c)\leq^+y$. Then $x<^\sim\gamma(c)$
and $y\not<^\sim \gamma(a)$, i.e. $\ua x\neq\ua y$, contrary to
the supposition. Thus $x\perp^\sim y$ cannot hold true.

Ad 2. We have $(\gamma\gamma(a),-,\ua \gamma(a))=(\gamma(x),-,\ua
x)$.  As $\gamma\gamma(a)=\gamma(x)$, we have either $\gamma(a)=x$
or $\gamma(a)\perp^+x$ or $\gamma(a)\perp^\sim x$. If
$\gamma(a)\leq^+x$ then
$\overline{a}<^+(\gamma(a),\emptyset,-)\leq^+(x,L,U)$.  The other
conditions are impossible. The condition $x<^+\gamma(a)$ is
impossible by Lemma \ref{fact16new0}, and the condition
$\gamma(a)\perp^\sim x$ is impossible as it is easily seen that we
were to have $\ua x\neq\ua \gamma(a)$.

Ad 3. We shall show that this case, i.e. $(\gamma\gamma(a),-,\ua
\gamma(a))=(\gamma\gamma(b),-,\ua \gamma(b))$ is impossible. As
$a,b\in\cI_{\gamma\gamma(a)}$ then $\gamma(a)\perp^\sim\gamma(b)$.
Suppose $\gamma(a)<^\sim\gamma(b)$. Then $b\in \gamma(a)$ and
$b\not\in\gamma(b)$. So we cannot have $\ua \gamma(a))=\ua
\gamma(b)$. This ends the proof of $\gamma$-linearity.

Finally, to verify $\delta$-linearity we need to consider the
following four cases:
\begin{enumerate}
  \item $\overline{z}\in\delta(x,L,U)\cap\delta(y,L',U')$,
  \item $(t,L'',U'')\in\delta(x,L,U)\cap\delta(y,L',U')$,
  \item $\delta(\overline{a})\in\delta(x,L,U)$,
  \item $\delta(\overline{a})=\delta(\overline{b})$.
\end{enumerate}

Ad 1. In this case we have $z\leq^+x$ and $z\leq^+y$. Thus by
Lemma \ref{le_linearity}, $x\perp^+ y$ or $x=y$. In both cases
$(x,L,U)\perp^+(y,L',U')$.

Ad 2. In this case $t\in\delta(x)\cap\delta(y)$ and
$(t,L'',U'')=(t,\da x,-)=(t,\da y,-)$. Thus either $x\perp^+y$ or
$x\perp^\sim y$. In the former case we have
$(x,L,U)\perp^+(y,L',U')$. We shall show that the later case is
impossible. Suppose $x<^\sim y$. Then $x\in S^\lambda$ and hence
there is $a\in\cI$ such that $\gamma(a)\leq^+x$. So
$a\in(\da_ty-\da_tx)$ and $(t,\da x,-)\neq(t,\da y,-)$ contrary to
the supposition.

Ad 3. In this case $\gamma\gamma(a)\in\delta(x)$ and
$(\gamma\gamma(a),\da \gamma(a),-)=(\gamma\gamma(a),\da x,-)$.
Thus either $\gamma(a)\perp^+x$ or $\gamma(a)\perp^\sim x$. If
$\gamma(a)\leq^+x$ then $\overline{a}<^+(x,L,U)$. The remaining
cases are impossible. $x<^+\gamma(a)$ is impossible by Lemma
\ref{fact16new0}, and if we were to have $\gamma(a)\perp^\sim x$
we would have
$\da_{\gamma\gamma(a)}\gamma(a)\neq\da_{\gamma\gamma(a)}x$.

Ad 4. We shall show that this case $(\gamma\gamma(a),\da
\gamma(a),-)=(\gamma\gamma(b),\da \gamma(b),-)$ is impossible. As
$a,b\in\cI_{\gamma\gamma(a)}$ we have
$\gamma(a)\perp^\sim\gamma(b)$. Say $\gamma(a)<^\sim\gamma(b)$.
Then $a\in\da\gamma(b)-\da\gamma(a)$ and
$\delta(\overline{a})\neq\delta(\overline{b})$ after all. This
ends the proof of $\delta$-linearity.

The fact that $q_S:S^\dag\lra S$ is a positive cover with the
kernel $\overline{\cI}$. $~\Box$

The theorem below show that if we take a positive cover of a
quotient by an unary ideal then we get the ordered face structure
back. Thus it shows that if we deal with unary ideals only the
construction of taking quotient of a positive face structure and
taking a positive cover of an ordered face structure are mutually
inverse.


\section{$k$-domains and $k$-codomains of ordered face
structures}

For any $k\in\o$, we introduce two operations
\[ \bd^{(k)}, \bc^{(k)}: Ob(\ofs) \lra Ob(\okfs) \]
of the $k$-th {\em domain} and the $k$-th {\em codomain}.

For a given ordered face structure $T$ the we shall define
$\bd^{(k)} T$ and $\bc^{(k)} T$ via convex subhypergraphs $d^{(k)}
T$ and $c^{(k)} T$ of $T$.  Then we shall put
\[ \bd^{(k)} T=[d^{(k)} T],\;\;\;\;\;\; \bc^{(k)} T=[c^{(k)} T]. \]

The operations $d^{(k)}X$ and $c^{(k)}X$ are defined for any
convex subset of any ordered face structure $T$. We put, for
$l\in\o$,
\[ (d^{(k)}X)_l = \left\{ \begin{array}{ll}
            \emptyset &  \mbox{if $l> k$,}  \\
            X_k-\gamma(X^{-\lambda}_{k+1}) &\mbox{if $l=k$,} \\
            X_l &\mbox{if $l<k$,}
                                    \end{array}
                \right. \]
and
\[ (c^{(k)}X)_l = \left\{ \begin{array}{ll}
            \emptyset &  \mbox{if $l> k$,}  \\
            X_k-\delta(X^{-\lambda}_{k+1})&\mbox{if $l=k$,} \\
            X_{k-1}-\iota(X_{k+1}) &\mbox{if $l=k-1$,} \\
            X_l &\mbox{if $l<k-1$.}
                                    \end{array}
                \right. \]
{\em Example.} Here is an example of an ordered face structure $T$
and its $1$-domain and $1$-codomain:
\begin{center} \xext=2400 \yext=660
\begin{picture}(\xext,\yext)(\xoff,\yoff)
\put(0,500){$T$}
  \put(50,130){\oval(100,100)[b]}
  \put(0,130){\line(1,4){40}}
  \put(100,130){\vector(-1,4){40}}
  \put(30,120){$^\Da$}
   \put(30,60){$^{a}$}
   \put(70,0){$^{x}$}
    \put(30,330){$s$}
\put(500,500){$\bd^{(1)} T$}   \put(600,330){$s$}

  \put(1050,500){$c^{(1)} T$}
  \put(1100,130){\oval(100,100)[b]}
  \put(1050,130){\line(1,4){40}}
  \put(1150,130){\vector(-1,4){40}}
  \put(1120,0){$^{x}$}
  \put(1080,330){$s$}

\put(1500,500){$\bc^{(1)} T$}
 \putmorphism(1700,300)(1,0)[(s,\emptyset,\{ x \})`(s,\{ x \},\emptyset)`_{x}]{700}{1}b
\end{picture}
\end{center}
The following is a more involved example. With ordered face
structure $S$ as below
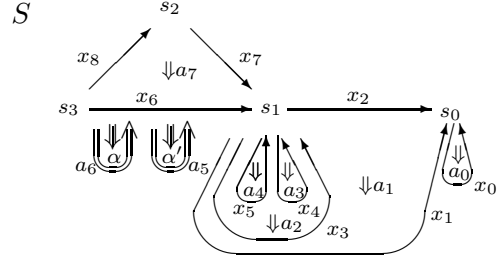
\begin{center} \xext=1500 \yext=840
\begin{picture}(\xext,\yext)(\xoff,\yoff)
\put(-200,800){$S$}
 \settriparms[-1`1`1;350]
 \putAtriangle(0,500)[_{s_2}`_{s_3}`_{s_1};_{x_8}`_{x_7}`]
 \put(230,490){$^{x_6}$}
  \put(320,580){$^{\Da a_7}$}
\putmorphism(700,500)(1,0)[\phantom{\bullet}`_{s_0}`_{x_2}]{620}{1}a

   \put(150,350){\oval(100,100)[b]}
  \put(150,350){\oval(130,130)[b]}
  \put(100,350){\line(0,1){80}}
  \put(85,350){\line(0,1){80}}
  \put(200,350){\line(0,1){80}}
  \put(215,350){\line(0,1){80}}
  \put(175,400){$\wedge$}
  \put(120,380){$\Da$}
  \put(130,280){$^\alpha$}
  \put(150,370){\line(0,1){70}}
  \put(20,250){$^{a_6}$}

      \put(350,350){\oval(100,100)[b]}
  \put(350,350){\oval(130,130)[b]}
  \put(300,350){\line(0,1){80}}
  \put(285,350){\line(0,1){80}}
  \put(400,350){\line(0,1){80}}
  \put(415,350){\line(0,1){80}}
  \put(375,400){$\wedge$}
  \put(320,380){$\Da$}
  \put(320,280){$^{\alpha'}$}
  \put(350,370){\line(0,1){70}}
  \put(410,250){$^{a_5}$}

 \put(830,150){\oval(800,300)[b]}
  \put(430,150){\line(1,2){125}}
  \put(1230,150){\vector(1,4){75}}
  \put(1000,180){$^\Da$}
  \put(1050,180){$^{a_1}$}

  \put(580,230){\line(1,2){85}}
  \put(680,230){\vector(0,1){180}}
  \put(620,230){$^\Da$}

  \put(700,200){\oval(400,300)[b]}
  \put(500,200){\line(1,2){100}}
  \put(900,200){\vector(-1,2){100}}
  \put(680,50){$^\Da$}
  \put(730,50){$^{a_2}$}

  \put(630,230){\oval(100,100)[b]}
  \put(580,230){\line(1,2){85}}
  \put(680,230){\vector(0,1){180}}
  \put(620,230){$^\Da$}
    \put(600,165){$^{a_4}$}
     \put(565,100){$^{x_5}$}

  \put(770,230){\oval(100,100)[b]}
   \put(720,230){\line(0,1){180}}
  \put(820,230){\vector(-1,2){85}}
   \put(730,230){$^\Da$}
   \put(740,165){$^{a_3}$}
   \put(790,100){$^{x_4}$}

  \put(1340,290){\oval(100,100)[b]}
  \put(1290,290){\line(1,4){40}}
  \put(1390,290){\vector(-1,4){40}}
  \put(1320,290){$^\Da$}
  \put(1310,230){$^{a_0}$}
  \put(1400,180){$^{x_0}$}
  \put(1250,80){$^{x_1}$}
   \put(890,30){$^{x_3}$}
\end{picture}
\end{center}
its $1$-domain is
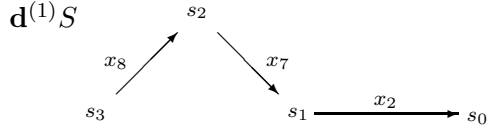
\begin{center} \xext=1500 \yext=400
\begin{picture}(\xext,\yext)(\xoff,\yoff)
\put(-300,300){$\bd^{(1)}S$}
 \settriparms[-1`1`0;350]
 \putAtriangle(0,0)[_{s_2}`_{s_3}`_{s_1};_{x_8}`_{x_7}`]
\putmorphism(700,0)(1,0)[\phantom{\bullet}`_{s_0}`_{x_2}]{620}{1}a
\end{picture}
\end{center}
the convex subset of $S$ defining $1$-codomain is
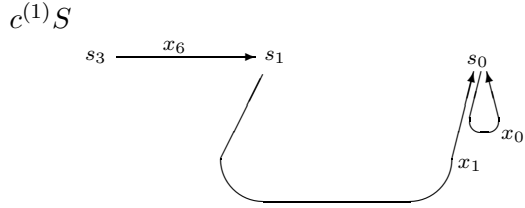
\begin{center} \xext=1500 \yext=700
\begin{picture}(\xext,\yext)(\xoff,\yoff)
\put(-300,600){$c^{(1)}S$}
 \settriparms[0`0`1;310]
 \putAtriangle(0,500)[`_{s_3}`_{s_1};``]
 \put(230,490){$^{x_6}$}
\putmorphism(700,500)(1,0)[\phantom{\bullet}`_{s_0}`]{620}{0}a

 \put(830,150){\oval(800,300)[b]}
  \put(430,150){\line(1,2){145}}
  \put(1230,150){\vector(1,4){75}}

  \put(1340,290){\oval(100,100)[b]}
  \put(1290,290){\line(1,4){40}}
  \put(1390,290){\vector(-1,4){40}}
  \put(1400,180){$^{x_0}$}
  \put(1250,80){$^{x_1}$}
\end{picture}
\end{center}
and finally the $1$-codomain of $S$ is
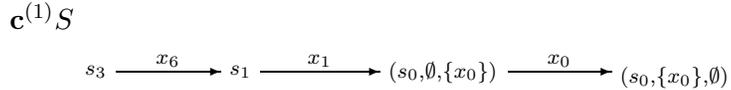
\begin{center} \xext=1500 \yext=700
\begin{picture}(\xext,\yext)(\xoff,\yoff)
\put(-300,150){$\bc^{(1)}S$}
\putmorphism(0,0)(1,0)[_{s_3}`_{s_1}`_{x_6}]{500}{1}a
 \putmorphism(500,0)(1,0)[\phantom{_{s_1}}`_{(s_0,\emptyset,\{ x_0 \})}`_{x_1}]{700}{1}a
  \putmorphism(1200,0)(1,0)[\phantom{_{(s_0,\emptyset,\{ x_0 \}})}`_{(s_0,\{ x_0 \},\emptyset)}`_{x_0}]{800}{1}a
\end{picture}
\end{center}

We have

\begin{lemma} Let  $T$ be an ordered face structure, $X$ a convex subset in $T$.
The subhypergraphs $d^{(k)} X$ and $c^{(k)} X$ of $T$ are convex.
Moreover, for $X=T$, $\cE^{d^{(k)} T}$ is empty, i.e. there are no
empty loops in $d^{(k)}T$ (hence $d^{(k)}T=\bd^{(k)}T$) and all
empty loops in $\cE^{c^{(k)} T}$ have dimension $k$.
\end{lemma}

{\it Proof.}~ The fact that $d^{(k)} X$ and $c^{(k)} X$ are convex
sets is an easy consequence Lemmas \ref{fact8} and \ref{fact13}.
$\cE^{d^{(k)} T}$ is empty by loop-filling.  The empty loops in
$\cE^{c^{(k)} T}$ have dimension $k$ by globularity. $~\Box$

Thus the ordered face structures $\bd^{(k)} T$ and $\bc^{(k)} T$
are well defined. We denote $\nu_{d^{(k)} T}$ by $\bd^{(k)}_T$ and
$\nu_{c^{(k)} T}$ by $\bc^{(k)}_T$. Thus we have defined a diagram
in $\ofs$:
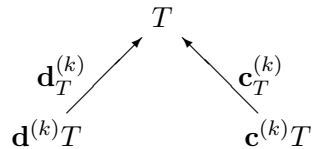
\begin{center} \xext=800 \yext=450
\begin{picture}(\xext,\yext)(\xoff,\yoff)
\settriparms[-1`-1`0;400]
 \putAtriangle(0,50)[T`\bd^{(k)}T`\bc^{(k)}T;\bd^{(k)}_T`\bc^{(k)}_T`]
\end{picture}
\end{center}

{\em Example.} Let $X\subset Y$ be convex subsets of an ordered
face structure $T$ as shown on the diagram below.
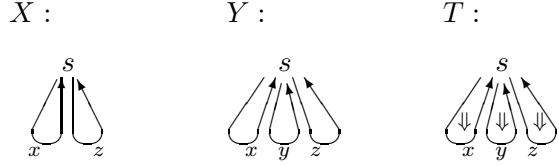
\begin{center} \xext=1500 \yext=660
\begin{picture}(\xext,\yext)(\xoff,\yoff)
\put(0,600){$X:$}
  \put(130,230){\oval(100,100)[b]}
  \put(80,230){\line(1,2){85}}
  \put(180,230){\vector(0,1){180}}
     \put(65,100){$^{x}$}
    \put(180,430){$s$}

  \put(270,230){\oval(100,100)[b]}
   \put(220,230){\line(0,1){180}}
  \put(320,230){\vector(-1,2){85}}
   \put(290,100){$^{z}$}

\put(750,600){$Y:$}
  \put(950,230){\oval(100,100)[b]}
  \put(900,230){\line(1,4){40}}
  \put(1000,230){\vector(-1,4){40}}
   \put(930,100){$^{y}$}

  \put(810,230){\oval(100,100)[b]}
  \put(760,230){\line(2,3){120}}
  \put(860,230){\vector(1,3){60}}
     \put(815,100){$^{x}$}
    \put(930,430){$s$}

  \put(1090,230){\oval(100,100)[b]}
   \put(1040,230){\line(-1,3){60}}
  \put(1140,230){\vector(-2,3){120}}
   \put(1040,100){$^{z}$}

\put(1500,600){$T:$}
  \put(1700,230){\oval(100,100)[b]}
  \put(1650,230){\line(1,4){40}}
  \put(1750,230){\vector(-1,4){40}}
  \put(1680,200){$^\Da$}
   \put(1680,100){$^{y}$}

  \put(1560,230){\oval(100,100)[b]}
  \put(1510,230){\line(2,3){120}}
  \put(1610,230){\vector(1,3){60}}
  \put(1550,200){$^\Da$}
     \put(1565,100){$^{x}$}
    \put(1680,430){$s$}

  \put(1840,230){\oval(100,100)[b]}
   \put(1790,230){\line(-1,3){60}}
  \put(1890,230){\vector(-2,3){120}}
   \put(1810,200){$^\Da$}
   \put(1790,100){$^{z}$}
\end{picture}
\end{center}
Clearly $X\subseteq Y$. And the stretching of $X$ and $Y$ gives
\begin{center}
\xext=1600 \yext=220
\begin{picture}(\xext,\yext)(\xoff,\yoff)
\put(0,200){\makebox(50,50){$[X]:$}}
\putmorphism(0,0)(1,0)[(s,\emptyset,-)`(s,\{ x \},-)`x]{800}{1}a
 \putmorphism(800,0)(1,0)[\phantom{(s,\{ x \},-)}`(s,\{ x,z\},-)`z]{800}{1}a
\end{picture}
\end{center}
and
\begin{center}
\xext=2800 \yext=220
\begin{picture}(\xext,\yext)(\xoff,\yoff)
\put(0,200){\makebox(50,50){$[Y]:$}}
\putmorphism(0,0)(1,0)[(s,\emptyset,-)`(s,\{ x \},-)`x]{800}{1}a
 \putmorphism(800,0)(1,0)[\phantom{(s,\{ x \},-)}`(s,\{ x,y\},-)`y]{900}{1}a
 \putmorphism(1700,0)(1,0)[\phantom{(s,\{ x,y \},-)}`(s,\{ x,y,z\},-)`z]{1000}{1}a
\end{picture}
\end{center}
respectively. Clearly there is no natural map from $[X]$ to $[Y]$.

This shows that there might be no natural comparison map between
stretchings even if one of the convex subset is contained in the
other. The Lemma below says that however in some important cases
we do have such comparison maps.

\begin{lemma} Let  $T$ be an ordered face structure, $X$ a convex subset in $T$.
The embeddings $d^{(k)}X\lra X$ and $c^{(k)}X\lra X$ induce
monotone morphisms $d^{(k)}_X :[d^{(k)}X]\lra [X]$ and $c^{(k)}_X
:[c^{(k)}X]\lra [X]$ so that the triangles
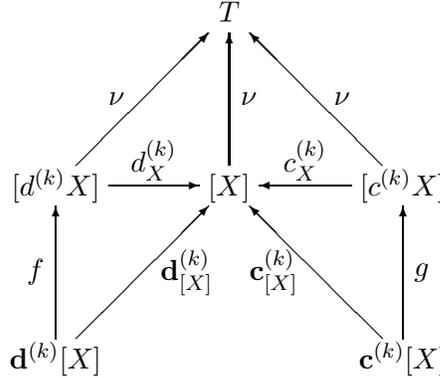
\begin{center}
\xext=1200 \yext=1200
\begin{picture}(\xext,\yext)(\xoff,\yoff)
\settriparms[-1`-1`1;600]
 \putdtriangle(0,600)[T`[d^{(k)}X]`[X];\nu`\nu`]
 \settriparms[0`-1`-1;600]
 \putbtriangle(600,600)[\phantom{T}`\phantom{[X]}`[c^{(k)}X];`\nu`]
 \settriparms[0`-1`-1;600]
 \putptriangle(0,0)[\phantom{[d^{(k)}X]}`\phantom{[X]}`\bd^{(k)}[X];d^{(k)}_X`f`\bd^{(k)}_{[X]}]
 \settriparms[0`-1`-1;600]
 \putqtriangle(600,0)[\phantom{[X]}`\phantom{[c^{(k)}X]}`\bc^{(k)}[X];c^{(k)}_X`\bc^{(k)}_{[X]}`g]
\end{picture}
\end{center}
commute, where $f$ and $g$ are monotone isomorphisms.
\end{lemma}

{\it Proof.}~The morphisms $\nu$ send cuts over a face to that face.
The commutation of the upper triangles comes to the observation
(see below) that both $d^{(k)}_X $ and $c^{(k)}_X$ sends cuts over
$a$ to cuts over $a$ for any $a$ in $d^{(k)}X $ and $c^{(k)}X$,
respectively.

Next we deal with the left lower triangle
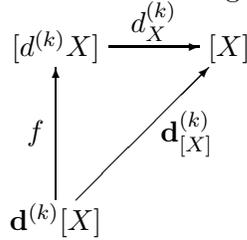
\begin{center}
\xext=1200 \yext=700
\begin{picture}(\xext,\yext)(\xoff,\yoff)
\settriparms[-1`-1`1;600]
 \settriparms[1`-1`-1;600]
 \putptriangle(0,0)[{[d^{(k)}X]}`[X]`\bd^{(k)}[X];d^{(k)}_X`f`\bd^{(k)}_{[X]}]
\end{picture}
\end{center}
 In dimensions $l<k$, we have
$[X]_l=[d^{(k)}X]_l=\bd^{(k)}[X]_l$ and
\[ f_l=(\bd^{(k)}_{[X]})_l=(d^{(k)}_X)_l=id_{[X]_l} \]
In dimension $k$, we have
\[ \bd^{(k)}[X]_k =\{ (a,\emptyset,\cE^X_a):a\in X_k-\gamma(X_{k+1}^{-\lambda}) \}, \]
\[ [d^{(k)}X]_k =\{  (a,\emptyset,\emptyset):a\in X_k-\gamma(X_{k+1}^{-\lambda})\} \]
and
\begin{center}
\xext=1200 \yext=600
\begin{picture}(\xext,\yext)(\xoff,\yoff)
 \settriparms[1`-1`0;600]
 \putptriangle(0,0)[{(a,\emptyset,\emptyset)}`(a,\emptyset,{\cal E}^X_a)`(a,\emptyset,{\cal E}^X_a);``]
 \put(-40,70){\line(1,0){80}}
  \put(190,560){\line(0,1){80}}
   \put(90,140){\line(1,-1){100}}
   \putmorphism(205,165)(1,1)[``]{300}{1}r
\end{picture}
\end{center}
From the description it is clear that both triangles commute and
that $f$ is an iso.

Now we shall describe the lower right triangle
\begin{center}
\xext=1200 \yext=700
\begin{picture}(\xext,\yext)(\xoff,\yoff)
 \settriparms[-1`-1`-1;600]
 \putqtriangle(600,0)[{[X]}`[c^{(k)}X]`\bc^{(k)}[X];c^{(k)}_X`\bc^{(k)}_{[X]}`g]
\end{picture}
\end{center}
In this case we need to look at the cells of both dimensions $k$
and $k-1$. In lower dimensions this triangle is, as in the
previous case, the triangle of identities. To describe the above
diagram, we shall describe the diagram
\begin{center} \xext=1400 \yext=550
\begin{picture}(\xext,\yext)(\xoff,\yoff)
\setsqparms[-1`-1`-1`-1;700`450]
 \putsquare(0,0)[X`c^{(k)}X`{[X]}`{c^{(k)}[X]};`\nu``]
 \setsqparms[-1`0`-1`-1;700`450]
 \putsquare(700,0)[\phantom{c^{(k)}X}`[c^{(k)}X]`\phantom{c^{(k)}[X]}`{[c^{(k)}[X]]};\nu``g`\nu]
\end{picture}
\end{center}
As the horizontal arrows in the left hand square are inclusions we
need to describe only the right hand square.  In dimension $k$, we
have
\[  (c^{(k)}X)_k = X_k  -\delta(X^{-\lambda}_{k+1}) \]
\[  [c^{(k)}X]_k =\{ (a,\emptyset,\emptyset)\in\cC(\cE^{c^{(k)}X}_a) : a\in X_k  -\delta(X^{-\lambda}_{k+1}) \} \]
\[ c^{(k)}[X]_k =[X]-\delta([X]^{-\lambda}_{k+1}) =
    \{(a,-,\emptyset)\in\cC(\cE^X_a):a\in X_k -\delta(X^{-\lambda}_{k+1})\}\]
\[  \bc^{(k)}[X]_k =[c^{(k)}[X]]_k =\{((a,-,\emptyset),\emptyset,\emptyset)\in \cC(\cE^{c^{(k)}[X]}_{(a,-,\emptyset)}) :
     (a,-,\emptyset)\in c^{(k)}[X]_k\} \]
and the commutation of the square is
\begin{center} \xext=1400 \yext=520
\begin{picture}(\xext,\yext)(\xoff,\yoff)
 \setsqparms[-1`-1`-1`-1;700`450]
 \putsquare(700,0)[a`(a,\emptyset,\emptyset)`(a,\emptyset,\emptyset)`((a,-,\emptyset),\emptyset,\emptyset);```]
  \put(660,70){\line(1,0){80}}
  \put(1360,70){\line(1,0){80}}
  \put(1210,410){\line(0,1){80}}
  \put(1060,40){\line(0,-1){80}}
\end{picture}
\end{center}
So the diagram in dimension $k$ commutes and $g_k$ is a bijection.

In dimension $k-1$ we have
\[  (c^{(k)}X)_{k-1} = X_{k-1}  -\iota(X_{k+1}) \]
\[  [c^{(k)}X]_{k-1} =\{ (x,L_0,U_0)\in\cC(\cE^{c^{(k)}X}_x) : x\in X_{k-1}  -\iota(X_{k+1}) \} \]
\[ c^{(k)}[X]_{k-1} = \{ (x,L_1,U_1)\in\cC(\cE^X_x) : \mbox{ there is no } \alpha\in X,\, {\rm such\; that\;} \]
\[\hskip 3cm \exists_{a,b\in\delta(\alpha)}\, \gamma(a)=x\in\delta(b),\, (a,\da
\alpha,-),(b,\da \alpha,-)\in [X]^{-\lambda}\,  \] \[\hskip 3cm
{\rm and\; } (\gamma(a),-,\ua a)=(x,L_1,U_1)=(x,\da b,-)  \}\]
\[  \bc^{(k)}[X]_{k-1} =[c^{(k)}[X]]_{k-1} = \hskip 6cm \]
\[ \hskip 1cm = \{((x,L_1,U_1),L_2,U_2)\in \cC(\cE^{c^{(k)}[X]}_{(x,L_1,U_1)}) :
     (x,L_1,U_1)\in c^{(k)}[X]_{k-1}\} \]
and the commutation of the square is
\begin{center} \xext=1400 \yext=520
\begin{picture}(\xext,\yext)(\xoff,\yoff)
 \setsqparms[-1`-1`-1`-1;1000`450]
 \putsquare(700,0)[x`(x,L_0,U_0)`(x,L_1,U_1)`((x,L_1,U_1),L_2,U_2);```]
  \put(660,70){\line(1,0){80}}
  \put(1660,70){\line(1,0){80}}
  \put(1445,410){\line(0,1){80}}
  \put(1250,40){\line(0,-1){80}}
\end{picture}
\end{center}
where the bijective correspondence between cuts $(x,L_0,U_0)$ in $
\cC(\cE^{c^{(k)}X}_x)$ and the cuts of cuts
$((x,L_1,U_1),L_2,U_2)$ in $[c^{(k)}[X]]_{k-1}$ is described
below.

First we introduce a piece of notation. We denote the faces
$a_{in}, a_{out}\in X^{-\lambda}_k-\delta(X^{-\lambda}_{k+1})$
such that $\gamma(a_{in})=x\in\delta(a_{out})$.  Such faces do not
need to exists but if they do they are unique. We have
\[ L_0=\{ l\in\cE^{c^{(k)}X}_x : (l,-,\emptyset)\in L_2 \mbox{ \rm  or }\exists_{a\in L_1}\, a\leq^+l
  \},  \]
\[ U_0=\{ l\in\cE^{c^{(k)}X}_x : (l,-,\emptyset)\in U_2 \mbox{ \rm  or }\exists_{a\in U_1}\, a\leq^+l
 \},  \]
\[ L_1=\{ a\in\cE^X_x : \exists_{l\in L_0}\, a\leq^+l \mbox{ \rm  or }
a_{in} \mbox{ \rm exists and } a<^+a_{in} \},  \]
\[ U_1=\{ a\in\cE^X_x : \exists_{l\in U_0}\, a\leq^+l \mbox{ \rm  or }
a_{out} \mbox{ \rm exists and } a<^+a_{out},  \}  \]
\[ L_2=\{ (l,-,\emptyset)\in\cE^{c^{(k)}[X]}_{(x,-,\ua l)} :  l\in L_0  \},  \]
\[ U_2=\{ (l,-,\emptyset)\in\cE^{c^{(k)}[X]}_{(x,-,\ua l)} : l\in U_0  \}.  \]
It is a matter of a check to see that this correspondence is
bijective and that $g$ is indeed an iso.  Note that in this
notation the map $c^{(k)}_X:[c^{(k)}X]\lra [X]$ is given by
\[ (x,L_0,U_0)\mapsto (x,L_1,U_1).\]
 $~\Box$

\begin{proposition}\label{dom and codom via covers}
Let $q:S\ra T$ be  a positive cover, $\cI=ker(q)$ be an ideal in
$T$ determining this cover. Then we have positive covers
$\bd^{(k)}(q):\bd^{(k)}S\ra \bd^{(k)}T$ and
$\bc^{(k)}(q):\bc^{(k)}S\ra \bc^{(k)}T$, with kernels
$\cI\cap\bd^{(k)}S$ and $\cI\cap\bc^{(k)}S$, respectively, making
both squares
\begin{center} \xext=1400 \yext=650
\begin{picture}(\xext,\yext)(\xoff,\yoff)
 \setsqparms[1`1`1`1;700`500]
 \putsquare(0,50)[\bd^{(k)}S`S`\bd^{(k)}T`T;\bd^{(k)}_S`\bd^{(k)}(q)`q`\bd^{(k)}_T]
 \setsqparms[-1`0`1`-1;700`500]
 \putsquare(700,50)[\phantom{S}`\bc^{(k)}S`\phantom{T}`\bc^{(k)}T;\bc^{(k)}_S``\bc^{(k)}(q)`\bc^{(k)}_T]
 \end{picture}
\end{center}
commute.
\end{proposition}
{\it Proof.}~To see that $d^{(k)}q$ exists, we shall show that if
$a\in S_k$ and $q(a)\in\gamma(T^{-\lambda}_{k+1})$ then
$a\in\gamma(S_{k+1})$.  So pick $\alpha\in T^{-\lambda}_{k+1}$
such that $\gamma(\alpha)=q(a)$.  As $q$ is a cover, there is
$\beta\in S_{k+1}$ such that $q(\beta)=\alpha$. Hence there is a
$\cI$-path from $a$ to $\gamma(\alpha)$ or from $\gamma(\alpha)$
to $a$.  In the latter case $a\in\gamma(S_{k+1})$ and we are done.
So assume that there is an upper $\cI$-path $a,\alpha_1,\ldots
,\alpha_n,\gamma(\beta)$. As $q(\beta)$ is not a loop and
$q(\alpha_i)=1_{\gamma(\beta)}$, we have $\alpha_i<^+\beta$ and
$a\in\iota(S)$. In particular, we have $a\in \gamma(S_{k+1})$, as
required.

Similarly we can show that we have a hypergraph morphism $q'$ as
in the diagram
\begin{center} \xext=1100 \yext=520
\begin{picture}(\xext,\yext)(\xoff,\yoff)
 \setsqparms[-1`1`1`-1;600`450]
 \putsquare(0,0)[S`\bc^{(k)}S`T`c^{(k)}T;`q``]
  \settriparms[0`1`-1;450]
 \putbtriangle(600,0)[\phantom{c^{(k)}S}`\phantom{c^{(k)}T}`[c^{(k)}T];q'`c^{(k)}q`\nu]
\end{picture}
\end{center}
making the square commutes.  We shall show that $q'$ can be lifted
to $c^{(k)}q : \bc^{(k)}S \lra [\bc^{(k)}T]$. As the only empty
loops in $c^{(k)}T$ have dimension $k$ we need to define the
function
\[ (c^{(k)}q)_{k-1} : S_{k-1}-\iota(S_{k+1}) \lra
    \bigcup \{ \cC(\cE^{c^{(k)}T}_x) : x\in T_{k-1}-\iota(T_{k+1}) \} \]
only.  For $x\in S_{k-1}-\iota(S_{k+1})$, we put

\[  (c^{(k)}q)_{k-1}(x) = \left\{ \begin{array}{ll}
            1_{(\gamma(q'(x)),\emptyset,\emptyset)} & \mbox{if $q'(x)\in 1_{T_{k-2}}$,}  \\
            (x,\emptyset,\emptyset) &
            \mbox{$\cE^{c^{(k)}T}_{q'(x)}=\emptyset$,}\\
             (x,\da a,\emptyset) &
            \mbox{$a\in S_k-\delta(S_{k+1})$ and $x\in\delta(a)$,}\\
            (x,\emptyset,\ua b) &
            \mbox{$b\in S_k-\delta(S_{k+1})$ and $x=\gamma(b)$.}
                                    \end{array}
                \right. \]
As for any $x\in S_{k-1}-\iota(S_{k+1})$ if $q'(x)\in T_{k-1}$ and
$\cE^{c^{(k)}T}_{q'(x)}\neq\emptyset$ then either there is a
unique $a\in S_k-\delta(S_{k+1})$ such that $x\in\delta(a)$ or
there is a unique $b\in S_k-\delta(S_{k+1})$ such that
$x=\gamma(b)$, $(c^{(k)}q)_{k-1}$ is well defined. The remaining
details are left for the reader. $~\Box$

In particular, from this Proposition and Theorem \ref{posivive
cover}, we have

\begin{corollary}\label{posivive cover d and c}
Let $S$ be an ordered face structure, $q_S:S^\dag\lra S$ it's
positive cover, with kernel $\overline{\cI}$, as defined in
section \ref{positive covers}. Then, $dim(\overline{\cI})<dim(S)$,
$\overline{\cI}\cap\bc(S^\dag) =\overline{\cI}_{\leq n-2}$  and
\[ \bc(S^\dag)_{/\overline{\cI}_{\leq n-2}}\cong \bc S,\;\;\;\;\;\;\;
\bd(S^\dag)_{/\overline{\cI}}\cong \bd S.  \]
\end{corollary}

The globularity equations for ordered face structures can be
deduced from the above Proposition.

\begin{proposition}\label{o-cat-dom-codom} Let $S$ be an ordered face
structure $k,l\in\o$, $k<l\leq dim(S)$. Then the diagram
  \begin{center}
\xext=1000 \yext=990
\begin{picture}(\xext,\yext)(\xoff,\yoff)
 \settriparms[-1`-1`0;500]
 \putAtriangle(0,500)[S`\bd^{(l)}S`\bc^{(l)}S;\bd^{(l)}_S`\bc^{(l)}_S`]
 \setsqparms[0`-1`-1`0;1000`500]
\putsquare(0,0)[\phantom{\bd^{(l)}S}`\phantom{\bc^{(l)}S}`\bd^{(k)}S`\bc^{(k)}S; `\bd^{(k)}_{\bd^{(l)}S}`\bc^{(k)}_{\bc^{(l)}S}`]

 \put(200,50){\vector(2,1){700}}
 \put(450,225){\vector(-2,1){350}}
 \put(800,50){\line(-2,1){250}}
 \put(640,380){\makebox(100,100){$\bd^{(k)}_{\bc^{(l)}S}$}}
 \put(240,380){\makebox(100,100){$\bc^{(k)}_{\bd^{(l)}S}$}}
\end{picture}
\end{center}
commutes.
\end{proposition}
{\it Proof.}~Having Theorem \ref{posivive cover} and  Proposition \ref{dom and codom via covers} we see that
the above diagram commutes as a consequence of the same diagram being commutative for the positive face structure.
$~\Box$

\section{$k$-tensor squares of ordered face structures}
Let $S$ and $T$ be ordered face structures such that $\bc^{(k)} S=\bd^{(k)} T$. In that case we define the
$k${\em-tensor}\index{tensor} $S\otimes_kT$ of $S$ and $T$ and the $k${\em-tensor square}\index{tensor!square} in $\ofs$
\begin{center}
\xext=850 \yext=550
\begin{picture}(\xext,\yext)(\xoff,\yoff)
 \setsqparms[1`-1`-1`1;850`450]
 \putsquare(0,50)[S`S\otimes_kT`\bc^{(k)} S=\bd^{(k)} T`T;\kappa_S`\bc^{(k)}_S`\kappa_T`\bd^{(k)}_T]
\end{picture}
\end{center}
The local part of $S\otimes_kT$ is defined so that the square
\begin{center}
\xext=850 \yext=550
\begin{picture}(\xext,\yext)(\xoff,\yoff)
 \setsqparms[1`-1`-1`1;850`450]
 \putsquare(0,50)[|S|`|S\otimes_kT|`|\bc^{(k)} S|`|T|;\kappa_S`\bc^{(k)}_S`\kappa_T`\bd^{(k)}_T]
\end{picture}
\end{center}
is a pushout in $\lfs$, so the faces of $S\otimes_kT$ are as in
the following table:
 \[ \begin{array}{|c|c|c|c|}\label{tensor faces}
    dim & \bc^{(k)}S      & \bd^{(k)}T & S\otimes_kT  \\ \hline
    l>k & \emptyset & \emptyset & S_l + T_l  \\ \hline
     k  & S_k-\delta(S^{-\lambda}_{k+1}) &
         T_k-\gamma(T^{-\lambda}_{k+1}) & S_k + \gamma(T^{-\lambda}_{k+1}) =T_k+\delta(S^{-\lambda}_{k+1} )\\ \hline
    k-1 & \cC(S_{k-1}-\iota(S_{k+1}))  & T_{k-1} & S_{k-1}   \\ \hline
    l<k-1 & S_l & T_l & S_l
  \end{array}\]
By the assumption the first and the second columns are equal and
the third describes the faces of $S\otimes_kT$. To simplify the
description of $S\otimes_kT$, we assume that
\[ S_k-\delta(S^{-\lambda}_{k+1}) =
T_k-\gamma(T^{-\lambda}_{k+1})= S_k\cap T_k, \] and we introduce
the notation for the function
\[ [-]=(\bc^{(k)}_S)_{k-1}:(\bc^{(k)}S)_{k-1}=T_{k-1}\lra S_{k-1}.\]
that sends $t$-cuts in $S$, (with $t\in S_{k-1}$, i.e. elements of
$T_{k-1}$) to $t$. All the components of the maps $\kappa_S : S
\lra S\otimes_kT$ and $\kappa_T : T \lra S\otimes_kT$ are
inclusions except for $(\kappa_T)_{k-1}$ which is $[-]$. The
domain and codomain maps in $S\otimes_kT$, denoted
$\gamma^\otimes$ and $\delta^\otimes$ for short, are obvious
except for $k$-faces in $\gamma(T^{-\lambda}_{k+1})$. If
$t\in\gamma(T^{-\lambda}_{k+1})$ we put:
\[ \gamma^{\otimes}(t) = [\gamma^T(t)],\hskip10mm
\delta^{\otimes}(t) = \{ [u] : u\in\delta^T(t)\}.\] To finish off
the definition of $S\otimes_kT$, it is enough to define
$<^{\otimes,l,\sim}$, for $l\geq 1$. For $l<k$,
$<^{(S\otimes_kT)_l,\sim}$ is $<^{S_l,\sim}$  and
$<^{(S\otimes_kT)_l,\sim}$ is $<^{S_l,\sim}+<^{T_l,\sim}$, for
$l>k+1$. Thus, it remains to define the orders
$<^{(S\otimes_kT)_k,\sim}$ and $<^{(S\otimes_kT)_{k+1},\sim}$. The
order $<^{(S\otimes_kT)_{k+1},\sim}$ is defined for
$a,b\in(S\otimes_kT)_{k+1}=S_{k+1} + T_{k+1}$ we put
\[ a<^{(S\otimes_kT)_{k+1},\sim} b \mbox{ iff } \left\{ \begin{array}{ll}
             \mbox{either $a,b\in S_{k+1}$ and $a<^{S_{k+1},\sim} b$,}  \\
             \mbox{or $a,b\in T_{k+1}$ and $a<^{T_{k+1},\sim} b$,}  \\
             \mbox{or $a\in S_{k+1}$, $b\in T_{k+1}$  and $a<^{S\otimes_kT,-} b$.}
                                    \end{array}
                \right. \]
i.e. it is $<^{S,\sim}$ on $S_{k+1}$, $<^{T,\sim}$ on is
$T_{k+1}$, and moreover if the faces comes from different parts
and are $<^{S\otimes_kT,-}$ related, then faces from $S$ comes
before the faces from $T$. The last clause of this definition is
the only reason $S\otimes_kT$ is not a pushout in $\ofs$, in
general. It may cause a face $a$ from $S$ to be $<^\sim$-smaller
than a face $b$ from $T$ even if there is no $\sim$-relation
between $a$ and $b$, whatsoever. By Lemma \ref{local to global},
to define the order $<^{(S\otimes_kT)_k,\sim}$ it is enough to say
that it agrees with $<^{T_k,\sim}$ on the set
$(S\otimes_kT)_k-\delta((S\otimes_kT)^{-\lambda}_{k+1})=(T_k-\delta(T^{-\lambda}_{k+1}))$.
However we give below the full, but more involved, definition of
the order $<^{(S\otimes_kT)_k,\sim}$.  We have
$(S\otimes_kT)_k=S_k + \gamma(T^{-\lambda}_{k+1})
=T_k+\delta(S^{-\lambda}_{k+1} )$. We define
$<^{(S\otimes_kT)_k,\sim}$ to be $<^{S,\sim}$ on $S_{k}$, and to
be $<^{T,\sim}$ on $T_{k+1}$. The essential case  is if
$x\in\delta(S^{-\lambda}_{k+1})$ and
$y\in\gamma(T^{-\lambda}_{k+1})$. In that case there is a unique
$x'\in S_k\cap T_k$ that $x<^+x'$. We put
$x<^{(S\otimes_kT)_k,\sim}y$ iff $x'<^{T_k,\sim}y$ and
$y<^{(S\otimes_kT)_k,\sim}x$ iff $y<^{T_k,\sim}x'$ and
$y<^{(S\otimes_kT)_k,-}x$. In other words for $x,y\in
(S\otimes_kT)_{k}$, we have:
\[ x<^{(S\otimes_kT)_{k},\sim} y \mbox{ iff } \left\{ \begin{array}{l}
            \mbox{either $x,y\in S$ and $x<^{S,\sim}y$,}  \\
            \mbox{or $x,y\in T$ and $x<^{T,\sim}y$,}  \\
            \mbox{or $x\in \delta(S^{-\lambda}_{k+1})$, $y\in \gamma(T^{-\lambda}_{k+1})$ and
            $\exists_{z\in S_k\cap T_k}\; x\leq^{S,+}z$  and $z<^{T,\sim}y$,}\\
            \mbox{or $x\in \gamma(T^{-\lambda}_{k+1})$, $y\in \delta(S^{-\lambda}_{k+1})$ and
            $x<^{(S\otimes_kT)_{k},-} y$} \\
            \mbox{ \hskip 3mm and $\exists_{z\in S_k\cap T_k}\; x\leq ^{T,+}z$  and $x<^{S,\sim}z$.}
                                    \end{array}
                \right. \]

{\em Examples.} Before we prove some properties of the above
construction let us look at some examples of $k$-tensors:
\begin{center} \xext=1500 \yext=520
\begin{picture}(\xext,\yext)(\xoff,\yoff)
\putmorphism(50,470)(1,0)[S`S\otimes_0 T`\kappa_S]{650}{1}a
\putmorphism(700,470)(1,0)[\phantom{S\otimes_0T}`T`\kappa_T]{650}{-1}a
  \put(630,130){\oval(100,100)[b]}
  \put(580,130){\line(1,2){85}}
  \put(680,130){\vector(0,1){180}}
  \put(620,130){$^\Da$}
    \put(610,65){$^{a}$}
     \put(565,0){$^x$}
     \put(680,320){$s$}

  \put(770,130){\oval(100,100)[b]}
   \put(720,130){\line(0,1){180}}
  \put(820,130){\vector(-1,2){85}}
   \put(730,130){$^\Da$}
   \put(750,65){$^{b}$}
   \put(790,0){$^{y}$}

  \put(50,130){\oval(100,100)[b]}
  \put(0,130){\line(1,4){40}}
  \put(100,130){\vector(-1,4){40}}
  \put(30,130){$^\Da$}
    \put(20,65){$^{a}$}
     \put(115,40){$^{x}$}
     \put(30,310){$s$}

  \put(1340,130){\oval(100,100)[b]}
  \put(1290,130){\line(1,4){40}}
  \put(1390,130){\vector(-1,4){40}}
  \put(1320,130){$^\Da$}
    \put(1310,65){$^{b}$}
     \put(1405,40){$^{y}$}
     \put(1320,310){$s$}
\end{picture}
\end{center}
In this case the only relation that is not coming from the fact
that $S\otimes_0T$ is a pushout locally is $x<^\sim y$. We have
that $x$ comes before $y$ as in 'case of doubt' faces from $S$
comes before those from $T$.

The next example is a bit more involved. For the ordered face
structures
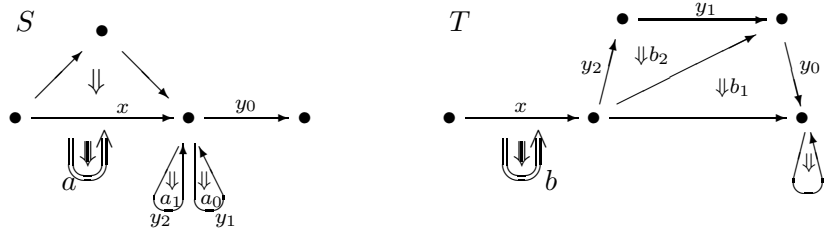
\begin{center} \xext=3000 \yext=750
\begin{picture}(\xext,\yext)(\xoff,\yoff)
\put(100,700){$S$}
 \settriparms[-1`1`1;300]
 \putAtriangle(100,400)[\bullet`\bullet`\bullet;``]
 \put(450,380){$^x$}
  \put(350,500){$\Da$}
\putmorphism(700,400)(1,0)[\phantom{\bullet}`\bullet`_{y_0}]{400}{1}a

   \put(350,250){\oval(100,100)[b]}
  \put(350,250){\oval(130,130)[b]}
  \put(300,250){\line(0,1){80}}
  \put(285,250){\line(0,1){80}}
  \put(400,250){\line(0,1){80}}
  \put(415,250){\line(0,1){80}}
  \put(375,300){$\wedge$}
  \put(320,250){$\Da$}
  \put(350,250){\line(0,1){70}}
  \put(260,150){$a$}

  \put(630,130){\oval(100,100)[b]}
  \put(580,130){\line(1,2){85}}
  \put(680,130){\vector(0,1){180}}
  \put(620,130){$^\Da$}
    \put(600,65){$^{a_1}$}
     \put(565,0){$^{y_2}$}

  \put(770,130){\oval(100,100)[b]}
   \put(720,130){\line(0,1){180}}
  \put(820,130){\vector(-1,2){85}}
   \put(730,130){$^\Da$}
   \put(740,65){$^{a_0}$}
   \put(790,0){$^{y_1}$}

\put(1600,700){$T$}
\putmorphism(1600,400)(1,0)[\bullet`\bullet`_x]{500}{1}a
\putmorphism(2100,400)(1,0)[\phantom{\bullet}`\bullet`]{720}{1}b

 \putmorphism(2200,740)(1,0)[\bullet`\bullet`]{550}{1}b
 \put(2120,450){\vector(1,4){55}}
 \put(2760,660){\vector(1,-4){55}}
 \put(2180,450){\vector(2,1){480}}
  \put(2520,450){$^\Da$} \put(2570,450){$^{b_1}$}
  \put(2240,560){$^\Da$} \put(2290,560){$^{b_2}$}
  \put(2810,530){$^{y_0}$}
  \put(2450,730){$^{y_1}$}
  \put(2055,530){$^{y_2}$}

      \put(1850,250){\oval(100,100)[b]}
  \put(1850,250){\oval(130,130)[b]}
  \put(1800,250){\line(0,1){80}}
  \put(1785,250){\line(0,1){80}}
  \put(1900,250){\line(0,1){80}}
  \put(1915,250){\line(0,1){80}}
  \put(1875,300){$\wedge$}
  \put(1820,250){$\Da$}
  \put(1850,250){\line(0,1){70}}
  \put(1930,150){$b$}

  \put(2840,190){\oval(100,100)[b]}
  \put(2790,190){\line(1,4){40}}
  \put(2890,190){\vector(-1,4){40}}
  \put(2820,190){$^\Da$}
\end{picture}
\end{center}
we have $\bc^{(1)}S=\bd^{(1)}T$
\begin{center} \xext=1500 \yext=40
\begin{picture}(\xext,\yext)(\xoff,\yoff)
\putmorphism(0,0)(1,0)[\bullet`\bullet`_{x}]{400}{1}a
\putmorphism(400,0)(1,0)[\phantom{\bullet}`\bullet`_{y_2}]{400}{1}a
\putmorphism(800,0)(1,0)[\phantom{\bullet}`\bullet`_{y_1}]{400}{1}a
\putmorphism(1200,0)(1,0)[\phantom{\bullet}`\bullet`_{y_0}]{400}{1}a
\end{picture}
\end{center}
and their $1$-tensor $S\otimes_1T$ square is
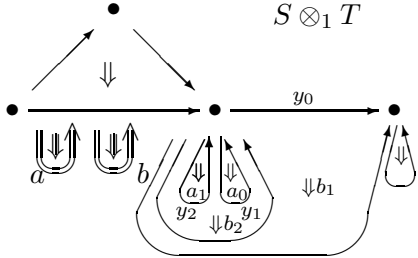
\begin{center} \xext=1500 \yext=870
\begin{picture}(\xext,\yext)(\xoff,\yoff)

\put(900,800){$S\otimes_1T$}
 \settriparms[-1`1`1;350]
 \putAtriangle(0,500)[\bullet`\bullet`\bullet;``]
  \put(300,600){$\Da$}
\putmorphism(700,500)(1,0)[\phantom{\bullet}`\bullet`_{y_0}]{620}{1}a

   \put(150,350){\oval(100,100)[b]}
  \put(150,350){\oval(130,130)[b]}
  \put(100,350){\line(0,1){80}}
  \put(85,350){\line(0,1){80}}
  \put(200,350){\line(0,1){80}}
  \put(215,350){\line(0,1){80}}
  \put(175,400){$\wedge$}
  \put(120,350){$\Da$}
  \put(150,350){\line(0,1){70}}
  \put(60,250){$a$}

      \put(350,350){\oval(100,100)[b]}
  \put(350,350){\oval(130,130)[b]}
  \put(300,350){\line(0,1){80}}
  \put(285,350){\line(0,1){80}}
  \put(400,350){\line(0,1){80}}
  \put(415,350){\line(0,1){80}}
  \put(375,400){$\wedge$}
  \put(320,350){$\Da$}
  \put(350,350){\line(0,1){70}}
  \put(430,250){$b$}

 \put(830,150){\oval(800,300)[b]}
  \put(430,150){\line(1,2){125}}
  \put(1230,150){\vector(1,4){75}}
  \put(1000,180){$^\Da$}
  \put(1050,180){$^{b_1}$}

  \put(580,230){\line(1,2){85}}
  \put(680,230){\vector(0,1){180}}
  \put(620,230){$^\Da$}

  \put(700,200){\oval(400,300)[b]}
  \put(500,200){\line(1,2){100}}
  \put(900,200){\vector(-1,2){100}}
  \put(680,50){$^\Da$}
  \put(730,50){$^{b_2}$}

  \put(630,230){\oval(100,100)[b]}
  \put(580,230){\line(1,2){85}}
  \put(680,230){\vector(0,1){180}}
  \put(620,230){$^\Da$}
    \put(600,165){$^{a_1}$}
     \put(565,100){$^{y_2}$}

  \put(770,230){\oval(100,100)[b]}
   \put(720,230){\line(0,1){180}}
  \put(820,230){\vector(-1,2){85}}
   \put(730,230){$^\Da$}
   \put(740,165){$^{a_0}$}
   \put(790,100){$^{y_1}$}

  \put(1340,290){\oval(100,100)[b]}
  \put(1290,290){\line(1,4){40}}
  \put(1390,290){\vector(-1,4){40}}
  \put(1320,290){$^\Da$}
\end{picture}
\end{center}
with $a<^\sim b$ as the only additional data not following from
the fact that $S\otimes_1T$ is a pushout locally.

\begin{proposition}\label{k-tensor} Let $S$ and $T$ be ordered face
  structures, $k\in\o$, and  $\bc^{(k)} S = \bd^{(k)} T$. Then
$S\otimes_kT$ is an ordered face structure, and the $k$-tensor
square
  \begin{center}
\xext=650 \yext=550
\begin{picture}(\xext,\yext)(\xoff,\yoff)
 \setsqparms[1`-1`-1`1;650`450]
 \putsquare(0,50)[S`S\otimes_kT`\bc^{(k)} S`T;\kappa_S`\bc^{(k)}_S`\kappa_T`\bd^{(k)}_T]
\end{picture}
\end{center}
  commutes in $\ofs$. Moreover the functor $|-|:\ofs\lra \lfs$
  sends the $k$-tensor squares to pushouts.
\end{proposition}
{\it Proof.}~The whole proof is a matter of a check. We shall
discuss globularity leaving the verification of other axioms of
ordered face structure for the reader.

The globularity condition for faces in $S\otimes_kT$ for other
faces than those in $T_k$ and $T_{k+1}$ holds as a direct
consequence of globularity for $S$ and $T$. A simple check shows
that in fact globularity for $T_k$ is also a consequence of
globularity for $T$. Thus we need to verify the globularity for
$a\in T_{k+1}\subseteq (S\otimes_kT)_{k+1}$. We will write
$\gamma$ for $\gamma^T$ and $\gamma^\otimes$ for
$\gamma^{S\otimes_kT}$. For empty-domain faces the globularity is
obvious so we assume that $a\in T_{k+1}^{-\varepsilon}$. Put
\[ L=\{ a\in \gamma^u(T^{-\lambda}_{k+1}) : {\rm there\; is\; a}\;
S^\lambda\cap T_k-{\rm path\;(possibly\; empty) \;
from\;}\delta(a)\;{\rm to\;} \gamma{a} \}. \] We have
\[ (S\otimes_kT)^\lambda_l= \left\{ \begin{array}{ll}
           S^\lambda_l+T^\lambda_l &  \mbox{for  $l>k$,}\\
           S^\lambda_k+L &  \mbox{for  $l=k$,}\\
           S^\lambda_l &  \mbox{for  $l<k$.}
                                    \end{array}
                \right. \]
 We shall describe the
sets involved in the globularity conditions:
\[ \delta^\otimes(a)=\delta(a),\;\;\; \dot{\delta}^{\otimes,-\lambda}(a)=\dot{\delta}^{-\lambda}(a)-L,\]
\[ \gamma^\otimes \gamma^\otimes(a)=[\gamma\gamma(a)],\;\;\;
\delta^\otimes \gamma^\otimes(a)=\{[t]: t\in\delta\gamma(a) ]\},\]
\[ \gamma^\otimes \delta^\otimes(a)=\{ [\gamma(x)] : x\in\delta(a)
\},\;\;\; \delta^\otimes \delta^\otimes(a)=\{[t]:
t\in\delta\delta(a) ] \},\]
\[ \gamma^\otimes \delta^{\otimes,-\lambda}(a)=\{ [\gamma(x)] :
x\in\delta(a)-L \},\;\;\; \delta^\otimes \delta^\otimes(a)=\{[t]:
\exists_{x\in\delta(a)-L}\;t\in\delta(x) ] \}.\] By assumption on
$T$ we have
$\gamma\gamma(a)=\gamma\delta(a)-\delta\dot{\delta}^{-\lambda}(a)$.
Thus to show the $\gamma$-globularity
\[ \gamma^\otimes\gamma^\otimes(a)=\gamma^\otimes\delta^\otimes(a)-
\delta^\otimes\dot{\delta}^{\otimes,-\lambda}(a)\] we need to show
\begin{enumerate}
  \item $\gamma^\otimes\gamma^\otimes(a)\not\in
\delta^\otimes\dot{\delta}^{\otimes,-\lambda}(a)$
  \item
  $\gamma^\otimes\delta^\otimes(a)\subseteq
  \gamma^\otimes\gamma^\otimes(a)\cup\delta^\otimes\dot{\delta}^{\otimes,-\lambda}(a)$.
\end{enumerate}

Ad 1. Suppose 1. does not hold and fix face $x\in\delta(a)-L$,
$t\in\dot{\delta}(x)$ such that there is an upper $(S^\lambda\cap
T)$-path $t,x_0,\ldots,x_k,\gamma\gamma(a)$ from $t$ to
$\gamma\gamma(a)$. In particular, this is
$T_k-\gamma(T^{-\lambda}_{k+1})$-path. As $x\in T_k$, we have
$x_0\leq^+x$.As $\gamma(x)\leq^+\gamma\gamma(a)$, by Path Lemma
and the definition of $L$, we have that $x\in L$, contrary to the
assumption.

Ad 2. Fix $x\in\delta(a)$.  Let
$\gamma(x),x_1,\ldots,x_k,\gamma\gamma(a)$ be the flat upper
$(T_k-\gamma(T^{-\lambda}_{k+1}))$-path. If this path is
$S^\lambda$-path then $[\gamma(x)]=[\gamma\gamma(a)]$, if it is
not then
$[\gamma\gamma(x)]\in\delta^\otimes\dot{\delta}^{\otimes,-\lambda}(a)$,
as required.

For $\delta$-globularity we consider only the case $\gamma(a)\in
T^{-\varepsilon}$. The other case is easy.  For empty-faces in $T$
we have
$\gamma\gamma\delta^\varepsilon(a)\subseteq\theta\delta\gamma(a)$
and hence passing to equivalence classes we also have
$\gamma^\otimes\gamma^\otimes\delta^{\otimes,\varepsilon}(a)\subseteq\theta^\otimes\delta^\otimes\gamma(a)$.
Moreover as
$\delta\gamma(a)=\dot{\delta}\delta(a)-\gamma\dot{\delta}^{-\lambda}(a)$
holds in $T$ to show $
\delta^\otimes\gamma^\otimes(a)=\dot{\delta}^\otimes\delta(a)-
\gamma^\otimes\dot{\delta}^{\otimes,-\lambda}(a)$ we need to show
again two things
\begin{enumerate}
  \item[3.] $\delta^\otimes\gamma^\otimes(a)\cap
  \gamma^\otimes\dot{\delta}^{\otimes,-\lambda}(a)=\emptyset$,
  \item[4.] $\delta^\otimes\delta^\otimes(a)\subseteq\delta^\otimes\gamma^\otimes(a)\cup
  \gamma^\otimes\dot{\delta}^{\otimes,-\lambda}(a)$
\end{enumerate}

Ad 3. Suppose contrary, that 3. does not hold. Fix
$t\in\delta\gamma(a)$ such that
$[t]\in\gamma^\otimes\dot{\delta}^{\otimes,-\lambda}(a)$. So there
is $x\in \delta^{-\lambda}(a)-L$ and upper $(S^\lambda_k\cap
T_k)$-path $t,x_1,\ldots,x_k,\gamma(x)$, with $k\geq 1$. If
$t\in\delta(x)$ or $\gamma(x_i)\in\delta(x)$, for some
$i=1,\ldots,k-1$ then $x\in L$, contrary to the supposition. If
$t\not\in\delta(x)$ and $x_i<^+x$, for some $i=1,\ldots,k$, then,
by Lemma \ref{fact9}.1 and Path Lemma, $\gamma(a)<^+x\in\delta(a)$
which is again a contradiction.  Thus 3. holds.

Ad 4. Fix $t\in\delta(x)$ such that $x\in\delta(a)$. Let
$x_1,\ldots,x_k,t$ be the maximal flat upper
$(T_k-\gamma(T^{-\lambda}_{k+1}))$-path ending at $t$. By Path
Lemma either there is $t'\in\delta\gamma(a)$ such that
$t'\in\delta(x_1)$ or $t'=\gamma(x_i)$, for some $i=1,\ldots,k-1$,
or $x_i<^+\gamma(a)$, for some $i=1,\ldots,k$, $x_i\in
T^\varepsilon$ and $\gamma\gamma(x)\in\theta\delta\gamma(a)$. In
the former case, if the path $t',x_j,\ldots,x_k,t$ is an
$S^\lambda$-path then $[t]\in\delta^\otimes\gamma^\otimes(a)$, if
not then using again Path Lemma we get that
$[t]\in\gamma^\otimes\dot{\delta}^{\otimes,-\lambda}(a)$. In the
later case we can also easily show that
$[t]\in\gamma^\otimes\dot{\delta}^{\otimes,-\lambda}(a)$, as
required. $~\Box$

The following propositions establish a connection between tensor
squares of ordered faces structures and special pushouts of
positive face structures.

\begin{proposition}\label{quotient of sp po}
Let $X$ and $Y$ be positive face structures, $k\in\o$,
$\bc^{(k)}(X)=\bd^{(k)}(Y)$, and $\cJ$ an ideal in the special
pushout\footnote{By this we mean the pushouts, in the category of
positive face structures $\posfs$, of a special kind that have
been described in \cite{Z}.} $X+_kY$. The quotient by ideal $\cJ$
of the special pushout being the top of the following cube
\begin{center}
\xext=1300 \yext=1100
\begin{picture}(\xext,\yext)(\xoff,\yoff)
\setsqparms[1`1`1`1;900`600]
 \putsquare(0,50)[\bc^{(k)}(X)`Y`\bc^{(k)}(X)/_{\cJ^k}`Y/_{\cJ^Y};
 \bd^{(k)}`p`p`\bd^{(k)}]

 \setsqparms[1`0`1`0;900`600]
 \putsquare(380,450)[X`X+_kY`X/_{\cJ^X}`(X+_kY)/_\cJ;
 \kappa``p`]

 \putmorphism(1050,200)(1,1)[``\kappa]{110}{1}r
 \putmorphism(1050,800)(1,1)[``\kappa]{110}{1}l

 \putmorphism(150,200)(1,1)[``]{110}{1}l
 \putmorphism(150,800)(1,1)[``\bc^{(k)}]{110}{1}l
\put(60,250){$\bc^{(k)}$}

\putmorphism(900,450)(1,0)[``]{120}{1}l
 \put(500,450){\line(1,0){360}}
 \put(630,360){$\kappa$}

\putmorphism(360,690)(0,-1)[``]{220}{1}l
 \put(360,980){\line(0,-1){280}}
 \put(240,730){$p$}
\end{picture}
\end{center}
is a $k$-tensor square on the bottom of the following cube, where
$\cJ^X$, $\cJ^Y$, and $\cJ^k$ are the ideals that arise  by
intersecting $\cJ$ with $X$, $Y$, and $\bc^{(k)}(X)$,
respectively. In the cube all squares commutes, and all vertical
maps are covers.
\end{proposition}
{\it Proof.}~   This is a matter of a simple check. $~\Box$

\begin{proposition}\label{cover of tensor sq}
Let $S$ and $T$ be ordered face structures, $k\in\o$, and
$\bc^{(k)} S = \bd^{(k)} T$ and a positive cover
$p:(S\otimes_kT)^\ddag\lra S\otimes_kT$ with the kernel $\cJ$.
Then there are covers $S^{\ddag}\ra S$ and $T^\ddag\ra T$, such
that the top square of the following cube
\begin{center}
\xext=1300 \yext=1100
\begin{picture}(\xext,\yext)(\xoff,\yoff)
\setsqparms[1`1`1`1;900`600]
 \putsquare(0,50)[\bc^{(k)}(S^\ddag)`T^\ddag`\bc^{(k)}(S)`T;\bd^{(k)}`p`p`\bd^{(k)}]

 \setsqparms[1`0`1`0;900`600]
 \putsquare(380,450)[S^\ddag`(S\otimes_kT)^\ddag`S`S\otimes_kT;
 \kappa``p`]

 \putmorphism(1050,200)(1,1)[``\kappa]{110}{1}r
 \putmorphism(1050,800)(1,1)[``\kappa]{110}{1}l

 \putmorphism(150,200)(1,1)[``]{110}{1}l
 \putmorphism(150,800)(1,1)[``\bc^{(k)}]{110}{1}l
\put(60,250){$\bc^{(k)}$}

\putmorphism(900,450)(1,0)[``]{220}{1}l
 \put(500,450){\line(1,0){360}}
 \put(630,360){$\kappa$}

\putmorphism(360,690)(0,-1)[``]{220}{1}l
 \put(360,980){\line(0,-1){280}}
 \put(240,730){$p$}
\end{picture}
\end{center}
is a special pushout in  $\posfs$, and the bottom square is the
quotient $k$-tensor square of the top by the kernel $\cJ$.
\end{proposition}
{\it Proof.}~ We denote $(S\otimes_kT)^\ddag$ by $P$.  We shall
define the positive face structures $S^\ddag$, $T^\ddag$, and the
morphisms from them in the diagram
\begin{center} \xext=1400 \yext=520
\begin{picture}(\xext,\yext)(\xoff,\yoff)
 \setsqparms[1`1`1`1;700`450]
 \putsquare(0,0)[S^\ddag`P`S`S\otimes_kT;\kappa^\ddag_S`q_S`q`\kappa_S]
 \setsqparms[-1`1`1`-1;700`450]
 \putsquare(700,0)[\phantom{P}`T^\ddag`\phantom{S\otimes_kT}`T;\kappa^\ddag_T``q_T`\kappa_S]
\end{picture}
\end{center}
$S$ can be identified with a subset of $S\otimes_kT$ (via
$\kappa_S$). We define $S^\ddag$ as the inverse image of $S$ i.e.
$S^\ddag=q^{-1}(S+1_S)$.  $S^\ddag$ is a positive face structure
as a convex subset of a positive face structure $P$. $q_S$ is the
restriction of $q$ to $S^\ddag$.  It is onto since $q$ is. It is
also easy to see that the kernel of $q_S$ is $\cJ\cap S^\ddag$.

The description of faces of $T^\ddag$ is more involved.
\begin{enumerate}
  \item $T^\ddag_{>k}=q^{-1}(T_{>k}\cup 1_{T_{\geq k}})$,
  \item $T^\ddag_{k}=P_k-\delta(S^\ddag_{k+1})$,
  \item $T^\ddag_{k-1}=P_k-\iota(S^\ddag_{k+1})$,
  \item $T^\ddag_{<k-1}=S^\ddag_{<k-1}(=P_{<k-1})$.
\end{enumerate}
$q_T$ is the restriction of $q$ to $T^\ddag$.

The verification that $q_T$ is a cover and
$\bc^{(k)}S^\ddag=\bd^{(k)}T^\ddag$,  which comes to verification
of two equalities
\[ \bc^{(k)}S^\ddag_{k}=q^{-1}(S_k\cup 1_{S_{k-1}})-\delta(q^{-1}(S_{k+1}\cup
1_{S_{k+1}}))=\] \[ = P_k-\delta(S^\ddag_{k+1})-\gamma(T_{k+1}\cup
1_{T_{k+1}})    = \bd^{(k)}T^\ddag_k\]

\[ \bc^{(k)}S^\ddag_{k-1}=S_{k-1}-\iota(S^\ddag_{k+1})=
P_{k-1}-\iota(S^\ddag_{k+1}) = \bd^{(k)}T^\ddag_{k-1}\]
 is left for the reader. $~\Box$


The following proposition describe explicitly the abstract properties of $k$-domain,
$k$-codomain, and $k$-tensor operations in $\ofs$. For more abstract treatment of these properties in terms
of the notion of a graded tensor category see \cite{Z1}.

\begin{proposition}\label{graded tensor} The $k$-tensor operation
$\ofs$ is functorial, compatible with the $k$-domain and $k$-codomain
operations, associative, and satisfy the middle exchange law.
\end{proposition}
{\it Proof.}~ In the course of the proof I will explain precisely what I mean by this statement in details.
Roughly speaking, it means that the all local morphisms form al objects of $\ofs$ into a single ordered faces structure $S$
has a natural structure of an $\o$-category $S^*$, with domains, codomains, and compositions in $S^*$ defined in terms of $k$-domain,
$k$-codomain, and $k$-tensor operations in $\ofs$.

The operations will be defined the operations on the skeleton of $\ofs$.
If $X$ and $Y$ are isomorphic ordered face structures there is a unique isomorphism between them
and in fact it is the only monotone morphism between them. We shall
identify two morphisms $f:X\ra Y$ and $f':X'\ra Y'$ in $\ofs$ iff
there are isomorphisms making the square
\begin{center} \xext=450 \yext=450
\begin{picture}(\xext,\yext)(\xoff,\yoff)
 \setsqparms[1`1`1`1;400`350]
 \putsquare(0,25)[X`Y`X'`Y';f`\cong`\cong`f']
\end{picture}
\end{center}
commutes. As these identifications are harmless we shall work in $\ofs$ recalling the identifications if needed.

To explain the functoriality of  $k$-tensor we define the category $\ofs\times_k\ofs$ as follows.
The objects of $\ofs\times_k\ofs$ are pairs of
ordered face structures $(S,S')$ such that
$\bc^{(k)}S=\bd^{(k)}S'$ and whose maps are pairs of monotone morphisms
$(f,f'):(S,S')\lra (T,T')$ such that the diagram
\begin{center} \xext=1400 \yext=600
\begin{picture}(\xext,\yext)(\xoff,\yoff)
   \setsqparms[-1`1`1`-1;700`400]
  \putsquare(0,70)[S`\bc^{(k)}S`T`\bc^{(k)}T;
 c^{(k)}_S`f`f''`\bc^{(k)}_{T}]
    \setsqparms[1`1`1`1;700`400]
   \putsquare(700,70)[\phantom{\bc^{(k)}S}`S'`\phantom{\bc^{(k)}T}`T';
 \bd^{(k)}_{S'}``f'`{\bd^{(k)}_{T'}}]
\end{picture}
\end{center}
commutes, where $f''$ is the restriction of $f'$ to $\bd^{(k)}S'$.
We have four functors
\[ \pi^0, \pi^1, \pi, \otimes_k :\ofs\times_k\ofs \lra \ofs. \]

The three first functors are defined on objects as follows
$\pi^0(S,S')=S$, $\pi^1(S,S')=S'$, $\pi(S,S')=\bc^{(k)}S$ for
$(S,S')$ in $\ofs\times_k\ofs$, and on morphisms they are defined
in the obvious way. The functor $\otimes_k$ is defined on object
and morphisms in the obvious way but we need to verify that the
local morphisms we get between local pushouts are in fact monotone.
This we leave for the reader. Moreover we have four obvious
natural transformations making the square
\begin{center} \xext=700 \yext=520
\begin{picture}(\xext,\yext)(\xoff,\yoff)
   \setsqparms[1`-1`-1`1;700`400]
  \putsquare(0,40)[\pi^0`\otimes_k`\pi`\pi^1;
 \kappa^0`\bc^{(k)}`\kappa^1`\bd^{(k)}]
\end{picture}
\end{center}
commutes, in $Nat(\ofs\times_k\ofs,\ofs)$.

By compatibility of $k$-tensor operation
 with the $k$-domain and $k$-codomain
operations, we mean that  for any ordered face structure $X$ the
squares
\begin{center} \xext=2000 \yext=580
\begin{picture}(\xext,\yext)(\xoff,\yoff)
 \setsqparms[1`-1`-1`1;700`450]
 \putsquare(0,50)[X`X`\bd^{(k)}X`\bd^{(k)}X;1_X`\bd^{(k)}_X`\bd^{(k)}_X`1_{\bd^{(k)}X}]
  \putsquare(1300,50)[X`X`\bc^{(k)}X`\bc^{(k)}X;1_X`\bc^{(k)}_X`\bc^{(k)}_X`1_{\bc^{(k)}X}]
\end{picture}
\end{center}
are $k$-tensor squares. Moreover, for $k>l$, there are isomorphism making the triangles
      \begin{center}
\xext=3000 \yext=550
\begin{picture}(\xext,\yext)(\xoff,\yoff)
 \settriparms[-1`-1`1;450]
 \putAtriangle(250,50)[X\otimes_l X'`\bd^{(k)}X\otimes_l \bd^{(k)}X'`\bd^{(k)}(X\otimes_l X');
 \bd^{(k)}_X\otimes_l \bd^{(k)}_{X'}`\bd^{(k)}_{X\otimes_l X'}`\cong]
  \putAtriangle(1900,50)[X\otimes_l X'`\bc^{(k)}X\otimes_l \bc^{(k)}X'`\bc^{(k)}(X\otimes_l X');
 \bc^{(k)}_X\otimes_l \bc^{(k)}_{X'}`\bc^{(k)}_{X\otimes_l X'}`\cong]
\end{picture}
\end{center}
commute, and for $k\leq l$, there are isomorphism making the triangles
\begin{center}
\xext=3000 \yext=550
\begin{picture}(\xext,\yext)(\xoff,\yoff)
 \settriparms[-1`-1`1;450]
 \putAtriangle(250,50)[X\otimes_l X'`\bd^{(k)}X`\bd^{(k)}(X\otimes_l X');
 \kappa^1_X\circ \bd^{(k)}_X`\bd^{(k)}_{X\otimes_l X'}`\cong]
  \putAtriangle(1900,50)[X\otimes_l X'`\bc^{(k)}X'`\bc^{(k)}(X\otimes_l X');
 \kappa^2_{X'}\circ \bc^{(k)}_X`\bc^{(k)}_{X\otimes_l X'}`\cong]
\end{picture}
\end{center}
commute.

The associativity isomorphisms come from the fact that for any
ordered face structure $X$, $Y$, $Z$ such that
$\bc^{(k)}X=\bd^{(k)}Y$ and $\bc^{(k)}Y=\bd^{(k)}Z$ both objects
\[ (X\otimes_kY)\otimes_kZ,\hskip 2cm X\otimes_k(Y\otimes_kZ) \]
are locally colimits of the diagram
\begin{center}
\xext=1600 \yext=450
\begin{picture}(\xext,\yext)(\xoff,\yoff)
\settriparms[0`-1`-1;400]
\putVtriangle(0,0)[X`Y`\bc^{(k)}X;`\bc^{(k)}_X`]
\putVtriangle(800,0)[\phantom{Y}`Z`\bc^{(k)}Y;``\bd^{(k)}_Z]
 \put(430,180){\mbox{$\bd^{(k)}_Y$}}
  \put(1050,180){\mbox{$\bc^{(k)}_Y$}}
\end{picture}
\end{center}
and the local isomorphism between them is in fact a monotone morphism.
This easily follows from Proposition \ref{loc map + are mono}.

Similarly, the interchange isomorphism between objects
\[ (X\otimes_kY)\otimes_l(Z\otimes_kT), \hskip 2cm (X\otimes_lZ)\otimes_k(Y\otimes_lT) \]
where $k<l$, is defined as the local isomorphism between two
colimits of the diagram
\begin{center}
\xext=1200 \yext=920
\begin{picture}(\xext,\yext)(\xoff,\yoff)
   \putmorphism(0,450)(1,0)[\bc^{(l)}X`\bc^{(k)}X`\bc^{(k)}_X]{700}{-1}a
   \putmorphism(700,450)(1,0)[\phantom{\bc^{(k)}X}`\bc^{(l)}Y`\bd^{(k)}_Y]{700}{1}a

   \putmorphism(0,450)(0,1)[\phantom{\bc^{(l)}X}`Z`\bd^{(l)}_Z]{450}{1}l
   \putmorphism(0,900)(0,1)[X`\phantom{\bc^{(l)}X}`\bc^{(l)}_X]{450}{-1}l

   \putmorphism(1400,450)(0,1)[\phantom{\bc^{(l)}Y}`T`\bd^{(l)}_T]{450}{1}r
   \putmorphism(1400,900)(0,1)[Y`\phantom{\bc^{(l)}X}`\bc^{(l)}_Y]{450}{-1}r
\end{picture}
\end{center}
which is in fact a monotone isomorphism.
 $~\Box$

{\em Remark.} It may seem that the $k$-tensor operation is a bit
arbitrary, as only part the order $<^\sim$ in the ordered face
structure $S\otimes_k S'$ is determined by the fact that it is a
pushout locally and that the embeddings $\kappa^S:S\ra S\otimes_k
S'$ and $\kappa^{S'}:S'\ra S\otimes_k S'$ are monotone morphisms. If
this structure determine uniquely the order $<^\sim$ in
$S\otimes_k S'$ (and hence the whole structure of $S\otimes_k S'$)
then we shall call such a $k$-tensor {\em locally
determined}\index{tensor!locally determined}. It is not hard to
see that the $k$-tensor $S\otimes_k S'$ is locally determined iff
there are no 'free' loops of dimension $k+1$ over the same
$k$-face $x\in \bc^{(k)}S$ that came from both $S$ and $S'$, i.e.
there are no $l\in S^{\lambda}_{k+1}-\delta(S^{-\lambda}_{k+2})$
and $l'\in S'^{\lambda}_{k+1}-\delta(S'^{-\lambda}_{k+2})$ such
that $\gamma(l)=\gamma(l')$ (as usually in such cases we assume
that $\bc^{(k)}S=\bd^{(k)}S'$). However if we ask for an operation
which is both pushout locally and functorial (in the sense
explained above) then the $k$-tensor operation is
the only possible one.

%
\begin{proposition}\label{tensor functorial}
The $k$-tensor operation is the unique functor $\otimes_k
:\ofs\times_k\ofs \lra \ofs$ which is a pushout functor locally,
i.e. the square
\begin{center} \xext=700 \yext=520
\begin{picture}(\xext,\yext)(\xoff,\yoff)
   \setsqparms[1`-1`-1`1;700`400]
  \putsquare(0,40)[\pi^0`\otimes_k`\pi`\pi^1;
 \kappa^0`\bc^{(k)}`\kappa^1`\bd^{(k)}]
\end{picture}
\end{center}
evaluated at any object of
$\ofs\times_k\ofs$ is a pushout in $\lfs$.
\end{proposition}

{\it Proof.}~Assume that for any $(X,X')\in\ofs\times_k\ofs$ the
square
\begin{center} \xext=700 \yext=520
\begin{picture}(\xext,\yext)(\xoff,\yoff)
   \setsqparms[1`-1`-1`1;700`400]
  \putsquare(0,40)[X`X\otimes_kX'`\bc^{(k)}X`X';
 \kappa^0_X`\bc^{(k)}_X`\kappa^1_{X'}`\bd^{(k)}_{X'}]
\end{picture}
\end{center}
is a pushout in $\lfs$.  This condition  determines the functor
$\otimes_k$ uniquely on all the objects $(Y,Y')$ of
$\ofs\times_k\ofs$ for which $k$-tensor $Y\otimes_kY'$ is locally
determined. However every object $(X,X')$ can be embedded in
$\ofs\times_k\ofs$ into a locally determined object $(Y,Y')$, i.e.
we have morphism $(f,f'):(X,X')\lra (Y,Y')$ in $\ofs\times_k\ofs$. As
the morphism $f\otimes f':X\otimes X'\lra Y\otimes Y'$ is monotone the
order $<^\sim$ in $X\otimes X'$ is uniquely determined by the
order $<^\sim$ in $Y\otimes Y'$, i.e. $\otimes_k$ is indeed the
unique functor satisfying the above requirements.  $~\Box$

Thus the above proposition says that $\otimes_k$ is the only
operation which is at the same time functorial and locally
determined as a pushout.

\section{$\o$-categories generated by local face structures}

Now we shall describe an $\o$-category $T^*$ generated by an
ordered face structure $T$, i.e. we shall describe a functor
\[ (-)^* : \ofs \lra \oC \]
however to prove some properties of $(-)^*$ it is more convenient
to describe this functor on a larger category $\lfs$, i.e. we
shall describe in fact a functor
\[ (-)^* : \lfs \lra \oC \]
We have forgetful functors, for $n\in\o$,
 \[ \bpi_n : \onfs \lra \lfs \]

For a local face structure  $T$, and for $n\in\o$, the set $T^*_n$
of $n$-cells of $T^*$ is the set of isomorphisms classes of
objects of  the comma category $\bpi_n \da T$.

For $k\leq n$, the {\em domain} and {\em codomain} operations in
$T^*$
\[ d^{(k),T^*}, c^{(k),T^*} :T^*_{n} \lra T^*_k\]
of the $k$-th {\em domain} and the $k$-th {\em codomain} are
defined by composition.  For an object $X$ of $\onfs$ and a cell
$x:X\lra T$ in $T^*_n$, we define
\[  d^{(k),T^*}(x)=\bd^{(k)}_X;x :\bd^{(k)}X\lra T  ,\hskip 10mm c^{(k),T^*}(x)=\bc^{(k)}_X;x : \bc^{(k)}X\lra
T\]
 The {\em identity} operation
\[ \bi^{(n)}:T^*_k \lra T^*_{n}\]
is an inclusion. The composition map
\[ \bm_{n,k,n} : T^*_{n}\times_{T^*_k} T^*_{n}\lra T^*_{n}\]
 is given by the $k$-tensor, i.e. for two $n$-cells $x:X\ra T$, $y:Y\ra T$
 in $T^*_n$ such that $\bc^{(k)}_X;x=\bd^{(k)}_Y;y$
then
 \[ \bm_{n,k,n}(x,y) = [x,y] \]
where $[x,y]$ is the unique map making the following diagram
\begin{center} \xext=1500 \yext=1050
\begin{picture}(\xext,\yext)(\xoff,\yoff)
 \setsqparms[1`-1`-1`1;900`600]
 \putsquare(0,50)[X`{X\otimes_kY}`\bc^{(k)}X=\bd^{(k)}Y`Y;`\bc^{(k)}_X``\bd^{(k)}_Y]
  \put(900,750){\makebox(100,100){$[x,y]$}}
  \put(500,800){\makebox(100,100){$x$}}
  \put(1200,350){\makebox(100,100){$y$}}
  \put(1500,975){\makebox(100,100){$T$}}
 \put(950,700){\vector(2,1){550}}
 \put(100,700){\vector(4,1){1400}}
 \put(1000,100){\vector(2,3){550}}
 \end{picture}
\end{center}
commutes. Note that $[x,y]$ exists and is unique since the
forgetful functor $|-|:\ofs\lra \lfs$ sends $X\otimes_kY$ to a
pushout. We often write $x \textbf{;}_ky$ for $\bm_{n,k,n}(x,y)$.

\begin{proposition}\label{Star_functor}
Let $T$ be a local face structure. Then $T^*$ is an $\o$-category.
In fact, we have a functor $(-)^* : \lfs \lra \oC.$
\end{proposition}
{\it Proof.}~ All the properties in question of $T^*$ follows more
or less in the same way from the fact that $\ofs$ is a monoidal
globular category and the the tensors in $\ofs$ are pushouts
locally.  To see how it goes we shall check the associativity of
the compositions.  So suppose we have local morphisms $x:X\ra T$,
$y:Y\ra T$, $z:Z\ra T$ such that $c^{(k)}(x)=d^{(k)}(y)$ and
$c^{(k)}(y)=d^{(k)}(z)$, i.e. the diagram
\begin{center}
\xext=1000 \yext=650
\begin{picture}(\xext,\yext)(\xoff,\yoff)
\settriparms[0`-1`-1;250] \putVtriangle(0,0)[X`Y`\bc^{(k)}X;``]
\putVtriangle(500,0)[\phantom{Y}`Z`\bc^{(k)}Y;``]
\putmorphism(500,550)(0,1)[T``y]{300}{-1}r

\putmorphism(200,380)(2,1)[`\phantom{T}`x]{100}{1}l
\putmorphism(430,570)(2,-1)[``z]{650}{-1}r
\end{picture}
\end{center}
commutes.  Hence the two compositions of these cells are
isomorphic via the canonical (local) isomorphism of pushouts
\begin{center}
\xext=1000 \yext=550
\begin{picture}(\xext,\yext)(\xoff,\yoff)
\settriparms[-1`-1`1;500]
\putAtriangle(0,20)[T`(X\otimes_kY)\otimes_kZ`X\otimes_k(Y\otimes_kZ);{[[x,y],z]}`{[x,[y,z]]}`\cong]
\end{picture}
\end{center}
But as we shown in Proposition \ref{graded tensor} these
isomorphisms are in fact monotone morphisms.  Thus the morphisms
$[[x,y],z]$ and $[x,[y,z]]$ represent the same cell in $T^*$.

The verification that $T^*$ satisfy also the remaining condition
of the definition of the $\o$-category is left for the reader. It
should be also obvious that any local morphism between local face
structures $f:S\ra T$ induces an $\o$-functor $f^*:S^*\ra T^*$ by
composition. $~\Box$

\vskip 2mm

The $k$-truncation $S_{\leq k}$ of an ordered face structure $S$
need not to be an ordered face structure, however it gives rise to
a local face structure of dimension $k$, i.e. for $k\in\o$ we have
a truncation functor
\[ tr_k : \ofs \lra \lfs_k \]
sending $(S,<^{S_k,\sim})_{k\in\o}$ to $(S,<^{\sim}_a)_{a\in
S_{>1,\leq k}}$, where $<^{\sim}_a$ is the restriction of $<^\sim$
to $\dot{\delta}(a)$, for $a\in S_{>1,\leq k}$. Here
$\lfs_k$\index{category!lfsk@$\lfs_k$} denotes the full
subcategory of $\lfs$ whose object have dimension at most $k$.
Clearly, we have a commuting square
\begin{center} \xext=700 \yext=600
\begin{picture}(\xext,\yext)(\xoff,\yoff)
 \setsqparms[1`1`1`1;700`500]
 \putsquare(0,50)[\ofs`\lfs_k`\oC`\kC;tr_k`(-)^*`(-)^*`tr_k]
 \end{picture}
\end{center}
Thus we have a functor \[(-)^*_{\leq k} : \ofs \lra \kC \] which
is defined as either of the above compositions.
$\kC$\index{category!kC@$\kC$} is the category of $k$-categories.

\section{Principal and Normal ordered face structures}

We recall few notions form section \ref{face stuctures}. Let $N$
be an ordered face structure.  $N$ is $k$-normal iff $dim(N)\leq
k$ and $size(N)_l=1$, for $l<k$. $N$ is $k$-principal iff
$size(N)_l=1$, for $l\leq k$. $N$ is principal iff
$size(N)_l\leq1$, for $l\leq \o$. $N$ is principal of dimension
$k$ iff $N$ is principal and $dim(N)=k$.

Notation for a $k$-normal $N$:
$\{\bp^N_l\}=\{\bp_l\}=N_l-\delta(N_{l+1})$, for $l<k$.

\begin{lemma}\label{principal ofs} Let $P$, $Q$, $N$, $T$ be an ordered face
  structures, $k\in\o$, $P$, $Q$ principal, $N$ $k$-normal.
\begin{enumerate}
 \item If the map $f:P\ra T$ is local then it is a monotone morphism.
 \item If the map $f:N\ra T$ is local then $f$ is monotone iff $f_k$ preserves $<^\sim$.
 \item If $dim(P)=dim(Q)$ and the maps $f:P\ra T$ and $g:Q\ra T$ are local such that $f(\bp^P)= g(\bp^Q)$
 then there is a unique monotone isomorphism $h:P\ra Q$ making the triangle
\begin{center}
\xext=600 \yext=350
\begin{picture}(\xext,\yext)(\xoff,\yoff)
\settriparms[1`1`1;300] \putVtriangle(0,0)[P`Q`T;h`f`g]
\end{picture}
\end{center}
 commutes, where $\bp^P$, $\bp^Q$ are the unique faces of dimension
 $k$ in $P$ and in $Q$, respectively.
 \item If $dim(P)=n>dim(Q)$ then any monotone morphism $x:Q\ra P$
factorizes either via $\bd_P: \bd P \ra P$ or $\bc_P: \bc P \ra
P$.
\end{enumerate}
\end{lemma}

{\it Proof.}~ 1. and 2. follows immediately from Lemma
 \ref{loc map + are mono}. 4. follows easily from 3.

Ad 3. Let $dim(P)=dim(Q)=k$. By 1. we need to construct a local
isomorphism only. The argument is by induction on $k$. For $k=0$
the claim is obvious. For $k>0$, we have by induction hypothesis
the local morphism $h':\bc^{(k-1)}P\lra \bc^{(k-1)}Q$. Then we note
that the bijections  $f: \delta(\bp^P) \lra \delta(f(\bp^P))$, $g:
\delta(\bp^Q) \lra \delta(g(\bp^Q))$ preserves order. As
$f(\bp^P)=g(\bp^Q)$ we get easily the local morphism $h:P\lra Q$.
 $~\Box$

{\em Example.} Note that in Lemma \ref{principal ofs}.4 it is
essential that $Q$ is principal and not any ordered face
structure. In the example as below
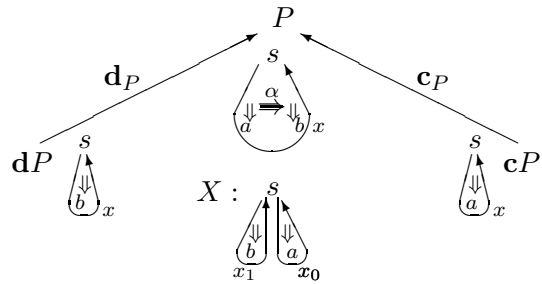
\begin{center} \xext=1500 \yext=1000
\begin{picture}(\xext,\yext)(\xoff,\yoff)
  \put(440,300){$X: $}
  \put(790,0){$^{x_0}$}
  \put(630,130){\oval(100,100)[b]}
  \put(580,130){\line(1,2){85}}
  \put(680,130){\vector(0,1){180}}
  \put(620,130){$^\Da$}
    \put(610,65){$^{b}$}
     \put(565,0){$^{x_1}$}
     \put(680,320){$s$}

  \put(770,130){\oval(100,100)[b]}
   \put(720,130){\line(0,1){180}}
  \put(820,130){\vector(-1,2){85}}
   \put(730,130){$^\Da$}
   \put(750,65){$^{a}$}
   \put(790,0){$^{x_0}$}

   \put(-200,400){$\bd P$}
  \put(50,300){\oval(100,100)[b]}
  \put(0,300){\line(1,4){40}}
  \put(100,300){\vector(-1,4){40}}
  \put(30,300){$^\Da$}
    \put(20,235){$^{b}$}
     \put(115,210){$^{x}$}
     \put(30,480){$s$}

  \put(1500,400){$\bc P$}
  \put(1400,300){\oval(100,100)[b]}
  \put(1350,300){\line(1,4){40}}
  \put(1450,300){\vector(-1,4){40}}
  \put(1380,300){$^\Da$}
    \put(1370,235){$^{a}$}
     \put(1465,210){$^{x}$}
     \put(1380,480){$s$}

 \put(700,900){$P$}
  \put(700,600){\oval(260,260)[b]}
  \put(570,600){\line(1,2){85}}
  \put(830,600){\vector(-1,2){85}}
  \put(600,550){$^\Da$}
  \put(750,550){$^\Da$}
  \put(650,600){$\Rightarrow$}
  \put(660,625){\line(1,0){75}}
  \put(675,620){$^{\alpha}$}
    \put(590,500){$^{a}$}
     \put(780,500){$^{b}$}
     \put(840,500){$^{x}$}
     \put(680,780){$s$}

  \put(120,700){$\bd_P$}
  \put(1200,700){$\bc_P$}
   \put(-100,500){\vector(2,1){750}}
    \put(1550,500){\vector(-2,1){750}}
\end{picture}
\end{center}
with morphism $f:X\ra P$ sending $x_i$ to $x$ and other cells to
the same cell we clearly cannot factor $f$ via neither $\bd_P$ nor
$\bc_P$.

\begin{lemma}\label{type of indet} Let $T$ be an ordered face
  structure, $l,k\in\o$, $l<k$, and  $\alpha\in T_k$. We have
\begin{enumerate}
  \item $\{\alpha\}$ is a convex set and  $[\alpha ]$ is a principal ordered face structure,
  \item $\delta(\alpha)$ is a convex set and  $[\delta(\alpha) ]$ is a $(k-1)$-normal ordered face structure,
 \item  Moreover, if $k>0$, then
\[ \bc^{(l)} [\alpha] = [\gamma^{(l)}(\alpha)] \hskip 15mm  \bd^{(l)} [\alpha] = [\delta^{(l)}(\alpha)]. \]
\end{enumerate}
\end{lemma}

{\it Proof.}~  We shall prove 1. The rest is left as an
exercise.

The proof goes by induction of the $dim(\alpha)=k$. If $k=0$ then
the thesis is obvious. So assume that $k>0$, the thesis holds for
$\gamma(\alpha)$ and we shall prove it for $\alpha$. If $\alpha\in
T^\lambda\cup T^\varepsilon$ then $<\alpha>$ has faces as follows
\[ \begin{array}{|c|c|}
    dim & {\rm faces}    \\ \hline
    k & \alpha   \\ \hline
     k-1  & \gamma(\alpha) \\ \hline
    k-2 &\dot{\delta}\gamma(\alpha)\cup\gamma\gamma(\alpha)
  \end{array}\]
and hence the thesis is obvious. So assume that $\alpha\in
T^{-\lambda\varepsilon}$.  Then $<\alpha>$ has faces as follows
\[ \begin{array}{|c|c|}
    dim & {\rm faces}    \\ \hline
    k & \alpha   \\ \hline
     k-1  & \delta(\alpha)\cup\gamma(\alpha) \\ \hline
    k-2 &\delta\dot{\delta}^{-\lambda}(\alpha)\cup\gamma\gamma(\alpha)   \\ \hline
    k-3 & \delta\dot{\delta}^{-\lambda}\gamma(\alpha)\cup\gamma\gamma\gamma(\alpha)
  \end{array}\]
For the faces of dimension $k$ and $k-1$ the thesis is obvious,
for dimensions $k-3$ and lower the thesis holds by inductive
assumption on $\gamma(\alpha)$. We need to check that $<^{S,+}$
and  $<^{<\alpha>,+}$ agree on $<\alpha>_{k-2}$. First note that
by Lemma \ref{fact2}.3, if
$x\in\delta\dot{\delta}^{-\lambda}(\alpha)$ then
$x<^{<\alpha>,+}\gamma\gamma(\alpha)$. So assume that
$x\in\delta(a)$ and $y\in\delta(b)$, $a,b\in\dot{\delta}(\alpha)$,
$x<^{S,+}y$. Let $x,a_1,\ldots, a_n,y$ be a flat upper
$T-\gamma(T^{-\lambda})$-path from $x$ to $y$ and $x,b_1,\ldots,
b_l,\gamma\gamma(\alpha)$ be a flat upper
$\dot{\delta}(\alpha)$-path from $x$ to $\gamma\gamma(\alpha)$. If
$y=\gamma(b_{l_0})$ for some $l_0\leq l$, we are done.  So assume
contrary. Then by Path Lemma, there is $l_1\leq l$ such that
$a_n<^+b_{l_1}$, $\gamma(a_n)=y\neq\gamma(b_{l_1})$.  Thus we have
a flat upper path $a_n,\beta_1,\ldots,\beta_r,b_{l_1}$ and, as
$y\neq\gamma(b_{l_1})$, there is $r_0\leq r$ such that
$y\in\iota(\beta_{r_0})$. Hence
$\delta(b)\cap\iota(\beta_{r_0})\neq\emptyset$ and by Lemma
\ref{fact9}.1, $b<^+\gamma(b_{r_0})\leq b_{l_1}$.  But $b,
b_{l_1}\in\delta(\alpha)$ and we get a contradiction with
discreetness. $~\Box$

Let $k\in\o$, $N$ be a $k$-normal ordered face structure. We
define a $(k+1)$-hypergraph $N^\bullet$, that contains two
additional faces: $\bp_{k+1}^{N^\bullet}$ of dimension $k+1$, and
$\bp_{k}^{N^\bullet}$ of dimension $k$. We shall drop superscripts
if it does not lead to confusions.  We also put
 \[ \gamma(\bp_{k+1})=\bp_k, \;\;\; \gamma(\bp_k)=\bp_{k-1},  \]

\[\delta(\bp_{k+1}) = \left\{ \begin{array}{ll}
            N_k &  \mbox{if $N_k\neq\emptyset$,}  \\
            1_{\bp_{k-1}} &\mbox{otherwise.}
                                    \end{array}
                \right.\;\;\;
 \delta(\bp_{k}) = \left\{ \begin{array}{ll}
            \delta(N_k)-\gamma(N_k)&  \mbox{if $N_k\neq\emptyset$,}  \\
            1_{\bp_{k-2}} &\mbox{otherwise.}
                                    \end{array}
                \right.
                 \]
Clearly,  $ \gamma(\bp_k)$  and $\delta(\bp_{k})$ are defined only
if $k>0$.

As $N$ is $k$-normal, $N_{k+1}=\emptyset$, so $N_k$ cannot contain
loops. Thus, if $N_k\neq\emptyset$ then
$\delta(N_k)-\gamma(N_k)\neq\emptyset$ and $\delta(\bp_k)$ is well
defined. This determines $N^\bullet$ uniquely. $N^\bullet$ is
called {\em the principal extension of}\index{principal extension}
$N$.

{\em Examples.} Here are some examples of $1$-normal ordered face
structures $N$ and their principal extensions $N^\bullet$:
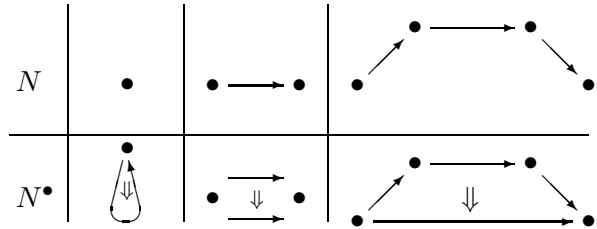
\begin{center} \xext=2000 \yext=750
\begin{picture}(\xext,\yext)(\xoff,\yoff)
  \put(0,300){\line(1,0){2000}}
  \put(200,0){\line(0,1){750}}
  \put(600,0){\line(0,1){750}}
  \put(1100,0){\line(0,1){750}}
 \put(20,450){$N$}
 \put(20,50){$N^\bullet$}

  \put(380,450){$\bullet$}
  \put(400,50){\oval(100,100)[b]}
  \put(350,50){\line(1,4){40}}
  \put(450,50){\vector(-1,4){40}}
  \put(380,50){$^\Da$}
   \put(380,230){$\bullet$}

 \putmorphism(700,470)(1,0)[\bullet`\bullet`]{300}{1}a
 \putmorphism(700,150)(1,0)[\phantom{\bullet}`\phantom{\bullet}`]{300}{1}a
 \putmorphism(700,80)(1,0)[\bullet`\bullet`]{300}{0}a
 \putmorphism(700,10)(1,0)[\phantom{\bullet}`\phantom{\bullet}`]{300}{1}a
 \put(830,20){$^\Da$}

 \putmorphism(1200,470)(1,0)[\bullet`\bullet`]{800}{0}a
 \putmorphism(1400,670)(1,0)[\bullet`\bullet`]{400}{1}a
 \put(1240,510){\vector(1,1){120}}
  \put(1840,630){\vector(1,-1){120}}
 \putmorphism(1200,0)(1,0)[\bullet`\bullet`]{800}{1}a
 \putmorphism(1400,200)(1,0)[\bullet`\bullet`]{400}{1}a
 \put(1240,40){\vector(1,1){120}}
 \put(1840,160){\vector(1,-1){120}}
 \put(1560,50){$\Da$}
\end{picture}
\end{center}
and some examples of $2$-normal ordered face structures $N$ and
their principal extensions $N^\bullet$:
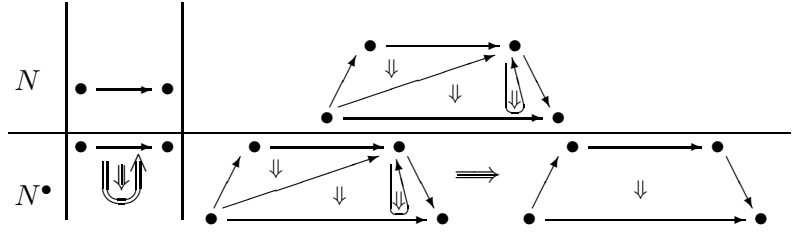
\begin{center} \xext=2700 \yext=750
\begin{picture}(\xext,\yext)(\xoff,\yoff)
  \put(0,300){\line(1,0){2700}}
  \put(200,0){\line(0,1){750}}
  \put(600,0){\line(0,1){750}}
 \put(20,450){$N$}
 \put(20,50){$N^\bullet$}

 \putmorphism(250,450)(1,0)[\bullet`\bullet`]{300}{1}a
 \putmorphism(250,250)(1,0)[\bullet`\bullet`]{300}{1}a
   \put(390,120){\oval(100,100)[b]}
  \put(390,120){\oval(130,130)[b]}
  \put(340,120){\line(0,1){80}}
  \put(325,120){\line(0,1){80}}
  \put(440,120){\line(0,1){80}}
  \put(455,120){\line(0,1){80}}
  \put(415,170){$\wedge$}
  \put(360,120){$\Da$}
  \put(390,120){\line(0,1){70}}

 \putmorphism(1100,350)(1,0)[\bullet`\bullet`]{800}{1}a
 \putmorphism(1250,600)(1,0)[\bullet`\bullet`]{500}{1}a
 \put(1110,390){\vector(1,2){90}}
 \put(1780,570){\vector(1,-2){90}}
 \put(1140,390){\vector(3,1){540}}
 \put(1520,380){$^\Da$}
 \put(1300,470){$^\Da$}

  \put(1750,400){\oval(60,60)[b]}
  \put(1720,400){\line(0,1){140}}
  \put(1780,400){\vector(-1,4){40}}
  \put(1725,360){$^\Da$}
 \putmorphism(700,0)(1,0)[\bullet`\bullet`]{800}{1}a
 \putmorphism(850,250)(1,0)[\bullet`\bullet`]{500}{1}a
 \put(710,40){\vector(1,2){90}}
 \put(1380,220){\vector(1,-2){90}}
 \put(740,40){\vector(3,1){540}}
 \put(1120,30){$^\Da$}
 \put(900,120){$^\Da$}

  \put(1350,50){\oval(60,60)[b]}
  \put(1320,50){\line(0,1){140}}
  \put(1380,50){\vector(-1,4){40}}
  \put(1325,10){$^\Da$}

 \putmorphism(1800,0)(1,0)[\bullet`\bullet`]{800}{1}a
 \putmorphism(1950,250)(1,0)[\bullet`\bullet`]{500}{1}a
 \put(1810,40){\vector(1,2){90}}
 \put(2480,220){\vector(1,-2){90}}
 \put(2160,50){$^\Da$}

\put(1540,130){$\Longrightarrow$}
  \put(1550,155){\line(1,0){135}}
\end{picture}
\end{center}
Clearly, $\bullet$ the 'empty' $1$-normal, and
$\bullet\longrightarrow\bullet$ is 'empty' $2$-normal ordered face
structure.

\begin{proposition}\label{bullet}
Let $N$ be a $k$-normal ordered face structure. Then
\begin{enumerate}
 \item $N^\bullet$ is a principal ordered face structure of dimension $k+1$.
 \item We have $\bd (N^\bullet)  \cong N $, $\bc (N^\bullet)  \cong (\bd N)^\bullet$.
 \item If $N$ is a principal, then $N\cong(\bd N)^\bullet$.
\end{enumerate}
\end{proposition}
{\it Proof.}~  Exercise. $~\Box$

\section{Decomposition of ordered face structures}

As positive face structures are easier and we understand well
their decompositions  we define decomposition of ordered face
structures via positive ones.  This will simplify the proof of
properties of the decompositions, as they will be easy
consequences of the analogous properties of decompositions of
positive face structures. However to get a better insight how the
ordered face structures are decomposed we shall characterize the
decompositions using convex subsets and stretching empty loops. We
decompose along an $\cI$-cut rather than a face.

NB. We write $\check{a}$ instead of $(a,L,U)$ if we don't need to
specify explicitly which cut over $a$ we consider.

The $k$-{\em decomposition}\index{decomposition} of $X$ is any
presentation of $X$ as a $k$-tensor $X=X_1\otimes_k X_2$ of two
other ordered face structures. $X_1$ is {\em the lower part of the
decomposition}\index{decomposition!lower part of -} and $X_2$ is
{\em the upper part of the
decomposition}\index{decomposition!upper part of -}.  The
$k$-decomposition of $X=X_1\otimes_k X_2$ is said to be {\em
proper}\index{decomposition!proper} iff
$size(X_1),size(X_2)<size(X)$.

Let $X$ be an ordered faces structure, $\check{a}\in X^\dag$,
$\cJ$ the kernel of the standard positive cover $q:X^\dag\ra X$.
We define the decomposition of $X$ along $\check{a}$ as the bottom
square of the following cube
\begin{center}
\xext=1300 \yext=1120
\begin{picture}(\xext,\yext)(\xoff,\yoff)
\setsqparms[1`1`1`1;900`600]
 \putsquare(0,50)[\bc^{(k)}(X^{\dag \da\check{a}})`X^{\dag \ua\check{a}}`
 \bc^{(k)}(X^{\da\check{a}})`X^{\ua\check{a}};
 \bd^{(k)}`p`p`\bd^{(k)}]

 \setsqparms[1`0`1`0;900`600]
 \putsquare(380,450)[X^{\dag \da\check{a}}`X^\dag`X^{\da\check{a}}`X;
 \kappa``p`]

 \putmorphism(1050,200)(1,1)[``\kappa]{110}{1}r
 \putmorphism(1050,800)(1,1)[``\kappa]{110}{1}l

 \putmorphism(150,200)(1,1)[``]{110}{1}l
 \putmorphism(150,800)(1,1)[``\bc^{(k)}]{110}{1}l
\put(60,250){$\bc^{(k)}$}

\putmorphism(910,450)(1,0)[``]{320}{1}l
 \put(500,450){\line(1,0){360}}
 \put(630,360){$\kappa$}

\putmorphism(360,690)(0,-1)[``]{220}{1}l
 \put(360,980){\line(0,-1){280}}
 \put(240,730){$p$}
\end{picture}
\end{center}
where the top square is the decomposition of the positive face
structure $X^\dag$ along $\check{a}$, and the bottom square
\begin{center}
\xext=1000 \yext=600
\begin{picture}(\xext,\yext)(\xoff,\yoff)
 \setsqparms[1`-1`-1`1;1000`500]
 \putsquare(0,40)[X^{\da \check{a}}`X`\bc^{(k)}(X^{\da \check{a}})`X^{\ua \check{a}};
 \kappa_X^{\da  \check{a}}`\bc^{(k)}`\kappa_X^{\ua \check{a}}`\bd^{(k)}]
\end{picture}
\end{center}
is obtained from the top square by dividing by $\cJ$.

\begin{lemma}\label{decomp_square lemma}
We note for the record
\[ \bd^{(k)}(X^{\da \check{a}})=\bd^{(k)}(X), \hskip 10mm \bc^{(k)}(X^{\ua
\check{a}})=\bc^{(k)}(X), \]
\[ \bc^{(k)}(X^{\da \check{a}})=\bd^{(k)}(X^{\ua
\check{a}}), \hskip 10mm X^{\da \check{a}}\otimes_k X^{\ua
\check{a}}= X.\] \end{lemma}

{\it Proof.}~Exercise. $~\Box$

\begin{lemma} \label{Sd}
Let $S$, $T$ be ordered face structures, $k\in\o$, and $\check{a}=(a,L,U)\in
(T^\dag_k -\iota(T^\dag_{k+2}))$. Then
\begin{enumerate}
  \item  $\check{a}\in Sd(T)$ iff there are $\alpha,\beta\in T_{k+1}$ such that
  $(\gamma^{(k)}(\alpha),-,\ua \gamma^{(k+1)}(\alpha))\leq^+\check{a}$ and
  $(\gamma^{(k)}(\beta),-,\ua \gamma^{(k+1)}(\beta))\not\leq^+\check{a}$.
  \item $Sd(T)=Sd(T^\dag)$.
  \item $size(T)=size(T^\dag)$.
  \item if $\bc^{(k)}(S)=\bd^{(k)}(T)$ then, for $l\in\o$,
\[  size(S\otimes_kT)_l = \left\{ \begin{array}{ll}
          size(S)_l+ size(T)_l  & \mbox{ if $l>k$,}  \\
       size(T)_l & \mbox{ if $l\leq k$.}
                                    \end{array}
                \right. \]
  \item $size(T)_k\geq 1$ iff $k\leq dim(T)$.
  \item $Sd(T)_k\neq\emptyset$ iff $size(T)_{k+1}\geq 2$.
  \item $T$ is principal iff $Sd(T)$ is empty.
\end{enumerate}
\end{lemma}
{\it Proof.}~ Easy.  $~\Box$

Before we shall establish the important properties of this
decomposition we shall show another way of constructing this
decomposition.
Let $Y$ be a convex subset of an ordered face structure $X$. We
define two subhypergraphs $Y^{\Da \check{a}}$ and $Y^{\Ua
\check{a}}$ of $X$:

\[ Y^{\Da \check{a}}_l = \left\{ \begin{array}{ll}
            \{ \alpha\in Y_l : (\gamma^{(k)}(\alpha),-,\ua \gamma^{(k+1)}(\alpha) )\leq ^+(a,L,U) \}
                                                                                        &  \mbox{for $l> k$,}  \\
            \{ b\in Y_k : b\leq^+a \;{\rm or}\;\; b\not\in\gamma(Y_{k+1}^{-\lambda}) \} & \mbox{for $l=k$} \\
            X_l &\mbox{for $l<k$.}
                                    \end{array}
                \right. \]

\[ Y^{\Ua \check{a}}_l = \left\{ \begin{array}{ll}
            \{ \alpha\in Y_l : (\gamma^{(k)}(\alpha),-,\ua \gamma^{(k+1)}(\alpha) )\not\leq^+(a,L,U)\}
                                                                                     &  \mbox{for $l> k$,}  \\
            \{ b\in Y_k : b\not<^+a \} &\mbox{for $l=k$} \\
            Y_{k-1} - \iota(Y^{\Da \check{a}}_{k+1}) &\mbox{for $l=k-1$} \\
            Y_l &\mbox{for $l<k-1$.}
                                    \end{array}
                \right. \]

\begin{lemma}\label{decomp}
With the notation as above $Y^{\Da \check{a}}$ and $Y^{\Ua
\check{a}}$ are convex subhypergraphs of $X$, $c^{(k)}(Y^{\Da
\check{a}})=d^{(k)}(Y^{\Ua \check{a}})$. Moreover $\cE^{Y^{\Da
\check{a}}}=\emptyset$ and
\[ \cE^{Y^{\Ua \check{a}}} = \left\{ \begin{array}{ll}
            \{ a \} &  \mbox{if $a$ is a loop,}\\
            \emptyset &\mbox{otherwise.}
                                    \end{array}
                \right. \]
\end{lemma}
{\it Proof.}~ Easy. $~\Box$

\begin{lemma}
With the notation as above there are monotone isomorphisms $h_\da$
and $h_\ua$ making the triangles
\begin{center} \xext=800 \yext=450
\begin{picture}(\xext,\yext)(\xoff,\yoff)
\settriparms[1`-1`-1;400]
 \putptriangle(0,0)[X^{\da \check{a}}`X`{[X^{\Da \check{a}}]}; \kappa`h_\da`\nu]
 \settriparms[1`-1`-1;400]
 \putqtriangle(400,0)[\phantom{X}`X^{\da \check{a}}`{[X^{\Ua \check{a}}]}; \kappa`\nu`h_\ua]
 \end{picture}
\end{center}
commute. $\kappa_X^{\ua \check{a}}$
\end{lemma}
{\it Proof.}~By Lemma \ref{subfs} it is enough to show that the
image of the monotone morphism $\kappa_X^{\ua \check{a}}: X^{\ua
\check{a}}\lra X$ is $X^{\Ua \check{a}}$ and the image of the
monotone morphism $\kappa_X^{\da \check{a}}: X^{\da \check{a}}\lra X$
is $X^{\Da \check{a}}$. The remaining details are left for the
reader. $~\Box$

Note that, by the above Lemma \ $X^{\Da \check{a}}$ is isomorphic
to $X^{\da \check{a}}$ and if $a$ is not a loop in $X$, $X^{\Ua
\check{a}}$ is isomorphic to $X^{\ua \check{a}}$. However the
ordered face structure $X^{\ua \check{a}}$ is not that complicated
even if $a$ is a loop. We shall describe it now. Thus $\cE^{X^{\Ua
\check{a}}} = \{ a\}$. In this case, up to isomorphism, the
underlying hypergraph of $X^{\ua \check{a}}$ can be describe as
follows.
\[ X^{\ua \check{a}}_l = \left\{ \begin{array}{ll}
           X^{\Ua \check{a}}_l &  \mbox{if $l\neq k-1$,}\\
           (X^{\Ua \check{a}}_{k-1}-\{ \gamma(a) \})\cup \{ \gamma(a)^-, \gamma(a)^+ \} &\mbox{if $l= k-1$.}
                                    \end{array}
                \right. \]
$\gamma^{X^{\ua \check{a}}}$ and $\delta^{X^{\ua \check{a}}}$ are
as in $X^{\Ua \check{a}}$ (and $X$) except for $\gamma^{X^{\ua
\check{a}}}_{k-1}$ and $\delta^{X^{\ua \check{a}}}_{k-1}$. For
$c\in X^{\ua \check{a}}_{k}$ we put ($\gamma$ and $\delta$ stands
for $\gamma^X$ and $\delta^X$, respectively)
\[ \gamma^{X^{\ua \check{a}}}(c) = \left\{ \begin{array}{ll}
            \gamma(c) &  \mbox{if $\gamma(c)\neq \gamma(a)$,}\\
           \gamma(a)^- &\mbox{if $\gamma(c)= \gamma(a)$ and $c<^\sim a$,}\\
           \gamma(a)^+ &\mbox{otherwise.}
                                    \end{array}
                \right. \]
\[ \delta^{X^{\ua \check{a}}}(c) = \left\{ \begin{array}{ll}
            \delta(c) &  \mbox{if $\gamma(b)\not\in \delta(c)$,}\\
           (\delta(c)-\{ \gamma(a)\})\cup \{ \gamma(a)^+ \} &\mbox{if $\gamma(a)\in\delta(c)$ and $a<^\sim c$,}\\
           (\delta(c)-\{ \gamma(a)\})\cup \{ \gamma(a)^- \} &\mbox{otherwise.}
                                    \end{array}
                \right. \]
The order $<^\sim$ in $X^{\da \check{a}}$ and $X^{\ua \check{a}}$
is uniquely determined by the fact that it is reflected from $X$
via $ \kappa^{\da \check{a}}_X$ and $\kappa^{\ua \check{a}}_X$.

{\em Examples.} For the ordered face structure $T$ as below
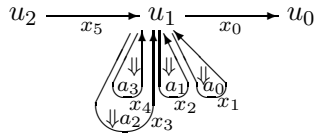
\begin{center} \xext=1320 \yext=700
\begin{picture}(\xext,\yext)(\xoff,\yoff)
\putmorphism(0,450)(1,0)[u_2`u_1`_{x_5}]{480}{1}b
\putmorphism(480,450)(1,0)[\phantom{u_1}`u_0`_{x_0}]{480}{1}b
 \put(-200,600){$T$}
 \put(350,150){\oval(200,200)[b]}
 \put(285,40){$^{\Da a_2}$}
 \put(250,150){\line(1,2){120}}
 \put(450,150){\vector(0,1){250}}
  \put(440,30){$^{x_3}$}

  \put(360,220){\oval(100,100)[b]}
  \put(310,220){\line(1,2){85}}
  \put(410,220){\vector(0,1){180}}
  \put(360,220){$^\Da$}
    \put(325,155){$^{a_3}$}
     \put(365,90){$^{x_4}$}
    \put(520,220){\oval(100,100)[b]}
  \put(470,220){\line(0,1){180}}
  \put(570,220){\vector(-1,2){85}}
  \put(480,220){$^\Da$}
    \put(490,155){$^{a_1}$}
     \put(520,90){$^{x_2}$}
     \put(650,220){\oval(100,100)[b]}
  \put(600,220){\line(-1,2){85}}
  \put(700,220){\vector(-1,1){180}}
  \put(600,190){$^\Da$}
    \put(620,155){$^{a_0}$}
     \put(670,100){$^{x_1}$}
\end{picture}
\end{center}
and a cut $\check{u_1}=(u_1,\{ a_3\},\{a_2,a_1\})$ we have the
following decomposition
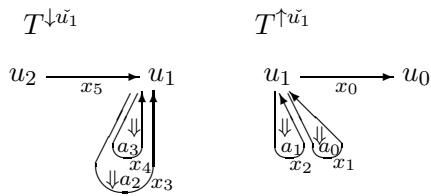
\begin{center} \xext=1320 \yext=700
\begin{picture}(\xext,\yext)(\xoff,\yoff)
 \put(0,600){$T^{\da \check{u_1}}$}
 \putmorphism(0,450)(1,0)[u_2`u_1`_{x_5}]{480}{1}b
  \put(350,150){\oval(200,200)[b]}
  \put(285,40){$^{\Da a_2}$}
  \put(250,150){\line(1,2){120}}
  \put(450,150){\vector(0,1){250}}
  \put(440,30){$^{x_3}$}

  \put(360,220){\oval(100,100)[b]}
  \put(310,220){\line(1,2){85}}
  \put(410,220){\vector(0,1){180}}
  \put(360,220){$^\Da$}
  \put(325,155){$^{a_3}$}
  \put(365,90){$^{x_4}$}
 \put(800,600){$T^{\ua \check{u_1}}$}
 \putmorphism(880,450)(1,0)[u_1`u_0`_{x_0}]{480}{1}b

    \put(920,220){\oval(100,100)[b]}
  \put(870,220){\line(0,1){180}}
  \put(970,220){\vector(-1,2){85}}
  \put(880,220){$^\Da$}
    \put(890,155){$^{a_1}$}
     \put(920,90){$^{x_2}$}
     \put(1050,220){\oval(100,100)[b]}
  \put(1000,220){\line(-1,2){85}}
  \put(1100,220){\vector(-1,1){180}}
  \put(1000,190){$^\Da$}
    \put(1020,155){$^{a_0}$}
     \put(1070,100){$^{x_1}$}
\end{picture}
\end{center}
and for the cut $\check{x_4}=(x_4,\emptyset,\emptyset)$ we have
the following decomposition
\begin{center} \xext=3000 \yext=650
\begin{picture}(\xext,\yext)(\xoff,\yoff)
 \put(0,550){$T^{\da \check{x_4}}$}
 \putmorphism(0,380)(1,0)[u_2`u_1`_{x_5}]{300}{1}b
 \putmorphism(300,380)(1,0)[\phantom{u_1}`u_0`_{x_0}]{300}{1}b

  \put(260,130){\oval(100,100)[b]}
  \put(210,130){\line(1,4){40}}
  \put(310,130){\vector(-1,4){40}}
  \put(240,130){$^\Da$}
  \put(225,65){$^{a_3}$}
  \put(265,0){$^{x_4}$}
  \put(800,550){$T^{\ua \check{x_4}}$}
  \putmorphism(900,380)(1,0)[u_2`(u_1,\emptyset,\{x_4 \})`_{x_0}]{600}{1}a
  \putmorphism(1500,430)(1,0)[\phantom{(u_1,\emptyset,\{x_4 \})}`\phantom{(u_1,\{x_4 \},\emptyset)}`_{x_4}]{900}{1}a
  \putmorphism(1500,380)(1,0)[\phantom{(u_1,\emptyset,\{x_4 \})}`(u_1,\{x_4 \},\emptyset)`]{900}{0}a
  \putmorphism(1500,330)(1,0)[\phantom{(u_1,\emptyset,\{x_4 \})}`\phantom{(u_1,\{x_3 \},\emptyset)}`_{x_3}]{900}{1}b
  \put(1900,320){$^{\Da a_2}$}
  \putmorphism(2400,380)(1,0)[\phantom{(u_1,\{x_4 \},\emptyset)}`x_0`u_0]{600}{1}a

 \put(2330,130){\oval(100,100)[b]}
  \put(2280,130){\line(1,2){85}}
  \put(2380,130){\vector(0,1){180}}
  \put(2320,130){$^\Da$}
  \put(2300,65){$^{a_1}$}
  \put(2265,0){$^{x_2}$}

  \put(2470,130){\oval(100,100)[b]}
  \put(2420,130){\line(0,1){180}}
  \put(2520,130){\vector(-1,2){85}}
  \put(2430,130){$^\Da$}
  \put(2440,65){$^{a_0}$}
  \put(2490,0){$^{x_1}$}
\end{picture}
\end{center}

The following Lemma establishes some properties of the double
decompositions. The double decomposition is meant in the sense of
convex set decomposition, i.e. when we write $X^{\da \check{x} \ua
\check{a}}$ we mean $[X^{\Da \check{x} \Ua \check{a}}]$.

\begin{lemma}\label{decomp1new}
Let $X$ be an ordered face structure,
$\check{a}=(a,L,U),\check{x}=(x,L',U')\in (X^\dag
-\iota(X^\dag))$, $k=dim(x)<dim(a)=m$.
\begin{enumerate}
  \item We have the following equations of ordered face structures:
\[ X^{\da \check{x} \da \check{a}}= X^{\da \check{a} \da \check{x}}, \;\;\;
   X^{\da \check{x} \ua \check{a}}= X^{\ua \check{a} \da \check{x}},
   \;\;\;
   X^{\ua \check{x} \da \check{a}}= X^{\da \check{a} \ua
\check{x}}, \;\;\;  X^{\ua \check{x} \ua \check{a}}= X^{\ua
\check{a} \ua \check{x}}, \]
   i.e. $x$-decompositions and $a$-decompositions commute.
   \item If $\check{x}\in Sd(X)$ then
   $\check{x}\in Sd_{\kappa^{\da \check{a}}}(X^{\da \check{a}})\cap Sd_{\kappa^{\ua \check{a}}}(X^{\ua \check{a}})$.
   \item Moreover, we have the following equations concerning
   domains and codomains
   \[ \bc^{(k)}(X^{\da \check{x} \da \check{a}})= \bc^{(k)}(X^{\da \check{x} \ua \check{a}})
    = \bd^{(k)}(X^{\ua \check{x} \da \check{a}})= \bd^{(k)}(X^{\ua \check{x} \ua \check{a}})\]
   \[ \bc^{(m)}(X^{\da \check{x} \da \check{a}})=  \bd^{(m)}(X^{\da \check{x} \ua \check{a}}), \;\;\;
     \bc^{(m)}(X^{\ua \check{x} \da \check{a}})= \bd^{(m)}(X^{\ua \check{x} \ua \check{a}}).\]
   \item Finally, we have the following equations concerning
   compositions
   \[ X^{\da \check{x} \da \check{a}}\otimes_m X^{\da \check{x} \ua \check{a}}= X^{\da \check{x}}, \;\;\;
    X^{\ua \check{x} \da \check{a}}\otimes_m X^{\ua \check{x} \ua \check{a}}=X^{\ua \check{x}},\]
    \[ X^{\da \check{x} \da \check{a}}\otimes_k X^{\ua \check{x} \da \check{a}}= X^{\da \check{a}}, \;\;\;
     X^{\da \check{x} \ua \check{a}}\otimes_k X^{\ua \check{x} \ua \check{a}}=X^{\ua \check{a}}.\]
\end{enumerate}
\end{lemma}

{\it Proof.}~ We need to verify the above equations for arrows
$\Da$ and $\Ua$ instead of $\da$ and $\ua$. $~\Box$ 

\begin{lemma}
\label{decomp2new} Let $T$ be ordered face structure,
 $X$ convex subhypergraph of $T$, and $a,b\in X$, $\check{a}=(a,L,U),\check{b}=(b,L',U')\in
T^\dag -\iota(T^\dag)$, $dim(a)=dim(b)=m$.
\begin{enumerate}
  \item We have the following equations of ordered face structures:
\[ X^{\da \check{a} \da \check{b}}= X^{\da \check{b} \da \check{a}}, \;\;\;
   X^{\ua \check{a} \ua \check{b}}= X^{\ua \check{b} \ua \check{a}},  \]
   i.e. the same direction $a$-decompositions and $b$-decompositions commute.
  \item Assume $\check{a}<^+\check{b}$. Then we have the following farther equations of ordered face structures:
  \[  X^{\ua \check{b}}=  X^{\ua \check{a} \ua \check{b}},\;\;\;
  X^{\da \check{a}}=  X^{\da \check{a} \da \check{b}},\;\;\;
 X^{\da \check{b} \ua \check{a}}=  X^{\ua \check{a} \da \check{b}}. \]
  Moreover, if $\check{a},\check{b}\in Sd(X)$ then $\check{a}\in Sd_{\kappa^{\ua \check{b}}}(X^{\da \check{b}})$
  and $\check{b}\in Sd_{\kappa^{\ua \check{a}}}(X^{\ua \check{a}})$.
  \item Assume $\check{a}<^-_l\check{b}$, for some $l<m$.
  Then $X^{\ua \check{b} \da \check{a}}$, $X^{\ua \check{a} \da \check{b}}$,
  are ordered face structures, and
  \[ X^{\da \check{a}}\otimes_m X^{\ua \check{a} \da \check{b}}=
  X^{\da \check{b}}\otimes_m X^{\ua \check{b} \da \check{a}}\]
  Moreover, if $a,b\in Sd(X)$ then either there is $k$ such that
  $l-1\leq k<m$ and $(\gamma^{(k)}(a),-,\ua \gamma^{(k+1)}(a))\in Sd(X)$
  or $\check{a}\in Sd_{\kappa^{\ua \check{b}}}(X^{\ua \check{b}})$ and
  $\check{b}\in Sd_{\kappa^{\ua \check{a}}}(X^{\ua \check{a}})$.
\end{enumerate}
\end{lemma}
{\it Proof.}~ Easy. $~\Box$

The following properties of ordered face structures are inherited
from the corresponding properties of positive face structures.

\begin{lemma}
\label{decomp4} Let $T$ be ordered face
structures of dimension $n$, $l< n-1$, $\check{a}=(a,L,U)\in
Sd(T)_l$. Then
\begin{enumerate}
  \item $\check{a}\in Sd(\bc T)\cap Sd(\bd T)$;
  \item $\bd(T^{\da \check{a}})=(\bd T)^{\da \check{a}}$;
  \item $\bd(T^{\ua \check{a}})=(\bd T)^{\ua \check{a}}$;
  \item $\bc(T^{\da \check{a}})=(\bc T)^{\da \check{a}}$;
  \item $\bc(T^{\ua \check{a}})=(\bc T)^{\ua \check{a}}$.
\end{enumerate}
\end{lemma}

{\it Proof.}~See the the corresponding properties of positive face
structures in \cite{Z}.$~\Box$

\begin{lemma}
\label{decomp3new} Let $T, T_1,T_2$ be ordered face structures,
$dim(T_1),dim(T_2)>k$, such that $\bc^{(k)}(T_1)=\bd^{(k)}(T_2)$
and $T=T_1\otimes_kT_2$, and let
$Z=\gamma((T_1)_{k+1})-\delta((T_1^{-\lambda})_{k+1})$. Then
$\emptyset\neq Z\subseteq\bc^{(k)}(T_1)_k$. For any face $a\in Z$,
the cut
$\check{a}=(a,\cI_a\cap(T_1)_{k+2},\cI_a\cap(T_2)_{k+2})\in Sd(T)$
and one of the following conditions holds:
\begin{enumerate}
  \item either $T_1=T^{\da \check{a}}$ and $T_2=T^{\ua \check{a}}$;
  \item or $\check{a}\in Sd(T_1)_k$, $T^{\da \check{a}}=T_1^{\da \check{a}}$ and
$T^{\ua \check{a}}=T_1^{\ua \check{a}}\otimes_kT_2$;
  \item or $\check{a}\in Sd(T_2)_k$, $T^{\ua \check{a}}=T_2^{\ua \check{a}}$ and
$T^{\da \check{a}}=T_1\otimes_kT_2^{\da \check{a}}$.
\end{enumerate}
\end{lemma}

{\it Proof.}~See the the corresponding properties of positive face
structures in \cite{Z}.$~\Box$

\section{$T^*$ is a many-to-one computad}\label{Sstar_is_pComp}

\begin{proposition}\label{Sstar} Let $T$ be an ordered face structure.
Then $T^*$ is a many-to-one computad, whose indets correspond to
the faces of $T$.
\end{proposition}

{\it Proof.}~In fact, to be able to carry on the induction we need
to prove more. Let $T$ be an ordered face structure, $n\in\o$.

{\em Inductive Hypothesis for $n$}. For any ordered face structure
$T$, the $n$-truncation $T^*_{\leq n}$ of $T^*$ is a many-to-one
computad whose $n$-indets are in the image of the embedding $\nu:
T_n\lra T^*_n$, sending $a\in T_n$ to the local morphism $\nu_a:[a]\lra
T$ in $T^*_{\leq n}$.

The proof proceeds by induction on $n$. The Inductive Hypothesis
for $n=0,1$ is obvious.

So assume that the Inductive Hypothesis holds already for some
$n\geq 1$. Suppose that $T$ is an ordered face structure. We shall
show that $T^*_{\leq n+1}$ is a many-to-one computad whose
$n+1$-indets are in the image of $\nu: T_{n+1}\lra T^*_{n+1}$.

We need to verify that for any $\o$-functor $f:T_{\leq n}^*\lra C$
to any $\o$-category $C$, and any function $|f|:T_{n+1}\lra
C_{n+1}$ such that for $a\in T_{n+1}$, and $\nu_a:[a]\ra T$
\[ d_C(|f|(a)) = f(\bd(\nu_a)), \;\;\;\;\; c_C(|f|(a)) = f(\bc(\nu_a)), \]
there is a unique $\o$-functor $F:T^*_{\leq n+1}\lra C$, such that
\[ F_{n+1}(\nu_a)=|f|(a), \;\;\;\;\; F_{\leq n}=f \]
as in the diagram
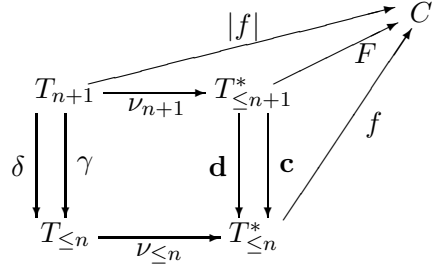
\begin{center}
\xext=650 \yext=1000
\begin{picture}(\xext,\yext)(\xoff,\yoff)
 \setsqparms[1`0`0`1;650`500]
 \putsquare(0,50)[T_{n+1}`T^*_{\leq n+1}`T_{\leq n}`T_{\leq n}^*;```{\nu_{\leq n}}]
  \putmorphism(-100,550)(0,-1)[``\delta ]{500}{1}l
  \putmorphism(0,550)(0,-1)[``\gamma ]{500}{1}r
  \putmorphism(600,550)(0,-1)[``\bd]{500}{1}l
  \putmorphism(700,550)(0,-1)[``\bc]{500}{1}r
\putmorphism(860,650)(2,1)[``F]{150}{1}r
\putmorphism(360,650)(4,1)[``]{500}{1}l
\putmorphism(800,180)(2,3)[``f]{350}{1}r
\put(1180,780){\makebox(100,100){$C$}}
\put(265,450){\makebox(100,100){$\nu_{n+1}$}}
\put(565,730){\makebox(100,100){$|f|$}}
\end{picture}
\end{center}

We need some notation for decompositions of cells in $T^*$. If
$\phi:X\ra T \in T^*$ and $\check{a}$ is a cut in $X^\dag$ then
$\phi^{\da \check{a}}=\kappa^{\da \check{a}};\phi : X^{\da
\check{a}}\lra T$, and $\phi^{\ua \check{a}}=\kappa^{\ua
\check{a}};\phi : X^{\da \check{a}}\lra T$.

We define $F_{n+1}$ as follows. For $\phi:X\ra T \in T^*_{n+1}$

\[ F_{n+1}(\phi) = \left\{ \begin{array}{ll}
                 id_{f(\phi)} & \mbox{ if $dim(X)\leq n$,} \\
                 |f|(a)  & \mbox{ if $\phi =\nu_a:[a]\ra T$, for some $a\in T_{n+1}$,}\\
                 F_{n+1}(\phi^{\da \check{a}});_l F_{n+1}(\phi^{\ua \check{a}}) &
                            \mbox{ if $dim(X)=n+1$, $\check{a}\in Sd(X)_l$.}
                 \end{array}  \right. \]
Clearly $F_k=f_k$, for $k\leq n$. The above morphism, if well
defined, clearly preserves identities. We need to verify, for any
$\phi:X\ra T$ in $T^*_{n+1}$, the following three conditions:
\begin{itemize}
  \item[\bf I] $F$ is well defined, i.e. for $\check{a},\check{b}\in Sd(X)$ we have
  $F_{n+1}(\phi^{\da \check{a}});_l F_{n+1}(\phi^{\ua \check{a}})=F_{n+1}(\phi^{\da \check{b}});_l F_{n+1}(\phi^{\ua
  \check{b}})$,
  \item[\bf II] $F$ preserves the domains and codomains i.e.we have
  $F(\bd \phi)=d(F(\phi))$ and $F(\bc \phi)=c(F(\phi))$,
  \item[\bf III] $F$ preserves compositions i.e.,  we have
  $F(\phi)=F(\phi_1);_kF(\phi_2)$
   whenever $\phi_i:X_i\ra T \in T^*_{n+1}$ for  $i=1,2$, $c^{(k)}(\phi_1)=d^{(k)}(\phi_2)$,
   and $\phi=\phi_1;_k\phi_2$.
\end{itemize}

Assume that $\phi:X\ra T\in T^*_{k+1}$, and for faces $y:Y\ra T$
of $T^*$ of size less than $size(X)$ the conditions [{\bf I}],
[{\bf II}], [{\bf III}] holds. We shall show that [{\bf I}], [{\bf
II}], [{\bf III}] hold for $\phi$, as well. For $X$ such that
$size(X)_{n+1}=0$ all three conditions are obvious.

If $X$ is principal of dimension $n+1$, [{\bf I}] is trivially
true as $Sd(X)=\emptyset$,  [{\bf III}] is true as if
$\phi=\phi_1;_k\phi_2$, with $X$ principal then either $\phi_1$ or
$\phi_2$ is an identity. So we need to check [{\bf II}]. We have
that $X_{n+1}=\{ m^X \}$ and $\phi(m^X)=a\in T_{k+1}$. By Lemma
\ref{principal ofs}.3, there is a unique isomorphism $h:[a]\ra X$
 making the
triangle
\begin{center}
\xext=600 \yext=350
\begin{picture}(\xext,\yext)(\xoff,\yoff)
\settriparms[1`1`1;300]
\putVtriangle(0,0)[{[\phi(m^X)]}`X`T;h`\nu_{\phi(m^X)}`\phi]
\end{picture}
\end{center}
commutes, i.e. $\nu_{\phi(m^X)}$ and $\phi$ represent the same
cell in $T^*$, and hence [{\bf II}] follows immediately from the
properties of $f$.

Now assume that $X$ is not principal and $dim(X)=n+1$.

Ad {\bf I}. First we will consider two saddle cuts
$\check{a},\check{x}\in Sd(X)$ of different dimension
$k=dim(x)<dim(a)=m$. Using Lemma \ref{decomp1new} we have
\begin{eqnarray*}
 F(\phi^{\da \check{a}});_m F(\phi^{\ua \check{a}}) = & ind.\; hyp.\; III    \\
 = (F(\phi^{\da \check{a}\da \check{x}});_k F(\phi^{\da \check{a}\ua \check{x}}));_m
 (F(\phi^{\ua \check{a} \da \check{x}});_k F(\phi^{\ua \check{a}\ua \check{x}})) = & MEL \\
 = (F(\phi^{\da \check{a}\da \check{x}});_m F(\phi^{\ua \check{a} \da \check{x}}));_k
 (F(\phi^{\da \check{a}\ua \check{x}});_m F(\phi^{\ua \check{a}\ua \check{x}})) = &  \\
 = (F(\phi^{\da \check{x}\da \check{a}});_m F(\phi^{\da \check{x} \ua \check{a}}));_k
 (F(\phi^{\ua \check{x}\da \check{a}});_m F(\phi^{\ua \check{x}\ua \check{a}})) =
 & ind.\; hyp.\; III  \\
 = F(\phi^{\da \check{x}});_m F(\phi^{\ua \check{x}}) &
\end{eqnarray*}
Now we will consider two saddle cuts $\check{a},\check{b}\in
Sd(X)$ of the same dimension $dim(a)=dim(b)=m$. We shall use Lemma
\ref{decomp2new}. Assume that $\check{a}<^-_l \check{b}$, for some
$l<m$. If $\check{x}=(\gamma^{(k)}(a),-,\ua\gamma^{(k+1)}(a))\in
Sd(X)$, for some $k<m$, then this case reduces to the previous one
for two pairs $\check{a},\check{x}\in Sd(X)$ and
$\check{b},\check{x}\in Sd(X)$. Otherwise $\check{a}\in Sd(X^{\ua
\check{b}})$, $\check{a}\in Sd(X^{\ua \check{b} })$, and we have
\begin{eqnarray*}
 F(\phi^{\da \check{a}});_k F(\phi^{\ua \check{a}}) = & ind.\; hyp\; III\\
 = F(\phi^{\da \check{a}});_k (F(\phi^{\ua \check{a}\da \check{b}});_k
 F(\phi^{\ua \check{a}\ua \check{b}})) = &  \\
 = (F(\phi^{\da \check{a}});_k F(\phi^{\ua \check{a}\da \check{b}}));_k
 F(\phi^{\ua \check{b}\ua \check{a}}) = & ind\; hyp\; III\\
  = F(\phi^{\da \check{a}};_k \phi^{\ua \check{a}\da \check{b}});_k
  F(\phi^{\ua \check{b}\ua \check{a}}) = &  \\
  = F(\phi^{\da \check{b}};_k \phi^{\ua \check{b}\da \check{a}});_k
  F(\phi^{\ua \check{b}\ua \check{a}}) =  & ind\; hyp\;  III\\
  = (F(\phi^{\da \check{b}});_k F(\phi^{\ua \check{b}\da \check{a}}));_k
  F(\phi^{\ua \check{b}\ua \check{a}}) = &   \\
  = F(\phi^{\da \check{b}});_k (F(\phi^{\ua \check{b}\da \check{a}});_k
  F(\phi^{\ua \check{b}\ua \check{a}})) = & ind\; hyp\;  III \\
 = F(\phi^{\da \check{b}});_k F(\phi^{\ua \check{b}})  &
\end{eqnarray*}
Finally, we consider the case $a<^+b$. We have
\begin{eqnarray*}
 F(\phi^{\da \check{a}});_k F(\phi^{\ua \check{a}}) = & ind.\; hyp\; III  \\
 = F(\phi^{\da \check{a}});_k (F(\phi^{\ua \check{a}\da \check{b}});_k F(\phi^{\ua \check{a}\ua \check{b}})) = & \\
 = (F(\phi^{\da \check{b}\da \check{a}});_k F(\phi^{\da \check{b}\ua \check{a}}));_k
 F(\phi^{\ua \check{b}}) = & ind\; hyp\; III \\
  = F(\phi^{\da \check{b}});_k F(\phi^{\ua \check{b}}) &
\end{eqnarray*}
This shows that $F(\phi)$ is well defined.

Ad {\bf II}. We shall show that the domains are preserved. The
proof that, the codomains are preserved, is similar.

The fact that if $Sd(X)=\emptyset$ then $F$ preserves domains and
codomains follows immediately from the assumption on $f$ and $|f|$
and Lemma \ref{Sd}. So assume that $Sd(X)\neq\emptyset$ and  let
$\check{a}\in Sd(X)_k$. We use Lemma \ref{decomp4}. We have to
consider two cases:  $k<n$, and $k=n$. If $k<n$ then
\begin{eqnarray*}
F_n(d(\phi))
 = F_n(d(\phi)^{\da \check{a}});_k F_n(d(\phi)^{\ua \check{a}})= & \\
 = F_n(d(\phi^{\da \check{a}}));_k F_n(d(\phi^{\ua \check{a}}))= & ind\; hyp\; II \\
 = d(F_{n+1}(\phi^{\da \check{a}}));_k d(F_{n+1}(\phi^{\ua \check{a}}))= &  \\
 = d(F_{n+1}(\phi^{\da \check{a}});_k F_{n+1}(\phi^{\ua \check{a}}))= & ind\; hyp\; I \\
 = d(F_{n+1}(\phi)) &
\end{eqnarray*}
If $k=n$ then
\begin{eqnarray*}
F_n(d(\phi))= F_n(d(\phi^{\da \check{a}};_n \phi^{\ua \check{a}}))= &   \\
 = F_n(d(\phi^{\da \check{a}}))= & ind\; hyp\; II \\
 = d(F_{n+1}(\phi^{\da \check{a}}))= &  ind\; hyp\; II\\
 = d(F_{n+1}(\phi^{\da \check{a}});_nF_{n+1}(\phi^{\ua \check{a}}) )= & ind\; hyp\; I \\
 = d(F_{n+1}(\phi)) &
\end{eqnarray*}

Ad {\bf III}. Suppose that $\phi=\phi_1;_k\phi_2$. We shall show
that $F$ preserves this composition.  If $dim(X_1)=k$ then
$\phi_2=\phi$, $\phi_1=\bd^{(k)}(\phi)$. We have
\[ F_{n+1}(\phi)=F_{n+1}(\phi_2)=1^{(n+1)}_{F_k(\bd^{(k)}(\phi_2))};_kF_{n+1}(\phi_2)= \]
\[ =1^{(n+1)}_{F_k(\phi_1)};_kF_{n+1}(\phi_2)= F_{n+1}(\phi_1);_kF_{n+1}(\phi_2) \]
The case $dim(X_2)=k$ is similar. So now assume that
$dim(X_1),dim(X_2)>k$. We shall use Lemma \ref{decomp3new}. Let us
fix a face
$a\in\gamma((X_1)_{k+1})-\delta((X_1^{-\lambda})_{k+1})$, and a
cut $\check{a}=(a,\cI_a\cap(X_1)_{k+2},\cI_a\cap(X_2)_{k+2})\in
Sd(X)$.

 If $X_1=X^{\da \check{a}}$ and $X_2=X^{\ua \check{a}}$
then we have
\[ F(\phi)=F(\phi^{\da \check{a}});_kF(\phi^{\ua \check{a}})= F(\phi_1);_kF(\phi_2). \]

If $\check{a}\in Sd(T_1)_k$, $T^{\da \check{a}}=T_1^{\da
\check{a}}$ and $T^{\ua \check{a}}=T_1^{\ua
\check{a}}\otimes_kT_2$
\begin{eqnarray*}
F(\phi)=F(\phi^{\da \check{a}});_kF(\phi^{\ua \check{a}})= & ind\; hyp\; III\\
 =F(\phi^{\da \check{a}});_k(F(\phi_1^{\ua \check{a}});_kF(\phi_2))= &  \\
 =(F(\phi_1^{\da \check{a}});_kF(\phi_1^{\ua \check{a}}));_kF(\phi_2)= & ind\; hyp\; III \\
 F(\phi_1);_kF(\phi_2) &
\end{eqnarray*}

If $\check{a}\in Sd(T_2)_k$, $T^{\ua \check{a}}=T_2^{\ua
\check{a}}$ and $T^{\da \check{a}}=T_1\otimes_kT_2^{\da
\check{a}}$

\begin{eqnarray*}
F(\phi)=F(\phi^{\da \check{a}});_kF(\phi^{\ua \check{a}})= & ind\; hyp\; III\\
 =(F(\phi_1);_kF(\phi_2^{\da \check{a}}));_kF(\phi_2^{\ua \check{a}}))= &  \\
 =F(\phi_1^{\da \check{a}});_k(F(\phi_2^{\ua \check{a}});_kF(\phi_2^{\ua \check{a}}))= & ind\; hyp\; III \\
 F(\phi_1);_kF(\phi_2) &
\end{eqnarray*}

So in any case the composition is preserved. This ends the proof
of the Lemma. $~~\Box$

For $n\in\o$, we have truncation functors
\[ (-)^{\sharp,n} : \ofs_{loc}\lra \mnComma,\hskip 10mm (-)^{*,n} : \ofs_{loc}\lra \mnComp \]
such that, for $S$ in $\ofs$
 \[ S^{\sharp,n}=(S_n,S^*_{<n}, [\delta],[\gamma]),\hskip 10mm (S)^{*,n}=S^*_{\leq n} \]
and for $f:S\ra T$ in $\ofs_{loc}$ we have
\[ f^{\natural,n}=(f_n,(f_{<n})^*),\hskip 10mm (f)^{*,n}=f^*_{\leq n}. \]

\begin{corollary}
For every $n\in\o$, the functors $(-)^{\sharp,n}$ and $(-)^{*,n}$
are well defined, full, faithful, and they send all tensor squares
to pushouts. Moreover, for $S$ in $\ofs$ we have
$S^*=\overline{S^{\sharp,n}}^n$.
\end{corollary}
 {\it Proof.}~ The functor $\overline{(-)}^n: \mnComma \lra
\mnComp$, which is an equivalence of categories, is described in
the Appendix.

Fullness and faithfulness of $(-)^{\sharp,n}$ is left for the
reader.  We shall show simultaneously that for every $n\in\o$,
both functors $(-)^{\sharp,n}$ and $(-)^{*,n}$ send
$n$-truncations of tensor squares to pushouts. For $n=0,1$ this is
obvious. So assume that $n\geq 1$ and that $(-)^{*,n}$ and
$(-)^{\sharp,n}$ send $n$-truncations of $k$-tensor squares to
pushouts. Let
\begin{center} \xext=500 \yext=400
\begin{picture}(\xext,\yext)(\xoff,\yoff)
 \setsqparms[1`-1`-1`1;500`400]
 \putsquare(0,0)[S`S\otimes_kT`\bc^{(k)} S`T;```]
\end{picture}
\end{center}
be a tensor squares in $\ofs$. The fact that the functor
$(-)^{\sharp,n+1}$ sends this square to a pushout in $\mnjComma$
can be verified in each dimension separately. In dimensions lower
or equal to $n$ this follows from the fact that the functor
$(-)^{*,n}$ sends $n$-truncations of tensor squares to pushouts.
In dimension $n+1$ we need to check that the square in $Set$
\begin{center} \xext=700 \yext=400
\begin{picture}(\xext,\yext)(\xoff,\yoff)
 \setsqparms[1`-1`-1`1;700`400]
 \putsquare(0,0)[S_{n+1}`(S\otimes_kT)_{n+1}`(\bc^{(k)} S)_{n+1}`T_{n+1};```]
\end{picture}
\end{center}
is a pushout. But this easily follows from the the description of
the tensor square given earlier. So the whole square
\begin{center} \xext=800 \yext=400
\begin{picture}(\xext,\yext)(\xoff,\yoff)
 \setsqparms[1`-1`-1`1;800`400]
 \putsquare(0,0)[S^{\sharp,n+1}`(S\otimes_kT)^{\sharp,n+1}`\bc^{(k)} S^{\sharp,n+1}`T^{\sharp,n+1};```]
\end{picture}
\end{center}
is a pushout in $\mnjComma$, i.e. $(-)^{\sharp,n+1}$ send
$(n+1)$-truncations of $k$-tensor squares to pushouts. As
$(-)^{*,n+1}$ is a composition of $(-)^{\sharp,n+1}$ with an
equivalence of categories it send $(n+1)$-truncations of
$k$-tensor squares to pushouts, as well. $~~\Box$

\begin{corollary}\label{ff-pres_sp} The functor
\[ (-)^* : \ofs_{loc} \lra \mComp \]
is full and faithful and sends tensor squares to pushouts.
\end{corollary}
{\it Proof.}~This follows from the previous Corollary and the fact
that the functor $\overline{(-)}^n: \mnComma \lra \mnComp$ (see
Appendix) is an equivalence of categories. $~~\Box$

Let $P$ be a many-to-one computad, $a$ a $k$-cell in $P$.  A {\em
description of the cell} $a$\index{description of a
cell}\index{cell!description of a -}  is a pair
  \[ <T_{a}, \tau_{a}: T_{a}^*\lra P> \]
where $T_{a}$ is an ordered face structure and $\tau_{a}$ is a
computad map such that \[ \tau_{a}(id_{T_{a}})=a. \]

\section{The terminal many-to-one computad}

In this section we shall describe the terminal many-to-one
computad $\cT$.

The set of $n$-cell $\cT_n$ consists of (isomorphisms classes of)
ordered face structures of dimension less than or equal to $n$.
For $n>0$, the operations of domain and codomain
$d^\cT,c^\cT:\cT_n\ra \cT_{n-1}$ are given, for $S\in\cT_n$  by

\[  d(S) = \left\{ \begin{array}{ll}
        S   & \mbox{if $dim(S)<n$,}  \\
        \bd S & \mbox{if $dim(S)=n$,}
                                    \end{array}
                \right. \]
and
\[  c(S) = \left\{ \begin{array}{ll}
        S   & \mbox{if $dim(S)<n$,}  \\
        \bc S & \mbox{if $dim(S)=n$.}
                                    \end{array}
                \right. \]
and, for $S,S'\in\cT_n$ such that $c^{(k)}(S)=d^{(k)}(S')$ the
composition is just the $k$-tensor of $S$ and $S'$ as ordered face
structures i.e. $S\otimes_kS'$

The identity $id_\cT : T_{n-1}\ra \cT_n$ is the inclusion map.

The $n$-indets in $\cT$ are the principal ordered $n$-face
structures.

\begin{proposition} $\cT$ described above is the terminal many-to-one
 computad.
\end{proposition}

{\it Proof.}~ The fact that $\cT$ is an $\o$-category is easy. The
fact that $\cT$ is free with free $n$-generators being principal
$n$-face structures can be shown much like the freeness of $S^*$
before. The fact that $\cT$ is terminal follows from the following
observation:

{\em Observation.} For every parallel pair of ordered face
structures $N$ and $B$ (i.e. $\bd N = \bd B$ and $\bc N = \bc B$)
such that $N$ is normal and $B$ is principal, there is a unique
(up to an iso) principal ordered face structure $N^\bullet$ such
that $\bd N^\bullet =N$ and $\bc N^\bullet =B$. $~~\Box$

\begin{lemma}\label{shape_pres} Let $S$ be an ordered face structure and
$!:S^*\lra \cT$ the unique computad map from $S^*$ to $\cT$. Then,
for $x:X\ra S\in S^*_k$ we have
\[ !_k(x)=X.\]
\end{lemma}

{\it Proof.}~ The proof is by induction on $k\in\o$ and the size
of $X$ in $S^*_k$. For $k=0,1$ the lemma is obvious. Let $k>1$ and
assume that lemma holds for $i<k$.

If $dim(X)=l<k$ then, using the inductive hypothesis and the fact
that $!$ is an $\o$-functor, we have
\[ !_k(x)=!_k(1^{(k)}_x) = 1^{(k)}_{!_l(x)} = 1^{(k)}_{X} = X \]

Suppose that $dim(X)=k$ and $X$ is principal. As $!$ is a computad
map $!_k(x)$ is an indet, i.e. it is principal, as well. We have,
using  again the inductive hypothesis and the fact that $!$ is an
$\o$-functor,
\[ d(!_k(x))=!_{k-1}(\bd x)= \bd X \]
\[ c(!_k(x))=!_{k-1}(\bc x)= \bc X \]
As $X$ is the only (up to a unique iso) ordered face structure
with the domain $\bd X$ and the codomain $\bc X$, it follows that
$!_k(x)=X$, as required.

Finally, suppose that $dim(X)=k$, $X$ is not principal, and for
the ordered face structures of size smaller than the size of $X$
the lemma holds. Thus there are $l\in\o$ and $\check{a}\in
Sd(X)_l$ so that
\[ !_k(x)= !_k(x^{\ua \check{a}};_lx^{\da \check{a}}) = !_k(x^{\ua \check{a}})\otimes_l!_k(x^{\da
\check{a}})= X^{\ua \check{a}}\otimes_lX^{\da\check{a}} = X, \]
 as required
 $~~\Box$

\section{A description of the many-to-one computads}

In this section we shall describe all the cells in many-to-one
computads using ordered face structures, in other words we shall
describe in concrete terms the functor:
\[ \overline{(-)}: \mnComma \lra \mnComp \]

More precisely, the many-to-one computads of dimension 1 (and all
computads as well) are free computads over graphs and they are
well understood. So suppose that $n>1$, and we are given an object
of $\mnComma$, i.e. a quadruple $(|\cP |_n,\cP,d,c)$ such that
\begin{enumerate}
  \item a many-to-one $(n-1)$-computad $\cP$;
  \item a set $|P|_n$ with two functions $c: |\cP|_n\lra |\cP|_{n-1}$
  and $d: |\cP|_n\lra \cP_{n-1}$ such that for $x\in |\cP|_n$, $cc(x)=cd(x)$ and
  $dc(x)=dd(x)$.
\end{enumerate}
If the maps $d$ and $c$ in the object $(|\cP|_n,\cP,d,c)$ are
understood from the context we can abbreviate notation to
$(|\cP|_n,\cP)$.

For an ordered face structure $S$, we denote by $S^{\sharp,n}$ the
object  $(S_n,(S_{<n})^*,[\delta],[\gamma])$ in $\mnComma$. In
fact, we have an obvious functor
\[ (-)^{\sharp,n}: \ofs_{loc} \lra \mnComma \]
such that
 \[ S\mapsto S^{\sharp,n}=(S_n,(S_{<n})^*,[\delta],[\gamma]) \]
Any many-to-one computad $\cP$ can be restricted to its part in
$\mnComma$.  So we have an obvious forgetful functor
\[ (-)^{\natural,n} :\mComp \lra \mnComma \]
such that
\[ \cP \mapsto \cP^{\natural,n}=(|\cP|_n,\cP_{<n},d,c) \]

We shall describe the many-to-one $n$-computad $\overline{\cP}$
whose $(n-1)$-truncation is $\cP$ and whose $n$-indets are
$|\cP|_n$ with the domains and codomains given by $c$ and $d$.

{\bf n-cells} of $\overline{\cP}$. An $n$-cell in
$\overline{\cP}_n$ is a(n equivalence class of) pair(s)
 $(S,f)$ where
\begin{enumerate}
  \item $S$ is an ordered face structure, $dim(S)\leq n$;
  \item $f: S^{\sharp,n}\lra \cP^{\natural,n}$
  is a morphism in $\mnComma$, i.e.
   \begin{center}
\xext=650 \yext=600
\begin{picture}(\xext,\yext)(\xoff,\yoff)
 \setsqparms[1`0`0`1;650`500]
 \putsquare(0,50)[S_n`|\cP|_n`S^*_{n-1}`\cP_{n-1};|f|_n```f_{n-1}]
  \putmorphism(-50,550)(0,-1)[``{[}\delta {]}]{500}{1}l
  \putmorphism(50,550)(0,-1)[``{[}\gamma {]}]{500}{1}r
  \putmorphism(600,550)(0,-1)[``d]{500}{1}l
  \putmorphism(700,550)(0,-1)[``c]{500}{1}r
\end{picture}
\end{center}
commutes.
\end{enumerate}

We identify two pairs $(S,f)$, $(S',f')$ if there is a monotone
isomorphism $h:S\lra S'$ such that the triangles of sets and of
$(n-1)$-computads
 \begin{center}
\xext=2200 \yext=500
\begin{picture}(\xext,\yext)(\xoff,\yoff)
 \settriparms[1`1`1;400]
 \putVtriangle(0,0)[S_n`S'_n`|\cP|_n;h_n`f_n`f'_n]
 \putVtriangle(1400,0)[(S_{< n})^*`(S'_{<n})^*`\cP;(h_{<n})^*`f_{<n}`f'_{<n}]
\end{picture}
\end{center}
commute. Clearly, such an $h$, if exists, is unique. Even if
formally cells in $\cP_n$ are equivalence classes of triples we
will work on triples themselves as if they were cells
understanding that equality between such cells is an isomorphism
in the sense defined above.

{\bf Domains and codomains} in $\overline{\cP}$. The domain and
codomain functions
\[ d^{(k)},c^{(k)}: \overline{\cP}_n \lra \overline{\cP}_{k} \]
are defined for an $n$-cell $(S,f)$ as follows:
\[ d^{(k)}(S,f) = \left\{ \begin{array}{ll}
                            (S,f) & \mbox{ if $dim(S)\leq k$,} \\
                            (\bd^{(k)}S,d^{(k)}f)  & \mbox{otherwise.}
                           \end{array}
                   \right. \]
where, for $x\in (\bd^{(k)}S)_k$ (and hence $\nu_x:[x]\ra
\bd^{(k)}S$),
\[ (d^{(k)}f)_k(x)= f_k(\nu_x)(x) \]
(i.e. we take the cell $\nu_x:[x]\ra \bd^{(k)}S$ of $S^*$, then
value of $f$ on it, and then we evaluate the map in $\mnComma$ on
$x$ the only element of $[x]_k$), and
\[ (\bd^{(k)}f)_{<k}= (\bd^{k}_S;f_{<n})_{<k}. \]

\[ c^{(k)}(S,f) = \left\{ \begin{array}{ll}
                            (S,f) & \mbox{ if $dim(S)\leq k$,} \\
                            (\bc^{(k)}S,c^{(k)}f)  & \mbox{otherwise.}
                           \end{array}
                   \right. \]
where, for $x\in (\bc^{(k)}S)_k$ (and hence $\nu_x:[x]\ra
 \bc^{(k)}S$),
\[ (c^{(k)}f)_k(x)= f_k(\nu_x)(x) \]
and
\[ (c^{(k)}f)_{<k}= (\bc^{k}_S;f_{<n})_{<k}. \]
i.e. we calculate the $k$-th domain and $k$-th codomain
of an $n$-cell $(S, f)$ by taking $\bd^{(k)}$ and $\bc^{(k)}$ of
the domain $S$ of the cell $f$, respectively, and by restricting
the maps $f$ accordingly.

{\bf Identities} in $\overline{\cP}$. The identity function
\[  \bi : \overline{\cP}_{n-1}\lra \overline{\cP}_n \]
is defined for an $(n-1)$-cell $((S,f)$ in $\cP_{n-1}$, as
follows:
\[ \bi(S, f) = \left\{ \begin{array}{ll}
                            (S,f) & \mbox{ if $dim(S)<n-1$,} \\
                            (S,\overline{f})  & \mbox{ if dim(S)=n-1}
                           \end{array}
                   \right. \]
Note that $\overline{f}$ is the map $\mnmjComp$ which is the value
of the functor $\overline{(-)}$ on a map $f$ from $\mnmjComma$. So
it is in fact defined as 'the same $(n-1)$-cell' but considered as
an $n$-cell.

{\bf Compositions} in $\overline{\cP}$. Suppose that $(S^i,f^i)$
are $n$-cells for $i=0,1$, such that
\[ c^{(k)}(S^0, f^0)= d^{(k)}(S^1, f^1). \]
Then their $k$-composition in $\mnComma$ is defined as
\[ (S^0, f^0);_k(S^1, f^1)= (S^0\otimes_k S^1,[f^0,f^1]) \]
i.e.
\begin{center}
\xext=1500 \yext=700
\begin{picture}(\xext,\yext)(\xoff,\yoff)
 \setsqparms[1`0`0`1;1450`500]
 \putsquare(0,100)[(S^0\otimes_kS^1)_n`|\cP|_n`
 (S^0\otimes_kS^1)^*_{n-1}`\cP_{n-1};
 [f^0_n,f^1_n]```{[f^0_{n-1},f^1_{n-1}]}]
  \putmorphism(-50,600)(0,-1)[``{[}\delta {]}]{500}{1}l
  \putmorphism(50,600)(0,-1)[``{[}\gamma {]}]{500}{1}r
  \putmorphism(1400,600)(0,-1)[``d]{500}{1}l
  \putmorphism(1500,600)(0,-1)[``c]{500}{1}r
\end{picture}
\end{center}
This ends the description of the computad $\overline{\cP}$.

Now let $h:\cP\ra \cQ$ be a morphism in $\mnComma$, i.e. a
function $h_n:|\cP|_n\lra |\cQ|_n$ and a $(n-1)$-computad morphism
$h_{<n}: \cP_{<n}\lra \cQ_{<n}$ such that the square
\begin{center} \xext=1500 \yext=700
\begin{picture}(\xext,\yext)(\xoff,\yoff)
 \setsqparms[1`0`0`1;950`500]
 \putsquare(0,100)[|\cP|_n`|Q|_n`\cP_{n-1}`\cQ_{n-1};
 h_n```h_{n-1}]
  \putmorphism(-50,600)(0,-1)[``d]{500}{1}l
  \putmorphism(50,600)(0,-1)[``c]{500}{1}r
  \putmorphism(900,600)(0,-1)[``d]{500}{1}l
  \putmorphism(1000,600)(0,-1)[``c]{500}{1}r
\end{picture}
\end{center}
commutes serially. We define
\[ \overline{h}: \overline{\cP}\lra \overline{Q} \]
by putting $\overline{h}_k=h_k$, for $k<n$, and for $(S,f)\in
\overline{\cP}_n$, we put
\[ \overline{h}(S,f)= (S,h\circ f).\]

{\bf Embedding} $ \eta_\cP: |\cP|_n\lra \overline{\cP}_n$ is
defined in the Proposition below.

{\em Notation.} Let $x=(X,f:X^{\sharp,n}\ra \cP^{\natural,n})$ be
a cell in $\overline{\cP}_n$ as above, and $\check{a}\in Sd(X)_k$.
Then by $x^{\ua \check{a}}=(X^{\ua \check{a}},f^{\ua \check{a}})$
and $x^{\da \check{a}}=(X^{\da \check{a}},f^{\da \check{a}})$ we
denote the cells in $\overline{\cP}_n$ that are the obvious
restrictions of $x$. Clearly, we have $c^{(k)}(x^{\ua \check{a}})=
d^{(k)}(x^{\da \check{a}})$ and that $x=x^{\ua \check{a}};_kx^{\da
\check{a}}$, where $k=dim(a)$.

The following Proposition contains several statements.  We have
put all of the together since they have to be proved
simultaneously, i.e. to prove them for $n$ we need to know all of
them for $n-1$.

\begin{proposition}\label{type} Let $n\in\o$. We have
\begin{enumerate}
  \item Let $\cP$ be an object of $\mnComma$. We define
  the function \[ \eta_\cP: |\cP|_n\lra \overline{\cP}_n\] as follows.
  Let $x\in |\cP|_n$. As $c(x)$ is an indet $d(x)$ is a normal cell
  of dimension $n-1$.  Thus there is a unique descriptions of the cells
  $d(x)$ and $c(x)$
  \[ <T_{d(x)}, \tau_{d(x)}: T_{d(x)}^*\lra \cP_{<n}>,\;\;\;\;
  <T_{c(x)}, \tau_{c(x)}: T_{c(x)}^*\lra \cP_{<n}> \]
   with $T_{d(x)}$ being $(n-1)$-normal ordered face structure and $T_{c(x)}$ being principal
   ordered face structure of dimension $n-1$. Then we
   have a unique $n$-cell in $\overline{\cP}$:
   \[ \bar{x} = <T_{d(x)}^{\bullet},\;\; |\overline{\tau}_x|_n :
   \{ 1_{T_{d(x)}^{\bullet}} \} \ra |\cP|_n,\;\;
    (\overline{\tau}_x)_{<n} : (T_{d(x)}^{\bullet})_{<n}^*\ra \cP_{<n}> \]
   (note: $|T_{d(x)}^{\bullet}|_n=\{ 1_{T_{d(x)}^{\bullet}} \}$) such that
   \[ |\overline{\tau}_x|_n(1_{T_{d(x)}^{\bullet}})=x\]
   and, for $y:Y\ra T_{d(x)}^{\bullet}\in (T_{d(x)}^{\bullet})^*_{<n}$
   \[  (\overline{\tau}_x)_{n-1}(y) = \left\{ \begin{array}{ll}
                         (\tau_{c(x)})_{n-1}(y') & \mbox{if $Y$ is principal} \\
                                                  & \mbox{ and $y=y';\bc_{(T_{d(x)}^{\bullet})}$,} \\
                         (\tau_{d(x)})_{n-1}(y'')  & \mbox{if $Y$ is
                         principal} \\
                         & \mbox{ and $y=y'';\bd_{(T_{d(x)}^{\bullet}),}$}\\
                         (\overline{\tau}_x)_{n-1}(y^{\da \check{a}});_k(\overline{\tau}_x)_{n-1}(y^{\ua
                            \check{a}})& \mbox{if $\check{a}\in Sd(Y)_k$}
                           \end{array}
                   \right. \]
  and $(\overline{\tau}_x)_{<(n-1)}=(\tau_{dx})_{<(n-1)}$. We put
  $\eta_\cP(x)=\bar{x}$.

Then $\overline{\cP}$ is a many-to-one computad with $\eta_\cP$
the inclusion of $n$-indeterminates. Moreover, any many-to-one
$n$-computad $Q$ is equivalent to a computad $\overline{\cP}$, for
some $\cP$ in $\mnComma$.

  \item  Let $\cP$ be an object of $\mnComma$, $!:\overline{\cP}\lra
  \cT$ the unique morphism into the terminal object $\cT$ and
  $f:S^{\sharp,n}\ra \cP$ a cell in $\overline{\cP}_n$. Then
  \[ !_n(f:S^{\sharp,n}\ra \cP  )=S. \]

  \item Let $h:\cP\ra \cQ$ be an object of $\mnComma$. Then
  $\overline{h}: \overline{\cP}\lra \overline{\cQ}$ is a computad morphism.

  \item Let $k\leq n$, $S$ be an ordered face structure of dimension at most
  $n$, $f:S^*\lra \cP$ a morphism in $\mnComp$ and $y:Y\ra S \in S^*_k$. We have that
  \[ \overline{f}_k(y)=(f\circ y^*)^{\natural,k}(=f^{\natural,k}\circ y^{\sharp,k}:Y^{\sharp,k}\lra \cP^{\natural,k}). \]

 \item Let $S$ be an ordered face structure of dimension $n$,
 $\cP$ many-to-one computad, $g,h:S^*\lra \cP$ computad maps. Then
 \[ g=h \hskip 10mm {\rm iff} \hskip 10mm g_n(1_S)=h_n(1_S). \]

 \item Let $S$ be an ordered face structure of dimension at most
 $n$, $\cP$ be an object in $\mnComma$. Then we have a bijective correspondence
  $$
\begin{array}{c}
f:S^{\sharp,n}\lra \cP \;\;\in \mnComma
\\ \hline
\overline{f}:S^*\lra \overline{\cP} \;\;\in \mnComp
\end{array}
$$
such that, $\overline{f}_n(1_S)=f$, and for $g:S^*\lra
\overline{\cP}$ we have $g=\overline{g_n(1_S)}$.
  \item The map
  \[ \kappa^\cP_n:\coprod_{S} \comp (S^*,\overline{\cP}) \lra \overline{\cP}_n \]
  \[ g:S^*\ra \overline{\cP}\;\;\; \mapsto\;\;\; g_n(1_S),\]
  where coproduct is taken over all (up to iso) ordered face structures $S$ of
  dimension at most $n$, is a bijection. In other words, any cell
  in $\overline{\cP}$ has a unique description.
\end{enumerate}
\end{proposition}

{\it Proof.} ~We prove all the statements simultaneously by
induction on $n$. For $n=0,1$ all of them are easy.

Ad 1. We have to verify that $\overline{\cP}$ satisfy the laws of
$\o$-categories and that it is free in the appropriate sense. Laws
$\o$-categories are left for the reader. We shall show that
$\overline{\cP}$ is free in the appropriate sense.

Let $C$ be an $\o$-category, $g_{<n}:\cP_{<n}\ra C_{<n}$ an
$(n-1)$-functor and $g_n:|\cP|_n\ra C_n$ a function so that the
diagram
\begin{center} \xext=1500 \yext=700
\begin{picture}(\xext,\yext)(\xoff,\yoff)
 \setsqparms[1`0`0`1;950`500]
 \putsquare(0,100)[|\cP|_n`C_n`\cP_{n-1}`C_{n-1};
 g_n```g_{n-1}]
  \putmorphism(-50,600)(0,-1)[``d]{500}{1}l
  \putmorphism(50,600)(0,-1)[``c]{500}{1}r
  \putmorphism(900,600)(0,-1)[``d]{500}{1}l
  \putmorphism(1000,600)(0,-1)[``c]{500}{1}r
\end{picture}
\end{center}
commutes serially. We shall define an $n$-functor
$\overline{g}:\overline{\cP}\ra C$ extending $g$ and $g_n$. For
$x=(X,f)\in \overline{\cP}_n$ we put
\[  \overline{g}_n(x)= \left\{ \begin{array}{ll}
               1_{g_{n-1}\circ f_{n-1}(x)} & \mbox{if $dim(X)<n$,} \\
               g_{n}\circ f_{n}(m_{X}) &
               \mbox{if $dim(X)=n$, $X$ is principal, $X_n=\{m_X\}$ } \\
               \overline{g}_n(x^{\ua \check{a}});_k
               \overline{g}_n(x^{\da \check{a}}) & \mbox{if $dim(X)=n$, $\check{a}\in Sd(S)_k$}
                           \end{array}
                   \right. \]
We need to check that $\overline{g}$ is well defined, unique one
that extends $g$, preserves domains, codomains, compositions and
identities.

All these calculations are similar, and they are very much like
those in the proof of Proposition \ref{Sstar}. We shall check,
assuming that we already know that  $\overline{g}$ is well
defined, and preserves identities that compositions are preserved.
So let $x=(X,f)$, $x_1=(X_1,f_1)$, $x_2=(X_2,f_2)$ be cells in
$\overline{\cP}_n$ such that $x=x_1\textbf{;}_kx_2$. Since
$\overline{g}$ preserves identities, we can assume that
$dim(X_1),dim(X_1)>k$. Let $l\in\o$ be minimal such that
$Sd(X)_l\neq\emptyset$. We have two cases:

{\em Case 1}. If $l<k$, then by Decomposition 3.2.a we have
$\check{a}\in Sd(T_2)_l$, and then
\begin{eqnarray*}
\overline{g}(x) =  \\
 \overline{g}(x^{\ua \check{a}});_l \overline{g}(x^{\da \check{a}}) =  \\
 \overline{g}((x_1;_kx_2)^{\ua \check{a}});_l \overline{g}((x_1;_kx_2)^{\da \check{a}}) =  \\
 \overline{g}(x_1^{\ua \check{a}};_kx_2^{\ua \check{a}});_l \overline{g}(x_1^{\da \check{a}};_kx_2^{\da \check{a}}) =  \\
 (\overline{g}(x_1^{\ua \check{a}});_k\overline{g}(x_2^{\ua \check{a}}));_l
 (\overline{g}(x_1^{\da \check{a}});_k\overline{g}(x_2^{\da \check{a}})) =  \\
 (\overline{g}(x_1^{\ua \check{a}});_l\overline{g}(x_1^{\da \check{a}}));_k
 (\overline{g}(x_2^{\ua \check{a}})   ;_l\overline{g}(x_2^{\da \check{a}})) =  \\
 = \overline{g}(x_1);_k \overline{g}(x_2)
\end{eqnarray*}

{\em Case 2}. If $l=k$ then by Decomposition 3.2.a we have
$\check{a}\in Sd(X_1)$ and

\begin{eqnarray*}
 \overline{g}(x) =  \\
 \overline{g}(x^{\ua \check{a}});_k \overline{g}(x^{\da \check{a}}) =  \\
 \overline{g}(x_1^{\ua \check{a}});_k \overline{g}(x_1^{\da \check{a}};_kx_2) =  \\
 \overline{g}(x_1^{\ua \check{a}});_k (\overline{g}(x_1^{\da \check{a}});_k\overline{g}(x_2)) = \\
 (\overline{g}(x_1^{\ua \check{a}});_k \overline{g}(x_1^{\da \check{a}}));_k\overline{g}(x_2) = \\
 = \overline{g}(x_1);_k \overline{g}(x_2)
\end{eqnarray*}

The remaining things are similar.

 Ad 2. Let $!:\overline{\cP}\lra \cT$ be the unique computad map
into the terminal object, $S$ an ordered face structure such that
$dim(S)=l\leq n$, $f:S^{\sharp,n}\lra \cP$ a cell in
$\overline{\cP}_n$.

If $l<n$ then by induction we have $!_n(f)=S$. If $l=n$ and $S$ is
principal then we have, by induction
 \[ !_n(d(f):(\bd S)^{\sharp,n}\ra \cP  )=\bd S, \hskip 10mm
 !_n(c(f):(\bc S)^{\sharp,n}\ra \cP  )=\bc S. \]
As $f$ is an indet in $\overline{\cP}$, $!_n(f)$ is a principal
ordered face structure. But the only (up to an iso) principal
ordered face structure $B$ such that
\[ \bd B =\bd S, \hskip 10mm \bd B =\bd S \]
is $S$ itself. Thus, in this case, $!_n(f)=S$.

Now assume that $l=n$, and $S$ is not principal, and that for
ordered face structures $T$ of smaller size than $S$ the statement
holds. Let $a\in Sd(S)_k$. We have
\[ !_n(f)= !_n(f^{\ua a};_k f^{\da a}) =  !_n(f^{\ua a});_k !_n
(f^{\da a}) = S^{\ua a};_k S^{\da a} = S \]
 where $f^{\ua a}=f\circ (\kappa^{\ua a})^{\sharp,n}$ and $f^{\da a}=f\circ
(\kappa^{\da a})^{\sharp,n}$ and $\kappa^{\ua a}$ and $\kappa^{\da
a}$ are the monotone morphisms as in the following tensor square
\begin{center}
\xext=650 \yext=400
\begin{picture}(\xext,\yext)(\xoff,\yoff)
 \setsqparms[1`-1`-1`1;550`400]
 \putsquare(0,0)[S^{\ua a}`S`\bc^{(k)} S`S^{\da a};\kappa^{\ua a}``\kappa^{\da a}`]
\end{picture}
\end{center}

 Ad 3. The main thing is to show that $\overline{h}$ preserves
compositions. This follows from the fact that the functor
\[ (-)^{\sharp,n} :\ofs \lra \mnComma \]
preserves special pushouts.

 Ad 4. This is an immediate consequence of 3.

 Ad 5. Let us fix ordered face structures $S$, $Y$, $dim(S)=n$,
 $\check{a}\in Sd(Y)$, and $f,g:S^*\lra \cP$. Clearly, if $f=g$
 then $f(1_S)=g(1_S)$.  We shall prove the converse.  As
\begin{center} \xext=800 \yext=550
\begin{picture}(\xext,\yext)(\xoff,\yoff)
 \setsqparms[1`-1`-1`1;800`500]
 \putsquare(0,0)[(Y^{\da \check{a}})^{\sharp,k}`Y^{\sharp,k}`
 (\bc^{(l)} Y^{\da \check{a}})^{\sharp,k}`(Y^{\ua \check{a}})^{\sharp,k}; ```]
\end{picture}
\end{center}
 is a pushout in $\mkComma$ we have that for any $y:Y\ra S\in S^*_k$
 \[ f^{\natural,k}\circ y^{\sharp,k}=g^{\natural,k}\circ
 y^{\sharp,k}\mbox{ iff }
  f^{\natural,k}\circ (y^{\da \check{a}})^{\sharp,k}=g^{\natural,k}\circ  (y^{\da \check{a}})^{\sharp,k} \mbox{ and }
  f^{\natural,k}\circ (y^{\ua \check{a}})^{\sharp,k}=g^{\natural,k}\circ  (y^{\ua \check{a}})^{\sharp,k} \]
From this observation it is easy to see that if for some $y:Y\ra
S\in S^*$ we have $f(y)\neq g(y)$ then we can assume that this $Y$
is principal. On the other hand, from the above observation, the
fact that both $f$ and $g$ are $\o$-functors and that
$f(1_S)=g(1_S)$ we can deduce that for any $y:Y\ra S\in S^*$ with
$Y$ principal we have $f(y)=g(y)$. This together shows 5.

 Ad 6. we shall use 5. Fix an ordered face structure $S$ of dimension $n$ and a many-to-one
 computad $\cP$. For $f:S^{\sharp,n}\ra \cP^{\natural,n}$ in $\mnComma$ we have
\[ \overline{f}_n(1_S)=(f\circ (1_S)^{\sharp,n})^{\natural,n}=f\circ
(1_S)^{\sharp,k}=f. \] On the other hand, for a computad map
$g:S^*\ra P$ we have
\[ \overline{g_n(1_S)}(1_S)= (g_n(1_S)\circ(1_S)^{\sharp,n})^{\natural,n} = \]
\[ = (g^{\natural,n}\circ(1_S)^{\sharp,n}\circ(1_S)^{\sharp,n})^{\natural,n} =
(g^{\natural,n}\circ(1_S)^{\sharp,n})^{\natural,n} = g(1_S).\]
Thus by 5. we have $\overline{g_n(1_S)}=g$.

Ad 7. It follows immediately from 6. $~\Box$

The following Proposition says a bit more about descriptions than
point 7. of the previous one.

\begin{proposition}\label{type1}
Let $\cP$ be a many-to-one computad, $n\in\o$, and $a\in \cP_n$.
Let $T_a=!^\cP_n(a)$ (where $!^\cP:\cP\lra \cT$ is the unique
morphism into the terminal many-to-one computad). Then there is a
unique computad map $\tau_a : T_a^* \lra \cP$ such that
$(\tau_a)_n(1_{T_a})=a$. Moreover, we have:
\begin{enumerate}
  \item for any $a\in \cP$ we have
  \[ \tau_{da}= d(\tau_a) = \tau_{da}= \tau_a\circ (\bd_{T_a})^*,\hskip 10mm
  \tau_{c(a)}= c(\tau_a) = \tau_{c(a)}= \tau_a\circ (\bc_{T_a})^*, \]
  \[ \tau_{1_a} = \tau_a \]
  \item for any $a,b\in \cP$ such that $c^{(k)}(a)=d^{(k)}(b)$ we have
  \[ \tau_{a;_kb} = [\tau_a,\tau_b]: T_a^*;_kT_b^* \lra \cP,\]
  \item for any ordered face structure $S$, for any computad map $f:S^*\lra \cP$,
  \[ \overline{\tau_{f_n(1_S)}}=f. \]
  \item for any ordered face structure $S$, any $\o$-functor $f:S^*\lra \cP$
  can be essentially uniquely factorized as
  \begin{center}
\xext=600 \yext=300 \adjust[`I;I`;I`;`I]
\begin{picture}(\xext,\yext)(\xoff,\yoff)
 \settriparms[1`1`-1;300]
 \putVtriangle(0,0)[S^*`\cP`T_{f(1_S)}^*;f`f^{in}`\tau_{f(1_S)}]
\end{picture}
\end{center}
  where $f^{in}$ is an inner map (i.e. $f^{in}(1_S)=1_{T_{f(1_S)}}$) and $(\tau_{f(1_S)},T_{f(1_S)})$ is the
  description of the cell $f(1_S)$.
  \end{enumerate}
\end{proposition}

{\it Proof.}~ Using the above description of the many-to-one
computad $\cP$ we have that $a:(T_a)^{\sharp,n}\lra
\cP^{\natural,n}$. We put $\tau_a=\overline{a}$. By Proposition
\ref{type} point 6, we have that
$(\tau_a)_n(1_{T_a})=\overline{a}_n(1_{T_a})=a$, as required. The
uniqueness of $(T_a,\tau_a)$ follows from Proposition \ref{type}
point 5. The remaining part is left for the reader. $~~\Box$

\newpage
\section{Appendix} 

{\bf A definition of the many-to-one computads and the comma
categories}

The notion of a computad was introduced by Ross Street. We repeat
this definition for a subcategory $\mComp$ of the category of all
computads $\comp$\index{category!comp@$\comp$} that have
indeterminates of a special shape, namely their codomains are
again indeterminates. We use
this opportunity to introduce the notation used in the paper. In
order to define $\mComp$\index{category!comp@$\mComp$} we define
three sequences of categories
$\mnComp$\index{category!compn@$\mnComp$},
$\mnComma$\index{category!comma@$\mnComma$}, and $\nComma$.

\begin{enumerate}
  \item For $n=0$, the categories $\mnComp$, $\mnComma$, and $\nComma$ are just
  $Set$, and the functor $\overline{(-)}^n: \mnComma \lra \mnComp$ is
  the identity.
  \item For $n=1$, the categories $\mnComma$ and $\nComma$ are the category of
  graphs (i.e. 1-graphs) and $\mnComp$ is the category of free
  $\o$-categories over graphs with morphisms being the functors sending
  indets (=indeterminates=generators) to indets.
  \item Let $n\geq 1$.  We define the following functor
  \[ \varpi^{m/1}_n : \mnComp \lra Set  \]
  such that
  \[ \varpi^{m/1}_n (\cP)=\{ (a,b): \; a\in
  \cP_n,\; b\in
  |\cP|_n,\; d(a)=d(b),\; c(a)=c(b)\} \]
 i.e. $\varpi^{m/1}_n (\cP)$ consists of those parallel pairs $(a,b)$ of $n$-cells of
 $\cP$ such that $b$ is an indet. On morphisms
 $\varpi_n$ is defined in the obvious way. We define $\mnjComma$ to be
 equal to the comma category $Set\da \varpi^{m/1}_n$. So we have a
 diagram
 \begin{center}
\xext=1200 \yext=600 \adjust[`I;I`;I`;`I]
\begin{picture}(\xext,\yext)(\xoff,\yoff)
 \settriparms[1`1`1;600]
 \putAtriangle(0,50)[\mnjComma`\mnComp`Set;(-)_{\leq n}`|-|_{n+1}`\varpi^{m/1}_n]
 \put(850,100){\makebox(100,100){$\mu\Downarrow$}}
\end{picture}
\end{center}

 \item For $n\geq 1$, we can define also a functor
 \[ \varpi_n : \nC \lra Set  \]
  such that
 \[ \varpi_n (C)=\{ (a,b): a,b \in C_n,\; d(a)=d(b),\; c(a)=c(b)\} \]
 i.e. $\varpi_n (C)$ consists of all parallel pairs $(a,b)$ of $n$-cells of
 the $n$-category $C$.
 We define $\njComma$ to be equal to the comma category $Set\da
 \varpi_n$. We often denote objects of $\njComma$ as quadruples
 $C=(|C|_{n+1},C_{\leq n},  d,c)$, where $C_{\leq n}$ is an
 $n$-category, $|C|_{n+1}$ is a set and $(d,c):|C|_{n+1}\lra
 \varpi_n(C_{\leq n})$ is a function. Clearly, the category $\mnjComma$
 is a full subcategory of $\njComma$, moreover we have a forgetful functor
 \[ {\cal U}_{n+1} : (n+1)\cat \lra \njComma\]
 such that for an $(n+1)$-category $A$
 \[ {\cal U}_{n+1}(A)=(A_{n+1},A_{\leq n}, d,c) \]
 i.e. ${\cal U}_{n+1}$ forgets the structure of compositions and identities at the
 top level. This functor has a left adjoint
 \[ {\cal F}_{n+1} : \njComma\lra (n+1)\cat  \]
 The category ${\cal F}_{n+1}(|B|_{n+1},B, d,c)$ is said to be
 a {\em free extension}\index{category!free extension of}\index{free extension}
 of the $n$-category $B$ by the indets $|B|_{n+1}$.
 The category of many-to-one $(n+1)$-computads $\mnjComp$
 is a subcategory of $(n+1)\cat$ whose objects are free extensions
 of objects from $\mnjComma$. The morphisms in $\mnjComp$ are
 $(n+1)$-functors that sends indets to indets. Thus the functor ${\cal
 F}_{n+1}$ restricts to an equivalence of categories
  \[ {\cal F}^{m/1}_{n+1} : \mnjComma\lra \mnjComp,  \]
 its essential inverse will be denoted by
 \[ \|-\|_{n+1} : \mnjComp\lra\mnjComma .\]
 Thus for an $(n+1)$-computad $\cP$ we have
 $\| \cP \|_{n+1} =(|\cP|_{n+1},\cP_{\leq n},d,c)$.

 \item The category $\mComp$ is the category of such
 $\o$-categories $\cP$, that for every $n\in\o$, $\cP_{\leq n}$ is a
 many-to-one $n$-computad, and whose morphisms are
 $\o$-functors sending indets to indets.

 For $n\in\o$, we have functors
 \[ |-|_n : \mComp \lra Set \]
  associating to computads their $n$-indets, i.e.
   \[ f:A\lra B \mapsto |f|_n:|A|_n\lra |B|_n, \]
  they all preserve colimits. Moreover we have a functor
 \[|-|: \mComp \lra Set \]
  associating to computads all their indets, i.e.
   \[ f:A\lra B \mapsto |f|:|A|\lra |B|, \]
   where
   \[ |A|=\coprod_{n\in\o}|A|_n. \]
  It also preserves colimits and moreover it is is
   faithful.

   \item We have a truncation functor \[(-)_{\leq n} : \oC \lra \nC \]
  such that
   \[ f:A\lra B \mapsto f_{\leq k}:A_{\leq k}\lra B_{\leq k} \]
   with $k\in\o$, it preserves limits and colimits.
\end{enumerate}


\begin{theindex}

  \item category
    \subitem $\hg$, 10
    \subitem $\lfs$, 13
    \subitem $\ofs$, 12
    \subitem $\mnComma$, 94
    \subitem $\comp$, 94
    \subitem $\mComp$, 94
    \subitem $\mnComp$, 94
    \subitem free extension of, 95
    \subitem $\kC$, 74
    \subitem $\lfs_k$, 74
    \subitem $\nfs$, 14
    \subitem $\pfs$, 14
  \item cell
    \subitem description of a -, 85
  \item convex subset, 23
  \item cut, 29, 50
    \subitem  lower description of -, 50
    \subitem  upper description of -, 50

  \indexspace

  \item decomposition, 77
    \subitem lower part of -, 77
    \subitem proper, 77
    \subitem upper part of -, 77
  \item depth, 23
  \item description of a cell, 85
  \item disjointness, 12

  \indexspace

  \item face
    \subitem -s based on $x$, 50
    \subitem depth of -, 23
    \subitem empty domain -, 11
    \subitem height of -, 23
    \subitem loop, 11
    \subitem non-empty domain, 11
    \subitem non-loop, 11
    \subitem unary, 11
    \subitem weight of -, 21
  \item face structure
    \subitem local -, 13
    \subitem ordered, 12
      \subsubitem $n$ -, 14
      \subsubitem normal -, 14
      \subsubitem principal -, 14
  \item free extension, 95

  \indexspace

  \item globularity, 12

  \indexspace

  \item height, 23
  \item hypergraph, 9
    \subitem convex sub-, 23, 29
    \subitem morphism, 10

  \indexspace

  \item ideal, 40
    \subitem unary -, 51

  \indexspace

  \item kernel, 39

  \indexspace

  \item linearity
    \subitem pencil -, 12
  \item local discreteness, 12
  \item $\cJ$-loop, 40
  \item loop-filling, 12

  \indexspace

  \item morphism
    \subitem  of hypergraphs, 10
    \subitem collapsing -, 39
    \subitem local, 13
    \subitem monotone, 12

  \indexspace

  \item path
    \subitem flat lower -, 11
    \subitem flat upper -, 11
    \subitem lower -, 11
    \subitem maximal, 22
    \subitem upper -, 11
  \item pencil linearity, 12
  \item positive cover, 49
  \item principal extension, 76

  \indexspace

  \item size
    \subitem of ordered face structure, 14
  \item strictness, 12

  \indexspace

  \item tensor, 65
    \subitem locally determined, 72
    \subitem square, 65
  \item truncation, 10

  \indexspace

  \item weight, 21

\end{theindex}

\end{document}